\documentclass{tac}

\usepackage[all]{xy}

\usepackage[breaklinks=true]{hyperref}%To get clickable links to papers

\thirdleveltheorems
\newtheorem{proposition}[theorem]{Proposition}
\newtheorem{conjecture}[theorem]{Conjecture}

\newcommand{\Set}{{\rm Set}}
\newcommand{\Cat}{{\rm Cat}}
\newcommand{\Span}{{\rm Span}}

\newcommand{\Bicat}{{\rm Bicat}}
\newcommand{\Tricat}{{\rm Tricat}}
\newcommand{\Gray}{{\rm Gray}}

\renewcommand{\hom}{{\rm hom}}
\newcommand{\To}{\Rightarrow}
\renewcommand{\to}{\rightarrow}

\newcommand{\maps}{\colon}
\newcommand{\Hom}{{\rm Hom}}
\newcommand{\B}{\mathcal{B}}
\newcommand{\T}{\mathcal{T}}
\newcommand{\op}{{\rm op}}

\newcommand{\C}{{\mathcal C}}

\newcommand{\Lim}{{\rm Lim}}
\newcommand{\Ob}{{\rm Ob}}

\newdir{ >}{{}*!/-7pt/@{>}}
\newdir{> >}{*!/0pt/@{>}*!/-5pt/@{>}}

\makeindex

\title{Spans in $2$-categories: a monoidal tricategory}

\author{Alexander E. Hoffnung}
\address{Temple University, Department of Mathematics\\ Rm $638$ Wachman Hall\\ $1805$ N. Broad Street \\ Philadelphia PA $19122$}
\eaddress{hoffnung@temple.edu}

\date{TODAY}


\begin{document}

\maketitle

\begin{abstract}
We present Trimble's definition of a tetracategory and prove that spans in (strict) $2$-categories with certain limits have the structure of a monoidal tricategory, defined as a one-object tetracategory.  We recall some notions of limits in $2$-categories for use in the construction of the monoidal tricategory of spans.
\end{abstract}

\section{Introduction}\label{sec-Introduction}

Spans have a central role in the development of higher category theory and its applications.  We prove a theorem on the structure of spans in strict $2$-categories admitting certain limits.  The statement of the theorem is approximately the following:
\[ \textit{Spans in a 2-category with finite limits are the morphisms of a monoidal tricategory.}\]
We make this statement precise in the context of Trimble's definition of \textit{weak $4$-category} or \textit{tetracategory} given by generators and relations.

Since we work with spans in strict $2$-categories rather than in the more general structure of a bicategory, we expect that the resulting monoidal tricategory will not be as weak as possible.  This can be seen explicitly in the data throughout the paper.  The structure of the underlying tricategory is an interesting example in that the modification components are trivial, yet the adjoint equivalences for associativity and units are not identities, and therefore we do not obtain a $\Gray$-category but rather a semi-strict tricategory.
  
We begin with some motivation for higher category theory and some history of key structures in its development.  We then review some key properties of spans found in the literature and give explicit definitions of the limits we intend to use.  An outline of the present article then concludes the introduction.

\subsection{The Categorical Ladder}

The language of category theory has, by its very nature, served as an impetus for the unification of structures across mathematical disciplines.  This has proved both practical and fruitful as the proliferation of new ideas in mathematics continues to explode.

Category theory is a natural extension of modern algebra.  In particular, as a mathematical structure, a category contains a set of objects, a set of morphisms, along with a notion of composition and unit.  This data should then satisfy appropriate coherence laws for associativity and unity.

While category-theoretic structures consist of a set of objects {\it and} a set of morphisms, one can easily find examples of ordinary set-theoretic structures by specifying to one-object categories, i.e., taking the set of objects to be the singleton set.  For example, a one-object category consists of a set of morphisms, an identity morphism for the single object, and an associative unital binary operation called composition --- that is, a one-object category is a monoid.

Passing from a category to a monoid is an example of a process called {\it looping} due to its relation to the loop space construction in topology.  Looping and its counterpart {\it delooping} have served as guiding princinples in the development of the structure of {\it low-dimensional higher categories}, especially in the iterated context of $k$-tuply monoidal $n$-categories.  See~\cite{BaDo98} for details and the ``periodic table of $n$-categories".

Categories were introduced by Eilenberg and Mac Lane as a step towards defining the notion of natural transformation in their study of axiomatic homology theory.  Further development of category theory, in large part through {\it internalization} and {\it enrichment}, vastly generalized the study of universal algebra allowing the proliferation throughout modern mathematics.

Researchers have developed aspects and applications of category theory extensively.  However, higher category theory --- the theory of weak $n$-categories --- is still very much a work in progress.  Numerous definitions of weak $n$-category have been proposed. See~\cite{Le04} and references therein.  The operadic approaches highlight some of the inherently combinatorial and topological aspects of higher categories.

The development of low-dimensional higher categories took strong cues from homotopy theory.  The axioms of tricategories are defined in terms of polyhedra which first appeared in Stasheff's thesis~\cite{Sta1}.  Further, Grothendieck's dream of developing $n$-category theory on top of homotopy theory by interpreting $n$-groupoids as homotopy $n$-types has continued to drive research since its appearance in his famous letter to Quillen~\cite{Gr}. 

A fundamental difficulty in comparing the definitions of an $n$-category is the need for a framework of $(n+1)$-categories in which to work.  The groupoidal approach to higher categories has lessened this difficulty to some extent by considering $\infty$-categories in which all higher morphisms (after a fixed level) are equivalences.  Recent successes in this direction include a rigorous statement and proof of the Stabilization Hypothesis~\cite{BaDo98, Lu}.

\subsection{Low-Dimensional Higher Categories}

While definitions of weak $n$-categories abound, explicit definitions by generators and relations have been given only for $n \leq 4$.  The prospect of writing down such definitions for higher-categorical dimensions is forboding, yet progress is possible.  Baez and Dolan developed a litmus test in the form of a list of key properties, including the Stablization Hypothesis, that a theory of weak $n$-categories should satisfy~\cite{BaDo98}.

\subsection{Bicategories}

Ehresmann introduced the technique of internalization in ambient categories such as differentiable manifolds and topological spaces to create a hybrid categorical setting for algebraic-geometric notions such as Lie groups and topological groups, respectively.

Internalization can be useful in climbing the so-called categorical ladder.  For example, a (strict) $2$-category can be defined by giving a diagramatic presentation of the data of a category and interpreting this in the category of (small) categories and functors.  However, this method promotes equations between functions to equations (rather than natural isomorphisms) between functors, and thus fails to capture important examples which satisfy more general coherence laws.

In $1968$, B\'enabou introduced bicategories (weak $2$-categories) and the three-dimensional structure of morphisms between these consisting of homomorphisms, transformations and modifications~\cite{Be}.  The definition of bicategory was designed, in part, to formalize the idea of treating mathematical structures as morphisms. 

One of the first examples is the bicategory of bimodules, which has (unital) rings as objects, bimodules of rings as morphisms, and bimodule maps as $2$-morphisms.  In this bicategory, composition is the tensor product of bimodules.  B\'enabou's next example is the bicategory of spans in a category with pullbacks.  He attributes the notion of span to Yoneda who first considered spans in the category of categories~\cite{Yo}.

The two most important aspects of the theory of bicategories in relation to this article and future directions are the statement of bicategorical coherence and the development of appropriate limits in the bicategorical setting.

The coherence theorem for bicategories states that every bicategory is biequivalent to a (strict) $2$-category.  This is proved in the abstract setting of Yoneda's lemma for bicategories in~\cite{GPS}.  The same article gives a more pedestrian proof following~\cite{MaPa},~\cite{Lew}, and~\cite{JoSt93}.  The development of tricategories and coherence for tricategories relies heavily on this result both for logical and aesthetic reasons.  Similarly, the tricategorical coherence statement of Gordon, Power and Street makes possible Trimble's definition of tetracategory.

The present article is concerned with spans as an example of a monoidal tricategory.  In particular, we consider spans in a $2$-category with certain limits.  The theory of $2$-dimensional limits has been successfully and carefully developed by a number of authors.  While the formulation of $2$-dimensional limits very much mirrors the traditional theory of limits for categories, the full story is considerably more complicated when developing completeness theorems for $2$-categories and bicategories.

Flexible limits, a certain type of weighted strict $2$-limit, were developed in~\cite{BiKePoSt} to obtain stronger completeness theorems for $2$-categories which admit all bilimits, but not all honest limits.  These authors found that completeness theorems could be formulated with the introduction of flexible limits.

In practice, bicategorical coherence allows for simpler computations when working in $2$-dimensional category theory.  However, in addition to such a coherence theorem, an additional coherence theorem for appropriately weak limits in bicategories was needed to remove many of the remaining imposing technicalities of fully weak $2$-dimensional category theory.  Power took an important next step providing the desired coherence theorem for bilimits~\cite{Po}.

\theorem [Power]
Every bicategory $\mathcal{B}$ with finite bilimits is biequivalent to a $2$-category with finite flexible limits.
\endtheorem

We briefly recall the details of Power's coherence theorem.  This result relies mainly on some nice features of the $2$-categorical Yoneda embedding
\[ Y\maps \mathcal{B}\to\hom(\mathcal{B}^{{\rm op}},\Cat),\]
\noindent which is defined as the contravariant $\hom$-functor.  This should not be surprising given the central role of representable presheaves in defining limits.

We first note that whenever a pseudo limit exists, then it also acts as a bilimit.  This follows from the fact that any strict $2$-functor is also a homomorphism of bicategories and both pseudo limits and bilimits are defined by strong transformations between $2$-functors and homomorphisms, respectively.  Further, the Yoneda $2$-functor preserves all pseudo limits and for any small $2$-category $\mathcal{B}$, the $2$-category $\hom(\mathcal{B}^{{\rm op}},\Cat)$ has all pseudo limits.

We will not need all finite bilimits in our present work; rather just products and pullbacks.  Our working definitions of products and pullbacks (see Section~\ref{explicit}) are both examples of finite flexible limits.  Recall that the terminal object is the `empty product'.  Although we do not consider spans in the fully weak setting of bicategories here, Power's result suggests the following.

\conjecture
Let $\mathcal{B}$ be a bicategory with pullbacks and finite products that is biequivalent to a $2$-category $\mathcal{B}'$ with pullbacks and finite strict products, then there is a monoidal tricategory $\Span(\mathcal{B})$ that is monoidally triequivalent (or tetraequivalent) to our construction $\Span(\mathcal{B}')$.
\endconjecture

\subsection{Tricategories}

The importance of a robust theory of bicategories is at this point apparant to many working mathematicians.  Gordon, Power and Street motivate the need for the introduction of tricategories both from category theory and from applications to representation theory, low-dimensional topology, and other areas.  We briefly recall some of their motivating examples.

The main push from within category theory comes from the need to consider monoidal structures on bicategories.  Recall, that a monoidal bicategory is defined to be a one-object tricategory.  Walter's theory of categories enriched in bicategories calls for monoidal structures on the enriching bicategory in order to define the domain of the composition operation in the enriched structure~\cite{Wa}.  Further, Carboni and Walter's work on bicategories of relations in a regular category naturally extends to the monoidal bicategory setting with the composition and monoidal product induced by finite limits~\cite{CaWa}.  The present work is a natural analogue of the work of Carboni and Walter, however this work considers finite limits in $2$-categories, and necessitates the introduction of monoidal structures on tricategories.

Gordon, Power, and Street note the appearance of tricategorical data and coherence outside of category theory as well.  In particular, in the study of algebraic homotopy $3$-types of Joyal-Tierney, the interplay of Zamolodchikov equations and monoidal bicategories as discussed in the strict setting by Kapranov and Voevodsky, the appearance of the $4$-cocycle condition in Drinfel'd's study of equivalences of representations of quasi-Hopf algebras, and the theory of operads developed to analyze iterated loop spaces which becomes relevant through examples in conformal field theory and string field theory.

The main theorem of~\cite{GPS} is the coherence theorem for tricategories.  The statement of tricategorical coherence is not as clean as the analogous statement for bicategories.  In particular, one does not expect that every tricategory is triequivalent to a (strict) $3$-category and the theorem requires a detailed analysis of the local structure of tricategories as well as techniques of enriched category theory involving the $\Gray$ tensor product.  The category $\Gray$ is closely related to the category of $2$-categories with the important distinction that the monoidal structure on $\Gray$ not be cartesian.

We recall the main coherence theorem.

\theorem [Gordon, Power, Street]
Every tricategory is triequivalent to a $\Gray$-category.
\endtheorem

Since the appearance of this coherence theorem, there has been continued investigation into the categorical structures formed by tricategories.  In~\cite{Le04}, it is pointed out that the tricategory definition in~\cite{GPS} is not algebraic.  This indicates that the definition is not amenable to certain operadic approaches to higher categories.  In particular, tricategories as defined therein are not governed by (that is, are not the algebras of) an appropriate operad.  The thesis of Gurski is devoted to the study of algebraic tricategories.  In particular, he is able to amend the original definition to produce a fully algebraic definition~\cite{Gu}.

\subsection{Tetracategories}

The coherence theorem for tricategories is essential to Trimble's definition of tetracategories as was the coherence theorem for bicategories in Gordon, Power, and Street's work on tricategories.  One of the fascinating aspects of Trimble's definition is the use of combinatorial structures to obtain a nearly algorithmic definition of weak $n$-categories given by generators and relations.  There is much to be done however as we do not have corresponding definitions for morphisms between tetracategories, a theory of limits, nor a coherence theorem for tetracategories.  Given the sheer size of the data and axioms for tetracategories, it remains unclear whether researchers will continue to pursue the development of higher category theory in this fashion.  Nonetheless, Trimble's work is intriguing, illuminating, and allows for a concrete example of a monoidal tricategory to be given in the present work.  To the best of our knowledge the span construction given here is the first explicit example of a monoidal tricategory in the literature.

Trimble includes remarks on his definition which have been presented on a web page of John Baez devoted to Trimble's work~\cite{Tr}.  Before presenting the definition in Section~\ref{tetra-definition}, we provide some exposition in an attempt elucidate some of the key ideas behind Trimble's work.  However, we also include Trimble's own remarks in Appendix~\ref{appendix}.

\subsection{Explicit Definitions of Pullbacks and Products}\label{explicit}

We intend to construct a composition operation and monoidal product on spans by assuming the existence of certain limits in our $2$-categories.  We give explicit definitions of these limits here and a more detailed exposition on $2$-dimensional limits in Appendix~\ref{sec-Limits}.

Composition of spans is defined by pullback of cospan diagrams, but there are several closely related limits of cospan diagrams.  The construction we work with at present is called the {\em iso-comma object}, however, we will usually refer to it simply as the {\em pullback}. The monoidal product is constructed by {\em products}.

\subsection*{Definition of Pullbacks}\label{pullback}
We give an explicit definition of the pullback or iso-comma object and its universal property as a limit. This is used to define certain composites of spans and higher morphisms as part of the structure of the tricategory $\Span(\B)$.

\begin{definition}\label{iso_comma}
Given a cospan
\[
   \xy
   (-15,0)*+{B}="2"; (15,0)*+{A}="3";
   (0,-15)*+{C}="4";
        {\ar_{g} "2";"4"};
        {\ar^{f} "3";"4"};
\endxy
\]
in a strict $2$-category $\mathcal{B}$, the {\bf pullback} is an object $BA$, equipped with projections $\pi_{A}^{B}$, $\pi_{B}^{A}$, and a $2$-cell $\kappa_{C}^{g,f}$
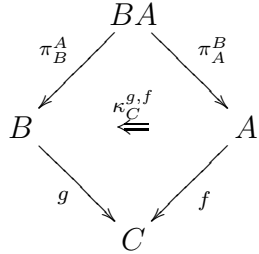
\[
   \xy
   (0,15)*+{BA}="1";
   (-15,0)*+{B}="2"; (15,0)*+{A}="3";
   (0,-15)*+{C}="4";
        {\ar_{g} "2";"4"};
        {\ar^{f} "3";"4"};
        {\ar_{\pi_{B}^{A}} "1";"2"};
        {\ar^{\pi_{A}^{B}} "1";"3"};
        {\ar@{=>}_{\scriptstyle \kappa_{C}^{g,f}} (2,0); (-2,0)};
\endxy
\]
such that

$\bullet$ for any pair of maps $p\maps X\to A$ and $q\maps X\to B$ and $2$-cell 
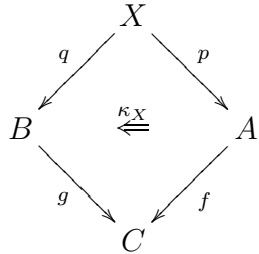
\[
   \xy
   (0,15)*+{X}="1";
   (-15,0)*+{B}="2"; (15,0)*+{A}="3";
   (0,-15)*+{C}="4";
        {\ar_{g} "2";"4"};
        {\ar^{f} "3";"4"};
        {\ar_{q} "1";"2"};
        {\ar^{p} "1";"3"};
        {\ar@{=>}_{\scriptstyle \kappa_{X}} (2,0); (-2,0)};
\endxy
\]
there exists a unique $1$-cell $h\maps X\to BA$ such that
\[ p=\pi_{A}^{B}h,\; q=\pi_{B}^{A}h, \textrm{ and } \kappa_{C}^{g,f}\cdot h=\kappa_{X},\]

$\bullet$ for any pair of $1$-cells $j,k\maps X\to BA$ and $2$-cells
\[\varpi\maps \pi_{A}^{B}j\To\pi_{A}^{B}k \textrm{ and } \varrho\maps \pi_{B}^{A}j\To \pi_{B}^{A}k,\]
such that
\[ (g\cdot\varrho)(\kappa_{C}^{g,f}\cdot j) = (\kappa_{C}^{g,f}\cdot k)(f\cdot\varpi),\]
there exists a unique $2$-cell $\gamma\maps j\To k$ such that
\[\pi_{A}^{B}\cdot\gamma = \varpi\textrm{ and } \pi_{B}^{A}\cdot\gamma = \varrho.\]
\end{definition}

\subsection*{Definition of Finite Products}

We define the product in $\mathcal{B}$ for use in constructing the monoidal structure on our tricategory of spans.

\begin{definition}\label{strictproduct}
Given a pair of objects in a strict $2$-category $\mathcal{B}$, the {\bf product} is an object $A\times B$ equipped with projections $A\stackrel{\pi_A}\leftarrow A\times B\stackrel{\pi_B}\rightarrow B$ such that

$\bullet$ for each pair of maps $p\maps X\to A$ and $q\maps X\to B$, there exists a unique $1$-cell $h\maps X\to A\times B$ such that
\[ p=\pi_{A}h \textrm{ and } q=\pi_{B}h,\]

$\bullet$ for each pair of $1$-cells $j,k\maps X\to A\times B$ and $2$-cells
\[\varpi\maps \pi_{A}j\To\pi_{A}k \textrm{ and } \varrho\maps \pi_{B}j\To \pi_{B}k,\]
\noindent there exists a unique $2$-cell $\gamma\maps j\To k$ such that
\[\pi_{A}\cdot\gamma = \varpi \textrm{ and } \pi_{B}\cdot\gamma = \varrho.\]
\end{definition}

\noindent We also include the nullary product, or terminal object as part of the structure of $\mathcal{B}$.  In fact, a $2$-category with iso-comma objects and a terminal object automatically has finite products, which are obtained by the obvious cospan with arrows into the terminal object.

\begin{definition}
We call an object $1\in\mathcal{B}$ the {\bf terminal object} if for every object $A\in\mathcal{B}$, there is a unique $1$-cell from $A$ to $1$.
\end{definition}

\subsection{Basics on Spans}

Motivating examples for the development of bicategories often do not satisfy the associative law and require the introduction of an {\it associator} natural transformation satisfying a coherence equation given by Mac Lane's pentagon.  Further, as discussed in the motivation for tricategories, it is often the case that given a theory of $n$-categories (for a fixed $n$) the need to develop a theory of ($n+1$)-categories arises from examples of $n$-categories which seem to naturally carry a monoidal structure. The span construction is a well-known example of both of these phenomena.  

Recall that a {\it span} in a category $\mathcal{C}$ is a pair of morphisms with a common domain.  This is often drawn in a shape reminiscent of a bridge or roof
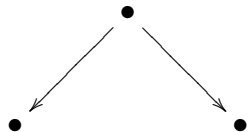
\[   \xy
   (0,15)*+{\bullet}="1";
   (-15,0)*+{\bullet}="2";
   (15,0)*+{\bullet}="3";
        {\ar_{} "1";"2"};
        {\ar^{} "1";"3"};
\endxy
\]
\noindent thus spawning a number of different names for this diagram, including the very common name {\it correspondence}.  Spans $R$ and $S$ are composable if they have common codomain and domain, respectively
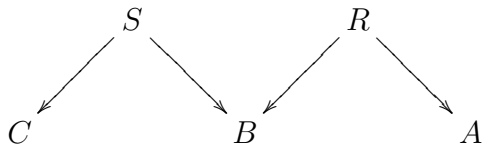
\[   \xy
   (-15,15)*+{S}="1";
   (-30,0)*+{C}="2";
   (0,0)*+{B}="3";
   (15,15)*+{R}="4";
   (30,0)*+{A}="5";
        {\ar_{} "1";"2"};
        {\ar^{} "1";"3"};
        {\ar_{} "4";"3"};
        {\ar^{} "4";"5"};
\endxy
\]
\noindent and there exists a reasonable notion of pullback, e.g., a limit $SR$
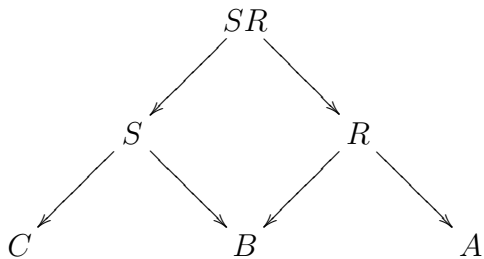
\[   \xy
   (0,30)*+{SR}="0";   
   (-15,15)*+{S}="1";
   (-30,0)*+{C}="2";
   (0,0)*+{B}="3";
   (15,15)*+{R}="4";
   (30,0)*+{A}="5";
        {\ar_{} "1";"2"};
        {\ar^{} "1";"3"};
        {\ar_{} "4";"3"};
        {\ar^{} "4";"5"};
        {\ar_{} "0";"1"};
        {\ar^{} "0";"4"};
\endxy
\]
\noindent exists.  These constructions are variously called `pullbacks', `fibered products', `homotopy pullbacks', `weak pullbacks', `pseudo pullbacks', `bipullbacks', `comma objects', `iso-comma objects', `lax pullbacks', `oplax pullbacks', etc.  In some cases, some of these names can be freely interchanged, however some are definitively different from each other.   See~\cite{KePr} for a detailed study of composition operations on spans.

This composition operation is {\em not} strictly associative.  As alluded to above this is essentially the original motivation for the generalization from (strict) $2$-categories to bicategories.  By introducing a suitable notion of `maps of spans', Benabou was able to introduce spans of sets, or, more generally, spans in a category with pullbacks, as an example of a bicategory~\cite{Be}.  In this paper, we categorify Benabou's work.  This requires us to introduce a suitable notion of `maps between maps of spans'.   These maps are the $3$-morphisms of the span construction.

It is often the case that a category with pullbacks also has finite products.  This is, in part, a consequence of the fact that finite products are limits of similar diagrams and can be obtained as a special case of pullbacks in the presence of a terminal object.  Now, given a category with pullbacks and finite products, applying the span construction yields not only a composition from pullbacks, but also a monoidal product induced by products.  It follows that applying the span construction to a category with pullbacks and finite products yields a monoidal bicategory.  Notice that we began with a $1$-category with some extra structure and after applying the span construction obtained a one-object tricategory.

Just as maps of monoidal bicategories include important examples such as bimodules and spans, the need for monoidal tricategories continues to increase as researchers begin to study spans in higher categories.  Common examples are spans of groupoids or stacks, which arise in geometric representation theory.  A letter entitled `Geometric function theory'~\cite{DBZ} written by David Ben-Zvi to the readers of the $n$-{\it Category Caf\`e} blog dicusses the appearance of spans in representation theory.

Moreover, spans are ubiquitous in mathematics for a very simple reason --- they are a straightforward generalization of relations, which can be used to to define partial maps, generalize aspects of quantum theory such as Heisenberg's matrix mechanics~\cite{MoVi}, give geometric constructions of convolution products in representation theory~\cite{ChGi}, and much more!

\subsection{Organization of Paper}

We now say a few words on the structure of this paper.  The span construction defined here utilizes particular examples of pseudolimits in $2$-categories, which we discuss in an exposition on $2$-dimensional limits in Appendix~\ref{sec-Limits}.  The definition of monoidal tricategory is given in Section~\ref{sec-definitions}.   The main theorem is the construction of a monoidal tricategory, or one-object tetracategory.  This comes in two parts.  The first is to construct a tricategory of spans denoted $\Span(\B)$ in Section~\ref{sec-tricategory} using pullbacks (by which we mean iso-comma objects).  The second is the construction of the monoidal structure on the tricategory of spans in Section~\ref{monoidal} using products.

The span construction yields a relatively weak monoidal tricategory, demanding a significant effort in verifying coherence axioms.  The techniques used in verifying these axioms all follow very similar reasoning as the components of the structural maps are defined almost invariably by the existence statements of the universal property of pseudolimits, leaving the uniqueness statements as the main tools used in verifying equations.  We work through a few of these arguments, but leave most out of the text.  Instead we include all of the structural data along with the equations satisfied by the components of this data.  In each instance, this is enough to routinely reproduce and verify the necessary coherence equations.

\section{Monoidal Tricategories as One-Object Tetracategories}\label{sec-definitions}
\subsection{Approaching a Definition}\label{subsec-approaching}

As the goal of this paper is to define a monoidal tricategory of spans and `monoidal tricategory' is not a well-defined notion in the literature, we need to first specify what structure we have in mind.  There is little doubt that a number of people either have or could write down a reasonable notion of monoidal structure on a tricategory if asked.  In 1995, Trimble went a step further and wrote down a definition of tetracategory with axioms sprawling over dozens of pages~\cite{Tr}.  Following the pattern of defining a monoidal category to be a one-object category one dimension above, we say
\[ \textrm{A monoidal tricategory is a one-object Trimble tetracategory.}\]
\noindent In recalling Trimble's definition we hopefully succeed in making tetracategories accessible to a wide audience.  We explain Trimble's notion of {\it product cells} for tritransformations and trimodifications, give a precise statement of the equivalence expected for structure cells at each level, and choose explicit $3$-cells (geometric $2$-cells in local tricategories) as the `interchange' cells appearing in the tetracategory axioms.  The choice of interchange cell is governed by coherence for tricategories, so all choices are suitably equivalent.

The definition of a tetracategory is largely straightforward.  In fact, Trimble's approach to defining tetracategories was to, as much as possible, formalize the process of drawing coherence axioms, at least up to coherent isomorphism.  Just as monoids, or one-object categories, have associativity and unit axioms, higher categories have generalized associativity and unit coherence axioms.  One starts by noting that the drawing of associativity axiom $K_{n}$ is nearly canonical at each level $n$.  These axioms first appeared as families of simplicial complexes called {\it associahedra} in work of Stasheff and called {\it orientals} in work of Street.  
These associator $K_{n+2}$ axioms can, in turn, be used to define unit axioms $U_{n+1,1} ,\ldots, U_{n+1,n+1}$ for weak $n$-categories.

It is useful to work up from the usual category axioms towards tetracategories developing intuition for the higher unit diagrams and building on successive steps.  It can also be useful to think of categorical structure as consisting of both associativity and unit operations and axioms.  The coherence axioms for categories include one associativity axiom $K_3$
\[ \otimes(\otimes\times 1) = \otimes(1\times\otimes)\]
\noindent and two unit axioms $U_{2,1}$
\[ \otimes (I\times 1) = 1\]
\noindent and $U_{2,2}$
\[ \otimes (1\times I) = 1.\]
\noindent Here the tensor product denotes the composition operation on the category and the three axioms often denoted $\alpha$, $\lambda$, and $\rho$, respectively, for obvious reasons.  We then recall that a category $\C$ also has a unit operation $I\in\Ob(\C)$.  Of course the unit object has an identity morphism, so we can write $I$ as a functor to make it appear more like an operation
\[ I\maps 1\to \C.\]
\noindent The unit operations and axioms are closely tied to those for associativity.  We know that $I$ is the unit for our composition operation $\otimes\maps \C\times\C\to\C$.  In fact, $\otimes$ is the associativity operation in this case and has the associativity axiom $K_3$ as noted above.  Although slightly awkward in many contexts it is useful here for combinatorial reasons to write $K_2$ for the composition operation.  Similarly, we can write $U_{1}$ as the unit operation.

Formally, bicategories include both of the category operations $K_2$ and $U_1$, which are now interpreted as functors between bicategories.  The category axioms become bicategory operations
\[ K_{3}\maps K_{2}(K_{2}, 1) \To K_{2}(1, K_{2}),\]
\[ U_{2,1}\maps K_{2}(U_{1}, 1) \To 1,\]
\noindent and 
\[ U_{2,2}\maps K_{2}(1, U_{1}) \To 1.\]
\noindent The bicategory associativity axiom is the MacLane pentagon, which we write algebraically as
\[
   \xy
   (-31,0)*+{K_{4}\maps K_{3}(1,K_{2})\circ K_{3}(K_{2},1)}="2"; (30,0)*+{K_{2}(1,K_{3})\circ K_{3}(1,K_{2})\circ K_{2}(K_{3},1).}="3";
        {\ar@3{->}_{} "2";"3"};
\endxy
\]
\noindent Notice that each possible $4$-ary operation with one occurrence of each $K_{2}$ and $K_{3}$ appears as a $1$-cell (edge) of the pentagon.  A similar pattern can be observed for the $0$-cells (vertices) of the pentagon.  We are left to derive the unit axioms from $K_3$.

We use $K_3$ to construct a template informing the general shape of the unit axioms $U_{3,1}$, $U_{3,2}$, and $U_{3,3}$.  The second index $i$ in $U_{3,i}$ tells us that the unit object should appear in the $i^{th}$ argument.  The unit axiom $U_{3,1}$ will have $1$-cells in the domain resembling those in the domain of $K_{3}$, except that a copy of the unit object will be placed in the first argument of each operation.  Similarly, the codomain will contain $1$-cells from the codomain of $K_{3}$ with units in the first argument.  There is an extra associativity term appearing in the unit operation, which we now describe.  

If we imagine an associativity operation for $0$-categories, $K_{1}\maps 1\to 1$, then we notice a pattern beginning with the bicategory unit operations.  The domain and codomain of $U_{2,1}$ contain the domain of and codomain of $K_{1}$, respectively.   We notice the appearance of the associativity operation $K_{2}$ in the domain.  It turns out that this is a very general phenomenon.  In each unit axiom $U_{n,i}$ there is a cell $K_{n}$ in the domain if $n$ is odd or in the codomain if $n$ is even.  Finally, there is a copy of $U_{1}$ appearing in the first and second arguments of $K_{2}$ in $U_{2,1}$ and $U_{2,2}$, respectively.

From the above considerations the pattern begins to become evident.  We have bicategory axioms
\[\textrm{Unit axiom } U_{3,1}\hspace{20pt} K_{2}(U_{2,1},1) = U_{2,1}(K_{2})\circ K_{3}(U_{1},1,1)\]
\[\textrm{Unit axiom } U_{3,2}\hspace{20pt} K_{2}(U_{2,2},1) = U_{2,2}(K_{2})\circ K_{3}(1,U_{1},1)\]
\[\textrm{Unit axiom } U_{3,3}\hspace{20pt} U_{2,2}(K_{2}) = K_{2}(1,U_{2,2})\circ K_{3}(1,1,U_{1}).\]
\noindent Trimble tames the complexity of his tetracategory axioms using operad-like trees to name the associativity and unit operations and axioms.  While the construction outlined in the previous paragraph may not be entirely transparent in the unit axioms above, a brief explanation of the tree diagrams should provide clarity.

The $K_{3}$ associativity operation is labelled by a ``$3$-sprout" with domain and codomain the expected pair  of rooted trees each having three leaves and one internal edge
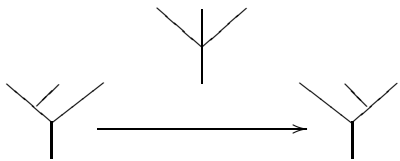
\[
   \xy
   (0,0)*+{}="-1";
   (0,7)*+{}="1"; 
   (0,5)*+{}="-2";
   (0,12)*+{}="2"; 
   (-7,12)*+{}="-3";  
   (1,5)*+{}="3";
   (7,12)*+{}="-4";
   (-1,5)*+{}="4";
   (-20,-10)*+{}="-1d";
   (-20,-3)*+{}="1d"; 
   (-23,-3)*+{}="-2d";
   (-18,2)*+{}="2d"; 
   (-27,2)*+{}="-3d";  
   (-19,-5)*+{}="3d";
   (-12,2)*+{}="-4d";
   (-21,-5)*+{}="4d";
   (-15,-5)*+{}="dom";
   (15,-5)*+{}="cod";
   (20,-10)*+{}="-1c";
   (20,-3)*+{}="1c"; 
   (23,-3)*+{}="-2c";
   (18,2)*+{}="2c"; 
   (27,2)*+{}="-3c";  
   (19,-5)*+{}="3c";
   (12,2)*+{}="-4c";
   (21,-5)*+{}="4c";
        {\ar@{->}_{} "dom";"cod"};
        {\ar@{-}_{} "-1";"1"};
        {\ar@{-}_{} "-2";"2"};
        {\ar@{-}_{} "-3";"3"};
        {\ar@{-}_{} "-4";"4"};
        {\ar@{-}_{} "-1d";"1d"};
        {\ar@{-}_{} "-2d";"2d"};
        {\ar@{-}_{} "-3d";"3d"};
        {\ar@{-}_{} "-4d";"4d"};
        {\ar@{-}_{} "-1c";"1c"};
        {\ar@{-}_{} "-2c";"2c"};
        {\ar@{-}_{} "-3c";"3c"};
        {\ar@{-}_{} "-4c";"4c"};
\endxy
\]
\noindent We can write down the unit axioms using trees this time.  We replace solid edges with dashed edges according to the indices of the unit axiom and include the $K_{3}$ factor as the $3$-sprout with a dashed edge.  We have, for example, the axiom $U_{3,1}$
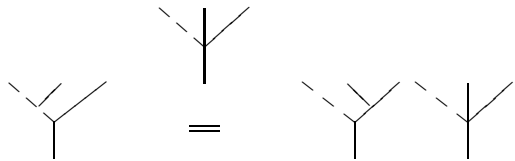
\[
   \xy
   (0,0)*+{}="-1";
   (0,7)*+{}="1"; 
   (0,5)*+{}="-2";
   (0,12)*+{}="2"; 
   (-7,12)*+{}="-3";  
   (1,5)*+{}="3";
   (7,12)*+{}="-4";
   (-1,5)*+{}="4";
   (-20,-10)*+{}="-1d";
   (-20,-3)*+{}="1d"; 
   (-23,-3)*+{}="-2d";
   (-18,2)*+{}="2d"; 
   (-27,2)*+{}="-3d";  
   (-19,-5)*+{}="3d";
   (-12,2)*+{}="-4d";
   (-21,-5)*+{}="4d";
   (-3,-5)*+{}="dom";
   (3,-5)*+{}="cod";
   (20,-10)*+{}="-1c";
   (20,-3)*+{}="1c"; 
   (23,-3)*+{}="-2c";
   (18,2)*+{}="2c"; 
   (27,2)*+{}="-3c";  
   (19,-5)*+{}="3c";
   (12,2)*+{}="-4c";
   (21,-5)*+{}="4c";
   (35,-10)*+{}="-1c2";
   (35,-3)*+{}="1c2"; 
   (35,-5)*+{}="-2c2";
   (35,2)*+{}="2c2"; 
   (42,2)*+{}="-3c2";  
   (34,-5)*+{}="3c2";
   (27,2)*+{}="-4c2";
   (36,-5)*+{}="4c2";
         {\ar@{-}_{} "-1";"1"};
        {\ar@{-}_{} "-2";"2"};
        {\ar@{--}_{} "-3";"3"};
        {\ar@{-}_{} "-4";"4"};
        {\ar@{=}_{} "dom";"cod"};
        {\ar@{-}_{} "-1d";"1d"};
        {\ar@{-}_{} "-2d";"2d"};
        {\ar@{--}_{} "-3d";"3d"};
        {\ar@{-}_{} "-4d";"4d"};
        {\ar@{-}_{} "-1c";"1c"};
        {\ar@{-}_{} "-2c";"2c"};
        {\ar@{-}_{} "-3c";"3c"};
        {\ar@{--}_{} "-4c";"4c"};
        {\ar@{-}_{} "-1c2";"1c2"};
        {\ar@{-}_{} "-2c2";"2c2"};
        {\ar@{-}_{} "-3c2";"3c2"};
        {\ar@{--}_{} "-4c2";"4c2"};
\endxy
\]
\noindent Now rewriting our unit operations as composites of $K$ operations and the unit object, e.g., $U_{2,1} = K_{2}(U_{1},1)$, we see the direct correspondence between the algebraic and tree descriptions of $U_{3,1}$
\[ U_{3,1}\maps K_{2}(K_{2}(U_{1},1),1)\to K_{2}(U_{1},K_{2})\circ K_{3}(U_{1},1,1). \]
\noindent Tree representations of the other unit axioms have the same underlying trees, so differ only by the dashed edge.
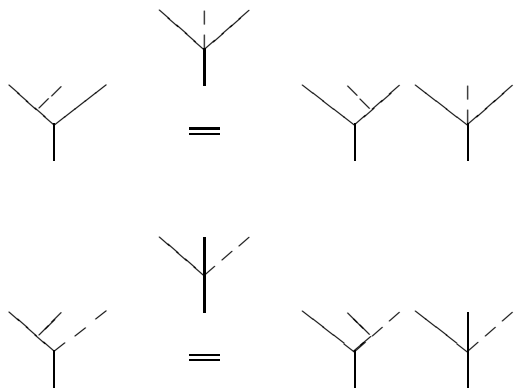
\[
   \xy
   (0,0)*+{}="-1";
   (0,7)*+{}="1"; 
   (0,5)*+{}="-2";
   (0,12)*+{}="2"; 
   (-7,12)*+{}="-3";  
   (1,5)*+{}="3";
   (7,12)*+{}="-4";
   (-1,5)*+{}="4";
   (-20,-10)*+{}="-1d";
   (-20,-3)*+{}="1d"; 
   (-23,-3)*+{}="-2d";
   (-18,2)*+{}="2d"; 
   (-27,2)*+{}="-3d";  
   (-19,-5)*+{}="3d";
   (-12,2)*+{}="-4d";
   (-21,-5)*+{}="4d";
   (-3,-5)*+{}="dom";
   (3,-5)*+{}="cod";
   (20,-10)*+{}="-1c";
   (20,-3)*+{}="1c"; 
   (23,-3)*+{}="-2c";
   (18,2)*+{}="2c"; 
   (27,2)*+{}="-3c";  
   (19,-5)*+{}="3c";
   (12,2)*+{}="-4c";
   (21,-5)*+{}="4c";
   (35,-10)*+{}="-1c2";
   (35,-3)*+{}="1c2"; 
   (35,-5)*+{}="-2c2";
   (35,2)*+{}="2c2"; 
   (42,2)*+{}="-3c2";  
   (34,-5)*+{}="3c2";
   (27,2)*+{}="-4c2";
   (36,-5)*+{}="4c2";
         {\ar@{-}_{} "-1";"1"};
        {\ar@{--}_{} "-2";"2"};
        {\ar@{-}_{} "-3";"3"};
        {\ar@{-}_{} "-4";"4"};
        {\ar@{=}_{} "dom";"cod"};
        {\ar@{-}_{} "-1d";"1d"};
        {\ar@{--}_{} "-2d";"2d"};
        {\ar@{-}_{} "-3d";"3d"};
        {\ar@{-}_{} "-4d";"4d"};
        {\ar@{-}_{} "-1c";"1c"};
        {\ar@{--}_{} "-2c";"2c"};
        {\ar@{-}_{} "-3c";"3c"};
        {\ar@{-}_{} "-4c";"4c"};
        {\ar@{-}_{} "-1c2";"1c2"};
        {\ar@{--}_{} "-2c2";"2c2"};
        {\ar@{-}_{} "-3c2";"3c2"};
        {\ar@{-}_{} "-4c2";"4c2"};
\endxy
\]
\[
   \xy
   (0,0)*+{}="-1";
   (0,7)*+{}="1"; 
   (0,5)*+{}="-2";
   (0,12)*+{}="2"; 
   (-7,12)*+{}="-3";  
   (1,5)*+{}="3";
   (7,12)*+{}="-4";
   (-1,5)*+{}="4";
   (-20,-10)*+{}="-1d";
   (-20,-3)*+{}="1d"; 
   (-23,-3)*+{}="-2d";
   (-18,2)*+{}="2d"; 
   (-27,2)*+{}="-3d";  
   (-19,-5)*+{}="3d";
   (-12,2)*+{}="-4d";
   (-21,-5)*+{}="4d";
   (-3,-5)*+{}="dom";
   (3,-5)*+{}="cod";
   (20,-10)*+{}="-1c";
   (20,-3)*+{}="1c"; 
   (23,-3)*+{}="-2c";
   (18,2)*+{}="2c"; 
   (27,2)*+{}="-3c";  
   (19,-5)*+{}="3c";
   (23.2,-1.3)*+{}="-5c";  
   (18.8,-5)*+{}="5c";
   (12,2)*+{}="-4c";
   (21,-5)*+{}="4c";
   (35,-10)*+{}="-1c2";
   (35,-3)*+{}="1c2"; 
   (35,-5)*+{}="-2c2";
   (35,2)*+{}="2c2"; 
   (42,2)*+{}="-3c2";  
   (34,-5)*+{}="3c2";
   (27,2)*+{}="-4c2";
   (36,-5)*+{}="4c2";
         {\ar@{-}_{} "-1";"1"};
        {\ar@{-}_{} "-2";"2"};
        {\ar@{-}_{} "-3";"3"};
        {\ar@{--}_{} "-4";"4"};
        {\ar@{=}_{} "dom";"cod"};
        {\ar@{-}_{} "-1d";"1d"};
        {\ar@{-}_{} "-2d";"2d"};
        {\ar@{-}_{} "-3d";"3d"};
        {\ar@{--}_{} "-4d";"4d"};
        {\ar@{-}_{} "-1c";"1c"};
        {\ar@{-}_{} "-2c";"2c"};
        {\ar@{--}_{} "-3c";"3c"};
        {\ar@{-}_{} "-5c";"5c"};
        {\ar@{-}_{} "-4c";"4c"};
        {\ar@{-}_{} "-1c2";"1c2"};
        {\ar@{-}_{} "-2c2";"2c2"};
        {\ar@{--}_{} "-3c2";"3c2"};
        {\ar@{-}_{} "-4c2";"4c2"};
\endxy
\]

We have now written down three distinct unit axioms for bicategories, while the usual definition only requires the associativity axiom $K_4$ and the single unit axiom $U_{3,2}$.  This highlights a bit of history in the early development of higher category theory.  When MacLane first defined monoidal categories, he included all three axioms (and some other axioms too).  In 1964 Max Kelly proved that one unit triangle was sufficient.

Thinking ahead one step we have the family of unit axioms for tricategories $U_{4,1}$, $U_{4,2}$, $U_{4,3}$, and $U_{4,4}$.  However, while the definition of tricategory contains all three structural unit operations $U_{3,1}$, $U_{3,2}$, $U_{3,3}$, only the two unit axioms $U_{4,2}$ and $U_{4,3}$ are required.  Gurski remarks on this appearance of cells that are not lifted from equations one-categorical rung down the ladder of $n$-categories as an example of the appearance of interesting higher structure replacing equality~\cite{Gu}.

Apart from the historical significance of the development of the unit axioms, it is also important in understanding Trimble's definition to highlight the relationship between the unit axioms.  In particular, Trimble only defined three unit axioms $U_{5,i}$, $2\leq i\leq 4$ conjecturing the continuation of the pattern for units.

\conjecture $\cite{TrPC}$
Given a non-negative integer $n$, the unit axioms $U_{n,1}$ and $U_{n,n}$ of weak $n$-categories follow from the associativity axiom $K_{n+1}$, the remaining unit axioms $U_{n,i},\; 2\leq i\leq n-1$, and the $U_{n+1,j}$ unit axioms of a weak $n+1$-category.
\endconjecture

Given the above conjecture, Trimble omits the unit axioms $U_{5,1}$ and $U_{5,5}$ in the tetracategory definition.  A proof of the conjecture would involve considering, at each categorical level, the structure of the unit operations and axioms one categorical dimension higher.

Having obtained all of the bicategory axioms we can move onto tricategories and tetracategories.  These have associativity axioms $K_5$ --- the Stasheff polytope --- and $K_6$, respectively.  At this point, the reader is urged to write down the unit axioms $U_{4,2}$ and $U_{4,3}$ for tricategories.  These diagrams appear in the definition below as perturbations along with $U_{4,1}$ and $U_{4,4}$.  Setting these operations to be identity perturbations, the tricategory unit axioms are recovered along with the equations $U_{4,1}$ and $U_{4,4}$, which are redundant axioms for tetracategories.

At this point we should be able to understand unit axioms for tetracategories.  These axioms are three-dimensional and can be rather intimidating at first glance, but can be understood very systematically with the aid of a few preliminary remarks.  The cells on either side of the equations are each components of perturbations, trimodifications, or tritransformations.  Trimble calls the component cells of trimodifications and tritransformations appearing in the axioms `product cells'.  The tetracategory axioms are presented as equations between composites of geometric $3$-cells ($4$-cells of the tetracategory) between surface diagrams.  

The following discussion will use the modification $m$ (see the diagram below) as an example.  There is a similar story for associativity in which triangles are replaced by pentagons.  A modification consists of a family of geometric $2$-cells indexed by objects and a family of invertible modifications indexed by morphisms, which are naturality cells for the geometric $2$-cells.  The components of these naturality modifications are geometric $3$-cells whose domains and codomains have factors corresponding to the domain and codomain of $m$, respectively, and modification $2$-cells corresponding to the domain and codomain of the indexing morphism.  In pictures, the `product' is manifest.  The unitor $2$-cell is a triangle
\[
   \xy
   (-25,15)*+{(X\otimes 1)\otimes Y}="1"; %DOMAIN
   (25,15)*+{X\otimes (1\otimes Y)}="2";
   (0,-10)*+{X\otimes Y}="3";
        {\ar^{\alpha_{X,1,Y}} "1";"2"}; %DOMAIN ARROWS
        {\ar^{1_{X}\otimes\lambda_{Y}} "2";"3"};
        {\ar_{\rho_{X}\otimes 1_{Y}} "1";"3"};
        {\ar@{=>}_{\scriptstyle m_{XY}} (-4,0); (4,6)}; %DOMAIN 2-CELLS
\endxy
\]

\noindent and, for morphisms $(f,f')\maps (X,Y)\to (X',Y')$, the modification $3$-cell, which is the desired product $3$-cell, fills a prism
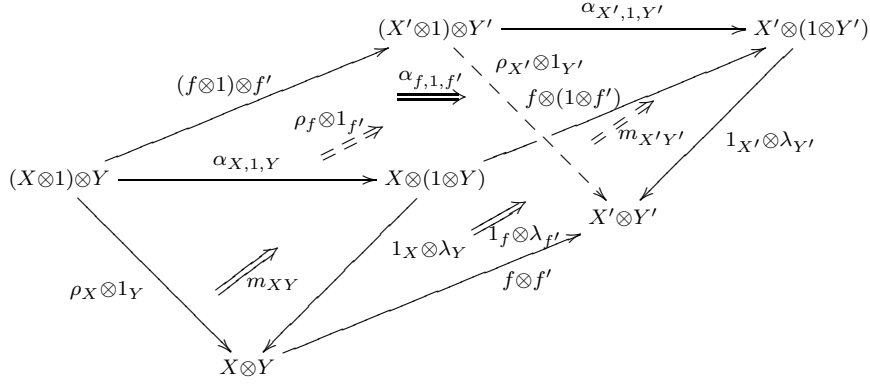
\[
\def\objectstyle{\scriptstyle}
  \def\labelstyle{\scriptstyle}
   \xy
   (-25,15)*+{(X\otimes 1)\otimes Y}="D1"; %DOMAIN
   (25,15)*+{X\otimes (1\otimes Y)}="D2";
   (0,-10)*+{X\otimes Y}="D3";
   (25,35)*+{(X'\otimes 1)\otimes Y'}="C1"; %CODOMAIN
   (75,35)*+{X'\otimes (1\otimes Y')}="C2";
   (50,10)*+{X'\otimes Y'}="C3";
   (39.5,20.5)*+{}="br";
        {\ar^{\alpha_{X,1,Y}} "D1";"D2"}; %DOMAIN ARROWS
        {\ar^<<<<<<<{1_{X}\otimes\lambda_{Y}} "D2";"D3"};
        {\ar_{\rho_{X}\otimes 1_{Y}} "D1";"D3"};
        {\ar^{\alpha_{X',1,Y'}} "C1";"C2"}; %CODOMAIN ARROWS
        {\ar^{1_{X'}\otimes\lambda_{Y'}} "C2";"C3"};
        {\ar@{-->}^<<<<<<<{\rho_{X'}\otimes 1_{Y'}} "C1";"C3"};
        {\ar^{(f\otimes 1)\otimes f'} "D1";"C1"}; %CROSS ARROWS
        {\ar@{-}_{} "D2";"br"};
        {\ar^<<<<<<<{f\otimes (1\otimes f')} "br";"C2"};
        {\ar_>>>>>>>>>>{f\otimes f'} "D3";"C3"};
        {\ar@{=>}_{\scriptstyle m_{XY}} (-4,0); (4,6)}; %DOMAIN 2-CELLS
        {\ar@{==>}_{\scriptstyle m_{X'Y'}} (46,20); (54,25.5)}; %CODOMAIN 2-CELLS
        {\ar@{==>}^{\scriptstyle \rho_{f}\otimes 1_{f'}} (10,18); (18,22)}; %DOMAIN CROSS 2-CELLS
        {\ar@{=>}^{\scriptstyle \alpha_{f,1,f'}} (20,26); (29,26)}; %CODOMAIN CROSS 2-CELLS
        {\ar@{=>}_{\scriptstyle 1_{f}\otimes\lambda_{f'}} (30,8); (37,12)};
\endxy
\]

\noindent We are left only to describe the three rectangular $2$-cells.  The three components of the domain and codomain of the triangle above correspond to the three rectangles.  Two of the cells are in the domain of our $3$-cell and the other is in the codomain.  Finally, we said the prism should be the unitor triangle cross a $K_{3}$-interval, so $(f,f')$ is either $(1,\alpha)$ or $(\alpha,1)$.  When writing down these product cells for the middle modification in the axioms, for example, we defer to coherence for tricategories.  This becomes manifest in the labelling of certain cells, which we now explain.

When we denote a geometric $2$-cell as the product of two structure cells, e.g., $\rho\times \alpha$, it is at times useful to employ tricategorical coherence.  Although, one may define these $2$-cells by various composites, tricategorical coherence assures us that the resulting $2$-cell diagrams are equivalent in the appropriate sense.

Consider, for example, the domain of the second cell in the composite of geometric $3$-cells forming the domain of the $U_{5,2}$ axiom.  The domain of the $3$-cell is the domain of a product cell of the middle mediator trimodification whiskered with various associativity tritransformation and trimodification cells.  The $2$-cell we now describe has the following form
\[ \rho\times\alpha^{-1}\maps (\rho_{A}\times 1_{B(CD)})*(1_{A1}\times \alpha_{BCD})\To (1_{A}\times \alpha_{BCD})*(\rho_{A}\times 1_{(BC)D}),\]
\noindent and it is evident that the composite is defined by a ``higher-dimensional Eckmann-Hilton argument".

The $2$-cell $\rho\times\alpha^{-1}$ is a square
\[
\def\objectstyle{\scriptstyle}
  \def\labelstyle{\scriptstyle}
   \xy
   (-15,15)*+{(xI)((yz)w)}="1";
   (15,15)*+{(xI)(y(zw))}="2";
   (-15,-15)*+{x((yz)w)}="3";
   (15,-15)*+{x(y(zw))}="4";
        {\ar^{1_{xI}*\alpha_{yzw}} "1";"2"};
        {\ar_{\rho_{x}*1_{(yz)w}} "1";"3"};
        {\ar^{\rho_{x}*1_{y(zw)}} "2";"4"};
        {\ar_{1_{x}*\alpha_{yzw}} "3";"4"};
\endxy
\]
Noticing that $\rho$ and $\alpha$ are acting, in some sense, independently of one another, we define this square as a composite of triangle $2$-cells by introducing the diagonal $1$-cell
\[ \rho_{x}*\alpha_{yzw}\maps (xI)((yz)w) \to x(y(zw)).\]
The $2$-cells comprising the triangles are coherence cells for interchange
and unit coherence cells.  The interchange coherence cells are structure cells of the strong transformation component $\chi_{\otimes}$ of the monoidal product trifunctor.  The unit coherence cells are the adjoint pairs of $1$-cells of adjoint equivalences in bicategories of $2$-functors, strong transformations, and modifications, i.e., the $1$-cells of the right and left unitor transformations of the local tricategories in which the diagrams of the $U_{5,i}$ axioms live.

The composite is
\[ \rho\times\alpha^{-1} := \chi^{-1}_{\otimes}(l^{-1}_{\rho}\otimes r^{-1}_{\alpha})(r_{\rho}\otimes l_{\alpha})\chi_{\otimes}\]
which, for computational purposes, is best drawn as a square whose interior has been sectioned into a pair of bigons and a pair of triangles.  We have
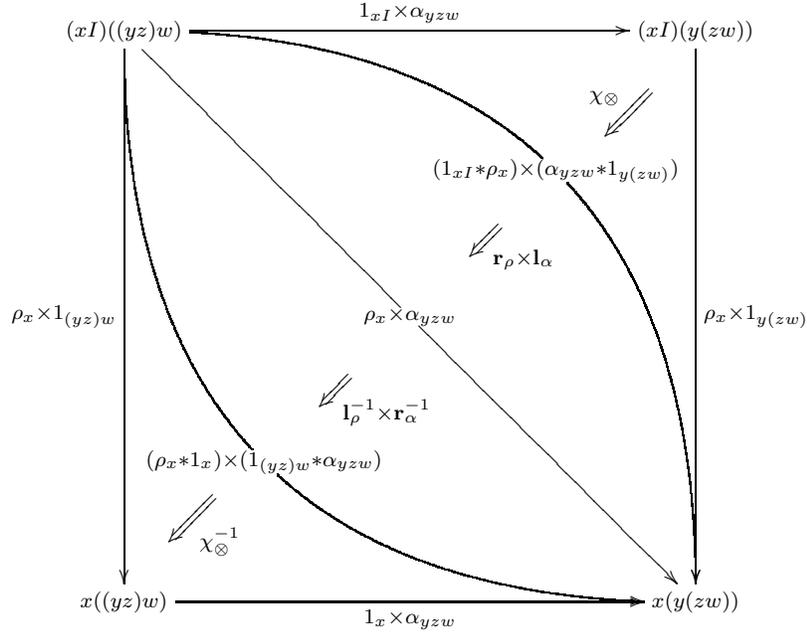
\[
\def\objectstyle{\scriptstyle}
  \def\labelstyle{\scriptstyle}
   \xy
   (-38,38)*+{(xI)((yz)w)}="1"; %DOMAIN
   (-38,-38)*+{x((yz)w)}="2";
   (38,38)*+{(xI)(y(zw))}="3"; %CODOMAIN
   (38,-38)*+{x(y(zw))}="4";
        {\ar_{\rho_{x}\times 1_{(yz)w}} "1";"2"};
        {\ar_{1_{x}\times\alpha_{yzw}} "2";"4"};
        {\ar^{1_{xI}\times \alpha_{yzw}} "1";"3"};
        {\ar^{\rho_{x}\times 1_{y(zw)}} "3";"4"};
        {\ar@/_6.5pc/|{(\rho_{x}* 1_{x})\times (1_{(yz)w}* \alpha_{yzw})} "1";"4"};
        {\ar@/^6.5pc/|{(1_{xI}*\rho_{x})\times (\alpha_{yzw}*1_{y(zw)})} "1";"4"};
        {\ar|{\rho_{x}\times \alpha_{yzw}} "1";"4"};
        {\ar@{=>}_{\chi_{\otimes}} (32,30); (26,24)};
        {\ar@{=>}^{\chi^{-1}_{\otimes}} (-26,-24); (-32,-30)};
        {\ar@{=>}^{{\bf r}_{\rho} \times {\bf l}_{\alpha}} (12,12); (8,8)};
        {\ar@{=>}^{{\bf l}^{-1}_{\rho} \times {\bf r}^{-1}_{\alpha}} (-8,-8); (-12,-12)};        
\endxy
\]
Similar interchange-type cells are used to define $\alpha^{-1}\widetilde{\times}\lambda$ and $\alpha^{-1}\widetilde{\times}\alpha$ explicitly.

\subsection{Trimble's Tetracategories}\label{tetra-definition}
The structure we work with in this paper is a monoidal tricategory.

\begin{definition}\label{monoidaltricategorydefinition}
A {\bf monoidal tricategory} is a one-object tetracategory in the sense of Trimble.
\end{definition}

\noindent We now give the definition of tetracategory following~\cite{Tr}.

\begin{definition}\label{tetracategorydefinition}
A {\bf tetracategory} $\T$ consists of
\begin{itemize}
\item a collection of objects $a,b,c,\ldots$,
\item for each pair of objects  $a,b$, a tricategory $\T(a,b)$ of $1$-, $2$-, and $3$-morphisms,
\item for each triple of objects  $a,b,c$,a trifunctor
\[\otimes\maps \T\times\T\to\T\]
\noindent called composition,
\item for each object $a$, a trifunctor
\[ I\maps 1\to \T\]
\noindent called a unit,
\item for each $4$-tuple of objects $a,b,c,d$, a biadjoint biequivalence (in the local tricategory of maps of tricategories)
\[ \alpha\maps \otimes(\otimes\times 1) \To \otimes(1\times \otimes)\]
\noindent called the associativity,
\item for each pair of objects  $a,b$, a biadjoint biequivalence (in the local tricategory of maps of tricategories)
\[ \lambda\maps \otimes(I \times 1)\To 1\]
\noindent  called the monoidal left unitor,
\item for each pair of objects  $a,b$, a biadjoint biequivalence (in the local tricategory of maps of tricategories)
\[ \rho\maps \otimes(1\times I)\To 1\]
\noindent called the monoidal right unitor,
\item for each $5$-tuple of objects $a,b,c,d,e$, an adjoint equivalence (in the local bicategory of maps of tricategories)
\[
   \xy
   (-13,0)*+{\pi\maps (1\times \alpha)\alpha(\alpha\times 1)}="2"; (12,0)*+{\alpha\alpha}="3";
        {\ar@3{->}_{} "2";"3"};
\endxy
\]
\noindent called the pentagonator,
\item for each triple of objects  $a,b,c$, an adjoint equivalence (in the local bicategory of maps of tricategories)
\[
   \xy
   (-8,0)*+{l \maps (1\times \lambda)\alpha}="2"; (9,0)*+{ \lambda}="3";
        {\ar@3{->}_{} "2";"3"};
\endxy
\]
\noindent called the left unit mediator,
\item for each triple of objects $a,b,c$, an adjoint equivalence (in the local bicategory of maps of tricategories)
\[
   \xy
   (-11,0)*+{m\maps (\rho\times 1)\alpha}="2"; (11,0)*+{\lambda\times 1}="3";
        {\ar@3{->}_{} "2";"3"};
\endxy
\]
\noindent called the middle unit mediator,
\item for each triple of objects $a,b,c$, an adjoint equivalence (in the local bicategory of maps of tricategories)
\[
   \xy
   (-8,0)*+{r\maps \rho\alpha}="2"; (7,0)*+{\rho\times 1}="3";
        {\ar@3{->}_{} "2";"3"};
\endxy
\]
\noindent called the right unit mediator,
\end{itemize}
\end{definition}
\begin{itemize}
\item for each $6$-tuple of objects $a,b,c,d,e,f$, a perturbation called a (non-abelian) $4$-cocycle (or $K_5$), consisting of, for each $5$-tuple of morphisms $A,B,C,D,E$, an invertible $4$-cell
\end{itemize}
\[
\def\objectstyle{\scriptstyle}
  \def\labelstyle{\scriptstyle}
   \xy
   (-35,32)*+{(((A\times B)\times C)\times D)\times E}="D1"; %DOMAIN
   (-60,22)*+{((A\times (B\times C))\times D)\times E}="D2";
   (-65,5)*+{(A\times ((B\times C)\times D))\times E}="D3";
   (-60,-15)*+{(A\times (B\times (C\times D)))\times E}="D4";
   (-55,-25)*+{A\times ((B\times (C\times D))\times E))}="D5";
   (-35,-34)*+{A\times (B\times ((C\times D)\times E))}="D6";
   (-9,-20)*+{A\times (B\times (C\times (D\times E)))}="D7";
   (-6,2.5)*+{(A\times B)\times (C\times (D\times E))}="D8";
   (-10,21)*+{((A\times B)\times C)\times(D\times E)}="D9";
   (-34,7.5)*+{((A\times B)\times (C\times D))\times E}="D10";
   (-31.5,-10)*+{(A\times B)\times ((C\times D)\times E)}="D11";
        {\ar_>>>{(\alpha\times 1)\times 1} "D1";"D2"}; %DOMAIN ARROWS
        {\ar_{\alpha\times 1} "D2";"D3"};
        {\ar_{(1 \times \alpha)\times 1} "D3";"D4"};
        {\ar_{\alpha} "D4";"D5"};
        {\ar_<<<<{1\times \alpha} "D5";"D6"};
        {\ar_{1\times (1\times \alpha)} "D6";"D7"};
        {\ar^{\alpha} "D1";"D9"};
        {\ar^{\alpha} "D9";"D8"};        
        {\ar^{\alpha} "D8";"D7"};       
        {\ar_{\alpha\times 1} "D1";"D10"}; 
        {\ar_{\alpha\times 1} "D10";"D4"};         
        {\ar_{\alpha} "D10";"D11"}; 
        {\ar_{1\times \alpha} "D11";"D8"}; 
        {\ar_{\alpha} "D11";"D6"};    
        {\ar@{=>}^<<{\scriptstyle \Pi} (-23,11); (-18,15)}; %DOMAIN 2-CELLS
        {\ar@{=>}^<<{\scriptstyle \alpha^{-1}_{1,1,\alpha}} (-20,-19); (-16,-15)};
        {\ar@{=>}_<<{\scriptstyle \Pi} (-45, -19); (-40, -17)};
        {\ar@{=>}^<<{\scriptstyle \Pi\times 1} (-53, 13); (-48, 10)};
  \endxy
 \] 
  \[
 \def\objectstyle{\scriptstyle}
  \def\labelstyle{\scriptstyle}
   \xy
         {\ar@{=>}_<<{\scriptstyle K_{5}} (0, 3); (0, -3)}; %PERTURBATION
   \endxy
   \]
 \[
 \def\objectstyle{\scriptstyle}
  \def\labelstyle{\scriptstyle}
   \xy
   (40,32)*+{(((A\times B)\times C)\times D)\times E}="C1"; %CODOMAIN
   (14,22)*+{((A\times (B\times C))\times D)\times E}="C2";
   (8,9.5)*+{(A\times ((B\times C)\times D))\times E}="C3";
   (8,-12)*+{(A\times (B\times (C\times D)))\times E}="C4";
   (17,-25)*+{A\times ((B\times (C\times D))\times E))}="C5";
   (41,-34)*+{A\times (B\times ((C\times D)\times E))}="C6";
   (65,-25)*+{A\times (B\times (C\times (D\times E)))}="C7";
   (70,-1)*+{(A\times B)\times (C\times (D\times E))}="C8";
   (64,22)*+{((A\times B)\times C)\times(D\times E)}="C9";
   (48,-11)*+{A\times ((B\times C)\times (D\times E))}="C10";
   (40,15)*+{(A\times (B\times C))\times (D\times E)}="C11";
   (26,-1)*+{A\times (((B\times C)\times D)\times E)}="C12";
        {\ar_>>{(\alpha\times 1)\times 1} "C1";"C2"}; %CODOMAIN ARROWS
        {\ar_{\alpha\times 1} "C2";"C3"};
        {\ar_{(1 \times \alpha)\times 1} "C3";"C4"};
        {\ar_{\alpha} "C4";"C5"};
        {\ar_{1\times\alpha} "C5";"C6"};
        {\ar_>>>{1\times (1\times \alpha)} "C6";"C7"};
        {\ar^{\alpha} "C1";"C9"};
        {\ar^{\alpha} "C9";"C8"};        
        {\ar^{\alpha} "C8";"C7"};        
        {\ar_{\alpha} "C2";"C11"}; 
        {\ar^<<<{\alpha\times 1} "C9";"C11"};         
        {\ar_{\alpha} "C11";"C10"}; 
        {\ar_{1\times \alpha} "C10";"C7"}; 
        {\ar^{\alpha} "C3";"C12"};   
        {\ar^{1\times (\alpha\times 1)} "C12";"C5"}; 
        {\ar_{1\times \alpha} "C12";"C10"};                 
        {\ar@{=>}_<<{\scriptstyle \Pi} (28, 5); (31, 10)};   %CODOMAIN 2-CELLS
        {\ar@{=>}_<<{\scriptstyle \alpha^{-1}_{\alpha,1,1}} (40,24); (45, 24)};
        {\ar@{=>}_<<{\scriptstyle \Pi} (55, 10); (59, 6)};
        {\ar@{=>}_<<{\scriptstyle 1\times \Pi} (35, -21); (39, -18)};
        {\ar@{=>}^<<{\scriptstyle \alpha^{-1}_{1,\alpha,1}} (16.5, -10); (20.5, -8)};
\endxy
\]

\begin{itemize}
\item for each $4$-tuple of objects $a,b,c,d$, a perturbation called the $U_{4,1}$ unit operation, consisting of, for each triple of morphisms $A,B,C$, an invertible $4$-cell
\end{itemize}
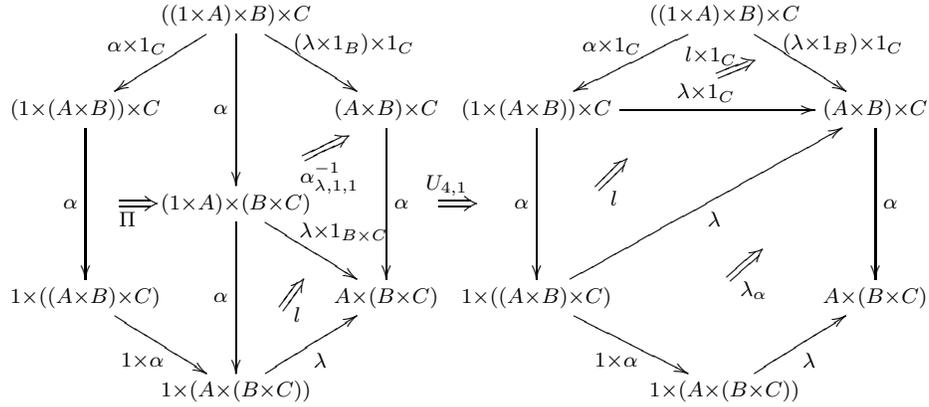
\[
\def\objectstyle{\scriptstyle}
  \def\labelstyle{\scriptstyle}
   \xy
   (-40,25)*+{((1\times A)\times B)\times C}="D1"; %DOMAIN
   (-60,12.5)*+{(1\times (A\times B))\times C}="D2";
   (-60,-12.5)*+{1\times ((A\times B)\times C)}="D3";
   (-40,0)*+{(1\times A)\times (B\times C)}="D4";
   (-40,-25)*+{1\times (A\times (B\times C))}="D5";
   (-20,-12.5)*+{A\times (B\times C)}="D6";
   (-20,12.5)*+{(A\times B)\times C}="D7";
   (25,25)*+{((1\times A)\times B)\times C}="C1"; %CODOMAIN
   (0,12.5)*+{(1\times (A\times B))\times C}="C2";
   (0,-12.5)*+{1\times ((A\times B)\times C)}="C3";
   (25,-25)*+{1\times (A\times (B\times C))}="C4";
   (45,-12.5)*+{A\times (B\times C)}="C5";
   (45,12.5)*+{(A\times B)\times C}="C6";
        {\ar_{\alpha\times 1_{C}} "D1";"D2"}; %DOMAIN ARROWS
        {\ar_{\alpha} "D2";"D3"};
        {\ar_{1 \times \alpha} "D3";"D5"};
        {\ar_{\alpha} "D1";"D4"};
        {\ar_{\alpha} "D4";"D5"};
        {\ar^{\lambda\times 1_{B\times C}} "D4";"D6"};
        {\ar^{(\lambda\times 1_{B})\times 1_{C}} "D1";"D7"};        
        {\ar^{\alpha} "D7";"D6"}; 
        {\ar_{\lambda} "D5";"D6"};         
        {\ar_{\alpha\times 1_{C}} "C1";"C2"}; %CODOMAIN ARROWS
        {\ar_{\alpha} "C2";"C3"};
        {\ar_{1 \times \alpha} "C3";"C4"};
        {\ar^{\alpha} "C6";"C5"};
        {\ar_{\lambda} "C4";"C5"};
        {\ar^{\lambda\times 1_{C}} "C2";"C6"};
        {\ar^{(\lambda\times 1_{B})\times 1_{C}} "C1";"C6"};         
        {\ar_{\lambda} "C3";"C6"};     
        {\ar@{=>}_<<{\scriptstyle l} (-34,-14); (-31,-10)}; %DOMAIN 2-CELLS
        {\ar@{=>}_<<{\scriptstyle \alpha^{-1}_{\lambda,1,1}} (-31, 6); (-25, 9)};
        {\ar@{=>}_<<{\scriptstyle \Pi} (-55.5, 0); (-50.5, 0)};
        {\ar@{=>}^<<{\scriptstyle l\times 1_{C}} (24,17); (29,19)}; %CODOMAIN 2-CELLS
        {\ar@{=>}_<<{\scriptstyle l} (8, 2); (12, 6)};
        {\ar@{=>}_<<{\scriptstyle \lambda_{\alpha}} (25.5, -10); (30, -6)};
        {\ar@{=>}^<<{\scriptstyle U_{4,1}} (-13, 0); (-8, 0)}; %PERTURBATION
\endxy
\]

\begin{itemize}
\item for each $4$-tuple of objects $a,b,c,d$, a perturbation called the $U_{4,2}$ unit operation, consisting of, for each triple of morphisms $A,B,C$, an invertible $4$-cell
\end{itemize}
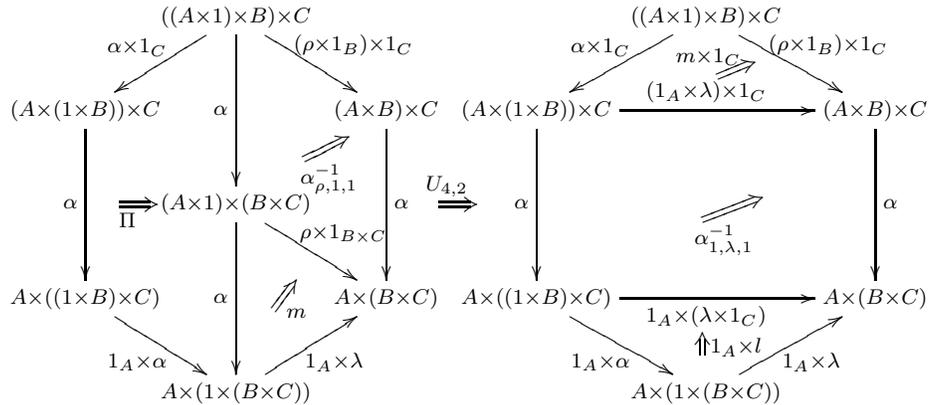
\[
\def\objectstyle{\scriptstyle}
  \def\labelstyle{\scriptstyle}
   \xy
   (-40,25)*+{((A\times 1)\times B)\times C}="D1"; %DOMAIN
   (-60,12.5)*+{(A\times (1\times B))\times C}="D2";
   (-60,-12.5)*+{A\times ((1\times B)\times C)}="D3";
   (-40,0)*+{(A\times 1)\times (B\times C)}="D4";
   (-40,-25)*+{A\times (1\times (B\times C))}="D5";
   (-20,-12.5)*+{A\times (B\times C)}="D6";
   (-20,12.5)*+{(A\times B)\times C}="D7";
   (22.5,25)*+{((A\times 1)\times B)\times C}="C1"; %CODOMAIN
   (0,12.5)*+{(A\times (1\times B))\times C}="C2";
   (0,-12.5)*+{A\times ((1\times B)\times C)}="C3";
   (22.5,-25)*+{A\times (1\times (B\times C))}="C4";
   (45,-12.5)*+{A\times (B\times C)}="C5";
   (45,12.5)*+{(A\times B)\times C}="C6";
        {\ar_{\alpha\times 1_{C}} "D1";"D2"}; %DOMAIN ARROWS
        {\ar_{\alpha} "D2";"D3"};
        {\ar_{1_{A} \times \alpha} "D3";"D5"};
        {\ar_{\alpha} "D1";"D4"};
        {\ar_{\alpha} "D4";"D5"};
        {\ar^{\alpha} "D7";"D6"};
        {\ar^{(\rho\times 1_{B})\times 1_{C}} "D1";"D7"};        
        {\ar^{\rho\times 1_{B\times C}} "D4";"D6"}; 
        {\ar_{1_{A}\times\lambda} "D5";"D6"};         
        {\ar_{\alpha\times 1_{C}} "C1";"C2"}; %CODOMAIN ARROWS
        {\ar_{\alpha} "C2";"C3"};
        {\ar_{1_{A} \times \alpha} "C3";"C4"};
        {\ar^{\alpha} "C6";"C5"};
        {\ar_{1_{A}\times \lambda} "C4";"C5"};
        {\ar^{(1_{A}\times \lambda)\times 1_{C}} "C2";"C6"};
        {\ar^{(\rho\times 1_{B})\times 1_{C}} "C1";"C6"};         
        {\ar_{1_{A}\times (\lambda\times 1_{C})} "C3";"C5"};     
        {\ar@{=>}_<<{\scriptstyle m} (-35,-14); (-32,-10)}; %DOMAIN 2-CELLS
        {\ar@{=>}_<<{\scriptstyle \alpha^{-1}_{\rho,1,1}} (-31, 6); (-25, 9)};
        {\ar@{=>}_<<{\scriptstyle \Pi} (-55.5, 0); (-50.5, 0)};
        {\ar@{=>}^<<{\scriptstyle m\times 1_{C}} (24,17); (29,19)}; %CODOMAIN 2-CELLS
        {\ar@{=>}_<<{\scriptstyle \alpha^{-1}_{1,\lambda,1}} (22, -2); (30, 1)};
        {\ar@{=>}_<<{\scriptstyle 1_{A}\times l} (22, -20); (22, -17)};
        {\ar@{=>}^<<{\scriptstyle U_{4,2}} (-13, 0); (-8, 0)}; %PERTURBATION
\endxy
\]
\begin{itemize}
\item for each $4$-tuple of objects $a,b,c,d$, a perturbation called the $U_{4,3}$ unit operation, consisting of, for each triple of morphisms $A,B,C$, an invertible $4$-cell
\end{itemize}
\[
\def\objectstyle{\scriptstyle}
  \def\labelstyle{\scriptstyle}
   \xy
   (-40,25)*+{((A\times B)\times 1)\times C}="D1"; %DOMAIN
   (-60,12.5)*+{(A\times (B\times 1))\times C}="D2";
   (-60,-12.5)*+{A\times ((B\times 1)\times C)}="D3";
   (-40,0)*+{(A\times B)\times (1\times C)}="D4";
   (-40,-25)*+{A\times (B\times (1\times C))}="D5";
   (-20,-12.5)*+{A\times (B\times C)}="D6";
   (-20,12.5)*+{(A\times B)\times C}="D7";
   (22.5,25)*+{((A\times B)\times 1)\times C}="C1"; %CODOMAIN
   (0,12.5)*+{(A\times (B\times 1))\times C}="C2";
   (0,-12.5)*+{A\times ((B\times 1)\times C)}="C3";
   (22.5,-25)*+{A\times (B\times (1\times C))}="C4";
   (45,-12.5)*+{A\times (B\times C)}="C5";
   (45,12.5)*+{(A\times B)\times C}="C6";
        {\ar_{\alpha\times 1_{C}} "D1";"D2"}; %DOMAIN ARROWS
        {\ar_{\alpha} "D2";"D3"};
        {\ar_{1_{A} \times \alpha} "D3";"D5"};
        {\ar_{\alpha} "D1";"D4"};
        {\ar_{\alpha} "D4";"D5"};
        {\ar^{\alpha} "D7";"D6"};
        {\ar^{\rho\times 1_{C}} "D1";"D7"};        
        {\ar_{1_{A\times B}\times\lambda} "D4";"D7"}; 
        {\ar_{1_{A}\times\lambda} "D5";"D6"};         
        {\ar_{\alpha\times 1_{C}} "C1";"C2"}; %CODOMAIN ARROWS
        {\ar_{\alpha} "C2";"C3"};
        {\ar_{1_{A} \times \alpha} "C3";"C4"};
        {\ar^{\alpha} "C6";"C5"};
        {\ar_{1_{A}\times (1_{B}\times\lambda)} "C4";"C5"};
        {\ar^{(1_{A}\times \rho)\times 1_{C}} "C2";"C6"};
        {\ar^{\rho\times 1_{C}} "C1";"C6"};         
        {\ar_{1_{A}\times (\rho\times 1_{C})} "C3";"C5"};     
        {\ar@{=>}^<<{\scriptstyle m} (-32,10); (-29,14)}; %DOMAIN 2-CELLS
        {\ar@{=>}^<<{\scriptstyle \alpha^{-1}_{1,1,\lambda}} (-31, -10); (-26, -4)};
        {\ar@{=>}_<<{\scriptstyle \Pi} (-55.5, 0); (-50.5, 0)};
        {\ar@{=>}^<<{\scriptstyle r\times 1_{C}} (24,17); (29,19)}; %CODOMAIN 2-CELLS
        {\ar@{=>}_<<{\scriptstyle \alpha^{-1}_{1,\rho,1}} (22, -2); (30, 1)};
        {\ar@{=>}_<<{\scriptstyle 1_{A}\times m} (22, -20); (22, -17)};
        {\ar@{=>}^<<{\scriptstyle U_{4,3}}  (-13, 0); (-8, 0)}; %PERTURBATION
\endxy
\]
\begin{itemize}
\item for each $4$-tuple of objects $a,b,c,d$, a perturbation called the $U_{4,4}$ unit operation, consisting of, for each triple of morphisms $A,B,C$, an invertible $4$-cell
\end{itemize}
\[
\def\objectstyle{\scriptstyle}
  \def\labelstyle{\scriptstyle}
   \xy
   (-40,25)*+{((A\times B)\times C)\times 1}="D1"; %DOMAIN
   (-60,12.5)*+{(A\times (B\times C))\times 1}="D2";
   (-60,-12.5)*+{A\times ((B\times C)\times 1)}="D3";
   (-40,0)*+{(A\times B)\times (C\times 1)}="D4";
   (-40,-25)*+{A\times (B\times (C\times 1))}="D5";
   (-20,-12.5)*+{A\times (B\times C)}="D6";
   (-20,12.5)*+{(A\times B)\times C}="D7";
   (22.5,25)*+{((A\times B)\times C)\times 1}="C1"; %CODOMAIN
   (0,12.5)*+{(A\times (B\times C))\times 1}="C2";
   (0,-12.5)*+{A\times ((B\times C)\times 1)}="C3";
   (22.5,-25)*+{A\times (1\times (B\times C))}="C4";
   (45,-12.5)*+{A\times (B\times C)}="C5";
   (45,12.5)*+{(A\times B)\times C}="C6";
        {\ar_{\alpha\times I} "D1";"D2"}; %DOMAIN ARROWS
        {\ar_{\alpha} "D2";"D3"};
        {\ar_{1_{A} \times \alpha} "D3";"D5"};
        {\ar_{\alpha} "D1";"D4"};
        {\ar_{\alpha} "D4";"D5"};
        {\ar^{\alpha} "D7";"D6"};
        {\ar^{\rho} "D1";"D7"};        
        {\ar_{1_{A\times B}\times\rho} "D4";"D7"}; 
        {\ar_{1_{A}\times (1_{B}\times \rho)} "D5";"D6"};         
        {\ar_{\alpha\times 1_{C}} "C1";"C2"}; %CODOMAIN ARROWS
        {\ar_{\alpha} "C2";"C3"};
        {\ar_{1_{A} \times \alpha} "C3";"C4"};
        {\ar^{\alpha} "C6";"C5"};
        {\ar_{1_{A}\times (1_{B}\times \rho)} "C4";"C5"};
        {\ar^{\rho} "C1";"C6"};
        {\ar^{\rho} "C2";"C5"};         
        {\ar_{1_{A}\times \rho} "C3";"C5"};     
        {\ar@{=>}^<<{\scriptstyle r} (-32,10); (-29,14)}; %DOMAIN 2-CELLS
        {\ar@{=>}^<<{\scriptstyle \alpha^{-1}_{1,1,\rho}} (-31, -10); (-26, -4)};
        {\ar@{=>}_<<{\scriptstyle \Pi} (-55.5, 0); (-50.5, 0)};
        {\ar@{=>}^<<{\scriptstyle \rho_{\alpha}} (24,12); (29,14)}; %CODOMAIN 2-CELLS
        {\ar@{=>}_<<{\scriptstyle \rho} (12, -7); (17, -4)};
        {\ar@{=>}_<<{\scriptstyle 1_{A}\times r} (20, -20); (20, -17)};
        {\ar@{=>}^<<{\scriptstyle U_{4,4}} (-13, 0); (-8, 0)}; %PERTURBATION
\endxy
\]
\begin{itemize}
\item all satisfying the $K_{6}$ associativity condition and the $U_{5,2}$, $U_{5,3}$, and $U_{5.4}$ unit conditions
\end{itemize}

\subsubsection*{$K_6$ Axiom}
\begin{itemize}
\item for each $7$-tuple of objects $a,b,c,d,e,f,g$, an equation called the $K_{6}$ associativity condition, consisting of, for each $6$-tuple of morphisms $A,B,C,D,E,F$, an equation of $4$-cells
\end{itemize}

\[
\def\objectstyle{\scriptstyle}
  \def\labelstyle{\scriptstyle}
   \xy
   (-22.5,20)*+{((((AB)C)D)E)F}="O1"; %BORDER VERTICES
   (10,20)*+{(((AB)C)D)(EF)}="O2";
   (45,0)*+{(AB)(C(D(EF)))}="O3";
    (-88,0)*+{((A(B(CD))E)F}="O4"; 
   (-88,-20)*+{(A((B(CD))E))F}="O5";   
   (-88, 20)*+{((A((BC)D))E)F}="O6";   
   (-55, 20)*+{(((A(BC))D)E)F}="O7";   
   (45,20)*+{((AB)C)(D(EF))}="O8";
   (45,-20)*+{A(B(C(D(EF))))}="O9"; 
   (-88 ,-40)*+{(A(B((CD)E)))F}="O10";   
   (-55,-40)*+{(A(B(C(DE))))F}="O11";   
   (10,-40)*+{A(B((C(DE))F))}="O12"; 
   (45,-40)*+{A(B(C((DE)F)))}="O13"; 
   (-22.5,-40)*+{A((B(C(DE)))F)}="O14";
   (-49,-10)*+{((AB)((CD)E))F}="I1"; %INNER VERTICES
   (0,-26)*+{(AB)((C(DE))F)}="I2";
   (-4,-5)*+{(((AB)C)(DE))F}="I3";
   (-35,5)*+{(((AB)(CD))E)F}="I4";
   (20,6)*+{((AB)C)((DE)F)}="I5";
   (-35,-25)*+{((AB)(C(DE)))F}="I6";
   (25,-10)*+{(AB)(C((DE)F))}="I7"; 
        {\ar@{-->}_{} "O5";"O10"}; % BORDER ARROWS
        {\ar@{-->}_{} "O10";"O11"};
        {\ar@{-->}_{} "O11";"O14"};
        {\ar@{-->}_{} "O13";"O9"};
        {\ar@{-->}_{} "O14";"O12"};
        {\ar@{-->}_{} "O12";"O13"};
        {\ar_{} "O1";"O2"};    
        {\ar@{-->}_{} "O1";"O7"};    
        {\ar@{-->}_{} "O7";"O6"};        
        {\ar@{-->}_{} "O6";"O4"};    
        {\ar@{-->}_{} "O4";"O5"};    
        {\ar^{} "O8";"O3"};
        {\ar^{} "O2";"O8"};                
        {\ar@{-->}_{} "O3";"O9"};       
        {\ar_{} "O1";"I3"}; % IN/OUT ARROWS
        {\ar_{} "O1";"I4"};       
        {\ar@{-->}_{} "I4";"O4"}; 
        {\ar_{} "I5";"O8"};
        {\ar@{-->}_{} "I1";"O10"}; 
        {\ar@{-->}_{} "I6";"O11"};
        {\ar@{-->}_{} "I2";"O12"};
        {\ar_{} "I7";"O3"};
        {\ar@{-->}_{} "I7";"O13"};
        {\ar_{} "I4";"I1"};  % INTERIOR ARROWS
        {\ar_{} "I1";"I6"};     
        {\ar_{} "I3";"I6"};     
        {\ar_{} "I3";"I5"};    
        {\ar_{} "I5";"I7"};
        {\ar_{} "I2";"I7"};
        {\ar_{} "I6";"I2"};
        {\ar@{=>}^{\scriptstyle \Pi} (7, 9); (10, 13)};  %DOMAIN 2-CELLS
        {\ar@{=>}_{\scriptstyle \alpha_{1,1,\alpha}^{-1}} (32, 4); (35, 8)};        
        {\ar@{=>}^{\scriptstyle \Pi} (5, -16); (10, -10)};
        {\ar@{=>}^{\scriptstyle \Pi\times 1} (-27, -8); (-23, -8)};
        {\ar@{=>}^{\scriptstyle (\Pi\times 1)\times 1} (-60, 14); (-56, 10)};
        {\ar@{=>}^{\scriptstyle \Pi\times 1} (-73, -7); (-68, -10)};
        {\ar@{=>}^{\scriptstyle \alpha_{1,1,\alpha}^{-1}\times 1} (-53, -27); (-49, -24)};
        {\ar@{=>}^{\scriptstyle \Pi} (-16, -35); (-11,-32)};
        {\ar@{=>}_{\scriptstyle \alpha^{-1}_{1,1,\alpha}} (23, -33); (27, -28)};
        {\ar@{=>}^{\scriptstyle \alpha^{-1}_{1,1,1\times\alpha}} (40, -17); (42, -13)};
\endxy
\]

\[
\def\objectstyle{\scriptstyle}
  \def\labelstyle{\scriptstyle}
   \xy  
    {\ar@{=>}_{\scriptstyle K_{5}} (-6, 2); (-6, -3)}; 
   \endxy
\]

\[
\def\objectstyle{\scriptstyle}
  \def\labelstyle{\scriptstyle}
   \xy
   (-22.5,20)*+{((((AB)C)D)E)F}="O1"; %BORDER VERTICES
   (10,20)*+{(((AB)C)D)(EF)}="O2";
   (45,0)*+{(AB)(C(D(EF)))}="O3";
    (-88,0)*+{((A(B(CD))E)F}="O4"; 
   (-88,-20)*+{(A((B(CD))E))F}="O5";   
   (-88, 20)*+{((A((BC)D))E)F}="O6";   
   (-55, 20)*+{(((A(BC))D)E)F}="O7";   
   (45,20)*+{((AB)C)(D(EF))}="O8";
   (45,-20)*+{A(B(C(D(EF))))}="O9"; 
   (-88 ,-40)*+{(A(B((CD)E)))F}="O10";   
   (-55,-40)*+{(A(B(C(DE))))F}="O11";   
   (10,-40)*+{A(B((C(DE))F))}="O12"; 
   (45,-40)*+{A(B(C((DE)F)))}="O13"; 
   (-22.5,-40)*+{A((B(C(DE)))F)}="O14";
   (-49,-10)*+{((AB)((CD)E))F}="I1"; %INNER VERTICES
   (0,-26)*+{(AB)((C(DE))F)}="I2";
   (-4,0)*+{((AB)(CD))(EF)}="I3";
   (-35,5)*+{(((AB)(CD))E)F}="I4";
   (25,7)*+{(AB)((CD)(EF))}="I5";
   (-35,-25)*+{((AB)(C(DE)))F}="I6";
   (25,-10)*+{(AB)(C((DE)F))}="I7";
   (-4,-15)*+{(AB)(((CD)E)F)}="I8";
        {\ar@{-->}_{} "O5";"O10"}; % BORDER ARROWS
        {\ar@{-->}_{} "O10";"O11"};
        {\ar@{-->}_{} "O11";"O14"};
        {\ar_{} "O13";"O9"};
        {\ar@{-->}_{} "O14";"O12"};
        {\ar_{} "O12";"O13"};
        {\ar@{-->}_{} "O1";"O2"};    
        {\ar@{-->}_{} "O1";"O7"};    
        {\ar@{-->}_{} "O7";"O6"};        
        {\ar@{-->}_{} "O6";"O4"};    
        {\ar@{-->}_{} "O4";"O5"};    
        {\ar@{-->}^{} "O8";"O3"};
        {\ar@{-->}^{} "O2";"O8"};                
        {\ar_{} "O3";"O9"};       
        {\ar@{-->}_{} "O2";"I3"};        % IN/OUT ARROWS
        {\ar@{-->}_{} "O1";"I4"}; 
        {\ar@{-->}_{} "I4";"O4"}; 
        {\ar@{-->}_{} "I1";"O10"}; 
        {\ar_{} "I5";"O3"}; 
        {\ar@{-->}_{} "I6";"O11"};
        {\ar_{} "I2";"O12"};
        {\ar_{} "I7";"O3"};
        {\ar_{} "I7";"O13"};
        {\ar@{-->}_{} "I4";"I1"};  % INTERIOR ARROWS
        {\ar@{-->}_{} "I1";"I6"};     
        {\ar@{-->}_{} "I1";"I8"};     
        {\ar@{-->}_{} "I3";"I5"};    
        {\ar_{} "I8";"I5"};   
        {\ar_{} "I2";"I7"};
        {\ar@{-->}_{} "I6";"I2"};
        {\ar_{} "I8";"I2"};
        {\ar@{-->}_{} "I4";"I3"};
        {\ar@{=>}^{\scriptstyle {\alpha_{\alpha,1,1}}} (-15, 7); (-11, 10)};  %DOMAIN 2-CELLS
        {\ar@{=>}^{\scriptstyle \Pi} (30,10); (34, 13)};        
        {\ar@{=>}^{\scriptstyle (1\times 1)\times \Pi} (22, -8); (25, -4)};
        {\ar@{=>}^{\scriptstyle \Pi} (-22, -8); (-18, -5)};
        {\ar@{=>}^{\scriptstyle  \alpha_{1,1,\alpha\times 1}} (-25, -21); (-20, -18)};
        {\ar@{=>}^{\scriptstyle (\Pi\times 1)\times 1} (-60, 14); (-56, 10)};
        {\ar@{=>}^{\scriptstyle \Pi\times 1} (-73, -7); (-68, -10)};
        {\ar@{=>}^{\scriptstyle \alpha_{1,1,\alpha}^{-1}\times 1} (-53, -27); (-49, -24)};
        {\ar@{=>}^{\scriptstyle \Pi} (-16, -35); (-11 ,-32)};
        {\ar@{=>}^{\scriptstyle \alpha^{-1}_{1,1,\alpha}} (23, -33); (27, -28)};
        {\ar@{=>}^{\scriptstyle \alpha^{-1}_{1,1,1\times\alpha}} (40, -17); (42, -13)};
\endxy
\]
\[
\def\objectstyle{\scriptstyle}
  \def\labelstyle{\scriptstyle}
   \xy  
    {\ar@{=>}_{\scriptstyle \alpha} (-6, 2); (-6, -3)};
   \endxy
\]

\[
\def\objectstyle{\scriptstyle}
  \def\labelstyle{\scriptstyle}
   \xy
   (-22.5,20)*+{((((AB)C)D)E)F}="O1"; %BORDER VERTICES
   (10,20)*+{(((AB)C)D)(EF)}="O2";
   (45,0)*+{(AB)(C(D(EF)))}="O3";
    (-88,0)*+{((A(B(CD))E)F}="O4"; 
   (-88,-20)*+{(A((B(CD))E))F}="O5";   
   (-88, 20)*+{((A((BC)D))E)F}="O6";   
   (-55, 20)*+{(((A(BC))D)E)F}="O7";   
   (45,20)*+{((AB)C)(D(EF))}="O8";
   (45,-20)*+{A(B(C(D(EF))))}="O9"; 
   (-88 ,-40)*+{(A(B((CD)E)))F}="O10";   
   (-55,-40)*+{(A(B(C(DE))))F}="O11";   
   (10,-40)*+{A(B((C(DE))F))}="O12"; 
   (45,-40)*+{A(B(C((DE)F)))}="O13"; 
   (-22.5,-40)*+{A((B(C(DE)))F)}="O14";
   (-49,-10)*+{((AB)((CD)E))F}="I1"; %INNER VERTICES
   (0,-26)*+{(AB)((C(DE))F)}="I2";
   (-4,0)*+{((AB)(CD))(EF)}="I3";
   (-35,5)*+{(((AB)(CD))E)F}="I4";
   (25,7)*+{(AB)((CD)(EF))}="I5";
   (-35,-25)*+{((AB)(C(DE)))F}="I6";
   (21,-21)*+{A(B(((CD)E)F))}="I7";
   (-4,-15)*+{(AB)(((CD)E)F)}="I8";
   (24,-12)*+{A(B((CD)(EF)))}="I9";
        {\ar@{-->}_{} "O5";"O10"}; % BORDER ARROWS
        {\ar_{} "O10";"O11"};
        {\ar_{} "O11";"O14"};
        {\ar@{-->}_{} "O13";"O9"};
        {\ar_{} "O14";"O12"};
        {\ar@{-->}_{} "O12";"O13"};
        {\ar@{-->}_{} "O1";"O2"};    
        {\ar@{-->}_{} "O1";"O7"};    
        {\ar@{-->}_{} "O7";"O6"};        
        {\ar@{-->}_{} "O6";"O4"};    
        {\ar@{-->}_{} "O4";"O5"};    
        {\ar@{-->}^{} "O8";"O3"};
        {\ar@{-->}^{} "O2";"O8"};                
        {\ar@{-->}_{} "O3";"O9"};       
        {\ar@{-->}_{} "I9";"O9"}; % IN/OUT ARROWS
        {\ar@{-->}_{} "O2";"I3"};       
        {\ar@{-->}_{} "O1";"I4"}; 
        {\ar@{-->}_{} "I4";"O4"}; 
        {\ar_{} "I1";"O10"}; 
        {\ar@{-->}_{} "I5";"O3"}; 
        {\ar_{} "I6";"O11"};
        {\ar_{} "I2";"O12"};
        {\ar_{} "I7";"O12"};
        {\ar@{-->}_{} "I4";"I1"};  % INTERIOR ARROWS
        {\ar_{} "I1";"I6"};     
        {\ar_{} "I1";"I8"};     
        {\ar@{-->}_{} "I3";"I5"};    
        {\ar@{-->}_{} "I8";"I5"};   
        {\ar_{} "I6";"I2"};
        {\ar_{} "I8";"I2"};
        {\ar@{-->}_{} "I4";"I3"};
        {\ar@{-->}_{} "I5";"I9"};
        {\ar@{-->}_{} "I7";"I9"};
        {\ar_{} "I8";"I7"};
        {\ar@{=>}^{\scriptstyle {\alpha_{\alpha,1,1}}} (-15, 7); (-11, 10)};  %DOMAIN 2-CELLS
        {\ar@{=>}^{\scriptstyle \Pi} (30,10); (34, 13)};          
        {\ar@{=>}^{\scriptstyle \alpha^{-1}_{1,1,\alpha}} (17, -11); (20, -8)};
        {\ar@{=>}^{\scriptstyle \Pi} (-22, -8); (-18, -5)};
        {\ar@{=>}^{\scriptstyle  \alpha_{1,1,\alpha\times 1}} (-25, -19); (-20, -16)};
        {\ar@{=>}^{\scriptstyle (\Pi\times 1)\times 1} (-60, 14); (-56, 10)};
        {\ar@{=>}^{\scriptstyle \Pi\times 1} (-68, -7); (-63, -10)};
        {\ar@{=>}^{\scriptstyle \alpha_{1,1,\alpha}^{-1}\times 1} (-53, -27); (-49, -24)};
        {\ar@{=>}^{\scriptstyle \Pi} (-16, -35); (-11,-32)};
        {\ar@{=>}^{\scriptstyle 1\times (1\times \Pi)} (33, -33); (37, -30)};
        {\ar@{=>}^{\scriptstyle \alpha_{1,1,\alpha}} (15, -26); (17, -23)};
        {\ar@{=>}^{\scriptstyle \alpha^{-1}_{1,1,\alpha}} (36, -9); (40, -6)};
\endxy
\]
\[
\def\objectstyle{\scriptstyle}
  \def\labelstyle{\scriptstyle}
   \xy  
    {\ar@{=>}_{\scriptstyle \Pi} (-6, 2); (-6, -3)};
   \endxy
\]

\[
\def\objectstyle{\scriptstyle}
  \def\labelstyle{\scriptstyle}
   \xy
   (-22.5,20)*+{((((AB)C)D)E)F}="O1"; %BORDER VERTICES
   (10,20)*+{(((AB)C)D)(EF)}="O2";
   (45,0)*+{(AB)(C(D(EF)))}="O3";
    (-88,0)*+{((A(B(CD))E)F}="O4"; 
   (-88,-20)*+{(A((B(CD))E))F}="O5";   
   (-88, 20)*+{((A((BC)D))E)F}="O6";   
   (-55, 20)*+{(((A(BC))D)E)F}="O7";   
   (45,20)*+{((AB)C)(D(EF))}="O8";
   (45,-20)*+{A(B(C(D(EF))))}="O9"; 
   (-88 ,-40)*+{(A(B((CD)E)))F}="O10";   
   (-55,-40)*+{(A(B(C(DE))))F}="O11";   
   (10,-40)*+{A(B((C(DE))F))}="O12"; 
   (45,-40)*+{A(B(C((DE)F)))}="O13"; 
   (-22.5,-40)*+{A((B(C(DE)))F)}="O14";
   (-49,-10)*+{((AB)((CD)E))F}="I1"; %INNER VERTICES
   (-10,-25)*+{A((B((CD)E))F)}="I2";
   (-4,0)*+{((AB)(CD))(EF)}="I3";
   (-35,5)*+{(((AB)(CD))E)F}="I4";
   (25,7)*+{(AB)((CD)(EF))}="I5";
   (21,-21)*+{A(B(((CD)E)F))}="I7";
   (-4,-15)*+{(AB)(((CD)E)F)}="I8";
   (24,-12)*+{A(B((CD)(EF)))}="I9";
        {\ar_{} "O5";"O10"}; % BORDER ARROWS
        {\ar@{-->}_{} "O10";"O11"};
        {\ar@{-->}_{} "O11";"O14"};
        {\ar@{-->}_{} "O13";"O9"};
        {\ar@{-->}_{} "O14";"O12"};
        {\ar@{-->}_{} "O12";"O13"};
        {\ar@{-->}_{} "O1";"O2"};    
        {\ar@{-->}_{} "O1";"O7"};    
        {\ar@{-->}_{} "O7";"O6"};        
        {\ar@{-->}_{} "O6";"O4"};    
        {\ar_{} "O4";"O5"};    
        {\ar@{-->}^{} "O8";"O3"};
        {\ar@{-->}^{} "O2";"O8"};                
        {\ar@{-->}_{} "O3";"O9"};       
        {\ar@{-->}_{} "I9";"O9"}; % IN/OUT ARROWS
        {\ar@{-->}_{} "O2";"I3"};       
        {\ar@{-->}_{} "O1";"I4"}; 
        {\ar_{} "I4";"O4"}; 
        {\ar_{} "I1";"O10"}; 
        {\ar@{-->}_{} "I5";"O3"}; 
        {\ar@{-->}_{} "I2";"O14"};
        {\ar@{-->}_{} "I7";"O12"};
        {\ar_{} "O10";"I2"};
        {\ar_{} "I4";"I1"};  % INTERIOR ARROWS
        {\ar_{} "I1";"I8"};     
        {\ar_{} "I3";"I5"};    
        {\ar_{} "I8";"I5"};   
        {\ar_{} "I2";"I7"};
        {\ar_{} "I4";"I3"};
        {\ar_{} "I5";"I9"};
        {\ar_{} "I7";"I9"};
        {\ar_{} "I8";"I7"};
        {\ar@{=>}^{\scriptstyle {\alpha_{\alpha,1,1}}} (-15, 7); (-11, 10)};  %DOMAIN 2-CELLS
        {\ar@{=>}^{\scriptstyle \Pi} (30,10); (34, 13)};       
        {\ar@{=>}^{\scriptstyle \alpha^{-1}_{1,1,\alpha}} (17, -11); (20, -8)};
        {\ar@{=>}^{\scriptstyle \Pi} (-22, -8); (-18, -5)};
        {\ar@{=>}^{\scriptstyle (\Pi\times 1)\times 1} (-60, 14); (-56, 10)};
        {\ar@{=>}^{\scriptstyle \Pi\times 1} (-68, -7); (-63, -10)};
        {\ar@{=>}^{\scriptstyle \Pi} (-43, -22); (-39, -18)};
        {\ar@{=>}_{\scriptstyle \alpha_{1,1\times\alpha,1}} (-35, -37); (-29,-33)};
        {\ar@{=>}^{\scriptstyle 1\times (1\times \Pi)} (33, -33); (37, -30)};
        {\ar@{=>}^{\scriptstyle 1\times \alpha_{1,\alpha,1}} (3,-33); (6, -30)};
        {\ar@{=>}^{\scriptstyle \alpha^{-1}_{1,1,\alpha}} (36, -9); (40, -6)};
\endxy
\]

\[
\def\objectstyle{\scriptstyle}
  \def\labelstyle{\scriptstyle}
   \xy  
    {\ar@{=>}_{\scriptstyle K_{5}} (-6, 2); (-6, -3)};
   \endxy
\]
\[
\def\objectstyle{\scriptstyle}
  \def\labelstyle{\scriptstyle}
   \xy
   (-22.5,20)*+{((((AB)C)D)E)F}="O1"; %BORDER VERTICES
   (10,20)*+{(((AB)C)D)(EF)}="O2";
   (45,0)*+{(AB)(C(D(EF)))}="O3";
    (-88,0)*+{((A(B(CD))E)F}="O4"; 
   (-88,-20)*+{(A((B(CD))E))F}="O5";   
   (-88, 20)*+{((A((BC)D))E)F}="O6";   
   (-55, 20)*+{(((A(BC))D)E)F}="O7";   
   (45,20)*+{((AB)C)(D(EF))}="O8";
   (45,-20)*+{A(B(C(D(EF))))}="O9"; 
   (-88 ,-40)*+{(A(B((CD)E)))F}="O10";   
   (-55,-40)*+{(A(B(C(DE))))F}="O11";   
   (10,-40)*+{A(B((C(DE))F))}="O12"; 
   (45,-40)*+{A(B(C((DE)F)))}="O13"; 
   (-22.5,-40)*+{A((B(C(DE)))F)}="O14";
   (-10,-25)*+{A((B((CD)E))F)}="I2"; %INNER VERTICES
   (-4,0)*+{((AB)(CD))(EF)}="I3";
   (-35,5)*+{(((AB)(CD))E)F}="I4";
   (25,7)*+{(AB)((CD)(EF))}="I5";
   (21,-21)*+{A(B(((CD)E)F))}="I7";
   (-39,-12)*+{(A(B(CD)))(EF)}="I8";
   (24,-12)*+{A(B((CD)(EF)))}="I9";
   (-45,-25)*+{A(((B(CD)E))F)}="I10";
   (-6,-12)*+{A((B(CD))(EF))}="I11";
        {\ar@{-->}_{} "O5";"O10"}; % BORDER ARROWS
        {\ar@{-->}_{} "O10";"O11"};
        {\ar@{-->}_{} "O11";"O14"};
        {\ar@{-->}_{} "O13";"O9"};
        {\ar@{-->}_{} "O14";"O12"};
        {\ar@{-->}_{} "O12";"O13"};
        {\ar_{} "O1";"O2"};    
        {\ar_{} "O1";"O7"};    
        {\ar_{} "O7";"O6"};        
        {\ar_{} "O6";"O4"};    
        {\ar@{-->}_{} "O4";"O5"};    
        {\ar@{-->}^{} "O8";"O3"};
        {\ar@{-->}^{} "O2";"O8"};                
        {\ar@{-->}_{} "O3";"O9"};       
        {\ar@{-->}_{} "I9";"O9"}; % IN/OUT ARROWS
        {\ar_{} "O2";"I3"};       
        {\ar_{} "O1";"I4"}; 
        {\ar_{} "I4";"O4"}; 
       {\ar_{} "O4";"I8"};  
        {\ar@{-->}_{} "I5";"O3"}; 
        {\ar@{-->}_{} "I2";"O14"};
        {\ar@{-->}_{} "I7";"O12"};
        {\ar@{-->}_{} "O10";"I2"};
        {\ar@{-->}_{} "O5";"I10"}; 
        {\ar@{-->}_{} "I3";"I5"};     % INTERIOR ARROWS    
        {\ar_{} "I3";"I8"}; 
        {\ar@{-->}_{} "I2";"I7"};
        {\ar_{} "I4";"I3"};
        {\ar@{-->}_{} "I5";"I9"};
        {\ar@{-->}_{} "I7";"I9"};
        {\ar@{-->}_{} "I10";"I2"};
        {\ar@{-->}_{} "I8";"I11"};
        {\ar@{-->}_{} "I11";"I9"};
        {\ar@{-->}_{} "I10";"I11"};
        {\ar@{=>}^{\scriptstyle {\alpha_{\alpha,1,1}}} (-15, 7); (-11, 10)};  %DOMAIN 2-CELLS
        {\ar@{=>}^{\scriptstyle \Pi} (30,10); (34, 13)};        
        {\ar@{=>}^{\scriptstyle \Pi} (10, -7); (14, -3)};
        {\ar@{=>}^{\scriptstyle 1\times \Pi} (-5, -20); (1, -17)};
        {\ar@{=>}^{\scriptstyle \alpha_{\alpha,1,1}} (-40, -4); (-36, -4)};
        {\ar@{=>}^{\scriptstyle (\Pi\times 1)\times 1} (-60, 14); (-56, 10)};
        {\ar@{=>}^{\scriptstyle \Pi} (-54, -20); (-50, -16)};
        {\ar@{=>}_{\scriptstyle \alpha_{1,1\times\alpha,1}} (-35, -37); (-29,-33)};
        {\ar@{=>}^{\scriptstyle \alpha_{1,\alpha,1}} (-67, -33); (-61,-29)};
        {\ar@{=>}^{\scriptstyle 1\times (1\times \Pi)} (33, -33); (37, -30)};
        {\ar@{=>}^{\scriptstyle 1\times \alpha_{1,\alpha,1}} (3,-33); (6, -30)};
        {\ar@{=>}^{\scriptstyle \alpha^{-1}_{1,1,\alpha}} (36, -9); (40, -6)};
\endxy
\]
\[
\def\objectstyle{\scriptstyle}
  \def\labelstyle{\scriptstyle}
   \xy  
    {\ar@{=>}_{\scriptstyle \alpha} (-6, 2); (-6, -3)};
   \endxy
\]
\[
\def\objectstyle{\scriptstyle}
  \def\labelstyle{\scriptstyle}
   \xy
   (-22.5,20)*+{((((AB)C)D)E)F}="O1"; %BORDER VERTICES
   (10,20)*+{(((AB)C)D)(EF)}="O2";
   (45,0)*+{(AB)(C(D(EF)))}="O3";
    (-88,0)*+{((A(B(CD))E)F}="O4"; 
   (-88,-20)*+{(A((B(CD))E))F}="O5";   
   (-88, 20)*+{((A((BC)D))E)F}="O6";   
   (-55, 20)*+{(((A(BC))D)E)F}="O7";   
   (45,20)*+{((AB)C)(D(EF))}="O8";
   (45,-20)*+{A(B(C(D(EF))))}="O9"; 
   (-88 ,-40)*+{(A(B((CD)E)))F}="O10";   
   (-55,-40)*+{(A(B(C(DE))))F}="O11";   
   (10,-40)*+{A(B((C(DE))F))}="O12"; 
   (45,-40)*+{A(B(C((DE)F)))}="O13"; 
   (-22.5,-40)*+{A((B(C(DE)))F)}="O14";
   (-10,-25)*+{A((B((CD)E))F)}="I2"; %INNER VERTICES
   (-4,0)*+{((AB)(CD))(EF)}="I3";
   (-30,5)*+{((A(BC))D)(EF)}="I4";
   (25,7)*+{(AB)((CD)(EF))}="I5";
   (21,-21)*+{A(B(((CD)E)F))}="I7";
   (-39,-12)*+{(A(B(CD)))(EF)}="I8";
   (24,-12)*+{A(B((CD)(EF)))}="I9";
   (-45,-25)*+{A(((B(CD)E))F)}="I10";
   (-6,-12)*+{A((B(CD))(EF))}="I12";
   (-60,5)*+{(A((BC)D))(EF)}="I13";
        {\ar@{-->}_{} "O5";"O10"}; % BORDER ARROWS
        {\ar@{-->}_{} "O10";"O11"};
        {\ar@{-->}_{} "O11";"O14"};
        {\ar@{-->}_{} "O13";"O9"};
        {\ar@{-->}_{} "O14";"O12"};
        {\ar@{-->}_{} "O12";"O13"};
        {\ar@{-->}_{} "O1";"O2"};    
        {\ar@{-->}_{} "O1";"O7"};    
        {\ar@{-->}_{} "O7";"O6"};        
        {\ar@{-->}_{} "O6";"O4"};    
        {\ar@{-->}_{} "O4";"O5"};    
        {\ar^{} "O8";"O3"};
        {\ar^{} "O2";"O8"};                
        {\ar_{} "O3";"O9"};       
        {\ar_{} "I9";"O9"}; % IN/OUT ARROWS
        {\ar_{} "O2";"I4"};       
        {\ar@{-->}_{} "O7";"I4"}; 
        {\ar_{} "I5";"O3"}; 
        {\ar@{-->}_{} "I2";"O14"};
        {\ar@{-->}_{} "I7";"O12"};
        {\ar@{-->}_{} "O10";"I2"};
        {\ar@{-->}_{} "O5";"I10"}; 
        {\ar@{-->}_{} "O6";"I13"}; 
        {\ar@{-->}_{} "O4";"I8"}; 
        {\ar_{} "I3";"I5"};     % INTERIOR ARROWS  
        {\ar_{} "I4";"I13"};   
        {\ar_{} "I3";"I8"}; 
        {\ar@{-->}_{} "I2";"I7"};
        {\ar_{} "O2";"I3"};
        {\ar_{} "I5";"I9"};
        {\ar@{-->}_{} "I7";"I9"};
        {\ar@{-->}_{} "I10";"I2"};
        {\ar@{-->}_{} "I10";"I12"};
       {\ar_{} "I13";"I8"};  
        {\ar_{} "I8";"I12"};
        {\ar_{} "I12";"I9"};
        {\ar@{=>}^{\scriptstyle {\alpha_{\alpha\times 1,1,1}}} (-25, 13); (-21, 13)};  %DOMAIN 2-CELLS
        {\ar@{=>}^{\scriptstyle \Pi} (30,10); (34, 13)};         
        {\ar@{=>}^{\scriptstyle \Pi} (10, -7); (14, -3)};
        {\ar@{=>}^{\scriptstyle 1\times \Pi} (-5, -20); (1, -17)};
        {\ar@{=>}^{\scriptstyle \Pi\times (1\times 1)} (-32, -2); (-27, -6)};
        {\ar@{=>}^{\scriptstyle \alpha_{1\times\alpha,1,1}} (-63, -4); (-59, -2)};
        {\ar@{=>}^{\scriptstyle \alpha_{\alpha,1,1}} (-60, 14); (-56, 11)};
        {\ar@{=>}^{\scriptstyle \Pi} (-54, -20); (-50, -16)};
        {\ar@{=>}_{\scriptstyle \alpha_{1,1\times\alpha,1}} (-35, -37); (-29,-33)};
        {\ar@{=>}^{\scriptstyle \alpha_{1,\alpha,1}} (-67, -33); (-61,-29)};
        {\ar@{=>}^{\scriptstyle 1\times (1\times \Pi)} (33, -33); (37, -30)};
        {\ar@{=>}^{\scriptstyle 1\times \alpha_{1,\alpha,1}} (3,-33); (6, -30)};
        {\ar@{=>}^{\scriptstyle \alpha^{-1}_{1,1,\alpha}} (36, -9); (40, -6)};
\endxy
\]
\[
\def\objectstyle{\scriptstyle}
  \def\labelstyle{\scriptstyle}
   \xy  
    {\ar@{=>}_{\scriptstyle K_{5}} (-6, 2); (-6, -3)};
   \endxy
\]
\[
\def\objectstyle{\scriptstyle}
  \def\labelstyle{\scriptstyle}
   \xy
   (-22.5,20)*+{((((AB)C)D)E)F}="O1"; %BORDER VERTICES
   (10,20)*+{(((AB)C)D)(EF)}="O2";
   (45,0)*+{(AB)(C(D(EF)))}="O3";
    (-88,0)*+{((A(B(CD))E)F}="O4"; 
   (-88,-20)*+{(A((B(CD))E))F}="O5";   
   (-88, 20)*+{((A((BC)D))E)F}="O6";   
   (-55, 20)*+{(((A(BC))D)E)F}="O7";   
   (45,20)*+{((AB)C)(D(EF))}="O8";
   (45,-20)*+{A(B(C(D(EF))))}="O9"; 
   (-88 ,-40)*+{(A(B((CD)E)))F}="O10";   
   (-55,-40)*+{(A(B(C(DE))))F}="O11";   
   (10,-40)*+{A(B((C(DE))F))}="O12"; 
   (45,-40)*+{A(B(C((DE)F)))}="O13"; 
   (-22.5,-40)*+{A((B(C(DE)))F)}="O14";
   (-10,-25)*+{A((B((CD)E))F)}="I2"; %INNER VERTICES
   (-30,5)*+{((A(BC))D)(EF)}="I4";
   (21,-21)*+{A(B(((CD)E)F))}="I7";
   (-39,-12)*+{(A(B(CD)))(EF)}="I8";
   (24,-12)*+{A(B((CD)(EF)))}="I9";
   (-45,-25)*+{A(((B(CD))E)F)}="I10";
   (-6,-12)*+{A((B(CD))(EF))}="I12";
   (-60,5)*+{(A((BC)D))(EF)}="I13";
   (20,10)*+{(A(BC))(D(EF)))}="I14";
  %  *+{((AB)(CD))(EF)}="I15";
   (-14,-2)*+{A(((BC)D)(EF))}="I16";
   (20,2)*+{A((BC)(D(EF)))}="I17";
        {\ar@{-->}_{} "O5";"O10"}; % BORDER ARROWS
        {\ar@{-->}_{} "O10";"O11"};
        {\ar@{-->}_{} "O11";"O14"};
        {\ar@{-->}_{} "O13";"O9"};
        {\ar@{-->}_{} "O14";"O12"};
        {\ar@{-->}_{} "O12";"O13"};
        {\ar@{-->}_{} "O1";"O2"};    
        {\ar@{-->}_{} "O1";"O7"};    
        {\ar@{-->}_{} "O7";"O6"};        
        {\ar@{-->}_{} "O6";"O4"};    
        {\ar@{-->}_{} "O4";"O5"};    
        {\ar@{-->}^{} "O8";"O3"};
        {\ar@{-->}^{} "O2";"O8"};                
        {\ar@{-->}_{} "O3";"O9"};       
        {\ar@{-->}_{} "I9";"O9"}; % IN/OUT ARROWS
        {\ar@{-->}_{} "O2";"I4"};       
        {\ar@{-->}_{} "O7";"I4"}; 
        {\ar@{-->}_{} "I2";"O14"};
        {\ar@{-->}_{} "I7";"O12"};
        {\ar@{-->}_{} "O10";"I2"};
        {\ar@{-->}_{} "O5";"I10"}; 
        {\ar@{-->}_{} "O6";"I13"}; 
        {\ar@{-->}_{} "O4";"I8"}; 
        {\ar@{-->}_{} "O8";"I14"};   
        {\ar@{-->}_{} "I4";"I13"};     % INTERIOR ARROWS  
        {\ar@{-->}_{} "I2";"I7"};
        {\ar@{-->}_{} "I7";"I9"};
        {\ar@{-->}_{} "I10";"I2"};
        {\ar@{-->}_{} "I10";"I12"};
       {\ar@{-->}_{} "I13";"I8"};  
        {\ar@{-->}_{} "I8";"I12"};
        {\ar@{-->}_{} "I12";"I9"};
        {\ar@{-->}_{} "I14";"I17"};
        {\ar@{-->}_{} "I17";"O9"};
        {\ar@{-->}_{} "I4";"I14"};
        {\ar@{-->}_{} "I13";"I16"};
        {\ar@{-->}_{} "I16";"I12"};
        {\ar@{-->}_{} "I16";"I17"};
        {\ar@{=>}^{\scriptstyle {\alpha_{\alpha\times 1,1,1}}} (-25, 13); (-21, 13)};  %DOMAIN 2-CELLS
        {\ar@{=>}^{\scriptstyle {\alpha_{\alpha,1,1}}} (16, 13); (20, 15)};   
        {\ar@{=>}^{\scriptstyle 1\times\Pi} (6, -8); (9, -5)};
        {\ar@{=>}^{\scriptstyle 1\times \Pi} (-5, -20); (1, -17)};
        {\ar@{=>}^{\scriptstyle \alpha_{1\times\alpha,1,1}} (-30, -8); (-27, -5)};
        {\ar@{=>}^{\scriptstyle \Pi} (-10, 0); (-7, 3)};
        {\ar@{=>}^{\scriptstyle \alpha_{1\times\alpha,1,1}} (-63, -4); (-59, -2)};
        {\ar@{=>}^{\scriptstyle \alpha_{\alpha,1,1}} (-60, 14); (-56, 11)};
        {\ar@{=>}^{\scriptstyle \Pi} (-54, -20); (-50, -16)};
        {\ar@{=>}_{\scriptstyle \alpha_{1,1\times\alpha,1}} (-35, -37); (-29,-33)};
        {\ar@{=>}^{\scriptstyle \alpha_{1,\alpha,1}} (-67, -33); (-61,-29)};
        {\ar@{=>}^{\scriptstyle 1\times (1\times \Pi)} (33, -33); (37, -30)};
        {\ar@{=>}^{\scriptstyle 1\times \alpha_{1,\alpha,1}} (3,-33); (6, -30)};
        {\ar@{=>}^{\scriptstyle \Pi} (33, 4); (37,4)};
\endxy
\]
\[\textrm{\Huge{=}}\]
\[
\def\objectstyle{\scriptstyle}
  \def\labelstyle{\scriptstyle}
   \xy
   (-22.5,20)*+{((((AB)C)D)E)F}="O1"; %BORDER VERTICES
   (10,20)*+{(((AB)C)D)(EF)}="O2";
   (45,0)*+{(AB)(C(D(EF)))}="O3";
    (-88,0)*+{((A(B(CD))E)F}="O4"; 
   (-88,-20)*+{(A((B(CD))E))F}="O5";   
   (-88, 20)*+{((A((BC)D))E)F}="O6";   
   (-55, 20)*+{(((A(BC))D)E)F}="O7";   
   (45,20)*+{((AB)C)(D(EF))}="O8";
   (45,-20)*+{A(B(C(D(EF))))}="O9"; 
   (-88 ,-40)*+{(A(B((CD)E)))F}="O10";   
   (-55,-40)*+{(A(B(C(DE))))F}="O11";   
   (10,-40)*+{A(B((C(DE))F))}="O12"; 
   (45,-40)*+{A(B(C((DE)F)))}="O13"; 
   (-22.5,-40)*+{A((B(C(DE)))F)}="O14";
   (-49,-10)*+{((AB)((CD)E))F}="I1"; %INNER VERTICES
   (0,-26)*+{(AB)((C(DE))F)}="I2";
   (-4,-5)*+{(((AB)C)(DE))F}="I3";
   (-35,5)*+{(((AB)(CD))E)F}="I4";
   (20,6)*+{((AB)C)((DE)F)}="I5";
   (-35,-25)*+{((AB)(C(DE)))F}="I6";
   (25,-10)*+{(AB)(C((DE)F))}="I7"; 
        {\ar_{} "O5";"O10"}; % BORDER ARROWS
        {\ar_{} "O10";"O11"};
        {\ar@{-->}_{} "O11";"O14"};
        {\ar@{-->}_{} "O13";"O9"};
        {\ar_{} "O14";"O12"};
        {\ar_{} "O12";"O13"};
        {\ar@{-->}_{} "O1";"O2"};    
        {\ar_{} "O1";"O7"};    
        {\ar_{} "O7";"O6"};        
        {\ar_{} "O6";"O4"};    
        {\ar_{} "O4";"O5"};    
        {\ar@{-->}^{} "O8";"O3"};
        {\ar@{-->}^{} "O2";"O8"};                
        {\ar@{-->}_{} "O3";"O9"};       
        {\ar_{} "O1";"I3"}; % IN/OUT ARROWS
        {\ar_{} "O1";"I4"};       
        {\ar_{} "I4";"O4"}; 
        {\ar@{-->}_{} "I5";"O8"};
        {\ar_{} "I1";"O10"}; 
        {\ar_{} "I6";"O11"};
        {\ar_{} "I2";"O12"};
        {\ar@{-->}_{} "I7";"O3"};
        {\ar@{-->}_{} "I7";"O13"};
        {\ar_{} "I4";"I1"};  % INTERIOR ARROWS
        {\ar_{} "I1";"I6"};     
        {\ar_{} "I3";"I6"};     
        {\ar@{-->}_{} "I3";"I5"};    
        {\ar@{-->}_{} "I5";"I7"};
        {\ar@{-->}_{} "I2";"I7"};
        {\ar@{-->}_{} "I6";"I2"};
        {\ar@{=>}^{\scriptstyle \Pi} (7, 9); (10, 13)};    %DOMAIN 2-CELLS
        {\ar@{=>}_{\scriptstyle \alpha_{1,1,\alpha}^{-1}} (32, 4); (35, 8)};               
        {\ar@{=>}^{\scriptstyle \Pi} (5, -16); (10, -10)};
        {\ar@{=>}^{\scriptstyle \Pi\times 1} (-27, -8); (-23, -8)};
        {\ar@{=>}^{\scriptstyle (\Pi\times 1)\times 1} (-60, 14); (-56, 10)};
        {\ar@{=>}^{\scriptstyle \Pi\times 1} (-73, -7); (-68, -10)};
        {\ar@{=>}^{\scriptstyle \alpha_{1,1,\alpha}^{-1}\times 1} (-53, -27); (-49, -24)};
        {\ar@{=>}^{\scriptstyle \Pi} (-16, -35); (-11,-32)};
        {\ar@{=>}_{\scriptstyle \alpha^{-1}_{1,1,\alpha}} (23, -33); (27, -28)};
        {\ar@{=>}^{\scriptstyle \alpha^{-1}_{1,1,1\times\alpha}} (40, -17); (42, -13)};
\endxy
\]
\[
\def\objectstyle{\scriptstyle}
  \def\labelstyle{\scriptstyle}
   \xy  
    {\ar@{=>}_{\scriptstyle K_{5}\times 1} (-6, 2); (-6, -3)}; 
   \endxy
\]
\[
\def\objectstyle{\scriptstyle}
  \def\labelstyle{\scriptstyle}
   \xy
   (-22.5,20)*+{((((AB)C)D)E)F}="O1"; %BORDER VERTICES
   (10,20)*+{(((AB)C)D)(EF)}="O2";
   (45,0)*+{(AB)(C(D(EF)))}="O3";
    (-88,0)*+{((A(B(CD))E)F}="O4"; 
   (-88,-20)*+{(A((B(CD))E))F}="O5";   
   (-88, 20)*+{((A((BC)D))E)F}="O6";   
   (-55, 20)*+{(((A(BC))D)E)F}="O7";   
   (45,20)*+{((AB)C)(D(EF))}="O8";
   (45,-20)*+{A(B(C(D(EF))))}="O9"; 
   (-88 ,-40)*+{(A(B((CD)E)))F}="O10";   
   (-55,-40)*+{(A(B(C(DE))))F}="O11";   
   (10,-40)*+{A(B((C(DE))F))}="O12"; 
   (45,-40)*+{A(B(C((DE)F)))}="O13"; 
   (-22.5,-40)*+{A((B(C(DE)))F)}="O14";
   (0,-26)*+{(AB)((C(DE))F)}="I2"; %INNER VERTICES
   (-4,-5)*+{(((AB)C)(DE))F}="I3";
   (20,6)*+{((AB)C)((DE)F)}="I5";
   (-30,-25)*+{((AB)(C(DE)))F}="I6";
   (25,-10)*+{(AB)(C((DE)F))}="I7"; 
   (-30,5)*+{((A(BC))(DE))F}="I8";
   (-45,-15)*+{(A((BC)(DE)))F}="I9"; 
   (-60,5)*+{(A(((BC)D)E))F}="I10";
        {\ar@{-->}_{} "O5";"O10"}; % BORDER ARROWS
        {\ar@{-->}_{} "O10";"O11"};
        {\ar_{} "O11";"O14"};
        {\ar@{-->}_{} "O13";"O9"};
        {\ar_{} "O14";"O12"};
        {\ar_{} "O12";"O13"};
        {\ar@{-->}_{} "O1";"O2"};    
        {\ar@{-->}_{} "O1";"O7"};    
        {\ar@{-->}_{} "O7";"O6"};        
        {\ar@{-->}_{} "O6";"O4"};    
        {\ar@{-->}_{} "O4";"O5"};    
        {\ar@{-->}^{} "O8";"O3"};
        {\ar@{-->}^{} "O2";"O8"};                
        {\ar@{-->}_{} "O3";"O9"};       
        {\ar@{-->}_{} "O1";"I3"}; % IN/OUT ARROWS
        {\ar@{-->}_{} "O7";"I8"};       
        {\ar@{-->}_{} "O6";"I10"};  
        {\ar@{-->}_{} "I5";"O8"}; 
        {\ar_{} "I6";"O11"};
        {\ar_{} "I9";"O11"};
        {\ar_{} "I2";"O12"};
        {\ar@{-->}_{} "I7";"O3"};
        {\ar_{} "I7";"O13"};  
        {\ar@{-->}_{} "I10";"O5"};   
        {\ar_{} "I3";"I6"};      % INTERIOR ARROWS
        {\ar_{} "I3";"I5"};    
        {\ar_{} "I5";"I7"};
        {\ar_{} "I2";"I7"};
        {\ar_{} "I6";"I2"};
        {\ar_{} "I3";"I8"};
        {\ar_{} "I8";"I9"};
        {\ar@{-->}_{} "I10";"I9"};
        {\ar@{=>}^{\scriptstyle \Pi} (7, 9); (10, 13)};  %DOMAIN 2-CELLS
        {\ar@{=>}_{\scriptstyle \alpha_{1,1,\alpha}^{-1}} (32, 4); (35, 8)};                        
        {\ar@{=>}^{\scriptstyle \Pi} (5, -16); (10, -10)};
        {\ar@{=>}^{\scriptstyle \Pi\times 1} (-27, -8); (-23, -8)};
        {\ar@{=>}^{\scriptstyle \Pi\times 1} (-58, 13); (-54, 12)};
        {\ar@{=>}^{\scriptstyle \alpha_{\alpha,1,1}\times 1} (-30, 13); (-27, 12)};
        {\ar@{=>}^{\scriptstyle \alpha_{1,\alpha,1}\times 1} (-77, -4); (-74, -4)};
        {\ar@{=>}^{\scriptstyle (1\times \Pi)\times 1} (-67, -24); (-63, -24)};
        {\ar@{=>}^{\scriptstyle \Pi} (-16, -35); (-11,-32)};
        {\ar@{=>}_{\scriptstyle \alpha^{-1}_{1,1,\alpha}} (23, -33); (27, -28)};
        {\ar@{=>}^{\scriptstyle \alpha^{-1}_{1,1,1\times\alpha}} (40, -17); (42, -13)};
\endxy
\]
\[
\def\objectstyle{\scriptstyle}
  \def\labelstyle{\scriptstyle}
   \xy  
    {\ar@{=>}_{\scriptstyle K_{5}} (-6, 2); (-6, -3)}; 
   \endxy
\]
\[
\def\objectstyle{\scriptstyle}
  \def\labelstyle{\scriptstyle}
   \xy
   (-22.5,20)*+{((((AB)C)D)E)F}="O1"; %BORDER VERTICES
   (10,20)*+{(((AB)C)D)(EF)}="O2";
   (45,0)*+{(AB)(C(D(EF)))}="O3";
    (-88,0)*+{((A(B(CD))E)F}="O4"; 
   (-88,-20)*+{(A((B(CD))E))F}="O5";   
   (-88, 20)*+{((A((BC)D))E)F}="O6";   
   (-55, 20)*+{(((A(BC))D)E)F}="O7";   
   (45,20)*+{((AB)C)(D(EF))}="O8";
   (45,-20)*+{A(B(C(D(EF))))}="O9"; 
   (-88 ,-40)*+{(A(B((CD)E)))F}="O10";   
   (-55,-40)*+{(A(B(C(DE))))F}="O11";   
   (10,-40)*+{A(B((C(DE))F))}="O12"; 
   (45,-40)*+{A(B(C((DE)F)))}="O13"; 
   (-22.5,-40)*+{A((B(C(DE)))F)}="O14";
   (-4,-5)*+{(((AB)C)(DE))F}="I3"; %INNER VERTICES
   (20,6)*+{((AB)C)((DE)F)}="I5";
   (25,-10)*+{(AB)(C((DE)F))}="I7"; 
   (-30,5)*+{((A(BC))(DE))F}="I8";
   (-45,-15)*+{(A((BC)(DE)))F}="I9"; 
   (-60,5)*+{(A(((BC)D)E))F}="I10";
   (-5,-30)*+{A(((BC)(DE))F)}="I11"; 
   (23,-25)*+{A((BC)((DE)F))}="I12"; 
   (0,-17)*+{(A(BC))((DE)F)}="I13"; 
        {\ar@{-->}_{} "O5";"O10"}; % BORDER ARROWS
        {\ar@{-->}_{} "O10";"O11"};
        {\ar@{-->}_{} "O11";"O14"};
        {\ar_{} "O13";"O9"};
        {\ar@{-->}_{} "O14";"O12"};
        {\ar@{-->}_{} "O12";"O13"};
        {\ar@{-->}_{} "O1";"O2"};    
        {\ar@{-->}_{} "O1";"O7"};    
        {\ar@{-->}_{} "O7";"O6"};        
        {\ar@{-->}_{} "O6";"O4"};    
        {\ar@{-->}_{} "O4";"O5"};    
        {\ar^{} "O8";"O3"};
        {\ar@{-->}^{} "O2";"O8"};                
        {\ar_{} "O3";"O9"};       
        {\ar@{-->}_{} "O1";"I3"}; % IN/OUT ARROWS
        {\ar@{-->}_{} "O7";"I8"};       
        {\ar@{-->}_{} "O6";"I10"};  
        {\ar_{} "I5";"O8"}; 
        {\ar@{-->}_{} "I9";"O11"};
        {\ar_{} "I7";"O3"};
        {\ar_{} "I7";"O13"};  
        {\ar@{-->}_{} "I10";"O5"};   
        {\ar_{} "I12";"O13"}; 
        {\ar@{-->}_{} "I11";"O14"}; 
        {\ar@{-->}_{} "I3";"I5"};     % INTERIOR ARROWS
        {\ar_{} "I5";"I7"};
        {\ar@{-->}_{} "I3";"I8"};
        {\ar@{-->}_{} "I8";"I9"};
        {\ar@{-->}_{} "I10";"I9"};
        {\ar_{} "I5";"I13"};
        {\ar_{} "I13";"I12"};
        {\ar@{-->}_{} "I11";"I12"};
        {\ar@{-->}_{} "I9";"I11"};
        {\ar@{-->}_{} "I8";"I13"};
        {\ar@{=>}^{\scriptstyle \Pi} (7, 9); (10, 13)};   %DOMAIN 2-CELLS
        {\ar@{=>}_{\scriptstyle \alpha_{1,1,\alpha}^{-1}} (32, 4); (35, 8)};      
        {\ar@{=>}^{\scriptstyle \Pi} (15, -15); (20, -15)};
        {\ar@{=>}^{\scriptstyle \alpha_{\alpha,1,1}} (-1, -11); (3, -11)};
        {\ar@{=>}^{\scriptstyle \Pi} (-25, -15); (-22, -12)};
        {\ar@{=>}^{\scriptstyle \Pi\times 1} (-58, 13); (-54, 12)};
        {\ar@{=>}^{\scriptstyle \alpha_{\alpha,1,1}\times 1} (-30, 13); (-27, 12)};
        {\ar@{=>}^{\scriptstyle \alpha_{1,\alpha,1}\times 1} (-77, -4); (-74, -4)};
        {\ar@{=>}^{\scriptstyle (1\times \Pi)\times 1} (-67, -24); (-63, -24)};
        {\ar@{=>}^{\scriptstyle \alpha_{1,\alpha,1}} (-34, -35); (-30,-32)};
        {\ar@{=>}_{\scriptstyle 1\times\Pi} (19, -35); (23, -32)};
        {\ar@{=>}^{\scriptstyle \alpha^{-1}_{1,1,1\times\alpha}} (40, -17); (42, -13)};
\endxy
\]
\[
\def\objectstyle{\scriptstyle}
  \def\labelstyle{\scriptstyle}
   \xy  
    {\ar@{=>}_{\scriptstyle \Pi} (-6, 2); (-6, -3)}; 
   \endxy
\]
\[
\def\objectstyle{\scriptstyle}
  \def\labelstyle{\scriptstyle}
   \xy
   (-22.5,20)*+{((((AB)C)D)E)F}="O1"; %BORDER VERTICES
   (10,20)*+{(((AB)C)D)(EF)}="O2";
   (45,0)*+{(AB)(C(D(EF)))}="O3";
    (-88,0)*+{((A(B(CD))E)F}="O4"; 
   (-88,-20)*+{(A((B(CD))E))F}="O5";   
   (-88, 20)*+{((A((BC)D))E)F}="O6";   
   (-55, 20)*+{(((A(BC))D)E)F}="O7";   
   (45,20)*+{((AB)C)(D(EF))}="O8";
   (45,-20)*+{A(B(C(D(EF))))}="O9"; 
   (-88 ,-40)*+{(A(B((CD)E)))F}="O10";   
   (-55,-40)*+{(A(B(C(DE))))F}="O11";   
   (10,-40)*+{A(B((C(DE))F))}="O12"; 
   (45,-40)*+{A(B(C((DE)F)))}="O13"; 
   (-22.5,-40)*+{A((B(C(DE)))F)}="O14";
   (-4,-5)*+{(((AB)C)(DE))F}="I3"; %INNER VERTICES
   (20,6)*+{((AB)C)((DE)F)}="I5";
   (-30,5)*+{((A(BC))(DE))F}="I8";
   (-45,-15)*+{(A((BC)(DE)))F}="I9"; 
   (-60,5)*+{(A(((BC)D)E))F}="I10";
   (-5,-30)*+{A(((BC)(DE))F)}="I11"; 
   (23,-25)*+{A((BC)((DE)F))}="I12"; 
   (0,-17)*+{(A(BC))((DE)F)}="I13"; 
   (30,-5)*+{(A(BC))(D(EF))}="I14"; 
   (32,-14)*+{A((BC)(D(EF)))}="I15"; 
        {\ar@{-->}_{} "O5";"O10"}; % BORDER ARROWS
        {\ar@{-->}_{} "O10";"O11"};
        {\ar@{-->}_{} "O11";"O14"};
        {\ar@{-->}_{} "O13";"O9"};
        {\ar@{-->}_{} "O14";"O12"};
        {\ar@{-->}_{} "O12";"O13"};
        {\ar_{} "O1";"O2"};    
        {\ar_{} "O1";"O7"};    
        {\ar@{-->}_{} "O7";"O6"};        
        {\ar@{-->}_{} "O6";"O4"};    
        {\ar@{-->}_{} "O4";"O5"};    
        {\ar@{-->}^{} "O8";"O3"};
        {\ar^{} "O2";"O8"};                
        {\ar@{-->}_{} "O3";"O9"};       
        {\ar_{} "O1";"I3"}; % IN/OUT ARROWS
        {\ar_{} "O7";"I8"};       
        {\ar@{-->}_{} "O6";"I10"};  
        {\ar_{} "I5";"O8"}; 
        {\ar@{-->}_{} "I9";"O11"};  
        {\ar@{-->}_{} "I10";"O5"};   
        {\ar@{-->}_{} "I12";"O13"}; 
        {\ar@{-->}_{} "I11";"O14"}; 
        {\ar_{} "O8";"I14"}; 
        {\ar@{-->}_{} "I15";"O9"}; 
        {\ar_{} "I3";"I5"};     % INTERIOR ARROWS
        {\ar_{} "I3";"I8"};
        {\ar@{-->}_{} "I8";"I9"};
        {\ar@{-->}_{} "I10";"I9"};
        {\ar_{} "I5";"I13"};
        {\ar@{-->}_{} "I13";"I12"};
        {\ar@{-->}_{} "I11";"I12"};
        {\ar@{-->}_{} "I9";"I11"};
        {\ar_{} "I8";"I13"};
        {\ar_{} "I13";"I14"};
        {\ar@{-->}_{} "I14";"I15"};
        {\ar@{-->}_{} "I12";"I15"};
        {\ar@{=>}^{\scriptstyle \Pi} (7, 9); (10, 13)}; %DOMAIN 2-CELLS
        {\ar@{=>}^{\scriptstyle {\rm \alpha^{-1}\widetilde{\times}\alpha}} (28, -2); (32, 2)};        
        {\ar@{=>}^{\scriptstyle \alpha_{1,1,\alpha}^{-1}} (20, -20); (24, -17)};
        {\ar@{=>}^{\scriptstyle \alpha_{\alpha,1,1}} (-1, -11); (3, -11)};
        {\ar@{=>}^{\scriptstyle \Pi} (-25, -15); (-22, -12)};
        {\ar@{=>}^{\scriptstyle \Pi\times 1} (-58, 13); (-54, 12)};
        {\ar@{=>}^{\scriptstyle \alpha_{\alpha,1,1}\times 1} (-30, 13); (-27, 12)};
        {\ar@{=>}^{\scriptstyle \alpha_{1,\alpha,1}\times 1} (-77, -4); (-74, -4)};
        {\ar@{=>}^{\scriptstyle (1\times \Pi)\times 1} (-67, -24); (-63, -24)};
        {\ar@{=>}^{\scriptstyle \alpha_{1,\alpha,1}} (-34, -35); (-30,-32)};
        {\ar@{=>}_{\scriptstyle 1\times\Pi} (19, -35); (23, -32)};
        {\ar@{=>}^{\scriptstyle 1\times \alpha^{-1}_{1,1,\alpha}} (43, -27); (39, -24)};
        {\ar@{=>}^{\scriptstyle \Pi} (37, -10); (41, -10)};
\endxy
\]
\[
\def\objectstyle{\scriptstyle}
  \def\labelstyle{\scriptstyle}
   \xy  
    {\ar@{=>}_{\scriptstyle \Pi} (-6, 2); (-6, -3)}; 
   \endxy
\]
\[
\def\objectstyle{\scriptstyle}
  \def\labelstyle{\scriptstyle}
   \xy
   (-22.5,20)*+{((((AB)C)D)E)F}="O1"; %BORDER VERTICES
   (10,20)*+{(((AB)C)D)(EF)}="O2";
   (45,0)*+{(AB)(C(D(EF)))}="O3";
    (-88,0)*+{((A(B(CD))E)F}="O4"; 
   (-88,-20)*+{(A((B(CD))E))F}="O5";   
   (-88, 20)*+{((A((BC)D))E)F}="O6";   
   (-55, 20)*+{(((A(BC))D)E)F}="O7";   
   (45,20)*+{((AB)C)(D(EF))}="O8";
   (45,-20)*+{A(B(C(D(EF))))}="O9"; 
   (-88 ,-40)*+{(A(B((CD)E)))F}="O10";   
   (-55,-40)*+{(A(B(C(DE))))F}="O11";   
   (10,-40)*+{A(B((C(DE))F))}="O12"; 
   (45,-40)*+{A(B(C((DE)F)))}="O13"; 
   (-22.5,-40)*+{A((B(C(DE)))F)}="O14"; 
   (-30,5)*+{((A(BC))(DE))F}="I8";   %INNER VERTICES
   (-45,-15)*+{(A((BC)(DE)))F}="I9"; 
   (-60,5)*+{(A(((BC)D)E))F}="I10";
   (-5,-30)*+{A(((BC)(DE))F)}="I11"; 
   (25,-25)*+{A((BC)((DE)F))}="I12"; 
   (0,-17)*+{(A(BC))((DE)F)}="I13"; 
   (30,-5)*+{(A(BC))(D(EF))}="I14"; 
   (32,-14)*+{A((BC)(D(EF)))}="I15"; 
   (-4,10)*+{((A(BC))D)(EF)}="I16"; 
        {\ar@{-->}_{} "O5";"O10"}; % BORDER ARROWS
        {\ar@{-->}_{} "O10";"O11"};
        {\ar@{-->}_{} "O11";"O14"};
        {\ar@{-->}_{} "O13";"O9"};
        {\ar@{-->}_{} "O14";"O12"};
        {\ar@{-->}_{} "O12";"O13"};
        {\ar@{-->}_{} "O1";"O2"};    
        {\ar@{-->}_{} "O1";"O7"};    
        {\ar_{} "O7";"O6"};        
        {\ar@{-->}_{} "O6";"O4"};    
        {\ar@{-->}_{} "O4";"O5"};    
        {\ar@{-->}^{} "O8";"O3"};
        {\ar@{-->}^{} "O2";"O8"};                
        {\ar@{-->}_{} "O3";"O9"};       
        {\ar_{} "O7";"I8"};        % IN/OUT ARROWS  
        {\ar_{} "O6";"I10"};   
        {\ar@{-->}_{} "I9";"O11"};  
        {\ar@{-->}_{} "I10";"O5"};   
        {\ar@{-->}_{} "I12";"O13"}; 
        {\ar@{-->}_{} "I11";"O14"}; 
        {\ar@{-->}_{} "O8";"I14"}; 
        {\ar@{-->}_{} "I15";"O9"}; 
        {\ar@{-->}_{} "O2";"I16"}; 
        {\ar_{} "O7";"I16"};    
        {\ar_{} "I8";"I9"};  % INTERIOR ARROWS
        {\ar_{} "I10";"I9"};
        {\ar_{} "I13";"I12"};
        {\ar_{} "I11";"I12"};
        {\ar_{} "I9";"I11"};
        {\ar_{} "I8";"I13"};
        {\ar_{} "I13";"I14"};
        {\ar_{} "I14";"I15"};
        {\ar_{} "I12";"I15"};
        {\ar_{} "I16";"I14"};
        {\ar@{=>}^{\scriptstyle \alpha_{\alpha,1,1}} (25, 9); (29, 10)};  %DOMAIN 2-CELLS   
        {\ar@{=>}^{\scriptstyle \alpha_{1,1,\alpha}^{-1}} (20, -20); (24, -17)};
        {\ar@{=>}^{\scriptstyle \Pi} (-1, -6); (3, -3)};
        {\ar@{=>}^{\scriptstyle \Pi} (-25, -15); (-22, -12)};
        {\ar@{=>}^{\scriptstyle \Pi\times 1} (-58, 13); (-54, 12)};
        {\ar@{=>}^{\scriptstyle \alpha_{\alpha\times 1,1,1}} (-14, 14); (-8, 14)};
        {\ar@{=>}^{\scriptstyle \alpha_{1,\alpha,1}\times 1} (-77, -4); (-74, -4)};
        {\ar@{=>}^{\scriptstyle (1\times \Pi)\times 1} (-67, -24); (-63, -24)};
        {\ar@{=>}^{\scriptstyle \alpha_{1,\alpha,1}} (-28, -35); (-24,-32)};
        {\ar@{=>}_{\scriptstyle 1\times\Pi} (19, -35); (23, -32)};
        {\ar@{=>}^{\scriptstyle 1\times \alpha^{-1}_{1,1,\alpha}} (43, -27); (39, -24)};
        {\ar@{=>}^{\scriptstyle \Pi} (37, -10); (41, -10)};
\endxy
\]
\[
\def\objectstyle{\scriptstyle}
  \def\labelstyle{\scriptstyle}
   \xy  
    {\ar@{=>}_{\scriptstyle K_{5}} (-6, 2); (-6, -3)}; 
   \endxy
\]
\[
\def\objectstyle{\scriptstyle}
  \def\labelstyle{\scriptstyle}
   \xy
   (-22.5,20)*+{((((AB)C)D)E)F}="O1"; %BORDER VERTICES
   (10,20)*+{(((AB)C)D)(EF)}="O2";
   (45,0)*+{(AB)(C(D(EF)))}="O3";
    (-88,0)*+{((A(B(CD))E)F}="O4"; 
   (-88,-20)*+{(A((B(CD))E))F}="O5";   
   (-88, 20)*+{((A((BC)D))E)F}="O6";   
   (-55, 20)*+{(((A(BC))D)E)F}="O7";   
   (45,20)*+{((AB)C)(D(EF))}="O8";
   (45,-20)*+{A(B(C(D(EF))))}="O9"; 
   (-88 ,-40)*+{(A(B((CD)E)))F}="O10";   
   (-55,-40)*+{(A(B(C(DE))))F}="O11";   
   (10,-40)*+{A(B((C(DE))F))}="O12"; 
   (45,-40)*+{A(B(C((DE)F)))}="O13"; 
   (-22.5,-40)*+{A((B(C(DE)))F)}="O14"; 
   (-45,-15)*+{(A((BC)(DE)))F}="I9";  %INNER VERTICES
   (-5,-30)*+{A(((BC)(DE))F)}="I11"; 
   (25,-25)*+{A((BC)((DE)F))}="I12"; 
   (30,-5)*+{(A(BC))(D(EF))}="I14"; 
   (32,-14)*+{A((BC)(D(EF)))}="I15"; 
   (-4,10)*+{((A(BC))D)(EF)}="I16";
   (-45,10)*+{(A((BC)D))(EF)}="I17";
   (-60,5)*+{(A(((BC)D)E))F}="I18";
   (-10,-5)*+{A(((BC)D)(EF))}="I19"; 
   (-17,-15)*+{A(((((BC)D)E)F)}="I20"; 
        {\ar_{} "O5";"O10"}; % BORDER ARROWS
        {\ar_{} "O10";"O11"};
        {\ar_{} "O11";"O14"};
        {\ar@{-->}_{} "O13";"O9"};
        {\ar@{-->}_{} "O14";"O12"};
        {\ar@{-->}_{} "O12";"O13"};
        {\ar@{-->}_{} "O1";"O2"};    
        {\ar@{-->}_{} "O1";"O7"};    
        {\ar@{-->}_{} "O7";"O6"};        
        {\ar@{-->}_{} "O6";"O4"};    
        {\ar@{-->}_{} "O4";"O5"};    
        {\ar@{-->}^{} "O8";"O3"};
        {\ar@{-->}^{} "O2";"O8"};                
        {\ar@{-->}_{} "O3";"O9"};       
        {\ar_{} "I9";"O11"};    % IN/OUT ARROWS     
        {\ar@{-->}_{} "I12";"O13"}; 
        {\ar_{} "I11";"O14"}; 
        {\ar@{-->}_{} "O8";"I14"}; 
        {\ar@{-->}_{} "I15";"O9"}; 
        {\ar@{-->}_{} "O2";"I16"}; 
        {\ar@{-->}_{} "O7";"I16"};    
        {\ar@{-->}_{} "O6";"I17"};       
        {\ar@{-->}_{} "O6";"I18"};
        {\ar_{} "I18";"O5"}; 
        {\ar@{-->}_{} "I11";"I12"};  % INTERIOR ARROWS
        {\ar_{} "I9";"I11"};
        {\ar@{-->}_{} "I14";"I15"};
        {\ar@{-->}_{} "I12";"I15"};
        {\ar@{-->}_{} "I16";"I14"};
        {\ar@{-->}_{} "I16";"I17"};
        {\ar@{-->}_{} "I17";"I19"};
        {\ar@{-->}_{} "I19";"I15"};
        {\ar@{-->}_{} "I20";"I19"};
        {\ar_{} "I20";"I11"};
        {\ar_{} "I18";"I9"};
        {\ar_{} "I18";"I20"};
        {\ar@{=>}^{\scriptstyle \alpha_{\alpha,1,1}} (25, 9); (29, 10)};  %DOMAIN 2-CELLS   
        {\ar@{=>}^{\scriptstyle 1\times \Pi} (12, -19); (13, -15)};
        {\ar@{=>}^{\scriptstyle \Pi} (-1,0); (3, 1)};
        {\ar@{=>}^{\scriptstyle \alpha_{1,\alpha,1}} (-27 ,-19); (-24, -19)};
        {\ar@{=>}^{\scriptstyle \Pi} (-35, -4); (-31, -1)};
        {\ar@{=>}^{\scriptstyle \alpha_{\alpha,1,1}} (-54, 15); (-50, 14)};
        {\ar@{=>}^{\scriptstyle \alpha_{\alpha\times 1,1,1}} (-14, 14); (-8, 14)};
        {\ar@{=>}^{\scriptstyle \alpha_{1,\alpha,1}\times 1} (-77, -4); (-74, -4)};
        {\ar@{=>}^{\scriptstyle (1\times \Pi)\times 1} (-67, -24); (-63, -24)};
        {\ar@{=>}^{\scriptstyle \alpha_{1,\alpha,1}} (-28, -35); (-24,-32)};
        {\ar@{=>}_{\scriptstyle 1\times\Pi} (19, -35); (23, -32)};
        {\ar@{=>}^{\scriptstyle 1\times \alpha^{-1}_{1,1,\alpha}} (43, -27); (39, -24)};
        {\ar@{=>}^{\scriptstyle \Pi} (37, -10); (41, -10)};
\endxy
\]
\[
\def\objectstyle{\scriptstyle}
  \def\labelstyle{\scriptstyle}
   \xy  
    {\ar@{=>}_{\scriptstyle \alpha} (-6, 2); (-6, -3)}; 
   \endxy
\]
\[
\def\objectstyle{\scriptstyle}
  \def\labelstyle{\scriptstyle}
   \xy
   (-22.5,20)*+{((((AB)C)D)E)F}="O1"; %BORDER VERTICES
   (10,20)*+{(((AB)C)D)(EF)}="O2";
   (45,0)*+{(AB)(C(D(EF)))}="O3";
    (-88,0)*+{((A(B(CD))E)F}="O4"; 
   (-88,-20)*+{(A((B(CD))E))F}="O5";   
   (-88, 20)*+{((A((BC)D))E)F}="O6";   
   (-55, 20)*+{(((A(BC))D)E)F}="O7";   
   (45,20)*+{((AB)C)(D(EF))}="O8";
   (45,-20)*+{A(B(C(D(EF))))}="O9"; 
   (-88 ,-40)*+{(A(B((CD)E)))F}="O10";   
   (-55,-40)*+{(A(B(C(DE))))F}="O11";   
   (10,-40)*+{A(B((C(DE))F))}="O12"; 
   (45,-40)*+{A(B(C((DE)F)))}="O13"; 
   (-22.5,-40)*+{A((B(C(DE)))F)}="O14";  
   (-5,-30)*+{A(((BC)(DE))F)}="I11"; %INNER VERTICES
   (25,-25)*+{A((BC)((DE)F))}="I12";
   (30,-5)*+{(A(BC))(D(EF))}="I14"; 
   (32,-14)*+{A((BC)(D(EF)))}="I15"; 
   (-4,10)*+{((A(BC))D)(EF)}="I16";
   (-45,10)*+{(A((BC)D))(EF)}="I17";
   (-60,5)*+{(A(((BC)D)E))F}="I18";
   (-10,-5)*+{A(((BC)D)(EF))}="I19";
   (-17,-15)*+{A(((((BC)D)E)F)}="I20"; 
   (-50,-20)*+{A(((B(CD))E)F)}="I21";  
    (-32,-30)*+{A((B((CD)E))F)}="I22";  
        {\ar@{-->}_{} "O5";"O10"}; % BORDER ARROWS
        {\ar@{-->}_{} "O10";"O11"};
        {\ar@{-->}_{} "O11";"O14"};
        {\ar_{} "O13";"O9"};
        {\ar_{} "O14";"O12"};
        {\ar_{} "O12";"O13"};
        {\ar@{-->}_{} "O1";"O2"};    
        {\ar@{-->}_{} "O1";"O7"};    
        {\ar@{-->}_{} "O7";"O6"};        
        {\ar@{-->}_{} "O6";"O4"};    
        {\ar@{-->}_{} "O4";"O5"};    
        {\ar@{-->}^{} "O8";"O3"};
        {\ar@{-->}^{} "O2";"O8"};                
        {\ar@{-->}_{} "O3";"O9"};       
        {\ar_{} "I12";"O13"};  % IN/OUT ARROWS     
        {\ar_{} "I11";"O14"}; 
        {\ar@{-->}_{} "O8";"I14"}; 
        {\ar_{} "I15";"O9"}; 
        {\ar@{-->}_{} "O2";"I16"}; 
        {\ar@{-->}_{} "O7";"I16"};    
        {\ar@{-->}_{} "O6";"I17"};       
        {\ar@{-->}_{} "O6";"I18"};
        {\ar@{-->}_{} "I18";"O5"}; 
        {\ar@{-->}_{} "O5";"I21"};
        {\ar@{-->}_{} "O10";"I22"};
        {\ar_{} "I22";"O14"}; 
        {\ar_{} "I11";"I12"};  % INTERIOR ARROWS
        {\ar@{-->}_{} "I14";"I15"};
        {\ar_{} "I12";"I15"};
        {\ar@{-->}_{} "I16";"I14"};
        {\ar@{-->}_{} "I16";"I17"};
        {\ar@{-->}_{} "I17";"I19"};
        {\ar_{} "I19";"I15"};
        {\ar_{} "I20";"I19"};
        {\ar_{} "I20";"I11"};
        {\ar@{-->}_{} "I18";"I20"};
        {\ar_{} "I20";"I21"};
        {\ar_{} "I21";"I22"};
        {\ar@{=>}^{\scriptstyle \alpha_{\alpha,1,1}} (25, 9); (29, 10)};  %DOMAIN 2-CELLS   
        {\ar@{=>}^{\scriptstyle 1\times \Pi} (12, -19); (13, -15)};
        {\ar@{=>}^{\scriptstyle \Pi} (-1,0); (3, 1)};
        {\ar@{=>}^{\scriptstyle \alpha_{1,\alpha \times 1,1}} (-51 ,-12); (-48, -11)};
        {\ar@{=>}^{\scriptstyle \Pi} (-35, -4); (-31, -1)};
        {\ar@{=>}^{\scriptstyle \alpha_{\alpha,1,1}} (-54, 15); (-50, 14)};
        {\ar@{=>}^{\scriptstyle \alpha_{\alpha\times 1,1,1}} (-14, 14); (-8, 14)};
        {\ar@{=>}^{\scriptstyle \alpha_{1,\alpha,1}\times 1} (-77, -4); (-74, -4)};
        {\ar@{=>}^{\scriptstyle \alpha_{1,\alpha,1}} (-67, -30); (-63, -26)};
        {\ar@{=>}^{\scriptstyle \alpha_{1,1\times\alpha,1}} (-36, -38); (-33,-35)};
        {\ar@{=>}^{\scriptstyle 1\times (\Pi\times 1)} (-26, -25); (-22,-28)};
        {\ar@{=>}_{\scriptstyle 1\times\Pi} (19, -35); (23, -32)};
        {\ar@{=>}^{\scriptstyle 1\times \alpha^{-1}_{1,1,\alpha}} (43, -27); (39, -24)};
        {\ar@{=>}^{\scriptstyle \Pi} (37, -10); (41, -10)};
\endxy
\]
\[
\def\objectstyle{\scriptstyle}
  \def\labelstyle{\scriptstyle}
   \xy  
    {\ar@{=>}_{\scriptstyle 1\times K_{5}} (-6, 2); (-6, -3)}; 
   \endxy
\]
\[
\def\objectstyle{\scriptstyle}
  \def\labelstyle{\scriptstyle}
   \xy
   (-22.5,20)*+{((((AB)C)D)E)F}="O1"; %BORDER VERTICES
   (10,20)*+{(((AB)C)D)(EF)}="O2";
   (45,0)*+{(AB)(C(D(EF)))}="O3";
    (-88,0)*+{((A(B(CD))E)F}="O4"; 
   (-88,-20)*+{(A((B(CD))E))F}="O5";   
   (-88, 20)*+{((A((BC)D))E)F}="O6";   
   (-55, 20)*+{(((A(BC))D)E)F}="O7";   
   (45,20)*+{((AB)C)(D(EF))}="O8";
   (45,-20)*+{A(B(C(D(EF))))}="O9"; 
   (-88,-40)*+{(A(B((CD)E)))F}="O10";   
   (-55,-40)*+{(A(B(C(DE))))F}="O11";   
   (10,-40)*+{A(B((C(DE))F))}="O12"; 
   (45,-40)*+{A(B(C((DE)F)))}="O13"; 
   (-22.5,-40)*+{A((B(C(DE)))F)}="O14";  
   (30,-5)*+{(A(BC))(D(EF))}="I14";  %INNER VERTICES
   (32,-14)*+{A((BC)(D(EF)))}="I15";
   (-4,10)*+{((A(BC))D)(EF)}="I16";
   (-45,10)*+{(A((BC)D))(EF)}="I17";
   (-60,5)*+{(A(((BC)D)E))F}="I18";
   (-10,-5)*+{A(((BC)D)(EF))}="I19";
   (-17,-15)*+{A(((((BC)D)E)F)}="I20"; 
   (-50,-20)*+{A(((B(CD))E)F)}="I21"; 
    (-32,-30)*+{A((B((CD)E))F)}="I22"; 
   (10,-20)*+{A((B(CD))(EF))}="I23"; 
   (32,-30)*+{A(B((CD)(EF)))}="I24"; 
   (2,-30)*+{A(B(((CD)E)F))}="I25"; 
        {\ar@{-->}_{} "O5";"O10"}; % BORDER ARROWS
        {\ar@{-->}_{} "O10";"O11"};
        {\ar@{-->}_{} "O11";"O14"};
        {\ar@{-->}_{} "O13";"O9"};
        {\ar@{-->}_{} "O14";"O12"};
        {\ar_{} "O12";"O13"};
        {\ar@{-->}_{} "O1";"O2"};    
        {\ar@{-->}_{} "O1";"O7"};    
        {\ar@{-->}_{} "O7";"O6"};        
        {\ar_{} "O6";"O4"};    
        {\ar_{} "O4";"O5"};    
        {\ar@{-->}^{} "O8";"O3"};
        {\ar@{-->}^{} "O2";"O8"};                
        {\ar@{-->}_{} "O3";"O9"};       
        {\ar@{-->}_{} "O8";"I14"}; % IN/OUT ARROWS      
        {\ar@{-->}_{} "I15";"O9"}; 
        {\ar@{-->}_{} "O2";"I16"}; 
        {\ar@{-->}_{} "O7";"I16"};    
        {\ar_{} "O6";"I17"};       
        {\ar_{} "O6";"I18"};
        {\ar_{} "I18";"O5"}; 
        {\ar_{} "O5";"I21"};
        {\ar@{-->}_{} "O10";"I22"};
        {\ar@{-->}_{} "I22";"O14"};  
        {\ar@{-->}_{} "I25";"O12"};  
        {\ar@{-->}_{} "I24";"O9"};  
        {\ar@{-->}_{} "I14";"I15"}; % INTERIOR ARROWS
        {\ar@{-->}_{} "I16";"I14"};
        {\ar@{-->}_{} "I16";"I17"};
        {\ar_{} "I17";"I19"};
        {\ar@{-->}_{} "I19";"I15"};
        {\ar_{} "I20";"I19"};
        {\ar_{} "I18";"I20"};
        {\ar_{} "I20";"I21"};
        {\ar@{-->}_{} "I21";"I22"};
        {\ar_{} "I19";"I23"};
        {\ar@{-->}_{} "I23";"I24"};
        {\ar@{-->}_{} "I22";"I25"};
        {\ar_{} "I25";"I24"};
        {\ar_{} "I21";"I23"};
        {\ar@{=>}^{\scriptstyle \alpha_{\alpha,1,1}} (25, 9); (29, 10)};  %DOMAIN 2-CELLS   
        {\ar@{=>}^{\scriptstyle 1\times \alpha_{\alpha,1,1}} (0, -19); (1, -16)};
        {\ar@{=>}^{\scriptstyle 1\times \Pi} (27, -24); (30, -21)};
        {\ar@{=>}^{\scriptstyle \Pi} (-1,0); (3, 1)};
        {\ar@{=>}^{\scriptstyle \alpha_{1,\alpha \times 1,1}} (-51 ,-12); (-48, -11)};
        {\ar@{=>}^{\scriptstyle \Pi} (-35, -4); (-31, -1)};
        {\ar@{=>}^{\scriptstyle \alpha_{\alpha,1,1}} (-54, 15); (-50, 14)};
        {\ar@{=>}^{\scriptstyle \alpha_{\alpha\times 1,1,1}} (-14, 14); (-8, 14)};
        {\ar@{=>}^{\scriptstyle \alpha_{1,\alpha,1}\times 1} (-77, -4); (-74, -4)};
        {\ar@{=>}^{\scriptstyle \alpha_{1,\alpha,1}} (-67, -30); (-63, -26)};
        {\ar@{=>}^{\scriptstyle \alpha_{1,1\times\alpha,1}} (-36, -38); (-33,-35)};
        {\ar@{=>}^{\scriptstyle 1\times \Pi} (-10, -28); (-7,-25)};
        {\ar@{=>}_{\scriptstyle 1\times  \alpha_{1,\alpha,1}} (-12, -37); (-8, -34)};
        {\ar@{=>}^{\scriptstyle 1\times (1\times\Pi)} (40,-36); (36, -	32)};
        {\ar@{=>}^{\scriptstyle \Pi} (37, -10); (41, -10)};
\endxy
\]
\[
\def\objectstyle{\scriptstyle}
  \def\labelstyle{\scriptstyle}
   \xy  
    {\ar@{=>}_{\scriptstyle \Pi} (-6, 2); (-6, -3)}; 
   \endxy
\]

\[
\def\objectstyle{\scriptstyle}
  \def\labelstyle{\scriptstyle}
   \xy
   (-22.5,20)*+{((((AB)C)D)E)F}="O1"; %BORDER VERTICES
   (10,20)*+{(((AB)C)D)(EF)}="O2";
   (45,0)*+{(AB)(C(D(EF)))}="O3";
    (-88,0)*+{((A(B(CD))E)F}="O4"; 
   (-88,-20)*+{(A((B(CD))E))F}="O5";   
   (-88, 20)*+{((A((BC)D))E)F}="O6";   
   (-55, 20)*+{(((A(BC))D)E)F}="O7";   
   (45,20)*+{((AB)C)(D(EF))}="O8";
   (45,-20)*+{A(B(C(D(EF))))}="O9"; 
   (-88 ,-40)*+{(A(B((CD)E)))F}="O10";   
   (-55,-40)*+{(A(B(C(DE))))F}="O11";   
   (10,-40)*+{A(B((C(DE))F))}="O12"; 
   (45,-40)*+{A(B(C((DE)F)))}="O13"; 
   (-22.5,-40)*+{A((B(C(DE)))F)}="O14";
    (-32,-30)*+{A((B((CD)E))F)}="I2";  %INNER VERTICES
   (-4,10)*+{((A(BC))D)(EF)}="I4";
   (2,-30)*+{A(B(((CD)E)F))}="I7"; 
   (-39,-12)*+{(A(B(CD)))(EF)}="I8";
   (32,-30)*+{A(B((CD)(EF)))}="I9"; 
   (-50,-20)*+{A(((B(CD))E)F)}="I10"; 
   (10,-20)*+{A((B(CD))(EF))}="I12"; 
   (-45,10)*+{(A((BC)D))(EF)}="I13";
   (30,-5)*+{(A(BC))(D(EF))}="I14";
   (-10,-5)*+{A(((BC)D)(EF))}="I16";
   (32,-14)*+{A((BC)(D(EF)))}="I17";
        {\ar@{-->}_{} "O5";"O10"}; % BORDER ARROWS
        {\ar@{-->}_{} "O10";"O11"};
        {\ar@{-->}_{} "O11";"O14"};
        {\ar@{-->}_{} "O13";"O9"};
        {\ar@{-->}_{} "O14";"O12"};
        {\ar@{-->}_{} "O12";"O13"};
        {\ar@{-->}_{} "O1";"O2"};    
        {\ar@{-->}_{} "O1";"O7"};    
        {\ar@{-->}_{} "O7";"O6"};        
        {\ar@{-->}_{} "O6";"O4"};    
        {\ar@{-->}_{} "O4";"O5"};    
        {\ar@{-->}^{} "O8";"O3"};
        {\ar@{-->}^{} "O2";"O8"};                
        {\ar@{-->}_{} "O3";"O9"};       
        {\ar@{-->}_{} "I9";"O9"}; % IN/OUT ARROWS
        {\ar@{-->}_{} "O2";"I4"};       
        {\ar@{-->}_{} "O7";"I4"}; 
        {\ar@{-->}_{} "I2";"O14"};
        {\ar@{-->}_{} "I7";"O12"};
        {\ar@{-->}_{} "O10";"I2"};
        {\ar@{-->}_{} "O5";"I10"}; 
        {\ar@{-->}_{} "O6";"I13"}; 
        {\ar@{-->}_{} "O4";"I8"}; 
        {\ar@{-->}_{} "O8";"I14"};   
        {\ar@{-->}_{} "I4";"I13"};     % INTERIOR ARROWS  
        {\ar@{-->}_{} "I2";"I7"};
        {\ar@{-->}_{} "I7";"I9"};
        {\ar@{-->}_{} "I10";"I2"};
        {\ar@{-->}_{} "I10";"I12"};
       {\ar@{-->}_{} "I13";"I8"};  
        {\ar@{-->}_{} "I8";"I12"};
        {\ar@{-->}_{} "I12";"I9"};
        {\ar@{-->}_{} "I14";"I17"};
        {\ar@{-->}_{} "I17";"O9"};
        {\ar@{-->}_{} "I4";"I14"};
        {\ar@{-->}_{} "I13";"I16"};
        {\ar@{-->}_{} "I16";"I12"};
        {\ar@{-->}_{} "I16";"I17"};
        {\ar@{=>}^{\scriptstyle \alpha_{\alpha\times 1,1,1}} (-14, 14); (-8, 14)};  %DOMAIN 2-CELLS
        {\ar@{=>}^{\scriptstyle \alpha_{\alpha,1,1}} (25, 9); (29, 10)};  
        {\ar@{=>}^{\scriptstyle 1\times \Pi} (27, -24); (30, -21)};
        {\ar@{=>}^{\scriptstyle 1\times \Pi} (-10, -28); (-7,-25)};
        {\ar@{=>}^{\scriptstyle \alpha_{1\times\alpha,1,1}} (-27, -8); (-24, -5)};
        {\ar@{=>}^{\scriptstyle \Pi} (-1,0); (3, 1)};
        {\ar@{=>}^{\scriptstyle \alpha_{1\times\alpha,1,1}} (-63, 0); (-59, 3)};
        {\ar@{=>}^{\scriptstyle \alpha_{\alpha,1,1}} (-54, 15); (-50, 14)};
        {\ar@{=>}^{\scriptstyle \Pi} (-57, -17); (-53, -14)};
        {\ar@{=>}^{\scriptstyle \alpha_{1,1\times\alpha,1}} (-36, -38); (-33,-35)};
        {\ar@{=>}^{\scriptstyle \alpha_{1,\alpha,1}} (-67, -30); (-63, -26)};
        {\ar@{=>}^{\scriptstyle 1\times (1\times\Pi)} (40,-36); (36, -	32)};
        {\ar@{=>}_{\scriptstyle 1\times  \alpha_{1,\alpha,1}} (-12, -37); (-8, -34)};
        {\ar@{=>}^{\scriptstyle \Pi} (37, -10); (41, -10)};
\endxy
\]

\subsubsection*{$U_{5,2}$ Axiom}

\begin{itemize}
\item for each $5$-tuple of objects $a,b,c,d,e$, an equation called the $U_{5,2}$ associativity condition, consisting of, for each $4$-tuple of morphisms $A,B,C,D$, an equation of $4$-cells
\end{itemize}
\[
\def\objectstyle{\scriptstyle}
  \def\labelstyle{\scriptstyle}
   \xy
   (-52,20)*+{(((A1)B)C)D}="D1";
   (-65,0)*+{((A1)(BC))D}="D2";
   (-52,-10)*+{(A1)((BC)D)}="D3";
   (-39,0)*+{((A1)B)(CD)}="D4";
   (-39,-20)*+{(A1)(B(CD))}="D5";
   (-18,0)*+{(AB)(CD)}="D6";
   (-18,20)*+{((AB)C)D}="D7";
   (-18,-20)*+{A(B(CD))}="D8";
   (-52,-40)*+{A(1((BC)D))}="D9";
   (-18,-40)*+{A(1(B(CD)))}="D10";
   (-88,-20)*+{A(((1B)C)D)}="D15"; 
   (-88,-40)*+{A((1(BC))D)}="D19";   
   (-88, 0)*+{(A((1B)C))D}="D16";   
   (-88, 20)*+{((A(1B))C)D}="D17";   
   (-65,-20)*+{(A(1(BC)))D}="D18";   
        {\ar_{} "D1";"D2"}; %DOMAIN ARROWS
        {\ar_{} "D2";"D3"};
        {\ar_{} "D3";"D5"};
        {\ar_{} "D1";"D4"};
        {\ar_{} "D4";"D6"};
        {\ar@{-->}_{} "D5";"D10"};
        {\ar^{} "D7";"D6"};
        {\ar^{} "D1";"D7"};        
        {\ar_{} "D4";"D5"}; 
        {\ar@{-->}_{} "D3";"D9"}; 
        {\ar_{} "D5";"D8"};           
        {\ar@{-->}_{} "D9";"D10"};           
        {\ar@{-->}_{} "D1";"D17"};    
        {\ar@{-->}_{} "D17";"D16"};    
        {\ar@{-->}_{} "D18";"D19"};    
        {\ar@{-->}_{} "D16";"D18"};    
        {\ar@{-->}_{} "D16";"D15"};    
        {\ar@{-->}_{} "D15";"D19"};    
        {\ar@{-->}_{} "D19";"D9"};    
        {\ar@{-->}_{} "D10";"D8"};   
        {\ar_{} "D6";"D8"};   
        {\ar@{-->}_{} "D2";"D18"};     
        {\ar@{=>}^{\scriptstyle \Pi\times 1} (-78, 7); (-73, 4)};
        {\ar@{=>}^{\scriptstyle \Pi} (-54.5, -2); (-48.5, -2)};
        {\ar@{=>}_<<{\scriptstyle \alpha^{-1}_{\rho\times 1,1,1}} (-34, 10); (-30, 12)};
        {\ar@{=>}_<<{\scriptstyle \alpha^{-1}_{\rho,1,1}} (-31, -12); (-27, -7)};
        {\ar@{=>}_<<{\scriptstyle m} (-25, -30); (-26, -26)};
        {\ar@{=>}_<<{\scriptstyle \alpha^{-1}_{1,1,\alpha}} (-46, -31); (-42, -27)};
        {\ar@{=>}^{\scriptstyle \Pi} (-65, -27); (-61, -30)};
        {\ar@{=>}_<<{\scriptstyle \alpha^{-1}_{1,\alpha,1}} (-82, -23); (-79, -23)};
\endxy
\]
\[
\def\objectstyle{\scriptstyle}
  \def\labelstyle{\scriptstyle}
   \xy  
        {\ar@{=>}^{\scriptstyle \Pi} (0, 2); (0, -2)}; %PERTURBATION
   \endxy
\]
\[
\def\objectstyle{\scriptstyle}
  \def\labelstyle{\scriptstyle}
   \xy  
   (35,20)*+{(((A1)B)C)D}="C1";
   (21,0)*+{((A1)(BC))D}="C2";
   (35,-10)*+{(A1)((BC)D)}="C3";
   (49.5,0)*+{(A(BC))D}="C4";
   (59.5,-10)*+{A((BC)D)}="C5";
   (71,0)*+{(AB)(CD)}="C6";
   (71,20)*+{((AB)C)D}="C7";
   (71,-20)*+{A(B(CD))}="C8";
   (35,-40)*+{A(1((BC)D))}="C9";
   (71,-40)*+{A(1(B(CD)))}="C10";
   (-1,-20)*+{A(((1B)C)D)}="C15"; 
   (-1,-40)*+{A((1(BC))D)}="C19";   
   (-1, 0)*+{(A((1B)C))D}="C16";   
   (-1, 20)*+{((A(1B))C)D}="C17";   
   (21, -20)*+{(A(1(BC)))D}="C18"; 
   (49.5,-20)*+{(A1)(B(CD))}="C20";  
        {\ar@{-->}_{} "C1";"C2"}; %DOMAIN ARROWS
        {\ar@{-->}_{} "C2";"C3"};
        {\ar_{} "C3";"C5"};
        {\ar@{-->}_{} "C7";"C4"};
        {\ar@{-->}_{} "C2";"C4"};
        {\ar@{-->}^{} "C7";"C6"};
        {\ar@{-->}^{} "C1";"C7"};        
        {\ar@{-->}_{} "C4";"C5"}; 
        {\ar_{} "C3";"C9"}; 
        {\ar_{} "C5";"C8"};           
        {\ar_{} "C9";"C10"};           
        {\ar@{-->}_{} "C1";"C17"};    
        {\ar@{-->}_{} "C17";"C16"};    
        {\ar@{-->}_{} "C18";"C19"};    
        {\ar@{-->}_{} "C16";"C18"};    
        {\ar@{-->}_{} "C16";"C15"};    
        {\ar@{-->}_{} "C15";"C19"};    
        {\ar@{-->}_{} "C19";"C9"};    
        {\ar_{} "C10";"C8"};   
        {\ar_{} "C6";"C8"};   
        {\ar@{-->}_{} "C2";"C18"}; 
        {\ar_{} "C3";"C20"}; 
        {\ar_{} "C20";"C8"}; 
        {\ar_{} "C20";"C10"};     
        {\ar@{=>}^{\scriptstyle \Pi\times 1} (7, 6); (11, 3)};
        {\ar@{=>}_<<{\scriptstyle \alpha^{-1}_{\rho,1,1}\times 1} (39, 10); (43, 12)};
        {\ar@{=>}_{\scriptstyle \alpha^{-1}_{\rho,1,1}} (39, -5); (43, -2)};
        {\ar@{=>}_<<{\scriptstyle \rho\widetilde{\times} \alpha^{-1}} (51, -15); (55, -13)};
        {\ar@{=>}_<<{\scriptstyle \Pi} (59, -3); (64, -3)};
        {\ar@{=>}_<<{\scriptstyle m} (65.5, -30); (63.5, -26)};
        {\ar@{=>}_<<{\scriptstyle \alpha^{-1}_{1,1,\alpha}} (43, -31); (47, -27)};
        {\ar@{=>}^{\scriptstyle \Pi} (21, -27); (25, -30)};
        {\ar@{=>}_<<{\scriptstyle \alpha^{-1}_{1,\alpha,1}} (6, -23); (11, -23)};
\endxy
\]
\[
\def\objectstyle{\scriptstyle}
  \def\labelstyle{\scriptstyle}
   \xy  
    {\ar@{=>}_{\scriptstyle m} (0, 2); (0, -2)}; 
   \endxy
\]
\[
\def\objectstyle{\scriptstyle}
  \def\labelstyle{\scriptstyle}
   \xy  
   (-52,20)*+{(((A1)B)C)D}="D1";
   (-65,0)*+{((A1)(BC))D}="D2";
   (-52,-10)*+{(A1)((BC)D)}="D3";
   (-39.5,0)*+{(A(BC))D}="D4";
   (-29.5,-10)*+{A((BC)D)}="D5";
   (-18,0)*+{(AB)(CD)}="D6";
   (-18,20)*+{((AB)C)D}="D7";
   (-18,-20)*+{A(B(CD))}="D8";
   (-52,-40)*+{A(1((BC)D))}="D9";
   (-18,-40)*+{A(1(B(CD)))}="D10";
   (-88,-20)*+{A(((1B)C)D)}="D15"; 
   (-88,-40)*+{A((1(BC))D)}="D19";   
   (-88, 0)*+{(A((1B)C))D}="D16";   
   (-88, 20)*+{((A(1B))C)D}="D17";   
   (-65,-20)*+{(A(1(BC)))D}="D18"; 
        {\ar@{-->}_{} "D1";"D2"}; %DOMAIN ARROWS
        {\ar_{} "D2";"D3"};
        {\ar_{} "D3";"D5"};
        {\ar@{-->}_{} "D7";"D4"};
        {\ar_{} "D2";"D4"};
        {\ar@{-->}^{} "D7";"D6"};
        {\ar@{-->}^{} "D1";"D7"};        
        {\ar_{} "D4";"D5"}; 
        {\ar_{} "D3";"D9"}; 
        {\ar@{-->}_{} "D5";"D8"};           
        {\ar@{-->}_{} "D9";"D10"};           
        {\ar@{-->}_{} "D1";"D17"};    
        {\ar@{-->}_{} "D17";"D16"};    
        {\ar_{} "D18";"D19"};    
        {\ar@{-->}_{} "D16";"D18"};    
        {\ar@{-->}_{} "D16";"D15"};    
        {\ar@{-->}_{} "D15";"D19"};    
        {\ar_{} "D19";"D9"};    
        {\ar@{-->}_{} "D10";"D8"};   
        {\ar@{-->}_{} "D6";"D8"};   
        {\ar_{} "D2";"D18"}; 
        {\ar_{} "D9";"D5"}; 
        {\ar@{=>}^{\scriptstyle \Pi\times 1} (-78, 7); (-73, 4)};
        {\ar@{=>}_<<{\scriptstyle \alpha^{-1}_{\rho,1,1}\times 1} (-48, 10); (-44, 12)};
        {\ar@{=>}^<<{\scriptstyle \alpha^{-1}_{\rho,1,1}} (-47, -7); (-43, -4)};
        {\ar@{=>}_<<{\scriptstyle \Pi} (-28, -3); (-23, -3)};
        {\ar@{=>}^<<{\scriptstyle 1\times \lambda_{\alpha}} (-24, -30); (-27, -26)};
        {\ar@{=>}_<<{\scriptstyle m} (-46, -22); (-43, -18)};
        {\ar@{=>}^{\scriptstyle \Pi} (-65, -27); (-61, -30)};
        {\ar@{=>}_<<{\scriptstyle \alpha^{-1}_{1,\alpha,1}} (-82, -23); (-79, -23)};
\endxy
\]
\[
\def\objectstyle{\scriptstyle}
  \def\labelstyle{\scriptstyle}
   \xy  
        {\ar@{=>}^{\scriptstyle U_{4,2}} (0, 2); (0, -2)}; %PERTURBATION
   \endxy
\]

\[
\def\objectstyle{\scriptstyle}
  \def\labelstyle{\scriptstyle}
   \xy  
   (35,20)*+{(((A1)B)C)D}="C1";
   (21,0)*+{((A1)(BC))D}="C2";
   (49.5,0)*+{(A(BC))D}="C4";
   (59.5,-10)*+{A((BC)D)}="C5";
   (71,0)*+{(AB)(CD)}="C6";
   (71,20)*+{((AB)C)D}="C7";
   (71,-20)*+{A(B(CD))}="C8";
   (35,-40)*+{A(1((BC)D))}="C9";
   (71,-40)*+{A(1(B(CD)))}="C10";
   (-1,-20)*+{A(((1B)C)D)}="C15"; 
   (-1,-40)*+{A((1(BC))D)}="C19";   
   (-1, 0)*+{(A((1B)C))D}="C16";   
   (-1, 20)*+{((A(1B))C)D}="C17";   
   (21, -20)*+{(A(1(BC)))D}="C18"; 
        {\ar_{} "C1";"C2"}; %DOMAIN ARROWS
        {\ar_{} "C18";"C4"};
        {\ar_{} "C7";"C4"};
        {\ar_{} "C2";"C4"};
        {\ar@{-->}^{} "C7";"C6"};
        {\ar^{} "C1";"C7"};        
        {\ar@{-->}_{} "C4";"C5"}; 
        {\ar@{-->}_{} "C5";"C8"};           
        {\ar@{-->}_{} "C9";"C10"};           
        {\ar_{} "C1";"C17"};    
        {\ar_{} "C17";"C16"};    
        {\ar@{-->}_{} "C18";"C19"};    
        {\ar_{} "C16";"C18"};    
        {\ar@{-->}_{} "C16";"C15"};    
        {\ar@{-->}_{} "C15";"C19"};    
        {\ar@{-->}_{} "C19";"C9"};    
        {\ar@{-->}_{} "C10";"C8"};   
        {\ar@{-->}_{} "C6";"C8"};   
        {\ar_{} "C2";"C18"}; 
        {\ar@{-->}_{} "C9";"C5"}; 
        {\ar@{-->}_{} "C19";"C5"}; 
        {\ar@{=>}^{\scriptstyle \Pi\times 1} (7, 6); (11, 3)};
        {\ar@{=>}_<<{\scriptstyle \alpha^{-1}_{\rho,1,1}\times 1} (39, 10); (43, 12)};
        {\ar@{=>}^<<{\scriptstyle m\times 1} (32, -8); (35, -4)};
        {\ar@{=>}_<<{\scriptstyle \Pi} (59, -3); (64, -3)};
        {\ar@{=>}^<<{\scriptstyle 1\times \lambda_{\alpha}} (64, -30); (61, -26)};
        {\ar@{=>}_<<{\scriptstyle 1\times l} (33, -32); (31, -29)};
        {\ar@{=>}^{\scriptstyle \alpha^{-1}_{1,\lambda,1}} (39, -18); (43, -16)};
        {\ar@{=>}_<<{\scriptstyle \alpha^{-1}_{1,\alpha,1}} (6, -23); (11, -23)};
\endxy
\]
\[
\def\objectstyle{\scriptstyle}
  \def\labelstyle{\scriptstyle}
   \xy  
    {\ar@{=>}_{\scriptstyle U_{4,2}\times 1} (0, 2); (0, -2)}; 
   \endxy
\]
\[
\def\objectstyle{\scriptstyle}
  \def\labelstyle{\scriptstyle}
   \xy  
   (-52,20)*+{(((A1)B)C)D}="C1";
   (-39.5,0)*+{(A(BC))D}="C4";
   (-29.5,-10)*+{A((BC)D)}="C5";
   (-18,0)*+{(AB)(CD)}="C6";
   (-18,20)*+{((AB)C)D}="C7";
   (-18,-20)*+{A(B(CD))}="C8";
   (-52,-40)*+{A(1((BC)D))}="C9";
   (-18,-40)*+{A(1(B(CD)))}="C10";
   (-88,-20)*+{A(((1B)C)D)}="C15"; 
   (-88,-40)*+{A((1(BC))D)}="C19";   
   (-88, 0)*+{(A((1B)C))D}="C16";   
   (-88, 20)*+{((A(1B))C)D}="C17";   
   (-65, -20)*+{(A(1(BC)))D}="C18"; 
        {\ar@{-->}@/_2pc/_{} "C17";"C7"}; %DOMAIN ARROWS
        {\ar_{} "C18";"C4"};
        {\ar_{} "C16";"C4"};
        {\ar@{-->}_{} "C7";"C4"};
        {\ar@{-->}^{} "C7";"C6"};
        {\ar@{-->}^{} "C1";"C7"};        
        {\ar_{} "C4";"C5"}; 
        {\ar@{-->}_{} "C5";"C8"};           
        {\ar@{-->}_{} "C9";"C10"};           
        {\ar@{-->}_{} "C1";"C17"};    
        {\ar@{-->}_{} "C17";"C16"};    
        {\ar_{} "C18";"C19"};    
        {\ar_{} "C16";"C18"};    
        {\ar_{} "C16";"C15"};    
        {\ar_{} "C15";"C19"};    
        {\ar@{-->}_{} "C19";"C9"};    
        {\ar@{-->}_{} "C10";"C8"};   
        {\ar@{-->}_{} "C6";"C8"};   
        {\ar@{-->}_{} "C9";"C5"}; 
        {\ar_{} "C19";"C5"};   
        {\ar@{=>}^{\scriptstyle \alpha^{-1}_{1,\lambda,1}\times 1} (-52, 2); (-47, 4)};
        {\ar@{=>}^<<{\scriptstyle (m\times 1)\times 1} (-50, 13); (-47, 14)};
        {\ar@{=>}^<<{\scriptstyle (1\times l)\times 1} (-63, -10); (-63, -4)};
        {\ar@{=>}_<<{\scriptstyle \Pi} (-28, -3); (-23, -3)};
        {\ar@{=>}^<<{\scriptstyle 1\times \lambda_{\alpha}} (-24, -30); (-27, -26)};
        {\ar@{=>}_<<{\scriptstyle 1\times l} (-55, -31); (-57, -28)};
        {\ar@{=>}^{\scriptstyle \alpha^{-1}_{1,\lambda,1}} (-51.5, -19.5); (-47.5, -17.5)};
        {\ar@{=>}_<<{\scriptstyle \alpha^{-1}_{1,\alpha,1}} (-82, -23); (-79, -23)};
\endxy
\]
\[
\def\objectstyle{\scriptstyle}
  \def\labelstyle{\scriptstyle}
   \xy  
        {\ar@{=>}^{\scriptstyle \alpha} (0, 2); (0, -2)}; %PERTURBATION
   \endxy
\]        
\[
\def\objectstyle{\scriptstyle}
  \def\labelstyle{\scriptstyle}
   \xy  
   (35,20)*+{(((A1)B)C)D}="D1";
   (49.5,0)*+{(A(BC))D}="D4";
   (59.5,-10)*+{A((BC)D)}="D5";
   (71,0)*+{(AB)(CD)}="D6";
   (71,20)*+{((AB)C)D}="D7";
   (71,-20)*+{A(B(CD))}="D8";
   (35,-40)*+{A(1((BC)D))}="D9";
   (71,-40)*+{A(1(B(CD)))}="D10";
   (-1,-20)*+{A(((1B)C)D)}="D15"; 
   (-1,-40)*+{A((1(BC))D)}="D19";   
   (-1, 0)*+{(A((1B)C))D}="D16";   
   (-1, 20)*+{((A(1B))C)D}="D17";      
        {\ar@{-->}@/_2pc/_{} "D17";"D7"}; %DOMAIN ARROWS
        {\ar@{-->}_{} "D16";"D4"};
        {\ar@{-->}_{} "D7";"D4"};
        {\ar@{-->}^{} "D7";"D6"};
        {\ar@{-->}^{} "D1";"D7"};        
        {\ar@{-->}_{} "D4";"D5"}; 
        {\ar@{-->}_{} "D5";"D8"};           
        {\ar@{-->}_{} "D9";"D10"};           
        {\ar@{-->}_{} "D1";"D17"};    
        {\ar@{-->}_{} "D17";"D16"};    
        {\ar@{-->}_{} "D16";"D15"};    
        {\ar@{-->}_{} "D15";"D19"};    
        {\ar@{-->}_{} "D19";"D9"};    
        {\ar@{-->}_{} "D10";"D8"};   
        {\ar@{-->}_{} "D6";"D8"};   
        {\ar@{-->}_{} "D15";"D5"}; 
        {\ar@{-->}_{} "D9";"D5"}; 
        {\ar@{-->}_{} "D19";"D5"};   
        {\ar@{=>}^{\scriptstyle \alpha^{-1}_{1,\lambda,1}\times 1} (30, 3); (35, 5)};
        {\ar@{=>}^<<{\scriptstyle (m\times 1)\times 1} (38, 13); (41, 14)};
        {\ar@{=>}_<<{\scriptstyle \Pi} (59, -3); (64, -3)};
        {\ar@{=>}^<<{\scriptstyle 1\times \lambda_{\alpha}} (64, -30); (61, -26)};
        {\ar@{=>}_<<{\scriptstyle 1\times l} (34, -31); (32, -28)};
        {\ar@{=>}^{\scriptstyle \alpha^{-1}_{1,\lambda\times 1,1}} (32, -11); (36, -8)};
        {\ar@{=>}_{\scriptstyle 1\times(l\times 1)} (10, -25); (14, -22)};
\endxy
\]
\[\textrm{\Huge{=}}\]
\[
\def\objectstyle{\scriptstyle}
  \def\labelstyle{\scriptstyle}
   \xy
   (-52,20)*+{(((A1)B)C)D}="D1";
   (-65,0)*+{((A1)(BC))D}="D2";
   (-52,-10)*+{(A1)((BC)D)}="D3";
   (-39.5,0)*+{((A1)B)(CD)}="D4";
   (-39.5,-20)*+{(A1)(B(CD))}="D5";
   (-18,0)*+{(AB)(CD)}="D6";
   (-18,20)*+{((AB)C)D}="D7";
   (-18,-20)*+{A(B(CD))}="D8";
   (-52,-40)*+{A(1((BC)D))}="D9";
   (-18,-40)*+{A(1(B(CD)))}="D10";
   (-88,-20)*+{A(((1B)C)D)}="D15"; 
   (-88,-40)*+{A((1(BC))D)}="D19";   
   (-88, 0)*+{(A((1B)C))D}="D16";   
   (-88, 20)*+{((A(1B))C)D}="D17";   
   (-65,-20)*+{(A(1(BC)))D}="D18";   
        {\ar_{} "D1";"D2"}; %DOMAIN ARROWS
        {\ar_{} "D2";"D3"};
        {\ar_{} "D3";"D5"};
        {\ar_{} "D1";"D4"};
        {\ar@{-->}_{} "D4";"D6"};
        {\ar_{} "D5";"D10"};
        {\ar@{-->}^{} "D7";"D6"};
        {\ar@{-->}^{} "D1";"D7"};        
        {\ar_{} "D4";"D5"}; 
        {\ar_{} "D3";"D9"}; 
        {\ar@{-->}_{} "D5";"D8"};           
        {\ar_{} "D9";"D10"};           
        {\ar_{} "D1";"D17"};    
        {\ar_{} "D17";"D16"};    
        {\ar_{} "D18";"D19"};    
        {\ar_{} "D16";"D18"};    
        {\ar_{} "D16";"D15"};    
        {\ar_{} "D15";"D19"};    
        {\ar_{} "D19";"D9"};    
        {\ar@{-->}_{} "D10";"D8"};   
        {\ar@{-->}_{} "D6";"D8"};   
        {\ar_{} "D2";"D18"};   
        {\ar@{=>}^{\scriptstyle \Pi\times 1} (-78, 7); (-73, 4)};
        {\ar@{=>}^{\scriptstyle \Pi} (-54.5, -2); (-48.5, -2)};
        {\ar@{=>}_<<{\scriptstyle \alpha^{-1}_{\rho,1,1}} (-31, -12); (-27, -7)};
        {\ar@{=>}_<<{\scriptstyle m} (-25, -30); (-26, -26)};
        {\ar@{=>}_<<{\scriptstyle \alpha^{-1}_{1,1,\alpha}} (-46, -31); (-42, -27)};
        {\ar@{=>}^{\scriptstyle \Pi} (-65, -27); (-61, -30)};
        {\ar@{=>}_<<{\scriptstyle \alpha^{-1}_{1,\alpha,1}} (-82, -23); (-79, -23)};
        {\ar@{=>}_<<{\scriptstyle \alpha^{-1}_{\rho\times 1,1,1}} (-34, 10); (-30, 12)};
\endxy
\]
\[
\def\objectstyle{\scriptstyle}
  \def\labelstyle{\scriptstyle}
   \xy  
                {\ar@{=>}^{\scriptstyle K_5} (0, 2); (0, -2)}; %PERTURBATION
   \endxy
\]
\[
\def\objectstyle{\scriptstyle}
  \def\labelstyle{\scriptstyle}
   \xy  
   (35,20)*+{(((A1)B)C)D}="C1";
   (21,-20)*+{(A(1B))(CD)}="C2";
   (50,0)*+{((A1)B)(CD)}="C4";
   (71,0)*+{(AB)(CD)}="C6";
   (71,20)*+{((AB)C)D}="C7";
   (71,-20)*+{A(B(CD))}="C8";
   (35,-40)*+{A(1((BC)D))}="C9";
   (71,-40)*+{A(1(B(CD)))}="C10";
   (-1,-20)*+{A(((1B)C)D)}="C15"; 
   (-1,-40)*+{A((1(BC))D)}="C19";   
   (-1, 0)*+{(A((1B)C))D}="C16";   
   (-1, 20)*+{((A(1B))C)D}="C17";   
   (21, -30)*+{A((1B)(CD))}="C18"; 
   (50,-20)*+{(A1)(B(CD))}="C20"; 
        {\ar@{-->}_{} "C1";"C4"}; %DOMAIN ARROWS
        {\ar_{} "C4";"C6"};
        {\ar_{} "C4";"C2"};
        {\ar@{-->}^{} "C7";"C6"};
        {\ar@{-->}^{} "C1";"C7"};               
        {\ar@{-->}_{} "C9";"C10"};           
        {\ar@{-->}_{} "C1";"C17"};    
        {\ar@{-->}_{} "C17";"C16"};    
        {\ar@{-->}_{} "C17";"C2"};   
        {\ar@{-->}_{} "C16";"C15"};    
        {\ar@{-->}_{} "C15";"C19"};    
        {\ar@{-->}_{} "C19";"C9"};    
        {\ar_{} "C10";"C8"};   
        {\ar_{} "C6";"C8"};   
        {\ar_{} "C2";"C18"}; 
        {\ar_{} "C4";"C20"}; 
        {\ar_{} "C20";"C8"}; 
        {\ar_{} "C18";"C10"};     
        {\ar_{} "C20";"C10"};  
        {\ar@{=>}_<<{\scriptstyle \alpha^{-1}_{\rho\times 1,1,1}} (55, 10); (59, 12)};
        {\ar@{=>}^{\scriptstyle \alpha^{-1}_{\alpha,1,1}} (22, 6); (26, 3)};
        {\ar@{=>}_{\scriptstyle \alpha^{-1}_{\rho,1,1}} (58, -12); (62, -7)};
        {\ar@{=>}_<<{\scriptstyle m} (65, -30); (63, -26)};
        {\ar@{=>}_<<{\scriptstyle \Pi} (35, -22); (39, -22)};
        {\ar@{=>}^{\scriptstyle 1\times \Pi} (16, -37); (20, -34)};
        {\ar@{=>}_<<{\scriptstyle \Pi} (5, -10); (9, -13)};
\endxy
\]
\[
\def\objectstyle{\scriptstyle}
  \def\labelstyle{\scriptstyle}
   \xy  
    {\ar@{=>}_{\scriptstyle U_{4,2}} (0, 2); (0, -2)}; 
   \endxy
\]
\[
\def\objectstyle{\scriptstyle}
  \def\labelstyle{\scriptstyle}
   \xy  
   (-52,20)*+{(((A1)B)C)D}="D1";
   (-65,-20)*+{(A(1B))(CD)}="D3";
   (-39,0)*+{((A1)B)(CD)}="D4";
   (-18,0)*+{(AB)(CD)}="D6";
   (-18,20)*+{((AB)C)D}="D7";
   (-18,-20)*+{A(B(CD))}="D8";
   (-52,-40)*+{A(1((BC)D))}="D9";
   (-18,-40)*+{A(1(B(CD)))}="D10";
   (-88,-20)*+{A(((1B)C)D)}="D15"; 
   (-88,-40)*+{A((1(BC))D)}="D19";   
   (-88, 0)*+{(A((1B)C))D}="D16";   
   (-88, 20)*+{((A(1B))C)D}="D17";   
   (-65, -30)*+{A((1B)(CD))}="D18"; 
        {\ar_{} "D3";"D6"};
        {\ar_{} "D4";"D6"};
        {\ar^{} "D7";"D6"};
        {\ar^{} "D1";"D7"};        
        {\ar_{} "D1";"D4"}; 
        {\ar@{-->}_{} "D3";"D18"}; 
        {\ar@{-->}_{} "D18";"D8"};           
        {\ar@{-->}_{} "D9";"D10"};           
        {\ar_{} "D1";"D17"};    
        {\ar@{-->}_{} "D17";"D16"};    
        {\ar@{-->}_{} "D15";"D18"};    
        {\ar_{} "D17";"D3"};    
        {\ar@{-->}_{} "D16";"D15"};    
        {\ar@{-->}_{} "D15";"D19"};    
        {\ar@{-->}_{} "D19";"D9"};    
        {\ar@{-->}_{} "D10";"D8"};   
        {\ar@{-->}_{} "D6";"D8"};   
        {\ar_{} "D4";"D3"}; 
        {\ar@{-->}_{} "D18";"D10"}; 
        {\ar@{=>}^{\scriptstyle \Pi } (-82, -10); (-78, -13)};
        {\ar@{=>}_<<{\scriptstyle \alpha^{-1}_{\rho\times 1,1,1}} (-37, 10); (-34, 12)};
        {\ar@{=>}^<<{\scriptstyle \alpha^{-1}_{\alpha,1,1}} (-65, 6); (-61, 3)};
        {\ar@{=>}^<<{\scriptstyle m\times (1\times 1)} (-40, -7); (-36, -6)};
        {\ar@{=>}^<<{\scriptstyle 1\times l} (-28, -31); (-29, -27)};
        {\ar@{=>}^{\scriptstyle 1\times\Pi} (-75, -37); (-71, -34)};
        {\ar@{=>}^<<{\scriptstyle \alpha^{-1}_{1,\lambda,1}} (-41, -21); (-38, -19)};
\endxy
\]
\[
\def\objectstyle{\scriptstyle}
  \def\labelstyle{\scriptstyle}
   \xy  
        {\ar@{=>}^{\scriptstyle \alpha} (0, 2); (0, -2)}; %PERTURBATION
   \endxy
\]

\[
\def\objectstyle{\scriptstyle}
  \def\labelstyle{\scriptstyle}
   \xy  
   (35,20)*+{(((A1)B)C)D}="C1";
   (21,-20)*+{(A(1B))(CD)}="C2";
   (71,0)*+{(AB)(CD)}="C6";
   (71,20)*+{((AB)C)D}="C7";
   (71,-20)*+{A(B(CD))}="C8";
   (35,-40)*+{A(1((BC)D))}="C9";
   (71,-40)*+{A(1(B(CD)))}="C10";
   (-1,-20)*+{A(((1B)C)D)}="C15"; 
   (-1,-40)*+{A((1(BC))D)}="C19";   
   (-1, 0)*+{(A((1B)C))D}="C16";   
   (-1, 20)*+{((A(1B))C)D}="C17";   
   (21, -30)*+{A((1B)(CD))}="C18"; 
        {\ar_{} "C17";"C2"}; %DOMAIN ARROWS
        {\ar^{} "C7";"C6"};
        {\ar@{-->}^{} "C1";"C7"};           
        {\ar@{-->}_{} "C9";"C10"};           
        {\ar@{-->}_{} "C1";"C17"};    
        {\ar_{} "C17";"C16"};    
        {\ar@{-->}_{} "C18";"C10"};    
        {\ar_{} "C15";"C18"};    
        {\ar_{} "C16";"C15"};    
        {\ar@{-->}_{} "C15";"C19"};    
        {\ar@{-->}_{} "C19";"C9"};    
        {\ar@{-->}_{} "C10";"C8"};   
        {\ar_{} "C6";"C8"};   
        {\ar_{} "C2";"C18"}; 
        {\ar@/_2pc/_{} "C17";"C7"};
        {\ar_{} "C2";"C6"}; 
        {\ar_{} "C18";"C8"}; 
        {\ar@{=>}^{\scriptstyle \alpha^{-1}_{1,\lambda,1}} (48, -21); (51, -19)};
        {\ar@{=>}^{\scriptstyle \Pi} (5, -10); (9, -13)};
        {\ar@{=>}_<<{\scriptstyle (m\times 1)\times 1} (33, 15); (37, 17)};
        {\ar@{=>}^<<{\scriptstyle \alpha^{-1}_{1\times\lambda,1,1}} (38, -3); (43, 0)};
        {\ar@{=>}^<<{\scriptstyle 1\times l} (61, -31); (59, -27)};
        {\ar@{=>}^{\scriptstyle \alpha^{-1}_{1,\lambda,1}} (48, -21); (51, -19)};
        {\ar@{=>}^<<{\scriptstyle 1\times\Pi} (16, -37); (20, -34)};
\endxy
\]
\[
\def\objectstyle{\scriptstyle}
  \def\labelstyle{\scriptstyle}
   \xy  
    {\ar@{=>}_{\scriptstyle \Pi} (0, 2); (0, -2)}; 
   \endxy
\]
\[
\def\objectstyle{\scriptstyle}
  \def\labelstyle{\scriptstyle}
   \xy  
   (-52,20)*+{(((A1)B)C)D}="C1";
   (-39.5,0)*+{(A(BC))D}="C4";
   (-29.5,-10)*+{A((BC)D)}="C5";
   (-18,0)*+{(AB)(CD)}="C6";
   (-18,20)*+{((AB)C)D}="C7";
   (-18,-20)*+{A(B(CD))}="C8";
   (-52,-40)*+{A(1((BC)D))}="C9";
   (-18,-40)*+{A(1(B(CD)))}="C10";
   (-88,-20)*+{A(((1B)C)D)}="C15"; 
   (-88,-40)*+{A((1(BC))D)}="C19";   
   (-88, 0)*+{(A((1B)C))D}="C16";   
   (-88, 20)*+{((A(1B))C)D}="C17";   
   (-65, -30)*+{A((1B)(CD))}="C18"; 
        {\ar@{-->}@/_2pc/_{} "C17";"C7"}; %DOMAIN ARROWS
        {\ar@{-->}_{} "C16";"C4"};
        {\ar@{-->}_{} "C7";"C4"};
        {\ar@{-->}^{} "C7";"C6"};
        {\ar@{-->}^{} "C1";"C7"};        
        {\ar@{-->}_{} "C4";"C5"}; 
        {\ar_{} "C5";"C8"};           
        {\ar_{} "C9";"C10"};           
        {\ar@{-->}_{} "C1";"C17"};    
        {\ar@{-->}_{} "C17";"C16"};    
        {\ar_{} "C15";"C18"};    
        {\ar_{} "C18";"C8"};    
        {\ar@{-->}_{} "C16";"C15"};    
        {\ar_{} "C15";"C19"};    
        {\ar_{} "C19";"C9"};    
        {\ar_{} "C10";"C8"};   
        {\ar@{-->}_{} "C6";"C8"};   
        {\ar_{} "C18";"C10"}; 
        {\ar_{} "C15";"C5"}; 
        {\ar@{=>}^{\scriptstyle \alpha^{-1}_{1,\lambda,1}\times 1} (-52, 4); (-48, 5)};
        {\ar@{=>}^<<{\scriptstyle \alpha^{-1}_{1,\lambda\times 1,1}}  (-52, -10); (-48, -8)};
        {\ar@{=>}^<<{\scriptstyle 1\times \lambda} (-29, -32); (-31, -28)};
        {\ar@{=>}^<<{\scriptstyle 1\times \Pi} (-75, -37); (-71, -34)};
        {\ar@{=>}^{\scriptstyle 1\times \alpha^{-1}_{\lambda,1,1}} (-52, -23); (-48, -21)};
        {\ar@{=>}^<<{\scriptstyle (m\times 1)\times 1}  (-52, 13); (-48, 14)};
        {\ar@{=>}_<<{\scriptstyle \Pi} (-30, -1); (-26, -1)};
\endxy
\]
\[
\def\objectstyle{\scriptstyle}
  \def\labelstyle{\scriptstyle}
   \xy  
        {\ar@{=>}^{\scriptstyle 1\times U_{4,1}} (0, 2); (0, -2)}; %PERTURBATION
\endxy
\]

\[
\def\objectstyle{\scriptstyle}
  \def\labelstyle{\scriptstyle}
   \xy  
   (35,20)*+{(((A1)B)C)D}="D1";
   (49.5,0)*+{(A(BC))D}="D4";
   (59.5,-10)*+{A((BC)D)}="D5";
   (71,0)*+{(AB)(CD)}="D6";
   (71,20)*+{((AB)C)D}="D7";
   (71,-20)*+{A(B(CD))}="D8";
   (35,-40)*+{A(1((BC)D))}="D9";
   (71,-40)*+{A(1(B(CD)))}="D10";
   (-1,-20)*+{A(((1B)C)D)}="D15"; 
   (-1,-40)*+{A((1(BC))D)}="D19";   
   (-1, 0)*+{(A((1B)C))D}="D16";   
   (-1, 20)*+{((A(1B))C)D}="D17";   
        {\ar@{-->}@/_2pc/_{} "D17";"D7"}; %DOMAIN ARROWS
        {\ar@{-->}_{} "D16";"D4"};
        {\ar@{-->}_{} "D7";"D4"};
        {\ar@{-->}^{} "D7";"D6"};
        {\ar@{-->}^{} "D1";"D7"};        
        {\ar@{-->}_{} "D4";"D5"}; 
        {\ar@{-->}_{} "D5";"D8"};           
        {\ar@{-->}_{} "D9";"D10"};           
        {\ar@{-->}_{} "D1";"D17"};    
        {\ar@{-->}_{} "D17";"D16"};       
        {\ar@{-->}_{} "D16";"D15"};    
        {\ar@{-->}_{} "D15";"D19"};    
        {\ar@{-->}_{} "D19";"D9"};    
        {\ar@{-->}_{} "D10";"D8"};   
        {\ar@{-->}_{} "D6";"D8"};   
        {\ar@{-->}_{} "D15";"D5"}; 
        {\ar@{-->}_{} "D9";"D5"}; 
        {\ar@{-->}_{} "D19";"D5"};   
        {\ar@{=>}^<<{\scriptstyle (m\times 1)\times 1} (40, 13); (43, 14)};
        {\ar@{=>}_<<{\scriptstyle \Pi} (59, -1); (64, -1)};
        {\ar@{=>}^{\scriptstyle \alpha^{-1}_{1,\lambda,1}\times 1} (39, 3); (43, 5)};
        {\ar@{=>}^<<{\scriptstyle 1\times \lambda_{\alpha}} (64, -30); (61, -26)};
        {\ar@{=>}_<<{\scriptstyle 1\times \lambda} (34, -31); (32, -28)};
        {\ar@{=>}^{\scriptstyle \alpha^{-1}_{1,\lambda\times 1,1}} (39, -11); (43, -8)};
        {\ar@{=>}_{\scriptstyle 1\times(\lambda\times 1)} (14, -23); (17, -20)};
\endxy
\]

\subsubsection*{$U_{5,3}$ Axiom}
\begin{itemize}
\item for each $5$-tuple of objects $a,b,c,d,e$, an equation called the $U_{5,3}$ unit condition, consisting of, for each $4$-tuple of morphisms $A,B,C,D$, an equation of $4$-cells
\end{itemize}
\[
\def\objectstyle{\scriptstyle}
  \def\labelstyle{\scriptstyle}
   \xy
   (-52,20)*+{(((AB)1)C)D}="D1";
   (-65,0)*+{((AB)(1C))D}="D2";
   (-52,-10)*+{(AB)((1C)D)}="D3";
   (-39,0)*+{((AB)1)(CD)}="D4";
   (-39,-20)*+{(AB)(1(CD))}="D5";
   (-18,0)*+{(AB)(CD)}="D6";
   (-18,20)*+{((AB)C)D}="D7";
   (-18,-20)*+{A(B(CD))}="D8";
   (-52,-40)*+{A(B((1C)D))}="D9";
   (-18,-40)*+{A(B(1(CD)))}="D10";
   (-88,-20)*+{A(((B1)C)D)}="D15"; 
   (-88,-40)*+{A((B(1C))D)}="D19";   
   (-88, 0)*+{(A((B1)C))D}="D16";   
   (-88, 20)*+{((A(B1))C)D}="D17";   
   (-65,-20)*+{(A(B(1C)))D}="D18";   
        {\ar_{} "D1";"D2"}; %DOMAIN ARROWS
        {\ar_{} "D2";"D3"};
        {\ar_{} "D3";"D5"};
        {\ar_{} "D1";"D4"};
        {\ar_{} "D4";"D6"};
        {\ar@{-->}_{} "D5";"D10"};
        {\ar^{} "D7";"D6"};
        {\ar^{} "D1";"D7"};        
        {\ar_{} "D4";"D5"}; 
        {\ar@{-->}_{} "D3";"D9"}; 
        {\ar_{} "D5";"D6"};           
        {\ar@{-->}_{} "D9";"D10"};           
        {\ar@{-->}_{} "D1";"D17"};    
        {\ar@{-->}_{} "D17";"D16"};    
        {\ar@{-->}_{} "D18";"D19"};    
        {\ar@{-->}_{} "D16";"D18"};    
        {\ar@{-->}_{} "D16";"D15"};    
        {\ar@{-->}_{} "D15";"D19"};    
        {\ar@{-->}_{} "D19";"D9"};    
        {\ar@{-->}_{} "D10";"D8"};   
        {\ar@{-->}_{} "D6";"D8"};   
        {\ar@{-->}_{} "D2";"D18"};   
        {\ar@{=>}^{\scriptstyle \Pi\times 1} (-78, 7); (-73, 4)};
        {\ar@{=>}^{\scriptstyle \Pi} (-54.5, -2); (-48.5, -2)};
        {\ar@{=>}_<<{\scriptstyle \alpha^{-1}_{\rho,1,1}} (-34, 10); (-30, 12)};
        {\ar@{=>}^<<{\scriptstyle m} (-34, -10); (-31, -6)};
        {\ar@{=>}^<<{\scriptstyle \alpha^{-1}_{1,1,\lambda}} (-25, -19); (-22, -17)};
        {\ar@{=>}_<<{\scriptstyle \alpha^{-1}_{1,1,\alpha}} (-46, -31); (-43, -28)};
        {\ar@{=>}^{\scriptstyle \Pi} (-64, -27); (-58, -27)};
        {\ar@{=>}_<<{\scriptstyle \alpha^{-1}_{1,\alpha,1}} (-78, -20); (-74, -20)};
\endxy
\]
\[
\def\objectstyle{\scriptstyle}
  \def\labelstyle{\scriptstyle}
   \xy  
        {\ar@{=>}^{\scriptstyle 1\times U_{4,2}} (0, 2); (0, -2)}; %PERTURBATION
   \endxy
\]
\[
\def\objectstyle{\scriptstyle}
  \def\labelstyle{\scriptstyle}
   \xy 
   (35,20)*+{(((AB)1)C)D}="C1";
   (21,0)*+{((AB)(1C))D}="C2";
   (35,-10)*+{(AB)((1C)D)}="C3";
   (71,0)*+{(AB)(CD)}="C6";
   (71,20)*+{((AB)C)D}="C7";
   (71,-20)*+{A(B(CD))}="C8";
   (35,-40)*+{A(B((1C)D))}="C9";
   (71,-40)*+{A(B(1(CD)))}="C10";
   (-1,-20)*+{A(((B1)C)D)}="C15"; 
   (-1,-40)*+{A((B(1C))D)}="C19";   
   (-1, 0)*+{(A((B1)C))D}="C16";   
   (-1, 20)*+{((A(B1))C)D}="C17";   
   (21, -20)*+{(A(B(1C)))D}="C18"; 
   (50,-20)*+{(AB)(1(CD))}="C20"; 
        {\ar@{-->}_{} "C1";"C2"}; %DOMAIN ARROWS
        {\ar@{-->}_{} "C2";"C3"};
        {\ar_{} "C20";"C6"};
        {\ar_{} "C3";"C6"};
        {\ar@{-->}_{} "C2";"C7"};
        {\ar@{-->}^{} "C7";"C6"};
        {\ar@{-->}^{} "C1";"C7"};        
        {\ar_{} "C3";"C9"};     
        {\ar_{} "C9";"C10"};           
        {\ar@{-->}_{} "C1";"C17"};    
        {\ar@{-->}_{} "C17";"C16"};    
        {\ar@{-->}_{} "C18";"C19"};    
        {\ar@{-->}_{} "C16";"C18"};    
        {\ar@{-->}_{} "C16";"C15"};    
        {\ar@{-->}_{} "C15";"C19"};    
        {\ar@{-->}_{} "C19";"C9"};    
        {\ar_{} "C10";"C8"};   
        {\ar_{} "C6";"C8"};   
        {\ar@{-->}_{} "C2";"C18"}; 
        {\ar_{} "C3";"C20"}; 
        {\ar_{} "C20";"C10"};     
        {\ar@{=>}^{\scriptstyle \Pi\times 1} (10, 6); (14, 3)};
        {\ar@{=>}_<<{\scriptstyle m\times 1} (41, 14); (45, 16)};
        {\ar@{=>}_{\scriptstyle \alpha^{-1}_{1,\lambda,1}} (47, 3); (51, 6)};
        {\ar@{=>}_<<{\scriptstyle (1\times 1)\times l} (48, -13); (50, -9)};
        {\ar@{=>}_<<{\scriptstyle \alpha^{-1}_{1,1,\lambda}} (59, -21); (61, -18)};
        {\ar@{=>}_<<{\scriptstyle \alpha^{-1}_{1,1,\alpha}} (44, -31); (48, -28)};
        {\ar@{=>}^{\scriptstyle \Pi} (25, -27); (29, -27)};
        {\ar@{=>}_<<{\scriptstyle \alpha^{-1}_{1,\alpha,1}} (9, -19); (12, -19)};
\endxy
\]
\[
\def\objectstyle{\scriptstyle}
  \def\labelstyle{\scriptstyle}
   \xy  
    {\ar@{=>}_{\scriptstyle \alpha} (0, 2); (0, -2)}; 
   \endxy
\]
\[
\def\objectstyle{\scriptstyle}
  \def\labelstyle{\scriptstyle}
   \xy  
   (-52,20)*+{(((AB)1)C)D}="D1";
   (-65,0)*+{((AB)(1C))D}="D2";
   (-52,-10)*+{(AB)((1C)D)}="D3";
   (-18,0)*+{(AB)(CD)}="D6";
   (-18,20)*+{((AB)C)D}="D7";
   (-18,-20)*+{A(B(CD))}="D8";
   (-52,-40)*+{A(B((1C)D))}="D9";
   (-18,-40)*+{A(B(1(CD)))}="D10";
   (-88,-20)*+{A(((B1)C)D)}="D15"; 
   (-88,-40)*+{A((B(1C))D)}="D19";   
   (-88, 0)*+{(A((B1)C))D}="D16";   
   (-88, 20)*+{((A(B1))C)D}="D17";   
   (-65, -20)*+{(A(B(1C)))D}="D18"; 
        {\ar@{-->}_{} "D1";"D2"}; %DOMAIN ARROWS
        {\ar_{} "D2";"D3"};
        {\ar_{} "D3";"D6"};
        {\ar_{} "D2";"D7"};
        {\ar^{} "D7";"D6"};
        {\ar@{-->}^{} "D1";"D7"};        
        {\ar_{} "D3";"D9"};    
        {\ar@{-->}_{} "D9";"D10"};           
        {\ar@{-->}_{} "D1";"D17"};    
        {\ar@{-->}_{} "D17";"D16"};    
        {\ar_{} "D18";"D19"};    
        {\ar@{-->}_{} "D16";"D18"};    
        {\ar@{-->}_{} "D16";"D15"};    
        {\ar@{-->}_{} "D15";"D19"};    
        {\ar_{} "D19";"D9"};    
        {\ar@{-->}_{} "D10";"D8"};   
        {\ar_{} "D6";"D8"};   
        {\ar_{} "D2";"D18"}; 
        {\ar_{} "D9";"D8"}; 
        {\ar@{=>}^{\scriptstyle \Pi\times 1} (-78, 7); (-73, 4)};
        {\ar@{=>}_<<{\scriptstyle m\times 1} (-50, 14); (-46, 16)};
        {\ar@{=>}^<<{\scriptstyle \alpha^{-1}_{1,\lambda,1}} (-42, -3); (-38, -1)};
        {\ar@{=>}_<<{\scriptstyle \alpha^{-1}_{1,1,\lambda\times 1}} (-38, -15); (-34, -11)};
        {\ar@{=>}^<<{\scriptstyle 1\times (1\times l)} (-25, -33); (-29, -29)};
        {\ar@{=>}^{\scriptstyle \Pi} (-64, -27); (-58, -27)};
        {\ar@{=>}_<<{\scriptstyle \alpha^{-1}_{1,\alpha,1}} (-78, -20); (-74, -20)};
\endxy
\]
\[
\def\objectstyle{\scriptstyle}
  \def\labelstyle{\scriptstyle}
   \xy  
        {\ar@{=>}^{\scriptstyle 1\times \Pi} (0, 2); (0, -2)}; %PERTURBATION
   \endxy
\]
\[
\def\objectstyle{\scriptstyle}
  \def\labelstyle{\scriptstyle}
   \xy  
   (35,20)*+{(((AB)1)C)D}="C1";
   (21,0)*+{((AB)(1C))D}="C2";
   (49.5,0)*+{(A(BC))D}="C4";
   (59.5,-10)*+{A((BC)D)}="C5";
   (71,0)*+{(AB)(CD)}="C6";
   (71,20)*+{((AB)C)D}="C7";
   (71,-20)*+{A(B(CD))}="C8";
   (35,-40)*+{A(B((1C)D))}="C9";
   (71,-40)*+{A(B(1(CD)))}="C10";
   (-1,-20)*+{A(((B1)C)D)}="C15"; 
   (-1,-40)*+{A((B(1C))D)}="C19";   
   (-1, 0)*+{(A((B1)C))D}="C16";   
   (-1, 20)*+{((A(B1))C)D}="C17";   
   (21, -20)*+{(A(B(1C)))D}="C18"; 
        {\ar_{} "C1";"C2"}; %DOMAIN ARROWS
        {\ar_{} "C18";"C4"};
        {\ar_{} "C7";"C4"};
        {\ar_{} "C2";"C7"};
        {\ar@{-->}^{} "C7";"C6"};
        {\ar^{} "C1";"C7"};        
        {\ar@{-->}_{} "C4";"C5"}; 
        {\ar@{-->}_{} "C5";"C8"};           
        {\ar@{-->}_{} "C9";"C10"};           
        {\ar_{} "C1";"C17"};    
        {\ar_{} "C17";"C16"};    
        {\ar@{-->}_{} "C18";"C19"};    
        {\ar_{} "C16";"C18"};    
        {\ar@{-->}_{} "C16";"C15"};    
        {\ar@{-->}_{} "C15";"C19"};    
        {\ar@{-->}_{} "C19";"C9"};    
        {\ar@{-->}_{} "C10";"C8"};   
        {\ar@{-->}_{} "C6";"C8"};   
        {\ar_{} "C2";"C18"}; 
        {\ar@{-->}_{} "C9";"C8"}; 
        {\ar@{-->}_{} "C19";"C5"}; 
        {\ar@{=>}^{\scriptstyle \Pi\times 1} (10, 6); (14, 3)};
        {\ar@{=>}_<<{\scriptstyle \alpha^{-1}_{1,1,\lambda}\times 1} (31, -1); (35, 2)};
        {\ar@{=>}_<<{\scriptstyle m\times 1} (41, 14); (45, 16)};
        {\ar@{=>}_<<{\scriptstyle \Pi} (59, -1); (63, -1)};
        {\ar@{=>}^<<{\scriptstyle 1\times (1\times l)} (68, -32); (65, -29)};
        {\ar@{=>}_<<{\scriptstyle 1\times \alpha^{-1}_{1,\lambda,1}} (34, -31); (36, -28)};
        {\ar@{=>}^{\scriptstyle \alpha^{-1}_{1,1\times\lambda,1}} (39, -19); (42, -17)};
        {\ar@{=>}_<<{\scriptstyle \alpha^{-1}_{1,\alpha,1}} (8, -20); (12, -20)};
\endxy
\]
\[
\def\objectstyle{\scriptstyle}
  \def\labelstyle{\scriptstyle}
   \xy  
    {\ar@{=>}_{\scriptstyle U_{4,3}\times 1} (0, 2); (0, -2)}; 
   \endxy
\]

\[
\def\objectstyle{\scriptstyle}
  \def\labelstyle{\scriptstyle}
   \xy  
   (-52,20)*+{(((AB)1)C)D}="C1";
   (-41.5,0)*+{(A(BC))D}="C4";
   (-35,-10)*+{A((BC)D)}="C5";
   (-18,0)*+{(AB)(CD)}="C6";
   (-18,20)*+{((AB)C)D}="C7";
   (-18,-20)*+{A(B(CD))}="C8";
   (-52,-40)*+{A(B((1C)D))}="C9";
   (-18,-40)*+{A(B(1(CD)))}="C10";
   (-88,-20)*+{A(((B1)C)D)}="C15"; 
   (-88,-40)*+{A((B(1C))D)}="C19";   
   (-88, 0)*+{(A((B1)C))D}="C16";   
   (-88, 20)*+{((A(B1))C)D}="C17";   
   (-65, -20)*+{(A(B(1C)))D}="C18"; 
        {\ar@{-->}@/_2pc/_{} "C17";"C7"}; %DOMAIN ARROWS
        {\ar_{} "C18";"C4"};
        {\ar_{} "C16";"C4"};
        {\ar@{-->}_{} "C7";"C4"};
        {\ar@{-->}^{} "C7";"C6"};
        {\ar@{-->}^{} "C1";"C7"};        
        {\ar_{} "C4";"C5"}; 
        {\ar@{-->}_{} "C5";"C8"};           
        {\ar@{-->}_{} "C9";"C10"};           
        {\ar@{-->}_{} "C1";"C17"};    
        {\ar@{-->}_{} "C17";"C16"};    
        {\ar_{} "C18";"C19"};    
        {\ar_{} "C16";"C18"};    
        {\ar_{} "C16";"C15"};    
        {\ar_{} "C15";"C19"};    
        {\ar_{} "C19";"C9"};    
        {\ar@{-->}_{} "C10";"C8"};   
        {\ar@{-->}_{} "C6";"C8"};   
        {\ar@{-->}_{} "C9";"C8"}; 
        {\ar_{} "C19";"C5"};      
        {\ar@{=>}^{\scriptstyle \alpha^{-1}_{1,\rho,1}\times 1} (-52, 4); (-48, 5)};
        {\ar@{=>}^<<{\scriptstyle (r\times 1)\times 1} (-52, 13); (-48, 14)};
        {\ar@{=>}^<<{\scriptstyle (1\times m)\times 1} (-63, -10); (-63, -4)};
        {\ar@{=>}_<<{\scriptstyle \Pi} (-33, -1); (-28, -1)};
        {\ar@{=>}^<<{\scriptstyle 1\times (1\times l)} (-25, -33); (-29, -29)};
        {\ar@{=>}_<<{\scriptstyle 1\times \alpha^{-1}_{1,\lambda,1}} (-57, -31); (-55, -28)};
        {\ar@{=>}^{\scriptstyle \alpha^{-1}_{1,1\times\lambda,1}} (-50.5, -18.5); (-46.5, -16.5)};
        {\ar@{=>}_<<{\scriptstyle \alpha^{-1}_{1,\alpha,1}} (-78, -20); (-74, -20)};
\endxy
\]
\[
\def\objectstyle{\scriptstyle}
  \def\labelstyle{\scriptstyle}
   \xy  
        {\ar@{=>}^{\scriptstyle \alpha} (0, 2); (0, -2)}; %PERTURBATION
\endxy
\]
\[
\def\objectstyle{\scriptstyle}
  \def\labelstyle{\scriptstyle}
   \xy  
   (35,20)*+{(((AB)1)C)D}="D1";
   (49.5,0)*+{(A(BC))D}="D4";
   (59.5,-10)*+{A((BC)D)}="D5";
   (71,0)*+{(AB)(CD)}="D6";
   (71,20)*+{((AB)C)D}="D7";
   (71,-20)*+{A(B(CD))}="D8";
   (35,-40)*+{A(B((1C)D))}="D9";
   (71,-40)*+{A(B(1(CD)))}="D10";
   (-1,-20)*+{A(((B1)C)D)}="D15"; 
   (-1,-40)*+{A((B(1C))D)}="D19";   
   (-1, 0)*+{(A((B1)C))D}="D16";   
   (-1, 20)*+{((A(B1))C)D}="D17";   
        {\ar@{-->}@/_2pc/_{} "D17";"D7"}; %DOMAIN ARROWS
        {\ar@{-->}_{} "D16";"D4"};
        {\ar@{-->}_{} "D7";"D4"};
        {\ar@{-->}^{} "D7";"D6"};
        {\ar@{-->}^{} "D1";"D7"};        
        {\ar@{-->}_{} "D4";"D5"}; 
        {\ar@{-->}_{} "D5";"D8"};           
        {\ar@{-->}_{} "D9";"D10"};           
        {\ar@{-->}_{} "D1";"D17"};    
        {\ar@{-->}_{} "D17";"D16"};    
        {\ar@{-->}_{} "D16";"D15"};    
        {\ar@{-->}_{} "D15";"D19"};    
        {\ar@{-->}_{} "D19";"D9"};    
        {\ar@{-->}_{} "D10";"D8"};   
        {\ar@{-->}_{} "D6";"D8"};   
        {\ar@{-->}_{} "D15";"D5"}; 
        {\ar@{-->}_{} "D9";"D8"}; 
        {\ar@{-->}_{} "D19";"D5"};   
        {\ar@{=>}^{\scriptstyle \alpha^{-1}_{1,\rho,1}\times 1} (32, 3); (37, 5)};
        {\ar@{=>}^<<{\scriptstyle (r\times 1)\times 1} (40, 13); (43, 14)};
        {\ar@{=>}_<<{\scriptstyle \Pi} (59, -1); (63, -1)};
        {\ar@{=>}^<<{\scriptstyle 1\times (1\times l)} (68, -32); (65, -29)};
        {\ar@{=>}_<<{\scriptstyle 1\times \alpha^{-1}_{1,\lambda,1}} (34, -31); (36, -28)};
        {\ar@{=>}^{\scriptstyle \alpha^{-1}_{1,\rho\times 1,1}} (32, -11); (36, -8)};
        {\ar@{=>}_{\scriptstyle 1\times(m\times 1)} (9, -25); (12, -22)};
\endxy
\]
\[ \textrm{\Huge{=}}\]
\[
\def\objectstyle{\scriptstyle}
  \def\labelstyle{\scriptstyle}
   \xy
   (-52,20)*+{(((AB)1)C)D}="D1";
   (-65,0)*+{((AB)(1C))D}="D2";
   (-52,-10)*+{(AB)((1C)D)}="D3";
   (-39,0)*+{((AB)1)(CD)}="D4";
   (-39,-20)*+{(AB)(1(CD))}="D5";
   (-18,0)*+{(AB)(CD)}="D6";
   (-18,20)*+{((AB)C)D}="D7";
   (-18,-20)*+{A(B(CD))}="D8";
   (-52,-40)*+{A(B((1C)D))}="D9";
   (-18,-40)*+{A(B(1(CD)))}="D10";
   (-88,-20)*+{A(((B1)C)D)}="D15"; 
   (-88,-40)*+{A((B(1C))D)}="D19";   
   (-88, 0)*+{(A((B1)C))D}="D16";   
   (-88, 20)*+{((A(B1))C)D}="D17";   
   (-65,-20)*+{(A(B(1C)))D}="D18";    
        {\ar_{} "D1";"D2"}; %DOMAIN ARROWS
        {\ar_{} "D2";"D3"};
        {\ar_{} "D3";"D5"};
        {\ar_{} "D1";"D4"};
        {\ar@{-->}_{} "D4";"D6"};
        {\ar_{} "D5";"D10"};
        {\ar@{-->}^{} "D7";"D6"};
        {\ar@{-->}^{} "D1";"D7"};        
        {\ar_{} "D4";"D5"}; 
        {\ar_{} "D3";"D9"}; 
        {\ar@{-->}_{} "D5";"D6"};           
        {\ar_{} "D9";"D10"};           
        {\ar_{} "D1";"D17"};    
        {\ar_{} "D17";"D16"};    
        {\ar_{} "D18";"D19"};    
        {\ar_{} "D16";"D18"};    
        {\ar_{} "D16";"D15"};    
        {\ar_{} "D15";"D19"};    
        {\ar_{} "D19";"D9"};    
        {\ar@{-->}_{} "D10";"D8"};   
        {\ar@{-->}_{} "D6";"D8"};   
        {\ar_{} "D2";"D18"};   
        {\ar@{=>}^{\scriptstyle \Pi\times 1} (-78, 7); (-73, 4)};
        {\ar@{=>}^{\scriptstyle \Pi} (-54.5, -2); (-48.5, -2)};
        {\ar@{=>}_<<{\scriptstyle \alpha^{-1}_{\rho,1,1}} (-34, 10); (-30, 12)};
        {\ar@{=>}^<<{\scriptstyle m} (-34, -10); (-31, -6)};
        {\ar@{=>}_<<{\scriptstyle \alpha^{-1}_{1,1,\lambda}} (-29, -21); (-26, -18)};
        {\ar@{=>}_<<{\scriptstyle \alpha^{-1}_{1,1,\alpha}} (-46, -31); (-43, -28)};
        {\ar@{=>}^{\scriptstyle \Pi} (-64, -27); (-58, -27)};
        {\ar@{=>}_<<{\scriptstyle \alpha^{-1}_{1,\alpha,1}} (-78, -20); (-74, -20)};
\endxy
\]
\[
\def\objectstyle{\scriptstyle}
  \def\labelstyle{\scriptstyle}
   \xy  
        {\ar@{=>}^{\scriptstyle K_5} (0, 2); (0, -2)}; %PERTURBATION
   \endxy
\]
\[
\def\objectstyle{\scriptstyle}
  \def\labelstyle{\scriptstyle}
   \xy 
   (35,20)*+{(((AB)1)C)D}="C1";
   (21,-20)*+{(A(B1))(CD)}="C2";
   (50,0)*+{((AB)1)(CD)}="C4";
   (71,0)*+{(AB)(CD)}="C6";
   (71,20)*+{((AB)C)D}="C7";
   (71,-20)*+{A(B(CD))}="C8";
   (35,-40)*+{A(B((1C)D))}="C9";
   (71,-40)*+{A(B(1(CD)))}="C10";
   (-1,-20)*+{A(((B1)C)D)}="C15"; 
   (-1,-40)*+{A((B(1C))D)}="C19";   
   (-1, 0)*+{(A((B1)C))D}="C16";   
   (-1, 20)*+{((A(B1))C)D}="C17";   
   (21, -30)*+{A((B1)(CD))}="C18"; 
   (50,-20)*+{(AB)(1(CD))}="C20";   
        {\ar@{-->}_{} "C1";"C4"}; %DOMAIN ARROWS
        {\ar_{} "C4";"C6"};
        {\ar_{} "C4";"C2"};
        {\ar@{-->}^{} "C7";"C6"};
        {\ar@{-->}^{} "C1";"C7"};                   
        {\ar@{-->}_{} "C9";"C10"};           
        {\ar@{-->}_{} "C1";"C17"};    
        {\ar@{-->}_{} "C17";"C16"};    
        {\ar@{-->}_{} "C17";"C2"};  
        {\ar@{-->}_{} "C15";"C18"};      
        {\ar@{-->}_{} "C16";"C15"};    
        {\ar@{-->}_{} "C15";"C19"};    
        {\ar@{-->}_{} "C19";"C9"};    
        {\ar_{} "C10";"C8"};   
        {\ar_{} "C6";"C8"};   
        {\ar_{} "C2";"C18"}; 
        {\ar_{} "C4";"C20"}; 
        {\ar_{} "C20";"C6"}; 
        {\ar_{} "C18";"C10"};     
        {\ar_{} "C20";"C10"};  
        {\ar@{=>}^{\scriptstyle \alpha^{-1}_{\alpha,1,1}} (17, 6); (20, 3)};
        {\ar@{=>}_<<{\scriptstyle \alpha^{-1}_{\rho,1,1}} (53, 10); (56, 12)};
        {\ar@{=>}^{\scriptstyle m} (55, -8); (58, -5)};
        {\ar@{=>}_<<{\scriptstyle \alpha_{1,1,\lambda}^{-1}} (59, -20); (62, -17)};
        {\ar@{=>}_<<{\scriptstyle 1\times \Pi} (16, -37); (20, -34)};
        {\ar@{=>}^{\scriptstyle \Pi} (33, -25); (38, -25)};
        {\ar@{=>}_<<{\scriptstyle \Pi} (7, -6); (10, -9)};
\endxy
\]
\[
\def\objectstyle{\scriptstyle}
  \def\labelstyle{\scriptstyle}
   \xy  
    {\ar@{=>}_{\scriptstyle U_{4,3}} (0, 2); (0, -2)}; 
   \endxy
\]

\[
\def\objectstyle{\scriptstyle}
  \def\labelstyle{\scriptstyle}
   \xy  
   (-52,20)*+{(((AB)1)C)D}="D1";
   (-65,-20)*+{(A(B1))(CD)}="D3";
   (-39,0)*+{((AB)1)(CD)}="D4";
   (-18,0)*+{(AB)(CD)}="D6";
   (-18,20)*+{((AB)C)D}="D7";
   (-18,-20)*+{A(B(CD))}="D8";
   (-52,-40)*+{A(B((1C)D))}="D9";
   (-18,-40)*+{A(B(1(CD)))}="D10";
   (-88,-20)*+{A(((B1)C)D)}="D15"; 
   (-88,-40)*+{A((B(1C))D)}="D19";   
   (-88, 0)*+{(A((B1)C))D}="D16";   
   (-88, 20)*+{((A(B1))C)D}="D17";   
   (-65, -30)*+{A((B1)(CD))}="D18"; 
        {\ar_{} "D3";"D6"};
        {\ar_{} "D4";"D6"};
        {\ar^{} "D7";"D6"};
        {\ar^{} "D1";"D7"};        
        {\ar_{} "D1";"D4"}; 
        {\ar@{-->}_{} "D3";"D18"}; 
        {\ar@{-->}_{} "D18";"D8"};           
        {\ar@{-->}_{} "D9";"D10"};           
        {\ar_{} "D1";"D17"};    
        {\ar@{-->}_{} "D17";"D16"};    
        {\ar@{-->}_{} "D15";"D18"};    
        {\ar_{} "D17";"D3"};    
        {\ar@{-->}_{} "D16";"D15"};    
        {\ar@{-->}_{} "D15";"D19"};    
        {\ar_{} "D19";"D9"};    
        {\ar@{-->}_{} "D10";"D8"};   
        {\ar@{-->}_{} "D6";"D8"};   
        {\ar_{} "D4";"D3"}; 
        {\ar@{-->}_{} "D18";"D10"}; 
        {\ar@{=>}^{\scriptstyle \Pi} (-83, -6); (-79, -9)};
        {\ar@{=>}_<<{\scriptstyle \alpha^{-1}_{\rho,1,1}} (-34, 10); (-30, 12)};
        {\ar@{=>}^<<{\scriptstyle \alpha^{-1}_{\alpha,1,1}} (-64, 6); (-61, 3)};
        {\ar@{=>}_<<{\scriptstyle r\times (1\times 1)} (-44, -6); (-40, -3)};
        {\ar@{=>}^<<{\scriptstyle 1\times m} (-28, -30); (-29, -26)};
        {\ar@{=>}_<<{\scriptstyle \alpha^{-1}_{1,\rho,1}} (-40, -18); (-37, -15)};
        {\ar@{=>}^{\scriptstyle 1\times\Pi} (-78, -34); (-74, -31)};
\endxy
\]
\[
\def\objectstyle{\scriptstyle}
  \def\labelstyle{\scriptstyle}
   \xy  
        {\ar@{=>}^{\scriptstyle \alpha} (0, 2); (0, -2)}; %PERTURBATION
   \endxy
\]
\[
\def\objectstyle{\scriptstyle}
  \def\labelstyle{\scriptstyle}
   \xy  
   (35,20)*+{(((AB)1)C)D}="C1";
   (21,-20)*+{(A(B1))(CD)}="C2";
   (71,0)*+{(AB)(CD)}="C6";
   (71,20)*+{((AB)C)D}="C7";
   (71,-20)*+{A(B(CD))}="C8";
   (35,-40)*+{A(B((1C)D))}="C9";
   (71,-40)*+{A(B(1(CD)))}="C10";
   (-1,-20)*+{A(((B1)C)D)}="C15"; 
   (-1,-40)*+{A((B(1C))D)}="C19";   
   (-1, 0)*+{(A((B1)C))D}="C16";   
   (-1, 20)*+{((A(B1))C)D}="C17";   
   (21, -30)*+{A((B1)(CD))}="C18"; 
        {\ar_{} "C17";"C2"}; %DOMAIN ARROWS
        {\ar^{} "C7";"C6"};
        {\ar@{-->}^{} "C1";"C7"};                 
        {\ar@{-->}_{} "C9";"C10"};           
        {\ar@{-->}_{} "C1";"C17"};    
        {\ar_{} "C17";"C16"};    
        {\ar@{-->}_{} "C18";"C10"};    
        {\ar_{} "C15";"C18"};    
        {\ar_{} "C16";"C15"};    
        {\ar@{-->}_{} "C15";"C19"};    
        {\ar@{-->}_{} "C19";"C9"};    
        {\ar@{-->}_{} "C10";"C8"};   
        {\ar_{} "C6";"C8"};   
        {\ar_{} "C2";"C18"}; 
        {\ar@/_2pc/_{} "C17";"C7"};
        {\ar_{} "C2";"C6"}; 
        {\ar_{} "C18";"C8"}; 
        {\ar@{=>}^{\scriptstyle \Pi} (5, -6); (8, -9)};
        {\ar@{=>}_<<{\scriptstyle \alpha^{-1}_{1\times \rho,1,1}} (37, -1); (41, 3)};
        {\ar@{=>}^<<{\scriptstyle (r\times 1)\times 1} (40, 13); (43, 14)};
        {\ar@{=>}_<<{\scriptstyle 1\times m} (54, -32); (52, -29)};
        {\ar@{=>}^{\scriptstyle \alpha^{-1}_{1,\rho,1}} (46, -20); (50, -17)};
        {\ar@{=>}^<<{\scriptstyle 1\times \Pi} (16, -37); (20, -34)};
\endxy
\]
\[
\def\objectstyle{\scriptstyle}
  \def\labelstyle{\scriptstyle}
   \xy  
    {\ar@{=>}_{\scriptstyle \Pi} (0, 2); (0, -2)}; 
   \endxy
\]

\[
\def\objectstyle{\scriptstyle}
  \def\labelstyle{\scriptstyle}
   \xy  
   (-52,20)*+{(((AB)1)C)D}="C1";
   (-41.5,0)*+{(A(BC))D}="C4";
   (-35,-10)*+{A((BC)D)}="C5";
   (-18,0)*+{(AB)(CD)}="C6";
   (-18,20)*+{((AB)C)D}="C7";
   (-18,-20)*+{A(B(CD))}="C8";
   (-52,-40)*+{A(B((1C)D))}="C9";
   (-18,-40)*+{A(B(1(CD)))}="C10";
   (-88,-20)*+{A(((B1)C)D)}="C15"; 
   (-88,-40)*+{A((B(1C))D)}="C19";   
   (-88, 0)*+{(A((B1)C))D}="C16";   
   (-88, 20)*+{((A(B1))C)D}="C17";   
   (-65, -30)*+{A((B1)(CD))}="C18"; 
   	{\ar@{-->}@/_2pc/_{} "C17";"C7"}; %DOMAIN ARROWS
        {\ar@{-->}_{} "C16";"C4"};
        {\ar@{-->}_{} "C7";"C4"};
        {\ar@{-->}^{} "C7";"C6"};
        {\ar@{-->}^{} "C1";"C7"};        
        {\ar@{-->}_{} "C4";"C5"}; 
        {\ar_{} "C5";"C8"};           
        {\ar@{-->}_{} "C9";"C10"};           
        {\ar@{-->}_{} "C1";"C17"};    
        {\ar@{-->}_{} "C17";"C16"};    
        {\ar_{} "C15";"C18"};    
        {\ar_{} "C18";"C8"};    
        {\ar@{-->}_{} "C16";"C15"};    
        {\ar_{} "C15";"C19"};    
        {\ar_{} "C19";"C9"};    
        {\ar_{} "C10";"C8"};   
        {\ar@{-->}_{} "C6";"C8"};   
        {\ar_{} "C18";"C10"}; 
        {\ar_{} "C15";"C5"};         
        {\ar@{=>}^{\scriptstyle \alpha^{-1}_{1,\rho,1}\times 1} (-52, 4); (-48, 5)};
        {\ar@{=>}^<<{\scriptstyle (r\times 1)\times 1} (-52, 13); (-48, 14)};
        {\ar@{=>}^<<{\scriptstyle \alpha^{-1}_{1,\rho\times 1,1}} (-52, -10); (-48, -8)};
        {\ar@{=>}_<<{\scriptstyle \Pi} (-30, -1); (-26, -1)};
        {\ar@{=>}^<<{\scriptstyle 1\times m} (-30, -30); (-31, -26)};
        {\ar@{=>}^{\scriptstyle 1\times \alpha^{-1}_{\rho,1,1}} (-52, -23); (-48, -21)};
        {\ar@{=>}^{\scriptstyle 1\times\Pi} (-78, -34); (-74, -31)};
\endxy
\]
\[
\def\objectstyle{\scriptstyle}
  \def\labelstyle{\scriptstyle}
   \xy  
        {\ar@{=>}^{\scriptstyle 1\times U_{4,2}} (0, 2); (0, -2)}; %PERTURBATION
\endxy
\]
\[
\def\objectstyle{\scriptstyle}
  \def\labelstyle{\scriptstyle}
   \xy  
   (35,20)*+{(((AB)1)C)D}="D1";
   (49.5,0)*+{(A(BC))D}="D4";
   (59.5,-10)*+{A((BC)D)}="D5";
   (71,0)*+{(AB)(CD)}="D6";
   (71,20)*+{((AB)C)D}="D7";
   (71,-20)*+{A(B(CD))}="D8";
   (35,-40)*+{A(B((1C)D))}="D9";
   (71,-40)*+{A(B(1(CD)))}="D10";
   (-1,-20)*+{A(((B1)C)D)}="D15"; 
   (-1,-40)*+{A((B(1C))D)}="D19";   
   (-1, 0)*+{(A((B1)C))D}="D16";   
   (-1, 20)*+{((A(B1))C)D}="D17";      
        {\ar@{-->}@/_2pc/_{} "D17";"D7"}; %DOMAIN ARROWS
        {\ar@{-->}_{} "D16";"D4"};
        {\ar@{-->}_{} "D7";"D4"};
        {\ar@{-->}^{} "D7";"D6"};
        {\ar@{-->}^{} "D1";"D7"};        
        {\ar@{-->}_{} "D4";"D5"}; 
        {\ar@{-->}_{} "D5";"D8"};           
        {\ar@{-->}_{} "D9";"D10"};           
        {\ar@{-->}_{} "D1";"D17"};    
        {\ar@{-->}_{} "D17";"D16"};     
        {\ar@{-->}_{} "D16";"D15"};    
        {\ar@{-->}_{} "D15";"D19"};    
        {\ar@{-->}_{} "D19";"D9"};    
        {\ar@{-->}_{} "D10";"D8"};   
        {\ar@{-->}_{} "D6";"D8"};   
        {\ar@{-->}_{} "D15";"D5"}; 
        {\ar@{-->}_{} "D9";"D8"}; 
        {\ar@{-->}_{} "D19";"D5"};          
        {\ar@{=>}^{\scriptstyle \alpha^{-1}_{1,\rho,1}\times 1} (32, 3); (37, 5)};
        {\ar@{=>}^<<{\scriptstyle (r\times 1)\times 1} (40, 13); (43, 14)};
        {\ar@{=>}_<<{\scriptstyle \Pi} (59, -1); (63, -1)};
        {\ar@{=>}^<<{\scriptstyle 1\times (1\times l)} (68, -32); (65, -29)};
        {\ar@{=>}_<<{\scriptstyle 1\times \alpha^{-1}_{1,\lambda,1}} (34, -31); (36, -28)};
        {\ar@{=>}^{\scriptstyle \alpha^{-1}_{1,\rho\times 1,1}} (32, -11); (36, -8)};
        {\ar@{=>}_{\scriptstyle 1\times(m\times 1)} (9, -25); (12, -22)};
\endxy
\]

\subsubsection*{$U_{5,4}$ Axiom}
\begin{itemize}
\item for each $5$-tuple of objects $a,b,c,d,e$, an equation called the $U_{5,4}$ unit condition, consisting of, for each $4$-tuple of morphisms $A,B,C,D$, an equation of $4$-cells
\end{itemize}
\[
\def\objectstyle{\scriptstyle}
  \def\labelstyle{\scriptstyle}
   \xy
   (-52,20)*+{(((AB)C)1)D}="D1";
   (-65,0)*+{((AB)(C1))D}="D2";
   (-52,-10)*+{(AB)((C1)D)}="D3";
   (-39,0)*+{((AB)C)(1D)}="D4";
   (-39,-20)*+{(AB)(C(1D))}="D5";
   (-18,0)*+{(AB)(CD)}="D6";
   (-18,20)*+{((AB)C)D}="D7";
   (-18,-20)*+{A(B(CD))}="D8";
   (-52,-40)*+{A(B((C1)D))}="D9";
   (-18,-40)*+{A(B(C(1D)))}="D10";
   (-88,-20)*+{A(((BC)1)D)}="D15"; 
   (-88,-40)*+{A((B(C1))D)}="D19";   
   (-88, 0)*+{(A((BC)1))D}="D16";   
   (-88, 20)*+{((A(BC))1)D}="D17";   
   (-65,-20)*+{(A(B(C1)))D}="D18";   
        {\ar_{} "D1";"D2"}; %DOMAIN ARROWS
        {\ar_{} "D2";"D3"};
        {\ar_{} "D3";"D5"};
        {\ar_{} "D1";"D4"};
        {\ar_{} "D4";"D7"};
        {\ar@{-->}_{} "D5";"D10"};
        {\ar^{} "D7";"D6"};
        {\ar^{} "D1";"D7"};        
        {\ar_{} "D4";"D5"}; 
        {\ar@{-->}_{} "D3";"D9"}; 
        {\ar_{} "D5";"D6"};           
        {\ar@{-->}_{} "D9";"D10"};           
        {\ar@{-->}_{} "D1";"D17"};    
        {\ar@{-->}_{} "D17";"D16"};    
        {\ar@{-->}_{} "D18";"D19"};    
        {\ar@{-->}_{} "D16";"D18"};    
        {\ar@{-->}_{} "D16";"D15"};    
        {\ar@{-->}_{} "D15";"D19"};    
        {\ar@{-->}_{} "D19";"D9"};    
        {\ar@{-->}_{} "D10";"D8"};   
        {\ar@{-->}_{} "D6";"D8"};   
        {\ar@{-->}_{} "D2";"D18"};   
        {\ar@{=>}^{\scriptstyle \Pi\times 1} (-78, 7); (-73, 4)};
        {\ar@{=>}^{\scriptstyle \Pi} (-54.5, -2); (-48.5, -2)};
        {\ar@{=>}_<<{\scriptstyle m} (-40, 10); (-40, 15)};
        {\ar@{=>}^<<{\scriptstyle \alpha^{-1}_{1,1,\lambda}} (-28, -1); (-24, 0)};
        {\ar@{=>}_<<{\scriptstyle \alpha^{-1}_{1,1,1\times\lambda}} (-30, -21); (-26, -19)};
        {\ar@{=>}_<<{\scriptstyle \alpha^{-1}_{1,1,\alpha}} (-45, -31); (-41, -27)};
        {\ar@{=>}^{\scriptstyle \Pi} (-64, -27); (-58, -27)};
        {\ar@{=>}_<<{\scriptstyle \alpha^{-1}_{1,\alpha,1}} (-78, -20); (-74, -20)};
\endxy
\]
\[
\def\objectstyle{\scriptstyle}
  \def\labelstyle{\scriptstyle}
   \xy  
        {\ar@{=>}^{\scriptstyle U_{4,3}} (0, 2); (0, -2)}; %PERTURBATION
   \endxy
\]

\[
\def\objectstyle{\scriptstyle}
  \def\labelstyle{\scriptstyle}
   \xy
   (35,20)*+{(((AB)C)1)D}="C1";
   (21,0)*+{((AB)(C1))D}="C2";
   (35,-10)*+{(AB)((C1)D)}="C3";
   (71,0)*+{(AB)(CD)}="C6";
   (71,20)*+{((AB)C)D}="C7";
   (71,-20)*+{A(B(CD))}="C8";
   (35,-40)*+{A(B((C1)D))}="C9";
   (71,-40)*+{A(B(C(1D)))}="C10";
   (-1,-20)*+{A(((BC)1)D)}="C15"; 
   (-1,-40)*+{A((B(C1))D)}="C19";   
   (-1, 0)*+{(A((BC)1))D}="C16";   
   (-1, 20)*+{((A(BC))1)D}="C17";   
   (21, -20)*+{(A(B(C1)))D}="C18"; 
   (50,-20)*+{(AB)(C(1D))}="C20";   
        {\ar@{-->}_{} "C1";"C2"}; %DOMAIN ARROWS
        {\ar@{-->}_{} "C2";"C3"};
        {\ar_{} "C20";"C6"};
        {\ar_{} "C3";"C6"};
        {\ar@{-->}_{} "C2";"C7"};
        {\ar@{-->}^{} "C7";"C6"};
        {\ar@{-->}^{} "C1";"C7"};        
        {\ar_{} "C3";"C9"};          
        {\ar_{} "C9";"C10"};           
        {\ar@{-->}_{} "C1";"C17"};    
        {\ar@{-->}_{} "C17";"C16"};    
        {\ar@{-->}_{} "C18";"C19"};    
        {\ar@{-->}_{} "C16";"C18"};    
        {\ar@{-->}_{} "C16";"C15"};    
        {\ar@{-->}_{} "C15";"C19"};    
        {\ar@{-->}_{} "C19";"C9"};    
        {\ar_{} "C10";"C8"};   
        {\ar_{} "C6";"C8"};   
        {\ar@{-->}_{} "C2";"C18"}; 
        {\ar_{} "C3";"C20"}; 
        {\ar_{} "C20";"C10"};     
        {\ar@{=>}^{\scriptstyle \Pi\times 1} (10, 6); (14, 3)};
        {\ar@{=>}_<<{\scriptstyle r\times 1} (37, 13); (40, 16)};
        {\ar@{=>}_{\scriptstyle \alpha^{-1}_{1,\rho,1}} (42, 2); (46, 5)};
        {\ar@{=>}_<<{\scriptstyle \alpha^{-1}_{1,1,1\times\lambda}} (60, -21); (64, -19)};
        {\ar@{=>}_<<{\scriptstyle (1\times 1)\times m} (50, -14); (48, -11)};
        {\ar@{=>}_<<{\scriptstyle \alpha^{-1}_{1,1,\alpha}} (43, -31); (46, -28)};        		
        {\ar@{=>}^{\scriptstyle \Pi} (25, -27); (29, -27)};
        {\ar@{=>}_<<{\scriptstyle \alpha^{-1}_{1,\alpha,1}} (8, -20); (12, -20)};
\endxy
\]
\[
\def\objectstyle{\scriptstyle}
  \def\labelstyle{\scriptstyle}
   \xy  
    {\ar@{=>}_{\scriptstyle \alpha} (0, 2); (0, -2)}; 
   \endxy
\]

\[
\def\objectstyle{\scriptstyle}
  \def\labelstyle{\scriptstyle}
   \xy  
   (-52,20)*+{(((AB)C)1)D}="D1";
   (-65,0)*+{((AB)(C1))D}="D2";
   (-52,-10)*+{(AB)((C1)D)}="D3";
   (-18,0)*+{(AB)(CD)}="D6";
   (-18,20)*+{((AB)C)D}="D7";
   (-18,-20)*+{A(B(CD))}="D8";
   (-52,-40)*+{A(B((C1)D))}="D9";
   (-18,-40)*+{A(B(C(1D)))}="D10";
   (-88,-20)*+{A(((BC)1)D)}="D15"; 
   (-88,-40)*+{A((B(C1))D)}="D19";   
   (-88, 0)*+{(A((BC)1))D}="D16";   
   (-88, 20)*+{((A(BC))1)D}="D17";   
   (-65, -20)*+{(A(B(C1)))D}="D18"; 
        {\ar@{-->}_{} "D1";"D2"}; %DOMAIN ARROWS
        {\ar_{} "D2";"D3"};
        {\ar_{} "D3";"D6"};
        {\ar_{} "D2";"D7"};
        {\ar^{} "D7";"D6"};
        {\ar@{-->}^{} "D1";"D7"};        
        {\ar_{} "D3";"D9"};          
        {\ar@{-->}_{} "D9";"D10"};           
        {\ar@{-->}_{} "D1";"D17"};    
        {\ar@{-->}_{} "D17";"D16"};    
        {\ar_{} "D18";"D19"};    
        {\ar@{-->}_{} "D16";"D18"};    
        {\ar@{-->}_{} "D16";"D15"};    
        {\ar@{-->}_{} "D15";"D19"};    
        {\ar_{} "D19";"D9"};    
        {\ar@{-->}_{} "D10";"D8"};   
        {\ar_{} "D6";"D8"};   
        {\ar_{} "D2";"D18"}; 
        {\ar_{} "D9";"D8"}; 
        {\ar@{=>}^{\scriptstyle \Pi\times 1} (-78, 7); (-73, 4)};
        {\ar@{=>}_<<{\scriptstyle r\times 1} (-51, 13); (-47, 15)};
        {\ar@{=>}^<<{\scriptstyle \alpha^{-1}_{1,\rho,1}} (-42, -3); (-38, -1)};
        {\ar@{=>}^<<{\scriptstyle 1\times (1\times m)} (-25, -33); (-29, -29)};
        {\ar@{=>}_<<{\scriptstyle \alpha^{-1}_{1,1,\rho\times 1}} (-40, -20); (-36, -15)};
        {\ar@{=>}^{\scriptstyle \Pi} (-64, -27); (-58, -27)};
        {\ar@{=>}_<<{\scriptstyle \alpha^{-1}_{1,\alpha,1}} (-78, -20); (-74, -20)};
\endxy
\]
\[
\def\objectstyle{\scriptstyle}
  \def\labelstyle{\scriptstyle}
   \xy  
        {\ar@{=>}^{\scriptstyle \Pi} (0, 2); (0, -2)}; %PERTURBATION
   \endxy
\]

\[
\def\objectstyle{\scriptstyle}
  \def\labelstyle{\scriptstyle}
   \xy  
   (35,20)*+{(((AB)C)1)D}="C1";
   (21,0)*+{((AB)(C1))D}="C2";
   (49.5,-10)*+{(A(BC))D}="C4";
   (49.5,-20)*+{A((BC)D)}="C5";
   (71,0)*+{(AB)(CD)}="C6";
   (71,20)*+{((AB)C)D}="C7";
   (71,-20)*+{A(B(CD))}="C8";
   (35,-40)*+{A(B((C1)D))}="C9";
   (71,-40)*+{A(B(C(1D)))}="C10";
   (-1,-20)*+{A(((BC)1)D)}="C15"; 
   (-1,-40)*+{A((B(C1))D)}="C19";   
   (-1, 0)*+{(A((BC)1))D}="C16";   
   (-1, 20)*+{((A(BC))1)D}="C17";   
   (21, -20)*+{(A(B(C1)))D}="C18"; 
        {\ar_{} "C1";"C2"}; %DOMAIN ARROWS
        {\ar_{} "C18";"C4"};
        {\ar_{} "C7";"C4"};
        {\ar_{} "C2";"C7"};
        {\ar@{-->}^{} "C7";"C6"};
        {\ar^{} "C1";"C7"};        
        {\ar@{-->}_{} "C4";"C5"}; 
        {\ar@{-->}_{} "C5";"C8"};           
        {\ar@{-->}_{} "C9";"C10"};           
        {\ar_{} "C1";"C17"};    
        {\ar_{} "C17";"C16"};    
        {\ar@{-->}_{} "C18";"C19"};    
        {\ar_{} "C16";"C18"};    
        {\ar@{-->}_{} "C16";"C15"};    
        {\ar@{-->}_{} "C15";"C19"};    
        {\ar@{-->}_{} "C19";"C9"};    
        {\ar@{-->}_{} "C10";"C8"};   
        {\ar@{-->}_{} "C6";"C8"};   
        {\ar_{} "C2";"C18"}; 
        {\ar@{-->}_{} "C9";"C8"}; 
        {\ar@{-->}_{} "C19";"C5"};    
        {\ar@{=>}^{\scriptstyle \Pi\times 1} (10, 6); (14, 3)};
        {\ar@{=>}_<<{\scriptstyle r\times 1} (37, 13); (40, 16)};
        {\ar@{=>}^{\scriptstyle \alpha^{-1}_{1,1,\rho}\times 1} (41, -6); (45, -2 )};
        {\ar@{=>}_<<{\scriptstyle \Pi} (60, -7); (65, -7)};
        {\ar@{=>}^<<{\scriptstyle 1\times (1\times m)} (68, -32); (65, -29)};
        {\ar@{=>}_<<{\scriptstyle 1\times \alpha^{-1}_{1,\rho,1}} (34, -30); (37, -27)};
        {\ar@{=>}^{\scriptstyle \alpha^{-1}_{1,1\times \rho,1}} (38, -23); (42, -19)};
        {\ar@{=>}_<<{\scriptstyle \alpha^{-1}_{1,\alpha,1}} (8, -20); (12, -20)};
\endxy
\]
\[
\def\objectstyle{\scriptstyle}
  \def\labelstyle{\scriptstyle}
   \xy  
      {\ar@{=>}^{\scriptstyle \Pi} (0, 2); (0, -2)}; %PERTURBATION
   \endxy
\]

\[
\def\objectstyle{\scriptstyle}
  \def\labelstyle{\scriptstyle}
   \xy  
   (-52,20)*+{(((AB)C)1)D}="C1";
   (-39,-10)*+{(A(BC))D}="C4";
   (-39,-20)*+{A((BC)D)}="C5";
   (-18,0)*+{(AB)(CD)}="C6";
   (-18,20)*+{((AB)C)D}="C7";
   (-18,-20)*+{A(B(CD))}="C8";
   (-52,-40)*+{A(B((C1)D))}="C9";
   (-18,-40)*+{A(B(C(1D)))}="C10";
   (-88,-20)*+{A(((BC)1)D)}="C15"; 
   (-88,-40)*+{A((B(C1))D)}="C19";   
   (-88, 0)*+{(A((BC)1))D}="C16";   
   (-88, 20)*+{((A(BC))1)D}="C17";   
   (-65, -20)*+{(A(B(C1)))D}="C18"; 
        {\ar@{-->}_{} "C17";"C4"}; %DOMAIN ARROWS
        {\ar_{} "C18";"C4"};
        {\ar_{} "C16";"C4"};
        {\ar@{-->}_{} "C7";"C4"};
        {\ar@{-->}^{} "C7";"C6"};
        {\ar@{-->}^{} "C1";"C7"};        
        {\ar_{} "C4";"C5"}; 
        {\ar@{-->}_{} "C5";"C8"};           
        {\ar@{-->}_{} "C9";"C10"};           
        {\ar@{-->}_{} "C1";"C17"};    
        {\ar@{-->}_{} "C17";"C16"};    
        {\ar_{} "C18";"C19"};    
        {\ar_{} "C16";"C18"};    
        {\ar_{} "C16";"C15"};    
        {\ar_{} "C15";"C19"};    
        {\ar@{-->}_{} "C19";"C9"};    
        {\ar@{-->}_{} "C10";"C8"};   
        {\ar@{-->}_{} "C6";"C8"};   
        {\ar@{-->}_{} "C9";"C8"}; 
        {\ar_{} "C19";"C5"};  
        {\ar@{=>}^{\scriptstyle r\times 1} (-73, 0); (-70, 3)};
        {\ar@{=>}^<<{\scriptstyle \rho_{\alpha}\times 1} (-51, 9); (-46, 9)};
        {\ar@{=>}^<<{\scriptstyle (1\times r)\times 1} (-64, -12); (-61, -9)};
        {\ar@{=>}_<<{\scriptstyle \Pi} (-30, -5); (-25, -5)};
        {\ar@{=>}^<<{\scriptstyle 1\times (1\times m)} (-25, -33); (-29, -29)};
        {\ar@{=>}_<<{\scriptstyle 1\times \alpha_{1,\rho,1}} (-59, -34); (-56, -30)};
        {\ar@{=>}_{\scriptstyle \alpha^{-1}_{1,\lambda,1}} (-68, -25); (-64, -22)};
        {\ar@{=>}_<<{\scriptstyle \alpha^{-1}_{1,\alpha,1}} (-78, -20); (-74, -20)};
\endxy
\]
\[
\def\objectstyle{\scriptstyle}
  \def\labelstyle{\scriptstyle}
   \xy  
        {\ar@{=>}^{\scriptstyle \alpha} (0, 2); (0, -2)}; %PERTURBATION
\endxy
\]
\[
\def\objectstyle{\scriptstyle}
  \def\labelstyle{\scriptstyle}
   \xy   
   (35,20)*+{(((AB)C)1)D}="D1";
   (49.5,-10)*+{(A(BC))D}="D4";
   (49.5,-20)*+{A((BC)D)}="D5";
   (71,0)*+{(AB)(CD)}="D6";
   (71,20)*+{((AB)C)D}="D7";
   (71,-20)*+{A(B(CD))}="D8";
   (35,-40)*+{A(B((C1)D))}="D9";
   (71,-40)*+{A(B(C(1D)))}="D10";
   (-1,-20)*+{A(((BC)1)D)}="D15"; 
   (-1,-40)*+{A((B(C1))D)}="D19";   
   (-1, 0)*+{(A((BC)1))D}="D16";   
   (-1, 20)*+{((A(BC))1)D}="D17";   
        {\ar@{-->}_{} "D17";"D4"}; %DOMAIN ARROWS
        {\ar@{-->}_{} "D16";"D4"};
        {\ar@{-->}_{} "D7";"D4"};
        {\ar@{-->}^{} "D7";"D6"};
        {\ar@{-->}^{} "D1";"D7"};        
        {\ar@{-->}_{} "D4";"D5"}; 
        {\ar@{-->}_{} "D5";"D8"};           
        {\ar@{-->}_{} "D9";"D10"};           
        {\ar@{-->}_{} "D1";"D17"};    
        {\ar@{-->}_{} "D17";"D16"};     
        {\ar@{-->}_{} "D16";"D15"};    
        {\ar@{-->}_{} "D15";"D19"};    
        {\ar@{-->}_{} "D19";"D9"};    
        {\ar@{-->}_{} "D10";"D8"};   
        {\ar@{-->}_{} "D6";"D8"};   
        {\ar@{-->}_{} "D15";"D5"}; 
        {\ar@{-->}_{} "D9";"D8"}; 
        {\ar@{-->}_{} "D19";"D5"};   
        {\ar@{=>}^{\scriptstyle r\times 1} (25, -2); (28, 0)};
        {\ar@{=>}^<<{\scriptstyle \rho_{\alpha}\times 1} (43, 9); (47, 9)};
        {\ar@{=>}_<<{\scriptstyle \Pi} (60, -7); (65, -7)};
        {\ar@{=>}^<<{\scriptstyle 1\times (1\times m)} (68, -32); (65, -29)};
        {\ar@{=>}_<<{\scriptstyle 1\times \alpha_{1,\rho,1}} (39, -31); (37, -28)};
        {\ar@{=>}^{\scriptstyle \alpha^{-1}_{1,\rho,1}} (26, -16); (30, -13)};
        {\ar@{=>}_{\scriptstyle 1\times(r\times 1)} (10, -27); (14, -24)};
\endxy
\]

\[ \textrm{\Huge{=}}\]

\[
\def\objectstyle{\scriptstyle}
  \def\labelstyle{\scriptstyle}
   \xy
   (-52,20)*+{(((AB)C)1)D}="D1";
   (-65,0)*+{((AB)(C1))D}="D2";
   (-52,-10)*+{(AB)((C1)D)}="D3";
   (-39,0)*+{((AB)C)(1D)}="D4";
   (-39,-20)*+{(AB)(C(1D))}="D5";
   (-18,0)*+{(AB)(CD)}="D6";
   (-18,20)*+{((AB)C)D}="D7";
   (-18,-20)*+{A(B(CD))}="D8";
   (-52,-40)*+{A(B((C1)D))}="D9";
   (-18,-40)*+{A(B(C(1D)))}="D10";
   (-88,-20)*+{A(((BC)1)D)}="D15"; 
   (-88,-40)*+{A((B(C1))D)}="D19";   
   (-88, 0)*+{(A((BC)1))D}="D16";   
   (-88, 20)*+{((A(BC))1)D}="D17";   
   (-65,-20)*+{(A(B(C1)))D}="D18";   
        {\ar_{} "D1";"D2"}; %DOMAIN ARROWS
        {\ar_{} "D2";"D3"};
        {\ar_{} "D3";"D5"};
        {\ar_{} "D1";"D4"};
        {\ar@{-->}_{} "D4";"D7"};
        {\ar_{} "D5";"D10"};
        {\ar@{-->}^{} "D7";"D6"};
        {\ar@{-->}^{} "D1";"D7"};        
        {\ar_{} "D4";"D5"}; 
        {\ar_{} "D3";"D9"}; 
        {\ar@{-->}_{} "D5";"D6"};           
        {\ar_{} "D9";"D10"};           
        {\ar_{} "D1";"D17"};    
        {\ar_{} "D17";"D16"};    
        {\ar_{} "D18";"D19"};    
        {\ar_{} "D16";"D18"};    
        {\ar_{} "D16";"D15"};    
        {\ar_{} "D15";"D19"};    
        {\ar_{} "D19";"D9"};    
        {\ar@{-->}_{} "D10";"D8"};   
        {\ar@{-->}_{} "D6";"D8"};   
        {\ar_{} "D2";"D18"};   
        {\ar@{=>}^{\scriptstyle \Pi\times 1} (-78, 7); (-73, 4)};
        {\ar@{=>}^{\scriptstyle \Pi} (-54.5, -2); (-48.5, -2)};
        {\ar@{=>}_<<{\scriptstyle m} (-40, 10); (-40, 15)};
        {\ar@{=>}^<<{\scriptstyle \alpha^{-1}_{1,1,\lambda}} (-28, -1); (-24, 0)};
        {\ar@{=>}_<<{\scriptstyle \alpha^{-1}_{1,1,1\times\lambda}} (-30, -21); (-26, -19)};
        {\ar@{=>}_<<{\scriptstyle \alpha^{-1}_{1,1,\alpha}} (-45, -31); (-41, -27)};
        {\ar@{=>}^{\scriptstyle \Pi} (-64, -27); (-58, -27)};
        {\ar@{=>}_<<{\scriptstyle \alpha^{-1}_{1,\alpha,1}} (-78, -20); (-74, -20)};
\endxy
\]

\[
\def\objectstyle{\scriptstyle}
  \def\labelstyle{\scriptstyle}
   \xy  
        {\ar@{=>}^{\scriptstyle K_5} (0, 2); (0, -2)}; %PERTURBATION
   \endxy
\]

\[
\def\objectstyle{\scriptstyle}
  \def\labelstyle{\scriptstyle}
   \xy
   (35,20)*+{(((AB)C)1)D}="C1";
   (21,-20)*+{(A(BC))(1D)}="C2";
   (49.5,0)*+{((AB)C)(1D)}="C4";
   (71,0)*+{(AB)(CD)}="C6";
   (71,20)*+{((AB)C)D}="C7";
   (71,-20)*+{A(B(CD))}="C8";
   (35,-40)*+{A(B((C1)D))}="C9";
   (71,-40)*+{A(B(C(1D)))}="C10";
   (-1,-20)*+{A(((BC)1)D)}="C15"; 
   (-1,-40)*+{A((B(C1))D)}="C19";   
   (-1, 0)*+{(A((BC)1))D}="C16";   
   (-1, 20)*+{((A(BC))1)D}="C17";   
   (21, -30)*+{A((BC)(1D))}="C18"; 
   (49.5,-20)*+{(AB)(C(1D))}="C20"; 
        {\ar@{-->}_{} "C1";"C4"}; %DOMAIN ARROWS
        {\ar_{} "C4";"C7"};
        {\ar_{} "C4";"C2"};
        {\ar^{} "C7";"C6"};
        {\ar@{-->}^{} "C1";"C7"};        
        {\ar@{-->}_{} "C9";"C10"};           
        {\ar@{-->}_{} "C1";"C17"};    
        {\ar@{-->}_{} "C17";"C16"};    
        {\ar@{-->}_{} "C17";"C2"};  
        {\ar@{-->}_{} "C15";"C18"};    
        {\ar@{-->}_{} "C16";"C15"};    
        {\ar@{-->}_{} "C15";"C19"};    
        {\ar@{-->}_{} "C19";"C9"};    
        {\ar_{} "C10";"C8"};   
        {\ar_{} "C6";"C8"};   
        {\ar_{} "C2";"C18"}; 
        {\ar_{} "C4";"C20"}; 
        {\ar_{} "C20";"C6"}; 
        {\ar_{} "C18";"C10"};     
        {\ar_{} "C20";"C10"};  
        {\ar@{=>}^{\scriptstyle \Pi} (5, -6); (8, -9)};
        {\ar@{=>}_<<{\scriptstyle m} (51, 10); (51, 15)};
        {\ar@{=>}_{\scriptstyle \alpha_{\alpha,1,1}} (26, 4); (30, 1)};
        {\ar@{=>}_<<{\scriptstyle \alpha^{-1}_{1,1,1\times\lambda}} (60, -21); (64, -19)};
        {\ar@{=>}^<<{\scriptstyle 1\times \Pi} (16, -37); (20, -34)};
        {\ar@{=>}_<<{\scriptstyle \alpha^{-1}_{1,1,\lambda}} (58, -1); (62, 0)};
        {\ar@{=>}^{\scriptstyle \Pi} (33, -26); (39, -26)};
   \endxy
\]
\[
\def\objectstyle{\scriptstyle}
  \def\labelstyle{\scriptstyle}
   \xy  
    {\ar@{=>}_{\scriptstyle \Pi} (0, 2); (0, -2)}; 
   \endxy
\]

\[
\def\objectstyle{\scriptstyle}
  \def\labelstyle{\scriptstyle}
   \xy  
   (-52,20)*+{(((AB)C)1)D}="D1";
   (-35,-20)*+{A((BC)D)}="D2";
   (-65,-20)*+{(A(BC))(1D)}="D3";
   (-40,0)*+{((AB)C)(1D)}="D4";
   (-35,-10)*+{(A(BC))D}="D5";
   (-18,0)*+{(AB)(CD)}="D6";
   (-18,20)*+{((AB)C)D}="D7";
   (-18,-20)*+{A(B(CD))}="D8";
   (-52,-40)*+{A(B((C1)D))}="D9";
   (-18,-40)*+{A(B(C(1D)))}="D10";
   (-88,-20)*+{A(((BC)1)D)}="D15"; 
   (-88,-40)*+{A((B(C1))D)}="D19";   
   (-88, 0)*+{(A((BC)1))D}="D16";   
   (-88, 20)*+{((A(BC))1)D}="D17";   
   (-65, -30)*+{A((BC)(1D))}="D18"; 
       {\ar@{-->}_{} "D2";"D8"}; %DOMAIN ARROWS
        {\ar_{} "D3";"D5"};
        {\ar_{} "D4";"D7"};
        {\ar@{-->}_{} "D5";"D2"};
        {\ar_{} "D4";"D3"};
        {\ar@{-->}^{} "D7";"D6"};
        {\ar^{} "D1";"D7"};        
        {\ar_{} "D1";"D4"}; 
        {\ar@{-->}_{} "D3";"D18"}; 
        {\ar@{-->}_{} "D18";"D2"};           
        {\ar@{-->}_{} "D9";"D10"};           
        {\ar_{} "D1";"D17"};    
        {\ar@{-->}_{} "D17";"D16"};    
        {\ar@{-->}_{} "D15";"D18"};    
        {\ar_{} "D17";"D3"};    
        {\ar@{-->}_{} "D16";"D15"};    
        {\ar@{-->}_{} "D15";"D19"};    
        {\ar_{} "D19";"D9"};    
        {\ar@{-->}_{} "D10";"D8"};   
        {\ar@{-->}_{} "D6";"D8"};   
        {\ar_{} "D7";"D5"}; 
        {\ar@{-->}_{} "D18";"D10"}; 
        {\ar@{=>}^{\scriptstyle \Pi} (-83, -5); (-79, -8)};
        {\ar@{=>}_<<{\scriptstyle m} (-37, 10); (-37, 15)};
        {\ar@{=>}^<<{\scriptstyle \alpha_{\alpha,1,1}} (-62, 3); (-57, 0)};
        {\ar@{=>}_<<{\scriptstyle \Pi} (-27, -5); (-23, -5)};
                {\ar@{=>}_<<{\scriptstyle \alpha^{-1}\widetilde{\times} \lambda} (-45.5, -6); (-41.5, -3)};
        {\ar@{=>}^<<{\scriptstyle 1\times {\alpha}^{-1}_{1,1,\lambda}} (-28, -30); (-29, -26)};
        {\ar@{=>}_<<{\scriptstyle \alpha^{-1}_{1,1,\lambda}} (-55, -20); (-50, -18)};
        {\ar@{=>}^{\scriptstyle 1\times\Pi} (-78, -34); (-74, -32)};
   \endxy
\]
\[
\def\objectstyle{\scriptstyle}
  \def\labelstyle{\scriptstyle}
   \xy  
        {\ar@{=>}^{\scriptstyle m} (0, 2); (0, -2)}; %PERTURBATION
   \endxy
\]

\[
\def\objectstyle{\scriptstyle}
  \def\labelstyle{\scriptstyle}
   \xy  
   (49.5,-20)*+{A((BC)D)}="C20"; 
   (36,20)*+{(((AB)C)1)D}="C1";
   (21,-20)*+{(A(BC))(1D)}="C2";
   (49.5,-10)*+{(A(BC))D}="C3";
   (71,0)*+{(AB)(CD)}="C6";
   (71,20)*+{((AB)C)D}="C7";
   (71,-20)*+{A(B(CD))}="C8";
   (35,-40)*+{A(B((C1)D))}="C9";
   (71,-40)*+{A(B(C(1D)))}="C10";
   (-1,-20)*+{A(((BC)1)D)}="C15"; 
   (-1,-40)*+{A((B(C1))D)}="C19";   
   (-1, 0)*+{(A((BC)1))D}="C16";   
   (-1, 20)*+{((A(BC))1)D}="C17";   
   (21, -30)*+{A((BC)(1D))}="C18"; 
        {\ar_{} "C17";"C2"}; %DOMAIN ARROWS
        {\ar_{} "C17";"C3"};
        {\ar@{-->}^{} "C7";"C6"};
        {\ar@{-->}^{} "C1";"C7"};        
        {\ar_{} "C3";"C20"}; 
        {\ar@{-->}_{} "C20";"C8"};           
        {\ar@{-->}_{} "C9";"C10"};           
        {\ar@{-->}_{} "C1";"C17"};    
        {\ar_{} "C17";"C16"};    
        {\ar@{-->}_{} "C18";"C10"};    
        {\ar_{} "C15";"C18"};    
        {\ar_{} "C16";"C15"};    
        {\ar@{-->}_{} "C15";"C19"};    
        {\ar@{-->}_{} "C19";"C9"};    
        {\ar@{-->}_{} "C10";"C8"};   
        {\ar@{-->}_{} "C6";"C8"};   
        {\ar_{} "C2";"C18"}; 
        {\ar@{-->}_{} "C7";"C3"};
        {\ar_{} "C2";"C3"}; 
        {\ar_{} "C18";"C20"}; 
        {\ar@{=>}^{\scriptstyle \Pi} (5, -6); (8, -9)};
        {\ar@{=>}_<<{\scriptstyle \rho_{\alpha}\times 1} (44, 10); (48, 12)};
        {\ar@{=>}^<<{\scriptstyle m} (25, -8); (30, -5)};
        {\ar@{=>}_<<{\scriptstyle \Pi} (61, -3); (65, -3)};
        {\ar@{=>}^{\scriptstyle \alpha^{-1}_{1,1,\lambda}} (36, -23); (39, -22)};
        {\ar@{=>}_<<{\scriptstyle 1\times \alpha^{-1}_{1,1,\lambda}} (56, -30); (54, -26)};
        {\ar@{=>}^<<{\scriptstyle 1\times \Pi} (16, -37); (20, -34)};
   \endxy
\]
\[
\def\objectstyle{\scriptstyle}
  \def\labelstyle{\scriptstyle}
   \xy  
    {\ar@{=>}_{\scriptstyle U_{4,3}} (0, 2); (0, -2)}; 
   \endxy
\]

\[
\def\objectstyle{\scriptstyle}
  \def\labelstyle{\scriptstyle}
   \xy  
   (-52,20)*+{(((AB)C)1)D}="C1";
   (-40,-10)*+{(A(BC))D}="C4";
   (-40,-20)*+{A((BC)D)}="C5";
   (-20,0)*+{(AB)(CD)}="C6";
   (-20,20)*+{((AB)C)D}="C7";
   (-20,-20)*+{A(B(CD))}="C8";
   (-52,-40)*+{A(B((C1)D))}="C9";
   (-20,-40)*+{A(B(C(1D)))}="C10";
   (-88,-20)*+{A(((BC)1)D)}="C15"; 
   (-88,-40)*+{A((B(C1))D)}="C19";   
   (-88, 0)*+{(A((BC)1))D}="C16";   
   (-88, 20)*+{((A(BC))1)D}="C17";   
   (-65, -30)*+{A((BC)(1D))}="C18"; 
        {\ar@{-->}_{} "C17";"C4"};
        {\ar@{-->}_{} "C16";"C4"};
        {\ar@{-->}_{} "C7";"C4"};
        {\ar@{-->}^{} "C7";"C6"};
        {\ar@{-->}^{} "C1";"C7"};        
        {\ar@{-->}_{} "C4";"C5"}; 
        {\ar_{} "C5";"C8"};           
        {\ar_{} "C9";"C10"};           
        {\ar@{-->}_{} "C1";"C17"};    
        {\ar@{-->}_{} "C17";"C16"};    
        {\ar_{} "C15";"C18"};    
        {\ar_{} "C18";"C5"};    
        {\ar@{-->}_{} "C16";"C15"};    
        {\ar_{} "C15";"C19"};    
        {\ar_{} "C19";"C9"};    
        {\ar_{} "C10";"C8"};   
        {\ar@{-->}_{} "C6";"C8"};   
        {\ar_{} "C18";"C10"}; 
        {\ar_{} "C15";"C5"}; 
        {\ar@{=>}^{\scriptstyle r\times 1} (-74, 3); (-70, 5)};
        {\ar@{=>}^<<{\scriptstyle  \rho_{\alpha}\times 1} (-41, 8); (-37, 13)};
        {\ar@{=>}^<<{\scriptstyle 1 \times m} (-65, -26); (-65, -22)};
        {\ar@{=>}_<<{\scriptstyle \Pi} (-30, -5); (-25, -5)};
        {\ar@{=>}^{\scriptstyle 1\times\Pi} (-78, -34); (-74, -32)};
        {\ar@{=>}^<<{\scriptstyle 1\times \alpha_{1,1,\lambda}} (-30, -30); (-31, -26)};
        {\ar@{=>}^{\scriptstyle \alpha^{-1}_{1,\lambda,1}} (-63, -15); (-59, -12)};
   \endxy
\]
\[
\def\objectstyle{\scriptstyle}
  \def\labelstyle{\scriptstyle}
   \xy  
        {\ar@{=>}^{\scriptstyle 1\times U_{4,3}} (0, 2); (0, -2)}; %PERTURBATION
   \endxy
\]
\[
\def\objectstyle{\scriptstyle}
  \def\labelstyle{\scriptstyle}
   \xy  
   (35,20)*+{(((AB)C)1)D}="D1";
   (49.5,-10)*+{(A(BC))D}="D4";
   (49.5,-20)*+{A((BC)D)}="D5";
   (71,0)*+{(AB)(CD)}="D6";
   (71,20)*+{((AB)C)D}="D7";
   (71,-20)*+{A(B(CD))}="D8";
   (35,-40)*+{A(B((C1)D))}="D9";
   (71,-40)*+{A(B(C(1D)))}="D10";
   (-1,-20)*+{A(((BC)1)D)}="D15"; 
   (-1,-40)*+{A((B(C1))D)}="D19";   
   (-1, 0)*+{(A((BC)1))D}="D16";   
   (-1, 20)*+{((A(BC))1)D}="D17";   ;   
        {\ar@{-->}_{} "D17";"D4"}; %DOMAIN ARROWS
        {\ar@{-->}_{} "D16";"D4"};
        {\ar@{-->}_{} "D7";"D4"};
        {\ar@{-->}^{} "D7";"D6"};
        {\ar@{-->}^{} "D1";"D7"};        
        {\ar@{-->}_{} "D4";"D5"}; 
        {\ar@{-->}_{} "D5";"D8"};           
        {\ar@{-->}_{} "D9";"D10"};           
        {\ar@{-->}_{} "D1";"D17"};    
        {\ar@{-->}_{} "D17";"D16"};     
        {\ar@{-->}_{} "D16";"D15"};    
        {\ar@{-->}_{} "D15";"D19"};    
        {\ar@{-->}_{} "D19";"D9"};    
        {\ar@{-->}_{} "D10";"D8"};   
        {\ar@{-->}_{} "D6";"D8"};   
        {\ar@{-->}_{} "D15";"D5"}; 
        {\ar@{-->}_{} "D9";"D8"}; 
        {\ar@{-->}_{} "D19";"D5"}; 
        {\ar@{=>}^{\scriptstyle r\times 1} (25, -2); (28, 0)};
        {\ar@{=>}^<<{\scriptstyle \rho_{\alpha}\times 1} (43, 9); (47, 9)};
        {\ar@{=>}_<<{\scriptstyle \Pi} (60, -7); (65, -7)};
        {\ar@{=>}^<<{\scriptstyle 1\times (1\times m)} (68, -32); (65, -29)};
        {\ar@{=>}_<<{\scriptstyle 1\times \alpha_{1,\rho,1}} (39, -31); (37, -28)};
        {\ar@{=>}^{\scriptstyle \alpha^{-1}_{1,\rho,1}} (26, -16); (30, -13)};
        {\ar@{=>}_{\scriptstyle 1\times(r\times 1)} (10, -27); (14, -24)};
   \endxy
\]

%%%%%%%%%%%%%%%%%%%%%%%%%%%%%%%%%%%%%%%%%%%%%%%%%%%%%%%%%%%%%%%%%%%%%%%%

\section{A Tricategory of Spans}\label{sec-tricategory}

We now restate the main theorem of this section followed by a series of definitions and propositions, which together prove the existence of a tricategory structure on spans in a $2$-category with pullbacks.

\theorem \label{spantricategory}
Let $\B$ be a (strict) $2$-category with pullbacks, consisting of
\begin{itemize}
\item $0$-cells $A,B,C,\ldots$,
\item $1$-cells $p,q,r,\ldots$, and
\item $2$-cells $\varpi,\varrho,\varsigma,\ldots$.
\end{itemize}

There is a semi-strict cubical tricategory $\Span(\B)$ consisting of
\begin{itemize}
\item the objects of $\B$ as objects,
\item for each pair of objects $A, B$, the (strict) $2$-category $\Span(A,B)$, defined in Proposition~\ref{hom2category}, consisting of
\begin{itemize}
\item spans in $\B$ (see Definition~\ref{1morphisms}),
\item maps of spans in $\B$ (see Definition~\ref{2morphisms}), and
\item maps of maps of spans in $\B$ (see Definition~\ref{3morphisms}),
\end{itemize}
\noindent which are, respectively, the $1$-morphisms, $2$-morphisms, and $3$-morphisms of $\Span(\B)$,
\item for each triple of objects $A, B, C$, a strict $2$-functor
\[ *_{ABC}\maps \Span(A,B)\otimes \Span(B,C)\to \Span(A,C),\]
\noindent defined in Proposition~\ref{compositionfunctor},
\item for each object $A$, a strict $2$-functor
\[ I_{A}\maps 1\to \Span(A,A),\]
\noindent defined in Proposition~\ref{unitfunctor},
\item for each $4$-tuple $A, B, C, D$ of objects, a strict adjoint equivalence
\[ {\bf a}_{ABCD}\maps *_{(AB)CD}(*_{ABC}\times 1)\To *_{A(BC)D}(1\times *_{BCD})\]
\noindent in the $2$-category of strict $2$-functors, strict transformations, and modifications, defined in Proposition~\ref{associatortransformation},
\item for each pair of objects $A, B$, an adjoint equivalence
\[ {\bf l}_{AB}\maps *_{ABB}(I_B\times 1)\To  1,\]
\noindent in the $2$-category of strict $2$-functors, strong transformations, and modifications, defined in Proposition~\ref{leftunitortransformation},
\item for each pair of objects $A, B$, an adjoint equivalence
\[ {\bf r}_{AB}\maps *_{AAB}(1\times I_A)\To  1,\]
\noindent in the $2$-category of strict $2$-functors, strong transformations, and modifications, defined in Proposition~\ref{rightunitortransformation},
\item for each $5$-tuple of objects $A, B, C, D, E$, an identity modification
\[
   \xy
   (-30,0)*+{\Pi_{ABCDE}\maps (*\cdot(1\times {\bf a}))({\bf a}\cdot (1\times * \times 1))(*\cdot({\bf a}\times 1))}="2"; (-36,-5)*+{}="4"; (-28,-5)*+{}="5"; (5,-5)*+{({\bf a}\cdot(1\times 1\times *))(*\cdot 1)({\bf a}\cdot(1\times 1 \times *))}="3";
        {\ar@3{->}_{} "4";"5"};
\endxy
\]
defined in Proposition~\ref{pentagonatormodification},
\item for each triple of objects $A, B, C$, an identity modification
\[
   \xy
   (-24,0)*+{\Lambda_{ABC}\maps  1(*\cdot(1\times {\bf l}))}="2"; (24,0)*+{({\bf l}\cdot *)({\bf a}\cdot (1\times 1\times I))(*\cdot 1)}="3";
        {\ar@3{->}_{} "2";"3"};
\endxy
\]
defined Proposition~\ref{leftunitormodification},
\item for each triple of objects $A, B, C$, an identity modification
\[
   \xy
   (-30,0)*+{M_{ABC}\maps (*\cdot({\bf l}\times 1))({\bf a}^{-1}\cdot (1\times I \times 1))(*\cdot
(1\times {\bf r}^{-1}))}="2"; (25,0)*+{*\cdot 1}="3";
        {\ar@3{->}_{} "2";"3"};
\endxy
\]
defined in Proposition~\ref{middleunitormodification}, and
\item for each triple of objects $A, B, C$, an identity modification
\[
   \xy
   (-25,0)*+{P_{ABC}\maps  (*\cdot(1 \times {\bf r}))1}="2"; (25,0)*+{({\bf r}\cdot *)(*\cdot 1)({\bf a}\cdot (I\times 1\times 1))}="3";
        {\ar@3{->}_{} "2";"3"};
\endxy
\]
defined Proposition~\ref{rightunitormodification}.
\end{itemize}
\endtheorem

\proof
The structural components are given in the referenced definitions and propositions.  The tricategory is cubical since it is locally strict with strict, thus cubical, composition and unit $2$-functors, and all product cells and modification components are identity morphisms.  The tricategory axioms then follow and we have the desired result.
\endproof

\begin{remark}
If the adjoint equivalences {\bf a}, {\bf l}, {\bf r} were identities, then $\Span(\B)$ would have the structure of a $\Gray$-category.  We are using the terminology {\it semi-strict} to reflect the fact that the tricategory has only trivial modifications, but possibly non-trivial adjoint equivalences.
\end{remark}

\subsection{Span Morphisms}
We define the morphisms of the span tricategory and set notation for the pullbacks and products used in defining the structure cells of the tricategory and the monoidal structure on the tricategory.

Let $\B$ be a strict $2$-category.  The morphisms of $\B$ are called $1$-cells and $2$-cells.  The structure cells of the span tricategory $\Span(\B)$ are called ($1$-)morphisms, $2$-morphisms, and $3$-morphisms.  We define the morphisms of $\Span(\B)$ here.

\subsubsection*{Definitions of morphisms}\label{morphisms_defn}

\noindent The $1$-morphisms in $\Span(\B)$ are `spans'.

\begin{definition}\label{1morphisms}
A {\bf span} in a $2$-category $\B$ is a pair of $1$-cells
in $\B$ with a common source object
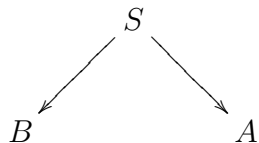
\[
   \xy
   (0,15)*+{S}="1";
   (-15,0)*+{B}="2"; 
   (15,0)*+{A}="3";
        {\ar_{} "1";"2"};
        {\ar^{} "1";"3"};
\endxy
\]
\end{definition}

\noindent The $2$-morphisms in $\Span(\B)$ are `maps of spans'.

\begin{definition}\label{2morphisms}
Given a pair of parallel spans in a $2$-category $\B$
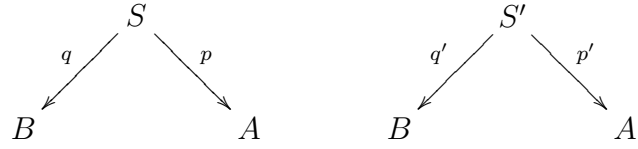
\[
   \xy
   (-25,15)*+{S}="1";
   (-40,0)*+{B}="2"; 
   (-10,0)*+{A}="3";
   (25,15)*+{S'}="1'";
   (10,0)*+{B}="2'"; 
   (40,0)*+{A}="3'";
        {\ar_{q} "1";"2"};
        {\ar^{p} "1";"3"};
        {\ar_{q'} "1'";"2'"};
        {\ar^{p'} "1'";"3'"};
\endxy
\]
\noindent a {\bf map of spans} is a triple $(\varpi,f,\varrho)$ consisting of a
$1$-cell $f\maps S \to S'$ together with a pair of invertible $2$-cells $\alpha\maps p\Rightarrow p'f$ and $\beta\maps q\Rightarrow q'f$
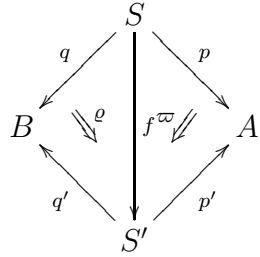
\[
   \xy
   (0,15)*+{S}="1";
   (-15,0)*+{B}="2"; (15,0)*+{A}="3";
   (0,-15)*+{S'}="4";
        {\ar_{q} "1";"2"};
        {\ar^{p} "1";"3"};
        {\ar^{q'} "4";"2"};
        {\ar_{p'} "4";"3"};
        {\ar^{f} "1";"4"};
        {\ar@{=>}^{\scriptstyle \varrho} (-8,2); (-5,-2)};
        {\ar@{=>}_{\scriptstyle \varpi} (8,2); (5,-2)};
\endxy
\]
\end{definition}

\noindent The $3$-morphisms in $\Span(\B)$ are `maps of maps of spans'.
\begin{definition}\label{3morphisms}
Given a parallel pair of maps of spans $(\varpi,f,\varrho)$ and $(\varpi',f',\varrho')$ in a $2$-category $\B$

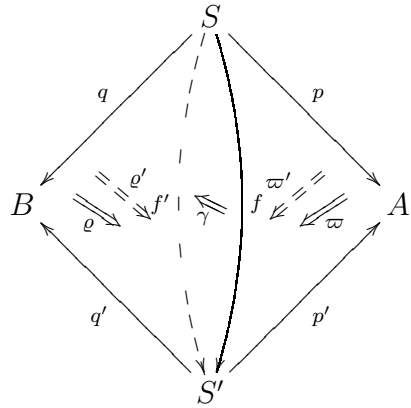
\[ 
   \xy
   (0,25)*+{S}="1";
   (-25,0)*+{B}="2";
   (25,0)*+{A}="3";
   (0,-25)*+{S'}="4";
        {\ar_{q} "1";"2"};
        {\ar^{p} "1";"3"};
        {\ar^{q'} "4";"2"};
        {\ar_{p'} "4";"3"};
        {\ar@/^1pc/^{f} "1";"4"};
	{\ar@{-->}@/_1pc/_{f'} "1";"4"};
        {\ar@{=>}_{\scriptstyle \varrho} (-18,1); (-12,-3)};
        {\ar@{=>}^{\scriptstyle \varpi} (18,1); (12,-3)};
        {\ar@{==>}^{\scriptstyle \varrho'} (-15,4); (-8,-2)};
        {\ar@{==>}_{\scriptstyle \varpi'} (15,4); (8,-2)};
        {\ar@{=>}^{\scriptstyle \gamma} (2,-1); (-2,1)};         		
\endxy
\]
\noindent a {\bf map of maps of spans} 
\[
   \xy
   (-12,0)*+{\gamma\maps (\varpi,f,\varrho)}="2"; (12,0)*+{(\varpi',f',\varrho')}="3";
        {\ar@3{->}_{} "2";"3"};
\endxy
\]
consists of a $2$-cell $\gamma\maps f \Rightarrow f'$ in $\B$ such that the following equations hold
\begin{equation}\label{3celleq1}
(p'\cdot\gamma)\varpi = \varpi'
\end{equation}

\noindent and

\begin{equation}\label{3celleq2}
(q'\cdot\gamma)\varrho = \varrho'
\end{equation}
\end{definition}

\noindent The notation $\varpi\cdot f$ denotes the {\em whiskering} of a $2$-cell $\varpi$ along a $1$-cell $f$.  This is defined by horizontal composition of the $2$-cell $\varpi$ with the identity $2$-cell $1_f$ in $\B$.

By definition a map of maps of spans $\gamma\maps (\varpi,f,\varrho)\To (\varpi',f',\varrho')$ satisfies  the equations
\[ (p'\cdot\gamma)\varpi = \varpi' \;\textrm{and}\; (q'\cdot\gamma)\varpi = \varrho'.\]
\noindent If there exists a $2$-cell $\gamma^{-1}$, we can compose on the left by $p'\cdot\gamma$ and $q'\cdot\gamma$, respectively, to obtain
\[ (p'\cdot\gamma^{-1})\varpi' = \varpi \;\textrm{and}\; (q'\cdot\gamma^{-1})\varrho' = \varrho.\]
 \noindent It follows that $\gamma^{-1} \maps (\varpi',f',\varrho')\To (\varpi,f,\varrho)$  is a map of maps of spans as well.

\begin{remark}
The maps of spans in Definition~\ref{2morphisms} are the $2$-morphisms in our construction.  Other occurrences of spans in the literatire require a more general definition of maps between spans in which, given parallel spans, a map between them consists of an object together with morphisms to each of the objects of the two spans, and invertible $2$-cells making the diagram commute up to isomorphism.  See, for example,~\cite{Mo13}.
\end{remark}

\begin{remark}\label{subsec-2cell note}
We make a brief remark on $2$-cells in maps of spans.  This paper originated from an interest in monoidal bicategories of spans of groupoids in the context of categorified representation theory in the {\em groupoidification program}~\cite{BaHoWa,Ho}.  At the time, we did not have a pressing need for the $3$-dimensional structure of spans.  Since the examples of interest were bicategories of spans of sets and bicategories of spans of groupoids, it was easy to believe that bicategories of spans could be constructed from any $2$-category with pullbacks.

We can view the category of sets as a $2$-category in which all $2$-cells are taken to be identities, and the $2$-cells in the $2$-category of groupoids are natural transformations, whose component $1$-cells are morphisms of groupoids, which are, of course, invertible.  This commonality turns out to be essential to the construction of spans in a bicategory in which composition is defined by the pullback; specifically, the iso-comma object.  If the $2$-cell components of maps of spans are not required to be {\em invertible}, then composition may not be defined.  We give a counterexample to illustrate the point.

Consider the $2$-category $\Cat$ of small categories, functors, and natural transformations.  Let $R$ be the category with objects $\{A,B\}$, morphisms $\{1_A,f\maps A\to B,1_B\}$, and define functors $a\maps R\to R$ and $b\maps R\to R$, which send all objects to $A$ and $B$, respectively.  Consider the following pair of maps of spans, where $1$ is the terminal category and $\sigma\maps a\to b1_R$ is the natural transformation that assigns $f$ to $A$ and $1_B$ to $B$.
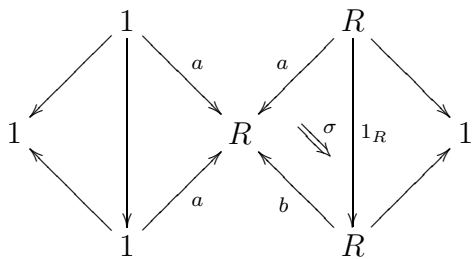
\[
   \xy
   (-30,0)*+{1}="1";
   (-15,15)*+{1}="2";
   (-15,-15)*+{1}="6";
   (0,0)*+{R}="3";
   (15,15)*+{R}="4";
   (15,-15)*+{R}="7";
   (30,0)*+{1}="5";
        {\ar_{} "2";"1"};
        {\ar^{a} "2";"3"};
        {\ar^{} "6";"1"};
        {\ar_{a} "6";"3"};
        {\ar_{a} "4";"3"};
        {\ar^{} "4";"5"};
        {\ar^{b} "7";"3"};
        {\ar_{} "7";"5"};
        {\ar^{} "2";"6"};
		{\ar^{1_R} "4";"7"};
        {\ar@{=>}^{\scriptstyle \sigma} (8,1); (12,-3)};
\endxy
\]

\noindent The $2$-cells in the pullback are required to be invertible, so they must have only invertible components, i.e., they must be identity natural transformations in this case.  Then the composite object on the top is the terminal category and the composite object on the bottom is the initial category.  There are no maps from the terminal category to the initial category, so it is not possible to construct a well-defined composition functor.  Then, to obtain a bicategory or a tricategory $\Span(\Cat)$, in the sense we desire, the definition of maps of spans should include an invertibility condition on the $2$-cells as in Definition~\ref{2morphisms}.
\end{remark}

\subsection{Composition Operations on Spans}\label{comp_spans}

We describe the composition operations in this section using the pullback construction in Definition~\ref{iso_comma}.  

\begin{definition}
For each cospan diagram
\[
   \xy
   (0,-15)*+{S}="1";
   (-15,0)*+{B}="2"; 
   (15,0)*+{R}="3";
        {\ar_{r} "2";"1"};
        {\ar^{q} "3";"1"};
\endxy
\]
\noindent in $B$, we choose (and denote) a limit object $SR$, projection morphisms
\[  \pi_R^S\maps SR\to R\;\;\;\; \pi_S^R\maps SR\to S,\]
\noindent and an invertible $2$-cell
\[\kappa_{B}^{r,q}\maps q\pi_{R}^{S}\To r\pi_{S}^{R}.\]
\noindent This data is called the {\bf pullback of the cospan}.
\end{definition}

\noindent  For convenience, we often refer to the object $SR$ itself as the pullback.  The chosen pullback data are fixed for the duration of the paper and the above notation will be reserved only for this data.

We think of a span as a $1$-morphism from an object $A$ to an object $B$ in $\Span(\B)$.  

\begin{definition}\label{horizcompspans}
Given a pair of composable spans
\[
   \xy
   (-15,0)*+{S}="2"; (15,0)*+{R}="3";
   (-30,-15)*+{C}="4";   (0,-15)*+{B}="5";   (30,-15)*+{A}="6";
        {\ar_{s} "2";"4"};
        {\ar^{r} "2";"5"};
        {\ar_{q} "3";"5"};
        {\ar^{p} "3";"6"};
\endxy
\]
\noindent there is a {\bf composite spans} $SR$ from $A$ to $C$
\[
   \xy
   (0,15)*+{SR}="1";
   (-15,0)*+{S}="2"; (15,0)*+{R}="3";
   (-30,-15)*+{C}="4";   (0,-15)*+{B}="5";   (30,-15)*+{A}="6";
        {\ar_{\pi_{S}^{R}} "1";"2"};
        {\ar^{\pi_{R}^{S}} "1";"3"};
        {\ar_{s} "2";"4"};
        {\ar_{r} "2";"5"};
        {\ar^{q} "3";"5"};
        {\ar^{p} "3";"6"};
        {\ar@{=>}_{\scriptstyle \kappa_{B}^{r,q}} (2,0); (-2,0)};
\endxy
\]
\noindent formed by the pullback.
\end{definition}

We use the pullback to define composition of maps of composable spans.

\begin{definition}\label{horizcompmaps}
Given a pair of maps of spans between composable pairs of spans, we define the {\bf horizontal composite of maps of spans}

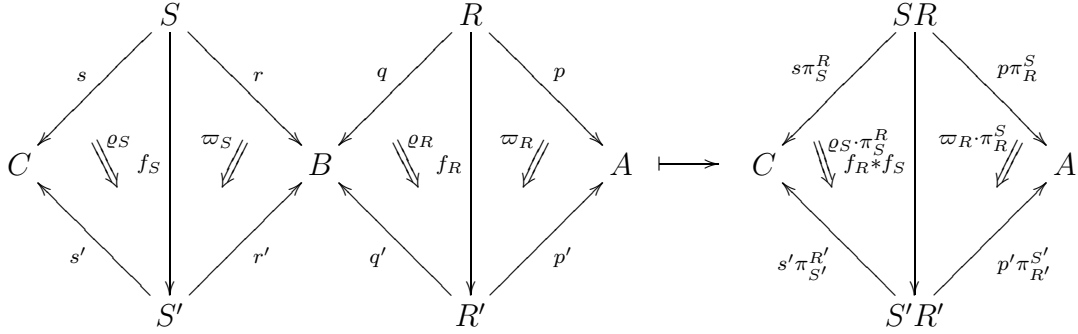
\[
   \xy
   (-50,20)*+{S}="1";
   (-70,0)*+{C}="2";
    (-30,0)*+{B}="4";
   (-50,-20)*+{S'}="5";
   (-10,20)*+{R}="-1";
    (10,0)*+{A}="-4";
   (-10,-20)*+{R'}="-5";
   (14,0)*+{}="6"; (24,0)*+{}="7";   
   (49,20)*+{SR}="8";
   (29,0)*+{C}="9"; (69,0)*+{A}="11";
   (49,-20)*+{S'R'}="12";
        {\ar@{|->}_{} "6";"7"}; 
        {\ar_{s} "1";"2"};
        {\ar^{r} "1";"4"};
        {\ar^{s'} "5";"2"};
        {\ar_{r'} "5";"4"}; 
        {\ar_{f_S} "1";"5"};
        {\ar_{q} "-1";"4"};
        {\ar^{p} "-1";"-4"};
        {\ar^{q'} "-5";"4"};
        {\ar_{p'} "-5";"-4"}; 
        {\ar_{f_R} "-1";"-5"};
        {\ar_{s\pi_S^R} "8";"9"};
        {\ar^{p\pi_R^S} "8";"11"};
        {\ar^{s'\pi_{S'}^{R'}} "12";"9"};
        {\ar_{p'\pi_{R'}^{S'}} "12";"11"}; 
        {\ar_{f_R*f_S} "8";"12"};
        {\ar@{=>}_<<{\scriptstyle \varpi_{R}} (0,3); (-3,-3)};
        {\ar@{=>}_<<{\scriptstyle \varpi_{S}} (-40,3); (-43,-3)};
        {\ar@{=>}^<<{\scriptstyle \varrho_{R}} (-20,3); (-17,-3)};
        {\ar@{=>}^<<{\scriptstyle \varrho_{S}} (-60,3); (-57,-3)};
        {\ar@{=>}_<<{\scriptstyle \varpi_{R}\cdot\pi_{R}^{S}} (63,3); (60,-3)};
        {\ar@{=>}^<<{\scriptstyle  \varrho_{S}\cdot\pi_{S}^{R}} (36,3); (38,-3)};
\endxy
\]
\noindent by the assignment
\[((\varpi_R,f_R,\varrho_R),(\varpi_S,f_S,\varrho_S)) \mapsto (\varpi_{R}\cdot\pi^{S}_{R}, \; f_R*f_S, \;\varrho_{S}\cdot \pi^{R}_{S}),\]
\noindent where
\[ f_R*f_S\maps SR\to S'R'\]
\noindent is the unique $1$-cell satisfying
\[ \pi_{R'}^{S'}(f_R*f_S) = f_R\pi_{R}^{S},\;\;\; \pi_{S'}^{R'}(f_R*f_S) = f_S\pi_{S}^{R},\;\textrm{ and }\;\kappa_{B}^{r'.q'}\cdot (f_R*f_S) =  (\varpi_{S}\cdot\pi_{S}^{R})\kappa_{B}^{r,q}(\varrho_{R}^{-1}\cdot \pi_{R}^{S}).\]
\end{definition}

The pullback also allows us to define the composite of maps of maps of spans.  (We omit the $2$-cell components of the maps of spans ($\varpi_S,f_S,\varrho_S$), etc., from the diagrams below for aesthetic purposes.)  In the following definition, we need to apply the universal property of pullbacks to define the necessary $2$-cell.  This requires that the equation
\[ (r'\cdot\gamma_{S}\cdot\pi_{S}^{R})(\kappa_{B}^{r',q'}\cdot(f_R*f_S)) = (\kappa_{B}^{r',q'}\cdot(f'_R*f'_S))(q'\cdot\gamma_R\cdot\pi_{R}^{S})\]
\noindent holds, which is straightforward to verify, although we leave the details to the reader.

\begin{definition}\label{horizcompmapsofmaps}
Given a pair of maps of maps of spans
\[
   \xy
   (-50,0)*+{\gamma_R\maps (\varpi_{R}, f_{R}, \varrho_{R})}="1"; (-18,0)*+{(\varpi'_{R}, f'_R, \varrho'_{R})}="2";
   (18,0)*+{\gamma\maps (\varpi,f,\varrho)}="3"; (43,0)*+{(\varpi',f',\varrho')}="4";
   (0,0)*+{\textrm{and}}="0";
        {\ar@3{->}_{} "1";"2"};
        {\ar@3{->}_{} "3";"4"};
\endxy
\]
\noindent between pairs of composable maps of spans, we define the {\bf horizontal composite of maps of maps of spans}
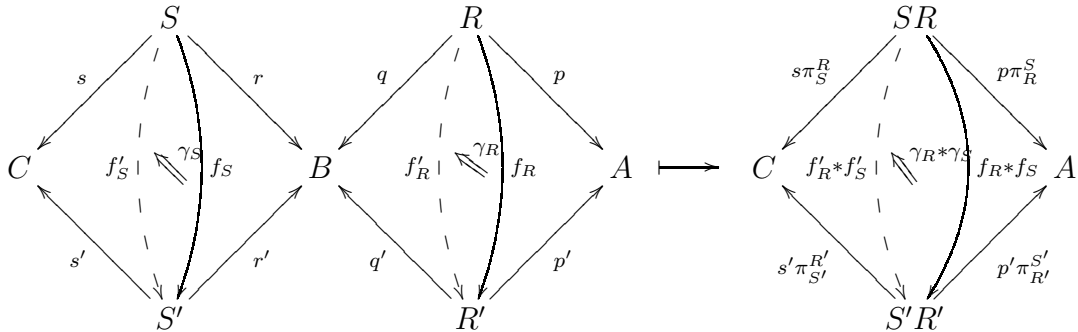
\[
   \xy
   (-50,20)*+{S}="1";
   (-70,0)*+{C}="2";
    (-30,0)*+{B}="4";
   (-50,-20)*+{S'}="5";
   (-10,20)*+{R}="-1";
    (10,0)*+{A}="-4";
   (-10,-20)*+{R'}="-5";
   (14,0)*+{}="6"; (24,0)*+{}="7";   
   (49,20)*+{SR}="8";
   (29,0)*+{C}="9"; (69,0)*+{A}="11";
   (49,-20)*+{S'R'}="12";
        {\ar@{|->}_{} "6";"7"}; 
        {\ar_{s} "1";"2"};
        {\ar^{r} "1";"4"};
        {\ar^{s'} "5";"2"};
        {\ar_{r'} "5";"4"}; 
        {\ar@/^1pc/^{f_S} "1";"5"};
         {\ar@{-->}@/_1pc/_{f'_S} "1";"5"};
        {\ar_{q} "-1";"4"};
        {\ar^{p} "-1";"-4"};
        {\ar^{q'} "-5";"4"};
        {\ar_{p'} "-5";"-4"}; 
        {\ar@/^1pc/^{f_R} "-1";"-5"};
         {\ar@{-->}@/_1pc/_{f'_R} "-1";"-5"};
        {\ar_{s\pi_S^R} "8";"9"};
        {\ar^{p\pi_R^S} "8";"11"};
        {\ar^{s'\pi_{S'}^{R'}} "12";"9"};
        {\ar_{p'\pi_{R'}^{S'}} "12";"11"}; 
        {\ar@/^1.7pc/^{f_R*f_S} "8";"12"};
         {\ar@{-->}@/_1.2pc/_{f'_R*f'_S} "8";"12"};
         {\ar@{=>}_{\scriptstyle \gamma_S} (-48,-2); (-52,2)};
        {\ar@{=>}_{\scriptstyle \gamma_R} (-8,-); (-12,2)};
        {\ar@{=>}_>>>{\scriptstyle \gamma_R *\gamma_S} (49,-2); (46,2)};
\endxy
\]
\noindent by the assignment
\[
   \xy
   (-35,0)*+{(\gamma_R, \gamma_S) \mapsto \gamma_R*\gamma_S\maps (\varpi_{R}\cdot \pi_{R}^{S}, f_R*f_S, \varrho_{S}\cdot \pi_{S}^{R})}="1"; (35,0)*+{(\varpi'_{R}\cdot \pi_{R}^{S}, f'_R*f'_S, \varrho'_{S}\cdot \pi_{S}^{R}),}="2";
        {\ar@3{->}_{} "1";"2"};
\endxy
\]
\noindent  where $\gamma_R*\gamma_S$ is the unique $2$-cell in $\B$ satisfying
\[ \pi_{R'}^{S'}\cdot (\gamma_R*\gamma_S) = \gamma_R\cdot\pi_{R}^{S} \;\textrm{ and }\;\pi^{R'}_{S'}\cdot(\gamma_R*\gamma_S) = \gamma_S\cdot\pi^{R}_{S}.\] 
\end{definition}

\noindent From the defining equations of the $2$-cell $\gamma_R*\gamma_S$, we have the equations
\[ (p'\pi^{S'}_{R'}\cdot(\gamma_R*\gamma_S))(\varpi_R\cdot\pi_R^S) = \varpi'_R\cdot\pi_R^S\]
\noindent and 
\[ (s'\pi_{S'}^{R'}\cdot(\gamma_R*\gamma_S))(\varrho_S\cdot\pi_S^R) = \varrho'_S\cdot\pi^R_S,\]

\noindent verifying that $\gamma_R*\gamma_S$, indeed, defines a map of maps of spans.

The three composition operations defined above were each induced by the universal property of the pullback.  There is a second set of composition operations obtained via the composition operations of $\B$, which we now define.

\begin{definition}\label{vertcompmaps}
Given a pair of composable maps of spans
\[ (\varpi_S, f_S, \varrho_S)\maps S\to S'\;\;\textrm{and}\;\;(\varpi'_{S}, f'_{S}, \varrho'_{S})\maps S'\to S''\]
\noindent between parallel spans, we define the {\bf vertical composite of maps of spans}
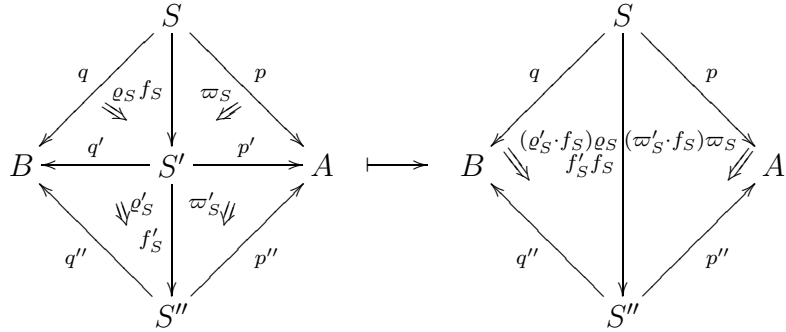
\[
   \xy
   (-30,20)*+{S}="1";
   (-50,0)*+{B}="2"; (-30,0)*+{S'}="3";  (-10,0)*+{A}="4";
   (-30,-20)*+{S''}="5";
   (-5,0)*+{}="6"; (5,0)*+{}="7";   
   (30,20)*+{S}="8";
   (10,0)*+{B}="9"; (50,0)*+{A}="11";
   (30,-20)*+{S''}="12";
        {\ar@{|->}_{} "6";"7"}; 
        {\ar_{q} "1";"2"};
        {\ar^{p} "1";"4"};
        {\ar_{q'} "3";"2"};
        {\ar^{p'} "3";"4"};
        {\ar^{q''} "5";"2"};
        {\ar_{p''} "5";"4"}; 
        {\ar_{f_S} "1";"3"};
        {\ar_{f'_{S}} "3";"5"};
        {\ar_{q} "8";"9"};
        {\ar^{p} "8";"11"};
        {\ar^{q''} "12";"9"};
        {\ar_{p''} "12";"11"}; 
        {\ar_{f'_{S}f_S} "8";"12"};
        {\ar@{=>}_<<{\scriptstyle \varpi_{S}} (-21,8); (-24,6)};
        {\ar@{=>}_<<{\scriptstyle \varpi'_{S}} (-22,-5); (-23,-8)};
        {\ar@{=>}^<<{\scriptstyle \varrho_{S}} (-39,8); (-36,6)};
        {\ar@{=>}^<<{\scriptstyle \varrho'_{S}} (-37,-5); (-36,-8)};
        {\ar@{=>}_<<{\scriptstyle (\varpi'_S\cdot f_{S})\varpi_{S}} (47,2); (44,-2)};
        {\ar@{=>}^<<{\scriptstyle (\varrho'_S\cdot f_{S})\varrho_{S}} (14.5,2); (17.5,-2)};
\endxy
\]

\noindent by the assignment
\[ ((\varpi_S,f_S,\varrho_S),(\varpi'_S,f'_S,\varrho'_S)) \mapsto ((\varpi'_S\cdot f_S)\varpi_S, f'_Sf_S, (\varrho'_S\cdot f_S)\varrho_S).\]
\end{definition}

\begin{definition}\label{vertcompmapsofmaps}
Given a pair of composable maps of maps of spans
\[
   \xy
   (-53,0)*+{\sigma\maps (\varpi,f,\varrho)}="1"; (-22,0)*+{(\varsigma,g,\varphi)\maps S\to S'}="2";
   (18,0)*+{\sigma'\maps (\varpi',f',\varrho')}="3"; (53,0)*+{(\varsigma',g',\varphi')\maps S'\to S'',}="4";
   (0,0)*+{\textrm{and}}="0";
        {\ar@3{->}_{} "1";"2"};
        {\ar@3{->}_{} "3";"4"};
\endxy
\]
\noindent we define the {\bf vertical composite of maps of maps of spans}
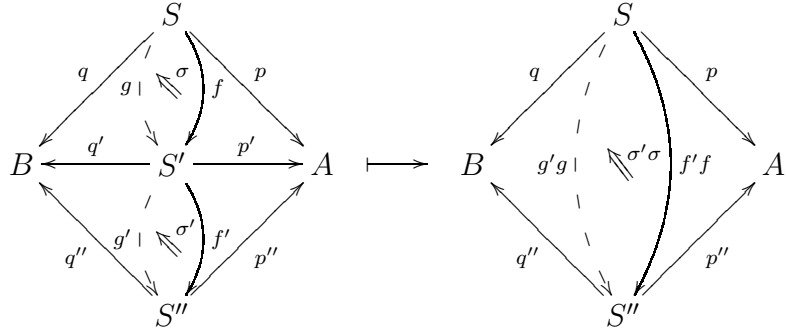
\[
   \xy
   (-30,20)*+{S}="1";
   (-50,0)*+{B}="2"; (-30,0)*+{S'}="3";  (-10,0)*+{A}="4";
   (-30,-20)*+{S''}="5";
   (-5,0)*+{}="6"; (5,0)*+{}="7";   
   (30,20)*+{S}="8";
   (10,0)*+{B}="9"; (50,0)*+{A}="11";
   (30,-20)*+{S''}="12";
        {\ar@{|->}_{} "6";"7"}; 
        {\ar_{q} "1";"2"};
        {\ar^{p} "1";"4"};
        {\ar_{q'} "3";"2"};
        {\ar^{p'} "3";"4"};
        {\ar^{q''} "5";"2"};
        {\ar_{p''} "5";"4"}; 
        {\ar@/^1pc/^{f} "1";"3"};
        {\ar@/^1pc/^{f'} "3";"5"};
        {\ar@{-->}@/_1pc/_{g} "1";"3"};
        {\ar@{-->}@/_1pc/_{g'} "3";"5"};
        {\ar_{q} "8";"9"};
        {\ar^{p} "8";"11"};
        {\ar^{q''} "12";"9"};
        {\ar_{p''} "12";"11"}; 
        {\ar@/^1.5pc/^{f'f} "8";"12"};
        {\ar@{-->}@/_1.5pc/_{g'g} "8";"12"};
        {\ar@{=>}_{\scriptstyle \sigma} (-29,9); (-32,12)};
        {\ar@{=>}_{\scriptstyle \sigma'} (-29,-12); (-32,-9)};
       {\ar@{=>}_{\scriptstyle \sigma'\sigma} (31,-2); (28,2)};
\endxy
\]
\noindent by the assignment
\[
   \xy
   (-31,0)*+{(\sigma, \sigma')\mapsto \sigma'\sigma\maps ((\varpi'\cdot f)\varpi, f'f, (\varrho'\cdot f)\varrho)}="1"; (31,0)*+{((\varsigma'\cdot g)\varsigma, g'g, (\varphi'\cdot g)\varphi),}="2";
        {\ar@3{->}_{} "1";"2"};
\endxy
\]
\noindent where $\sigma'\sigma\maps f'f\To g'g$ is the horizontal composite of $2$-cells in $\B$.  (We suppress the $2$-cell components of the maps of spans for aesthetic purposes.)
\end{definition}

\noindent That $\sigma'\sigma$ is a map of maps of spans follows from the interchange law and that $\sigma$ and $\sigma'$ are maps of maps of spans.  We have
\[ (\varsigma'\cdot f)\varsigma = ((p''\cdot\sigma')\varpi)\cdot f)((p'\cdot\sigma)\varpi) = p''\cdot(\sigma'\sigma)((\varpi'\cdot f)\varpi)\]
\noindent and
\[ (\varphi'\cdot g)\varphi = ((q''\cdot\sigma')\varrho)\cdot g)((q'\cdot\sigma)\varrho) = q''\cdot(\sigma'\sigma)((\varrho'\cdot g)\varrho).\]

\medskip

\begin{definition}\label{horizcompparallelmaps}
Given a pair of composable maps
\[\sigma\maps (\varpi,f,\varrho)\To (\varpi',f',\varrho')\;\;\textrm{ and } \;\;\sigma'\maps (\varpi',f',\varrho')\To (\varpi'',f'',\varrho'')\]
\noindent between parallel maps of spans from $S$ to $S'$, we define the {\bf horizontal composite of maps of parallel maps of spans}
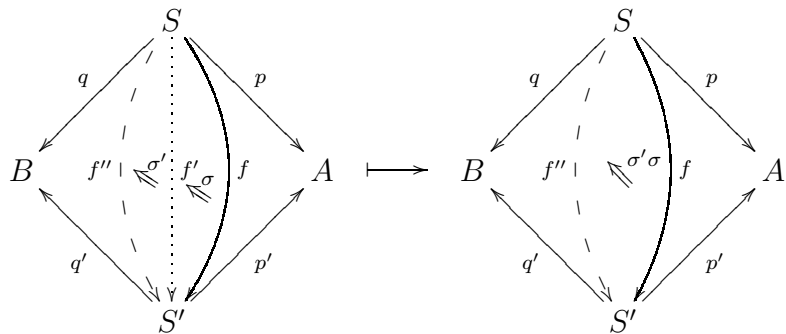
\[
   \xy
   (-30,20)*+{S}="1";
   (-50,0)*+{B}="2";  (-10,0)*+{A}="4";
   (-30,-20)*+{S'}="5";
   (-5,0)*+{}="6"; (5,0)*+{}="7";   
   (30,20)*+{S}="8";
   (10,0)*+{B}="9"; (50,0)*+{A}="11";
   (30,-20)*+{S'}="12";
        {\ar@{|->}_{} "6";"7"}; 
        {\ar_{q} "1";"2"};
        {\ar^{p} "1";"4"};
        {\ar^{q'} "5";"2"};
        {\ar_{p'} "5";"4"}; 
        {\ar@/^1.8pc/^{f} "1";"5"};
        {\ar@{..>}@/^0pc/^{f'} "1";"5"};
        {\ar@{-->}@/_1.6pc/_{f''} "1";"5"};
        {\ar_{q} "8";"9"};
        {\ar^{p} "8";"11"};
        {\ar^{q'} "12";"9"};
        {\ar_{p'} "12";"11"}; 
        {\ar@/^1.5pc/^{f} "8";"12"};
        {\ar@{-->}@/_1.5pc/_{f''} "8";"12"};
        {\ar@{=>}_{\scriptstyle \sigma} (-25,-4); (-28,-2)};
        {\ar@{=>}_{\scriptstyle \sigma'} (-32,-2); (-35,0)};
       {\ar@{=>}_{\scriptstyle \sigma'\sigma} (31,-2); (28,1)};
\endxy
\]

\noindent by the assignment
\[(\sigma,\sigma')\mapsto \sigma'\sigma\maps (\varpi,f,\varrho)\To(\varpi'',f'',\varrho''),\]
\noindent where $\sigma'\sigma\maps f\To f''$ is the vertical composite of $2$-cells in $\B$.  (We suppress the $2$-cell components of the maps of spans for aesthetic purposes.)
\end{definition}

\noindent It is straightforward to check that $\sigma'\sigma$ is a map of maps of spans.   We have
\[ (p'\cdot\sigma'\sigma)\varpi = (p'\cdot\sigma')(p'\cdot\sigma) \varpi = (p'\cdot\sigma')\varpi' = \varpi''\]
\noindent and
\[ (q'\cdot\sigma'\sigma)\varrho = (q'\cdot\sigma')(q'\cdot\sigma)\varrho = (q'\cdot\sigma')\varrho' = \varrho''.\]

\subsection{Identity Morphisms}

For each composition operation in the tricategory $\Span(\B)$, we define an identity morphism.  

\begin{definition}\label{identityspan}
Given an object $A\in\B$, the {\bf identity span} is the diagram
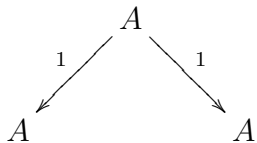
\[
   \xy
   (0,15)*+{A}="1";
   (-15,0)*+{A}="2"; 
   (15,0)*+{A}="3";
        {\ar_{1} "1";"2"};
        {\ar^{1} "1";"3"};
\endxy
\]
\end{definition}

\noindent This is the identity for the horizontal composition operation on spans in Definition~\ref{horizcompspans}.

\begin{definition}\label{identitymap}
Given a span 
\[
   \xy
   (0,15)*+{S}="1";
   (-15,0)*+{B}="2"; 
   (15,0)*+{A}="3";
        {\ar_{q} "1";"2"};
        {\ar^{p} "1";"3"};
\endxy
\]
\noindent in $\B$, the {\bf identity map of spans} consists of the identity $1$-cell for $S$ and the identity $2$-cells $p\To p1_S$ and $q\To q1_S$ in $\B$.
\end{definition}

\noindent This is the identity for both the horizontal and the vertical composition operations on maps of spans in Definition~\ref{horizcompmaps} and Definition~\ref{vertcompmaps}, respectively.  Note that $p1_S = p$ and $q1_S = q$, since $\B$ is a strict $2$-category.

\begin{definition}\label{identitymapofmaps}
Given a map of spans $(\varpi,f,\varrho)$, the {\bf identity map of maps of spans} consists of the identity $2$-cell $1_f\maps f\To f$ in $\B$.
\end{definition}

\noindent This is the identity for the horizontal and vertical composition operations on maps of maps of spans in Definition~\ref{horizcompmapsofmaps} and Definition~\ref{vertcompmapsofmaps}, respectively, and the horizontal composition operation on maps of parallel spans in Definition~\ref{horizcompparallelmaps}.

\subsection{Strict $\Hom$-$2$-Categories}

We describe the local structure of the tricategory of spans.

\proposition \label{hom2category}
 For each pair of objects $A,B\in\B$, there is a strict $2$-category $\Span(\B)(A,B)$ consisting of
\begin{itemize}
\item spans from $A$ to $B$ as objects (see Definition~\ref{1morphisms}), 
\item for each pair $R,S$ of spans from $A$ to $B$, a category $\Span(\B)(A,B)(R,S)$, consisting of
\begin{itemize}
\item maps of spans (see Definition~\ref{2morphisms}),
\item maps of maps of spans (see Definition~\ref{3morphisms}),
\item a composition operation on maps between parallel maps of spans (see Definition~\ref{horizcompparallelmaps}),
\item the identity map of maps of spans (see Definition~\ref{identitymapofmaps}),
\end{itemize}
\item for each triple $R,S,T$ of spans from $A$ to $B$, a composition functor
\[ *_{v}\maps \Span(\B)(A,B)(R,S)\times\Span(\B)(A,B)(S,T)\to\Span(\B)(A,B)(R,T),\]
\noindent consisting of
\begin{itemize}
\item a vertical composition operation on maps of spans (see Definition~\ref{vertcompmaps}),
\item a vertical composition operation on maps of maps of spans (see Definition~\ref{vertcompmapsofmaps}),
\end{itemize}
\item and, for each span $R$ from $A$ to $B$, a unit functor, consisting of the corresponding identity map of spans and identity map of maps of spans (see Definitions~\ref{identitymap} and~\ref{identitymapofmaps}).
\end{itemize}
\endproposition

\proof
Local composition is defined by vertical composition of $2$-cells in the strict $2$-category $\B$.  It follows that $\Span(\B)(A,B)(R,S)$ is a category.  Functoriality of composition of maps of maps of spans follows from the interchange law for horizontal and vertical composition of $2$-cells in the $2$-category $\B$.  Preservation of identities is immediate.  The axioms are straightforward from the associative and unital composition of $1$-cells and $2$-cells in $\B$.
\endproof

\begin{remark}
We will often write the vertical composites as concatenation suppressing the symbol $*_{v}$.
\end{remark}

\subsection{Strict $2$-Functors}

We define {\em composition} and {\em unit} strict $2$-functors between strict $\hom$-$2$-categories.

\subsubsection*{Composition $2$-Functor}

We first define horizontal composition.

\proposition \label{compositionfunctor}
For each triple of objects $A, B, C\in\B$, there is a strict $2$-functor
\[ *_{h}\maps \Span(A,B)\times \Span(B,C)\to\Span(A,C),\]
\noindent consisting of
\begin{itemize}
\item horizontal composition of spans (see Definition~\ref{horizcompspans}),
\item horizontal composition of maps of spans (see Definition~\ref{horizcompmaps}), and
\item horizontal composition of maps of maps of spans (see Definition~\ref{horizcompmapsofmaps}).
\end{itemize}
\endproposition

\proof
We need to check that horizontal composition operation on maps of maps of spans preserves vertical composition and identities, and that the naturality equations expressing functoriality for composition of maps of spans hold.  The axioms of a $2$-functor are immediate since $\B$ is strict, i.e, the associator natural isomorphisms are identities.

Functoriality of composition of maps of maps of spans is immediate, i.e., given horizontally composable pairs of vertically composable pairs of maps of maps of spans, $(\sigma_{R}, \tau_{R})$ and $(\sigma_{S}, \tau_{S})$, the equations
\[ (\sigma_{S}*_{v}\sigma_{R}) *_{h} (\tau_{S}*_{v}\tau_{R}) = (\tau_{R}*_{h}\tau_{S})*_{v}(\sigma_{R}*_{h}\sigma_{S})\]
\noindent and
\[ 1_{R}*_{h}1_{S} = 1_{RS}\]
\noindent hold.

We check that vertical composition of maps of spans is preserved by horizontal composition of maps of spans.  Consider the maps of spans
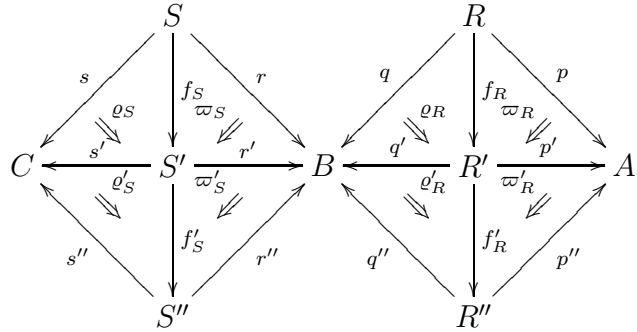
\[
   \xy
   (-20,20)*+{S}="1";    (20,20)*+{R}="2";
   (-40,0)*+{C}="3"; (0,0)*+{B}="4";  (40,0)*+{A}="5";
   (-20,0)*+{S'}="6";    (20,0)*+{R'}="7"; 
   (-20,-20)*+{S''}="8";    (20,-20)*+{R''}="9";
        {\ar_{s} "1";"3"};
        {\ar^{r} "1";"4"};
        {\ar_{q} "2";"4"};
        {\ar^{p} "2";"5"};
        {\ar_{s'} "6";"3"};
        {\ar^{r'} "6";"4"};
        {\ar_{q'} "7";"4"};
        {\ar^{p'} "7";"5"};
        {\ar^{s''} "8";"3"};
        {\ar_{r''} "8";"4"};
        {\ar^{q''} "9";"4"};
        {\ar_{p''} "9";"5"};
        {\ar^{f_{R}} "2";"7"};
        {\ar^{f'_{R}} "7";"9"};
        {\ar^{f_{S}} "1";"6"};
        {\ar^{f'_{S}} "6";"8"};
        {\ar@{=>}_<<{\scriptstyle \varpi_{R}} (30,6); (27,3)};
        {\ar@{=>}_<<{\scriptstyle \varpi'_{R}} (30,-4); (27,-7)};
        {\ar@{=>}^<<{\scriptstyle \varrho_{R}} (11,6); (14,3)};
        {\ar@{=>}^<<{\scriptstyle \varrho'_{R}} (11,-4); (14,-7)};
        {\ar@{=>}_<<{\scriptstyle \varpi_{S}} (-11,6); (-14,3)};
        {\ar@{=>}_<<{\scriptstyle \varpi'_{S}} (-11,-4); (-14,-7)};
        {\ar@{=>}^<<{\scriptstyle \varrho_{S}} (-30,6); (-27,3)};
        {\ar@{=>}^<<{\scriptstyle \varrho'_{S}} (-30,-4); (-27,-7)};
\endxy
\]
\noindent Beginning with the vertical composite followed by the horizontal composite, we have
\[ ((\varpi'_R\cdot f_R\pi_{R}^{S})(\varpi_R\cdot \pi_{R}^{S}),\; (f'_R*_{v}f_R) *_{h} (f'_S*_{v}f_S),\; (\varrho'_S\cdot f_S\pi_{S}^{R})(\varrho_S\cdot\pi_{S}^{R}))\]

\noindent Beginning with the horizontal composite followed by the vertical composite, we have
\[ ((\varpi'_R\cdot\pi_{R'}^{S'}(f_{R}*f_{S}))(\varpi_R\cdot\pi_{R}^{S}), \;(f'_R*_{h}f'_S)*_{v}(f_R*_{h}f_S),\; (\varrho'_S\cdot\pi_{S'}^{R'}(f_R*f_S))(\varrho_S\cdot\pi_{S}^{R})).\]

\noindent These maps of spans are equal, thus composition is preserved on the nose.  Similarly, identity maps of spans are preserved.  It follows from these naturality equations of composites of maps of spans that horizontal composition is a strict $2$-functor on $\hom$-$2$-categories.
\endproof

\begin{remark}
We will often write the horizontal composite $*_{h}$ simply as $*$.  This should not cause confusion with the vertical composite, which will usually be written as concatenation.
\end{remark}

\subsubsection*{Unit $2$-Functor}

\proposition \label{unitfunctor}
For each object $A\in\B$, there is a strict $2$-functor
\[ I_{A}\maps {\bf 1}\to \Span(A,A),\]
\noindent which consists only of the identity span, identity map of spans, and the identity map of maps of spans (see Definitions~\ref{identityspan}, ~\ref{identitymap}, and~\ref{identitymapofmaps}).
\endproposition

\proof
Straightforward.
\endproof

\subsection{Adjoint Equivalence Transformations}

We define associativity and unit adjoint equivalences in $2$-categories $[\C,\C']$ of maps between $2$-categories.  These $2$-categories of maps are strict if the codomain $2$-category $\C'$ is strict.  Since the $\hom$-$2$-categories of the span construction are all strict, then we consider only adjoint equivalences in strict $2$-categories, which simplifies the bicategorical triangle axioms of an adjoint equivalence.  Each internal adjunction consists of an adjoint pair of transformations together with both counit and unit modifications.  In case these modifications are trivial, then we say the adjoint equivalence is {\em strict}.  

\subsubsection*{Associator}

To define the components of the associator transformation, first consider the diagram
\[
   \xy
   (-45,0)*+{D}="0";
   (-30,15)*+{T}="1";
   (-30,-15)*+{T}="7";
   (-15,0)*+{C}="2";
   (0,15)*+{S}="3";
   (0,-15)*+{S}="8";
   (15,0)*+{B}="4";
   (30,15)*+{R}="5";
   (30,-15)*+{R}="9";
   (45,0)*+{A}="6";
   (15,30)*+{SR}="10";
   (0,45)*+{T(SR)}="11";
   (-15,-30)*+{TS}="12";
   (0,-45)*+{(TS)R}="13";
        {\ar_{s} "1";"0"};
        {\ar^{r} "1";"2"};
        {\ar_{q} "3";"2"};
        {\ar^{p} "3";"4"};
        {\ar_{n} "5";"4"};
        {\ar^{m} "5";"6"};
        {\ar^{s} "7";"0"};
        {\ar_{r} "7";"2"};
        {\ar^{q} "8";"2"};
        {\ar_{p} "8";"4"};
        {\ar^{n} "9";"4"};
        {\ar_{m} "9";"6"};
        {\ar_{\pi^{R}_{S}} "10";"3"};
        {\ar^{\pi^{S}_{R}} "10";"5"};
        {\ar_{\pi^{SR}_{T}} "11";"1"};
        {\ar^{\pi^{T}_{SR}} "11";"10"};
        {\ar^{\pi^{S}_{T}} "12";"7"};
        {\ar_{\pi^{T}_{S}} "12";"8"};
        {\ar^{\pi^{R}_{TS}} "13";"12"};
        {\ar_{\pi^{TS}_{R}} "13";"9"};
        {\ar@{=>}_<<{\scriptstyle \kappa_{B}^{p,n}} (17,15); (13,15)};
        {\ar@{=>}_<<{\scriptstyle \kappa_{C}^{r,q\pi_{S}^{R}}} (-13,15); (-17,15)};
        {\ar@{=>}_<<{\scriptstyle \kappa_{B}^{p\pi_{S}^{T},n}} (17,-15); (13,-15)};
        {\ar@{=>}_<<{\scriptstyle \kappa_{C}^{r,q}} (-13,-15); (-17,-15)};
\endxy
\]

\noindent Applying the universal property to induce the unique $1$-cell $a\maps T(SR)\to TS$ satisfying the equations
\[\pi^{T}_{S}a = \pi^{R}_{S}\pi^{T}_{SR},\;\;\pi^{S}_{T}a = \pi^{SR}_{T},\;\textrm{ and }\;\kappa_{C}^{r,q}\cdot a = \kappa_{C}^{r,q\pi_{S}^{R}}.\]
\noindent Similarly, applying the universal property to induce the unique $1$-cell $a^{-1}\maps (TS)R\to SR$ satisfying the equations
\[ \pi_{R}^{S}\cdot a^{-1} = \pi_{R}^{TS},\;\; \pi_{S}^{R}\cdot a^{-1} = \pi_{S}^{T}\pi_{TS}^{R},\;\;\textrm{ and }\;\;\kappa_{B}^{p,n}\cdot a^{-1} = \kappa_{B}^{p\pi_{S}^{T},n}.\]
\noindent Another application of the universal property yields the desired $1$-cell components of the transformation and its inverse in the following proposition.

\proposition \label{associatortransformation}
For each $4$-tuple of objects $A, B, C, D\in\B$ there is an strict adjoint equivalence $({\bf a},\; {\bf a}^{-1},\; 1,\; 1)$ in the strict $2$-category $[\Span(A,B)\times\Span(B,C)\times\Span(B,D),\Span(A,D)]$, with
\begin{itemize}
\item strict transformations 
\[ {\bf a}\maps *_{(AB)CD}(*_{ABC} \times 1_{D})\To *_{A(BC)D}(1_{A}\times *_{BCD}),\]
\noindent and
\[ {\bf a}^{-1}\maps *_{A(BC)D}(1_{A}\times *_{BCD})\To *_{(AB)CD}(*_{ABC}\times 1_{D})\]
\noindent consisting of, for each triple of composable spans
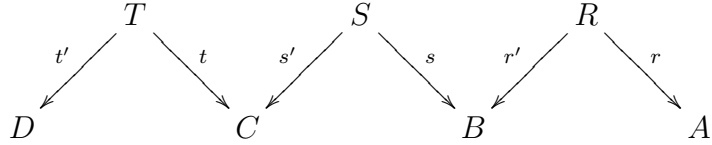
\[
   \xy
   (-45,0)*+{D}="2";
   (-30,15)*+{T}="3";
   (-15,0)*+{C}="4";
   (0,15)*+{S}="5";
   (15,0)*+{B}="6";
   (30,15)*+{R}="7";
   (45,0)*+{A}="8";
        {\ar_{t'} "3";"2"};
        {\ar^{t} "3";"4"};
        {\ar_{s'} "5";"4"};
        {\ar^{s} "5";"6"};
        {\ar_{r'} "7";"6"};
        {\ar^{r} "7";"8"};
\endxy
\]
\begin{itemize}
\item a map of spans ${\bf a}_{RST}$
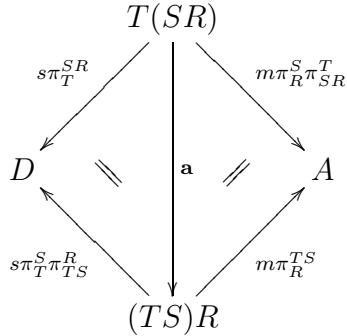
\[
   \xy
   (0,20)*+{T(SR)}="1";
   (-20,0)*+{D}="2"; (20,0)*+{A}="3";
  (0,-20)*+{(TS)R}="5"; 
        {\ar_{s\pi_{T}^{SR}} "1";"2"};
        {\ar^{m\pi_{R}^{S}\pi_{SR}^{T}} "1";"3"};
        {\ar^{s\pi_{T}^{S}\pi_{TS}^{R}} "5";"2"};
        {\ar_{m\pi_{R}^{TS}} "5";"3"};
        {\ar^{{\bf a}} "1";"5"};
        {\ar@{=}^{} (-10,1); (-7,-2)};
        {\ar@{=}_{} (10,1); (7,-2)};
\endxy
\]
\noindent where {\bf a} is a $1$-cell in $\B$ satisfying
\[ \pi_{R}^{TS}{\bf a} = \pi_{R}^{S}\pi_{SR}^{T},\;\; \pi_{S}^{T}\pi_{TS}^{R}{\bf a} = \pi_{S}^{R}\pi_{SR}^{T},\;\; \pi_{T}^{S}\pi_{TS}^{R}{\bf a} = \pi^{SR}_{T},\;\textrm{ and }\;\kappa_{B}^{r\pi_{S}^{T},q}\cdot{\bf a} = \kappa_{B}^{r,q}\cdot\pi_{SR}^{T},\]
\item and a map of spans ${\bf a}^{-1}_{RST}$
\[
   \xy
   (0,-20)*+{T(SR)}="5";
   (-20,0)*+{D}="2"; (20,0)*+{A}="3";
  (0,20)*+{(TS)R}="1"; 
        {\ar^{s\pi_{T}^{SR}} "5";"2"};
        {\ar_{m\pi_{R}^{S}\pi_{SR}^{T}} "5";"3"};
        {\ar_{s\pi_{T}^{S}\pi_{TS}^{R}} "1";"2"};
        {\ar^{m\pi_{R}^{TS}} "1";"3"};
        {\ar_{{\bf a}^{-1}} "1";"5"};
        {\ar@{=}^{} (-10,1); (-7,-2)};
        {\ar@{=}_{} (10,1); (7,-2)};
\endxy
\]
\noindent where ${\bf a}^{-1}$ is a $1$-cell in $\B$ satisfying
\[ \pi_{R}^{S}\pi_{SR}^{T}{\bf a}^{-1} = \pi_{R}^{TS},\;\; \pi_{S}^{R}\pi_{SR}^{T}{\bf a}^{-1} = \pi_{S}^{T}\pi_{TS}^{R},\;\; \pi_{T}^{SR}{\bf a}^{-1} = \pi_{T}^{S}\pi_{TS}^{R},\;\textrm{ and }\;\kappa_{C}^{r,q\pi_{S}^{R}}\cdot{\bf a}^{-1} = \kappa_{C}^{r,q}\cdot\pi_{TS}^{R},\]
\end{itemize}
\noindent respectively,
\item for each triple of maps of composable spans
\[
   \xy
   (-45,0)*+{D}="0";
   (-30,15)*+{T}="1";
   (-30,-15)*+{T'}="7";
   (-15,0)*+{C}="2";
   (0,15)*+{S}="3";
   (0,-15)*+{S'}="8";
   (15,0)*+{B}="4";
   (30,15)*+{R}="5";
   (30,-15)*+{R'}="9";
   (45,0)*+{A}="6";
        {\ar_{s} "1";"0"};
        {\ar^{r} "1";"2"};
        {\ar_{q} "3";"2"};
        {\ar^{p} "3";"4"};
        {\ar_{n} "5";"4"};
        {\ar^{m} "5";"6"};
        {\ar^{s'} "7";"0"};
        {\ar_{r'} "7";"2"};
        {\ar^{q'} "8";"2"};
        {\ar_{p'} "8";"4"};
        {\ar^{n'} "9";"4"};
        {\ar_{m'} "9";"6"};
        {\ar^{f_{T}} "1";"7"};
        {\ar^{f_{S}} "3";"8"};
        {\ar^{f_{R}} "5";"9"};
        {\ar@{=>}_<<{\scriptstyle \varpi_{R}} (38,2); (35,-2)};
        {\ar@{=>}_<<{\scriptstyle \varpi_{S}} (8,2); (5,-2)};
        {\ar@{=>}_<<{\scriptstyle \varpi_{T}} (-22,2); (-25,-2)};
        {\ar@{=>}^<<{\scriptstyle \varrho_{R}} (22,2); (25,-2)};
        {\ar@{=>}^<<{\scriptstyle \varrho_{S}} (-8,2); (-5,-2)};
        {\ar@{=>}^<<{\scriptstyle \varrho_{T}} (-38,2); (-35,-2)};
\endxy
\]
\noindent respective naturality equations between maps of spans
\[ ((\varpi_{R}\cdot\pi_{R}^{S})\cdot\pi_{SR}^{T},\;{\bf a}'((f_{R}*f_{S})*f_{T}),\;\varrho_{T}\cdot\pi_{T}^{SR}) = \]
\[((\varpi_{R}\cdot\pi_{R}^{TS}){\bf a},\;(f_{R}*(f_{S}*f_{T})){\bf a},\;((\varrho_{T}\cdot\pi_{T}^{S})\cdot\pi_{TS}^{R}){\bf a})\]
\noindent and
\[ (\varpi_{R}\cdot\pi_{R}^{TS},\;{\bf a}'^{-1}(f_{R}*(f_{S}*f_{T})),\;(\varrho_{T}\cdot\pi_{T}^{S})\cdot\pi_{TS}^{R}) = \] \[((\varpi_R\cdot\pi_{R}^{S})\cdot\pi_{SR}^{T}{\bf a}^{-1},\;((f_{R}*f_{S})*f_{T}){\bf a}^{-1},\;(\varrho_{T}\cdot\pi_{T}^{SR}){\bf a}^{-1}),\]
\item and identity modifications consisting of, for each span, equations
\[ (1,\; {\bf a}{\bf a}^{-1},\; 1) = (1,1,1)\;\textrm{ and }\;(1,1,1) = (1,\;{\bf a}^{-1}{\bf a},\; 1).\]
\end{itemize}
\endproposition

\proof
The equations for {\bf a} and ${\bf a}^{-1}$ are obtained by combining the uniqueness equations for these $1$-cells with the uniqueness equations for $a$ and $a^{-1}$, respectively.  The transformation and modification axioms are immediate, as are the adjoint equivalence axioms.  The result follows.
\endproof

\subsubsection*{Left Unitor}

We define, for each pair of objects $A, B\in \B$, the components of the left unitor adjoint equivalence $({\bf l},\;{\bf l}^{-1},\;\epsilon_{{\bf l}},\;\eta_{{\bf l}})$. The adjoint equivalence is a pair of strong transformations together with unit and counit modifications.

\proposition \label{leftunitortransformation}
For each pair of objects $A, B\in\B$ there is an adjoint equivalence
\[ ({\bf l},\; {\bf l}^{-1},\; \epsilon_{{\bf l}},\; \eta_{{\bf l}})\maps *_{AB1}(I_{B}\times 1)\To 1,\]
\noindent in the strict $2$-category $[\Span(A,B)\times\Span(B,B),\Span(A,B)]$, with
\begin{itemize}
\item strong transformations
\[ {\bf l}\maps *_{ABB}(I_{B} \times 1)\To 1\]
\noindent and
\[ {\bf l}^{-1}\maps 1\To *_{ABB}(I_{B}\times 1)\]
\noindent consisting of, for each span
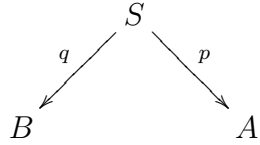
\[
   \xy
   (-15,0)*+{B}="0";
   (0,15)*+{S}="1";
   (15,0)*+{A}="2";
        {\ar^{p} "1";"2"};
        {\ar_{q} "1";"0"};
\endxy
\]
\begin{itemize}
\item a map of spans ${\bf l}_S$
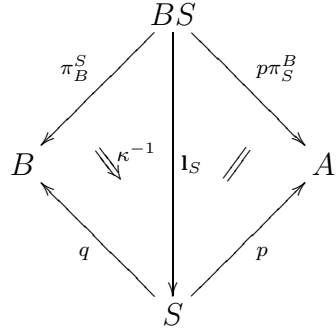
\[
   \xy
   (0,20)*+{BS}="1";
   (-20,0)*+{B}="2"; (20,0)*+{A}="3";
  (0,-20)*+{S}="5"; 
        {\ar_{\pi_{B}^{S}} "1";"2"};
        {\ar^{p\pi_{S}^{B}} "1";"3"};
        {\ar^{q} "5";"2"};
        {\ar_{p} "5";"3"};
        {\ar^{{\bf l}_{S}} "1";"5"};
        {\ar@{=>}^{\kappa^{-1}} (-10,2); (-7,-2)};
        {\ar@{=}_{} (10,2); (7,-2)};
\endxy
\]
\noindent where ${\bf l}_{S} := \pi_{S}^{B}$ and $\kappa := \kappa_{B}^{1,q}$, and
\item a map of spans ${\bf l}_{S}^{-1}$
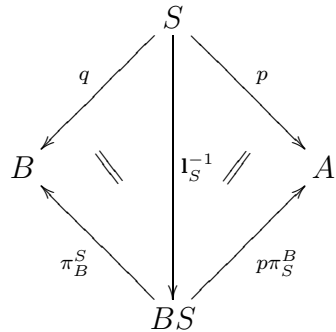
\[
   \xy
   (0,20)*+{S}="1";
   (-20,0)*+{B}="2"; (20,0)*+{A}="3";
  (0,-20)*+{BS}="5"; 
        {\ar^{\pi_{B}^{S}} "5";"2"};
        {\ar_{p\pi_{S}^{B}} "5";"3"};
        {\ar_{q} "1";"2"};
        {\ar^{p} "1";"3"};
        {\ar^{{\bf l}_{S}^{-1}} "1";"5"};
        {\ar@{=}^{} (-10,2); (-7,-2)};
        {\ar@{=}_{} (10,2); (7,-2)};
\endxy
\]
\noindent where ${\bf l}_{S}^{-1} := \pi_{BS}^{S}$, the unique $1$-cell in $\B$ satisfying
\[ \pi_{S}^{B}{\bf l}_{S}^{-1} = 1,\;\;\pi_{B}^{S}{\bf l}_{S}^{-1} = q,\;\textrm{ and }\; \kappa_{B}^{1,q}\cdot{\bf l}_{S}^{-1} = 1,\]
\end{itemize}
\noindent respectively, and
\item for each pair of parallel spans $S, S'$, a pair of natural isomorphisms
\[ {\bf l}\maps ({\bf l}_{S'})_*(*(I\times 1)) \To 1({\bf l}_{S})^*\]
\noindent and
\[ {\bf l}^{-1}\maps ({\bf l}^{-1}_{S'})_{*}1\To *(I\times 1)({\bf l}_{S}^{-1})^{*},\]
\noindent consisting of, for each map of spans
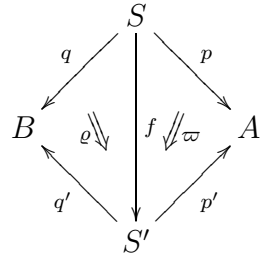
\[
   \xy
   (0,15)*+{S}="1";
   (-15,0)*+{B}="2"; (15,0)*+{A}="3";
   (0,-15)*+{S'}="4";
        {\ar_{q} "1";"2"};
        {\ar^{p} "1";"3"};
        {\ar^{q'} "4";"2"};
        {\ar_{p'} "4";"3"};
        {\ar^{f} "1";"4"};
        {\ar@{=>}_{\scriptstyle \varrho} (-6,2); (-4,-3)};
        {\ar@{=>}^{\scriptstyle \varpi} (6,2); (4,-3)};
\endxy
\]
\begin{itemize}
\item the equation of maps of spans
\[ {\bf l}_{f}\maps (\varpi\cdot\pi_{S}^{B},\;{\bf l}_{S'}(f*1),\; \kappa'^{-1}\cdot(f*1)) = (\varpi\cdot\pi_{S}^{B}{\bf l}_{S},\; f{\bf l}_{S},\;(\varrho\cdot\pi_{S}^{B}{\bf l}_{S})\kappa^{-1}),\]
\item and the isomorphism of maps of spans
\[ {\bf l}_{f}^{-1}\maps (\varpi,\; {\bf l}^{-1}_{S'}f,\;\varrho) \To ((\varpi\cdot\pi_{S}^{B})\cdot {\bf l}_{S}^{-1},\;(f*1){\bf l}^{-1}_{S},\; 1),\]
\noindent defined by the unique $2$-cell
\[ {\bf l}_{f}^{-1}\maps  {\bf l}_{S'}^{-1}f \To (f*1){\bf l}_{S},\]
\noindent in $\B$ satisfying
\[ \pi_{S'}^{B}\cdot{\bf l}_{f}^{-1} = 1 ;\;\textrm{ and }\;\;\pi_{B}^{S'}\cdot {\bf l}_{f}^{-1} = \varrho^{-1},\]
\end{itemize}
\noindent respectively, and, 
\item a pair of invertible modifications
\[
   \xy
   (-7,0)*+{\epsilon_{{\bf l}}\maps {\bf l}{\bf l}^{-1}}="1"; (8,0)*+{1,}="2";
        {\ar@3{->}_{} "1";"2"};
\endxy
\]
\noindent and
\[
   \xy
   (-7,0)*+{\eta_{{\bf l}}\maps 1}="1"; (7,0)*+{{\bf l}^{-1}{\bf l},}="2";
        {\ar@3{->}_{} "1";"2"};
\endxy
\]
\noindent consisting of, for each span
\[
   \xy
   (0,15)*+{S}="1";
   (-15,0)*+{B}="2"; (15,0)*+{A}="3";
        {\ar_{q} "1";"2"};
        {\ar^{p} "1";"3"};
\endxy
\]
\begin{itemize}
\item an equation of maps of spans
\[ \epsilon_{{\bf l}_{S}}\maps (1 ,\;{\bf l}_{S}{\bf l}_{S}^{-1},\; {\kappa_{B}^{1,q}}^{-1}\cdot{\bf l}_{S}^{-1}) = (1,\; 1,\; 1),  \]

\item and, an isomorphism of maps of spans
\[
   \xy
   (-17,0)*+{\eta_{{\bf l}_{S}}\maps (1,1,1)}="1"; (17,0)*+{(1,\;{\bf l}_{S}^{-1}{\bf l}_{S},\; {\kappa_{B}^{1,q}}^{-1}),}="2";
        {\ar@3{->}_{} "1";"2"};
\endxy
\]
\noindent defined by the unique $2$-cell
\[
   \xy
   (-8,0)*+{\eta_{{\bf l}_{S}}\maps 1}="1"; (8,0)*+{{\bf l}_{S}^{-1}{\bf l}_{S}}="2";
        {\ar@3{->}_{} "1";"2"};
\endxy
\]
\noindent in $\B$ satisfying
\[ \pi_{S}^{B}\cdot \eta_{{\bf l}_{S}} = 1\;\;\textrm{ and }\;\; \pi_{B}^{S}\cdot \eta_{{\bf l}_{S}} = {\kappa_{B}^{1,q}}^{-1}.\]

\end{itemize}
\end{itemize}
\endproposition

\proof
We need to verify naturality for the components of the transformations and and verify the axioms of a transformation.  It will then follow that the transformations are strong since all $2$-cell data is defined via the universal property and is therefore invertible.  

The equation of $2$-cells
\[ (1\cdot (\pi_{B}^{S'}\cdot {\bf l}_{f}^{-1}))(\kappa_{B}^{1,q'}\cdot {\bf l}_{S'}^{-1}f) = (\kappa_{B}^{1,q'}\cdot (f*1_B){\bf l}_{S}^{-1})(q'\cdot (\pi^{B}_{S'}\cdot {\bf l}_{f}^{-1}))\]
\noindent allows us to apply the universal property to obtain the $2$-cell ${\bf l}_{f}^{-1}$.

For each map of maps of spans
\[ \sigma\maps (\varpi_{f},\; f,\; \varrho_{f})\To (\varpi_{g},\; g,\; \varrho_{g})\]
\noindent it follows from the equation
\[ ((\sigma*1)\cdot{\bf l}_{S}^{-1}){\bf l}_{f}^{-1} = {\bf l}_{g}^{-1}({\bf l}_{S'}^{-1}\cdot\sigma)\]
\noindent that ${\bf l}^{-1}$ is a natural isomorphism.

Since composition of maps of maps of spans is strictly associative and unital, and the composition and unit $2$-functors of the span construction are strict, the transformation axioms reduce to the equation of maps of maps of spans
\[  {\bf l}_{g}^{-1}{\bf l}_{f}^{-1} = {\bf l}_{gf}^{-1}.\]
\noindent It follows that ${\bf l}^{-1}$ is a strong transformation as desired.  Since the natural isomorphism ${\bf l}$ is the identity, the transformation ${\bf l}$ is strict.

The equation of $2$-cells
\[ (1\cdot (\pi_{B}^{S}\cdot\eta_{{\bf l}_{S}}))(\kappa_{B}^{1,q}\cdot 1) = (\kappa_{B}^{1,q}\cdot {\bf l}_{S}^{-1}{\bf l}_{S})(q\cdot  (\pi^{B}_{S}\cdot\eta_{{\bf l}_{S}}))\]
\noindent is an equation of identity $2$-cells.  We can apply the universal property to define the component $2$-cells of the unit modification
\[
   \xy
   (-15,0)*+{\eta_{{\bf l}_{S}}\maps (1,1,1)}="1"; (15,0)*+{(1,\;{\bf l}_{S}^{-1}{\bf l}_{S},\; \kappa^{-1})}="2";
        {\ar@3{->}_{} "1";"2"};
\endxy
\]
\noindent as the unique $2$-cells
\[
   \xy
   (-8,0)*+{{\eta_{{\bf l}}}_{S}\maps 1}="1"; (8,0)*+{{\bf l}_{S}^{-1}{\bf l}_{S}.}="2";
        {\ar@3{->}_{} "1";"2"};
\endxy
\]

Finally, we have the bicategorical triangle identities.  The first identity reduces to the equation of isomorphisms of maps of spans
\[ {\eta_{{\bf l}}}_{S}1_{{\bf l}^{-1}_{S}} = 1.\]
\noindent The second identity reduces to the equation
\[ 1_{{\bf l}_{S}} {\eta_{{\bf l}}}_{S} = 1.\]
\noindent Recall that the vertical composite of maps of maps of spans is defined by horizontal composition of $2$-cells in $\B$.  Both equations follow and we have the desired adjoint equivalence.
\endproof

\subsubsection*{Right Unitor}

We define, for each pair of objects $A, B\in \B$, the components of the right unitor adjoint equivalence $({\bf r},\;{\bf r}^{-1},\;\epsilon_{{\bf r}},\;\eta_{{\bf r}})$. As in the case of the left unitor, this adjoint equivalence is a pair of strong transformations together with unit and counit modifications.

\proposition \label{rightunitortransformation}
For each pair of objects $A, B\in\B$ there is an adjoint equivalence
\[ ({\bf r},\; {\bf r}^{-1},\; \epsilon_{{\bf r}},\; \eta_{{\bf r}})\maps *_{AAB}(1\times I_{A})\To 1,\]
\noindent in the strict $2$-category $[\Span(A,A)\times\Span(A,B),\Span(A,B)]$, with
\begin{itemize}
\item strong transformations
\[ {\bf r}\maps *_{AAB}(1 \times I_{A})\To 1\]
\noindent and
\[ {\bf r}^{-1}\maps 1\To *_{AAB}(1\times I_{A})\]
\noindent consisting of, for each span
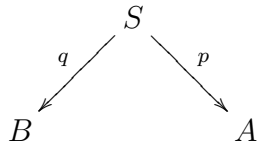
\[
   \xy
   (-15,0)*+{B}="0";
   (0,15)*+{S}="1";
   (15,0)*+{A}="2";
        {\ar^{p} "1";"2"};
        {\ar_{q} "1";"0"};
\endxy
\]
\begin{itemize}
\item a map of spans ${\bf r}_S$
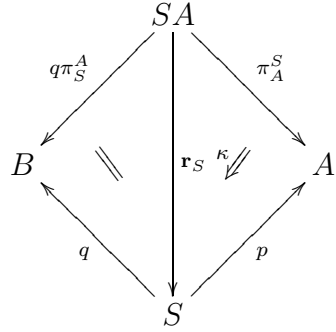
\[
   \xy
   (0,20)*+{SA}="1";
   (-20,0)*+{B}="2"; (20,0)*+{A}="3";
  (0,-20)*+{S}="5"; 
        {\ar_{q\pi_{S}^{A}} "1";"2"};
        {\ar^{\pi_{A}^{S}} "1";"3"};
        {\ar^{q} "5";"2"};
        {\ar_{p} "5";"3"};
        {\ar^{{\bf r}_{S}} "1";"5"};
        {\ar@{=}^{} (-10,2); (-7,-2)};
        {\ar@{=>}_{\kappa} (10,2); (7,-2)};
\endxy
\]
\noindent where ${\bf r} := \pi_{S}^{A}$ and $\kappa := \kappa_{A}^{p,1}$, and
\item a map of spans ${\bf r}_{S}^{-1}$
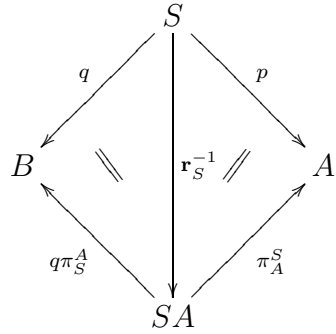
\[
   \xy
   (0,20)*+{S}="1";
   (-20,0)*+{B}="2"; (20,0)*+{A}="3";
  (0,-20)*+{SA}="5"; 
        {\ar^{q\pi_{S}^{A}} "5";"2"};
        {\ar_{\pi_{A}^{S}} "5";"3"};
        {\ar_{q} "1";"2"};
        {\ar^{p} "1";"3"};
        {\ar^{{\bf r}_{S}^{-1}} "1";"5"};
        {\ar@{=}^{} (-10,2); (-7,-2)};
        {\ar@{=}_{} (10,2); (7,-2)};
\endxy
\]
\noindent where ${\bf r}^{-1}_{S} := \pi_{SA}^{S}$, the unique $1$-cell in $\B$ satisfying
\[ \pi_{A}^{S}\pi_{SA}^{S} = p\;\textrm{ and }\;\pi_{S}^{A}\pi_{SA}^{S} = 1,\;\textrm{ and }\kappa_{A}^{p,1}\cdot{\bf r}_{S}^{-1} = 1,\]
\end{itemize}
\noindent respectively, and
\item for each pair of parallel spans $S, S'$, a pair of natural isomorphisms
\[ {\bf r}\maps ({\bf r}_{S'})_*(*(1\times I)) \To 1({\bf r}_{S})^*\]
\noindent and
\[ {\bf r}^{-1}\maps ({\bf r}_{S'}^{-1})_{*}1\To *(1\times I)({\bf r}_{S}^{-1})^{*},\]
\noindent consisting of, for each map of spans
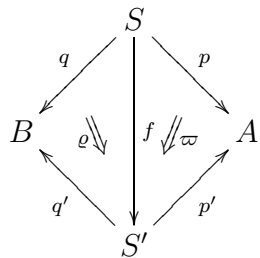
\[
   \xy
   (0,15)*+{S}="1";
   (-15,0)*+{B}="2"; (15,0)*+{A}="3";
   (0,-15)*+{S'}="4";
        {\ar_{q} "1";"2"};
        {\ar^{p} "1";"3"};
        {\ar^{q'} "4";"2"};
        {\ar_{p'} "4";"3"};
        {\ar^{f} "1";"4"};
        {\ar@{=>}_{\scriptstyle \varrho} (-6,2); (-4,-3)};
        {\ar@{=>}^{\scriptstyle \varpi} (6,2); (4,-3)};
\endxy
\]
\begin{itemize}
\item the equation of maps of spans
\[ {\bf r}_{f}\maps (\kappa\cdot(1*f),\;{\bf r}_{S'}(1*f),\; \varrho\cdot \pi_{S}^{A}) \To ((\varpi\cdot {\bf r}_{S})\cdot\kappa ,\; f{\bf r}_{S},\; \varrho\cdot{\bf r}_{S})\]
\item and an isomorphism of maps of spans
\[ {\bf r}_{f}^{-1}\maps (\varpi,\; {\bf r}^{-1}_{S'}f,\;\varrho) \To (1,\;(1*f){\bf r}^{-1}_{S},\; \varrho\cdot\pi_{S}^{A}{\bf r}_{S}^{-1}),\]
\noindent defined by the unique $2$-cell
\[ {\bf r}_{f}^{-1}\maps {\bf r}^{-1}_{S'}f \To (1*f){\bf r}^{-1}_{S}\]
\noindent in $\B$ satisfying
\[ \pi_{A}^{S'}\cdot {\bf r}_{f}^{-1} = \varpi\;\textrm{ and }\; \pi_{S}^{A}\cdot {\bf r}_{f}^{-1} = 1,\]
\end{itemize}
\noindent respectively, and, 
\item a pair of invertible modifications
\[
   \xy
   (-8,0)*+{\epsilon_{{\bf r}}\maps {\bf r}{\bf r}^{-1}}="1"; (8,0)*+{1}="2";
        {\ar@3{->}_{} "1";"2"};
\endxy
\]
\noindent and
\[
   \xy
   (-8,0)*+{\eta_{{\bf r}}\maps 1}="1"; (8,0)*+{{\bf r}^{-1}{\bf r},}="2";
        {\ar@3{->}_{} "1";"2"};
\endxy
\]
\noindent consisting of, for each span
\[
   \xy
   (0,15)*+{S}="1";
   (-15,0)*+{B}="2"; (15,0)*+{A}="3";
        {\ar_{q} "1";"2"};
        {\ar^{p} "1";"3"};
\endxy
\]
\begin{itemize}
\item an equation of maps of spans
\[ \epsilon_{{\bf r}}\maps(\kappa\cdot{\bf r}_{S}^{-1},\; {\bf r}_{S}{\bf r}_{S}^{-1},\; 1) = (1,\; 1,\; 1),\]

\item and, an isomorphism of maps of spans
\[
   \xy
   (-16,0)*+{\eta_{{\bf r}}\maps (1,1,1)}="1"; (15,0)*+{({\kappa_{A}^{p,1}},\;{\bf r}_{S}^{-1}{\bf r}_{S},\; 1).}="2";
        {\ar@3{->}_{} "1";"2"};
\endxy
\]
\end{itemize}
\end{itemize}
\endproposition

\proof
We need to verify naturality for the components of the transformations and and verify the axioms of a transformation.  It will then follow that the transformations are strong since all $2$-cell data is defined via the universal property and is therefore invertible.

The equation of $2$-cells
\[ (p'\cdot (\pi^{S'}_{A}\cdot {\bf r}_{f}^{-1}))(\kappa_{A}^{p',1}\cdot \pi_{S'A}^{S'}f) = (\kappa_{A}^{p',1}\cdot (1*f)\pi^{S}_{SA})(q'\cdot (\pi^{A}_{S'}\cdot {\bf r}_{f}^{-1}))\]
\noindent allows us to apply the universal property to obtain the $2$-cell ${\bf r}_{f}^{-1}$.

For each map of maps of spans
\[ \sigma\maps (\varpi_{f},\; f,\; \varrho_{f})\To (\varpi_{g},\; g,\; \varrho_{g})\]
\noindent it follows from the equation
\[ ({\bf r}_{S'}^{-1}\cdot(1*\sigma)){\bf r}_{f}^{-1} = {\bf r}_{g}^{-1}({\bf r}_{S'}^{-1}\cdot\sigma)\]
\noindent that ${\bf r}^{-1}$ is a natural isomorphism.

Since composition of maps of maps of spans is strictly associative and unital, and the composition and unit $2$-functors of the span construction are strict, the transformation axioms reduce to the equation of maps of maps of spans
\[  {\bf r}_{g}^{-1}{\bf r}_{f}^{-1} = {\bf r}_{gf}^{-1}.\]
\noindent It follows that ${\bf r}^{-1}$ is a strong transformation as desired.  Since the natural isomorphism ${\bf r}$ is the identity, the transformation ${\bf r}$ is strict.

The equation
\[ (p\cdot (\pi_{S}^{A}\cdot\eta_{{\bf r}}))(\kappa_{A}^{p,1}\cdot 1) = (\kappa_{A}^{p,1}\cdot {\bf r}^{-1}{\bf r})(1\cdot  (\pi^{B}_{S}\cdot\eta_{{\bf r}}))\]
\noindent is an equation of identity $2$-cells.  We can apply the universal property to define the component $2$-cells of the unit modification
\[
   \xy
   (-15,0)*+{\eta_{{\bf r}_{S}}\maps (1,1,1)}="1"; (15,0)*+{(\kappa,\;{\bf r}_{S}^{-1}{\bf r}_{S},\; 1)}="2";
        {\ar@3{->}_{} "1";"2"};
\endxy
\]
\noindent as the unique $2$-cells
\[
   \xy
   (-9,0)*+{\eta_{{\bf r}_{S}}\maps 1}="1"; (8,0)*+{{\bf r}_{S}^{-1}{\bf r}_{S}.}="2";
        {\ar@3{->}_{} "1";"2"};
\endxy
\]

Finally, we have the bicategorical triangle identities.  The first identity reduces to the equation of isomorphisms of maps of spans
\[ \eta_{{\bf r}_{S}}1_{{\bf r}^{-1}_{S}} = 1.\]
\noindent The second identity reduces to the equation
\[ 1_{{\bf r}_{S}} \eta_{{\bf r}_{S}} = 1.\]
\noindent Recall that the vertical composite of maps of maps of spans is defined by horizontal composition of $2$-cells in $\B$.  Both equations follow and we have the desired adjoint equivalence.
\endproof

\subsection{Invertible Modifications}

We now define four modifications each of which is an identity.

\subsubsection*{Pentagonator Modification}

\proposition \label{pentagonatormodification}
For objects $A, B, C, D, E\in \B$, there is an identity  modification
\[
   \xy
   (-37,0)*+{\Pi_{ABCDE}\maps (*\cdot(1\times {\bf a}))({\bf a}\cdot (1\times * \times 1))(*\cdot({\bf a}\times 1))}="1"; (-25,-5)*+{}="3"; (-17,-5)*+{}="4";  (17,-5)*+{({\bf a}\cdot(1\times 1\times *))(*\cdot 1)({\bf a}\cdot(1\times 1 \times *)),}="2";
        {\ar@3{->}_{} "3";"4"};
\endxy
\]
\noindent consisting of, for each four composable spans
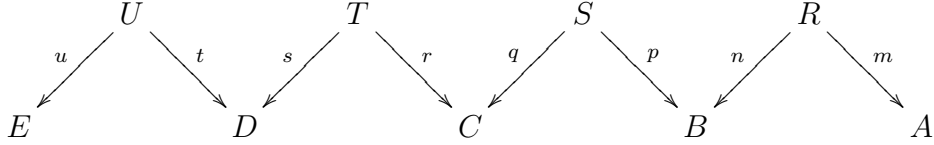
\[
   \xy
   (-60,0)*+{E}="0";
   (-45,15)*+{U}="1";
   (-30,0)*+{D}="2";
   (-15,15)*+{T}="3";
   (0,0)*+{C}="4";
   (15,15)*+{S}="5";
   (30,0)*+{B}="6";
   (45,15)*+{R}="7";
   (60,0)*+{A}="8";
        {\ar_{u} "1";"0"};
        {\ar^{t} "1";"2"};
        {\ar_{s} "3";"2"};
        {\ar^{r} "3";"4"};
        {\ar_{q} "5";"4"};
        {\ar^{p} "5";"6"};
        {\ar_{n} "7";"6"};
        {\ar^{m} "7";"8"};
\endxy
\]
\noindent an equation of maps of spans 
\[ (1,\;(1_{R}*{\bf a}_{STU}){\bf a}_{R(ST)U}({\bf a}_{RST}*1_U),\;1) = (1,\;{\bf a}_{RS(TU)}(1_{SR}*1_{UT}){\bf a}_{(RS)TU},\;1).\] 
\endproposition

\proof
Straightforward.
\endproof

\subsubsection*{Unit Modifications}

\proposition \label{leftunitormodification}
For each triple of objects $A, B, C\in\B$, there is an identity modification
\[
   \xy
   (-25,0)*+{\Lambda_{ABC}\maps  1(*\cdot(1\times {\bf l}))}="1"; (25,0)*+{({\bf l}\cdot *)({\bf a}\cdot (1\times 1\times I))(*\cdot 1),}="2";
        {\ar@3{->}_{} "1";"2"};
\endxy
\]
\noindent consisting of, for each pair of composable spans
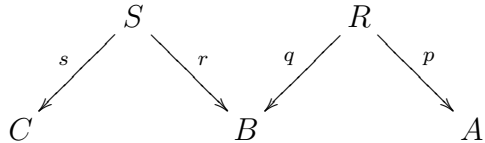
\[
   \xy
   (0,0)*+{C}="2";
   (15,15)*+{S}="3";
   (30,0)*+{B}="4";
   (45,15)*+{R}="5";
   (60,0)*+{A}="6";
        {\ar_{s} "3";"2"};
        {\ar^{r} "3";"4"};
        {\ar_{q} "5";"4"};
        {\ar^{p} "5";"6"};
\endxy
\]
\noindent an equation of maps of spans 
\[ (1,\; 1_{R}*{\bf l}_{S},\; {\kappa_{C}^{1,s}}^{-1}\cdot \pi_{CS}^{R}) = (1 ,\; {\bf l}_{SR}{\bf a}_{CSR},\; {\kappa_{C}^{1,s\pi_{S}^{R}}}^{-1}\cdot {\bf a}_{CSR}).\]
\endproposition

\proof
Straightforward.
\endproof

\proposition \label{middleunitormodification}
For each triple of objects $A, B, C\in\B$, there is an identity modification
\[
   \xy
   (-31,0)*+{M_{ABC}\maps (*\cdot({\bf l}\times 1))\circ ({\bf a}^{-1}\cdot (1\times I \times 1))\circ (*\cdot (1\times {\bf r}^{-1}))}="1"; (31,0)*+{*\cdot 1,}="2";
        {\ar@3{->}_{} "1";"2"};
\endxy
\]
\noindent consisting of, for each pair of composable spans
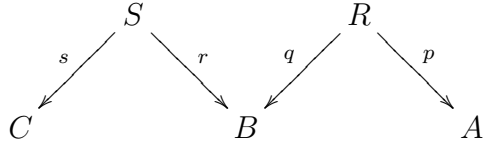
\[
   \xy
   (0,0)*+{C}="2";
   (15,15)*+{S}="3";
   (30,0)*+{B}="4";
   (45,15)*+{R}="5";
   (60,0)*+{A}="6";
        {\ar_{s} "3";"2"};
        {\ar^{r} "3";"4"};
        {\ar_{q} "5";"4"};
        {\ar^{p} "5";"6"};
\endxy
\]
\noindent an equation of maps of spans 
\[(1,({\bf l}_{R}*1_S){\bf a}^{-1}_{RBS}(1_R*{\bf r}_{S}^{-1}),1) = (1,1_{SR},1). \]
\endproposition

\proof
Straightforward.
\endproof

\proposition \label{rightunitormodification}
For each triple of objects $A, B, C\in\B$, there is an identity modification
\[
   \xy
   (-25,0)*+{P_{ABC}\maps  (*\cdot(1 \times {\bf r}))1}="1"; (25,0)*+{({\bf r}\cdot *)(*\cdot 1)({\bf a}\cdot (I\times 1\times 1)),}="2";
        {\ar@3{->}_{} "1";"2"};
\endxy
\]

\noindent consisting of, for each pair of composable spans
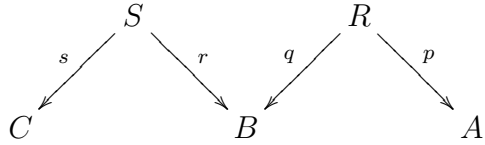
\[
   \xy
   (0,0)*+{C}="2";
   (15,15)*+{S}="3";
   (30,0)*+{B}="4";
   (45,15)*+{R}="5";
   (60,0)*+{A}="6";
        {\ar_{s} "3";"2"};
        {\ar^{r} "3";"4"};
        {\ar_{q} "5";"4"};
        {\ar^{p} "5";"6"};
\endxy
\]
\noindent an identity isomorphism of maps of spans 
\[(\kappa_{A}^{p,1}\cdot\pi_{RA}^{S},\; {\bf r}_{R}*1_{S},\; 1) = ({\kappa_{A}^{p\pi^{S}_{R},1}}\cdot {\bf a}_{SRA} ,\; {\bf r}_{SR}{\bf a}_{SRA},\; 1).\]
\endproposition

\proof
Straightforward.
\endproof

\section{Monoidal Structure on the Tricategory of Spans}\label{monoidal}

In the previous section, we gave an explicit construction of the tricategory $\Span(\B)$, where $\B$ is a strict $2$-category with pullbacks.  The main result of the present section is a construction of a monoidal structure on $\Span(\B)$, where $\B$ is again a strict $2$-category with pullbacks and, in addition, finite products.  Recall, pullback refers to the iso-comma object and products are pseudo (or, equivalently, strict) products.

\subsection{Product Operations on Spans}

We define the basic components of the monoidal structure using products on the objects and morphisms of $\Span(\B)$.

\begin{definition}
For each pair of objects $A, B\in \B$, we choose an object denoted $A\times B$ and projection $1$-cells
\[ \widetilde{\pi}_{A}^{B}\maps A\times B\to A\;\; \textrm{ and } \;\;\widetilde{\pi}^{A}_{B}\maps A\times B\to B\]
\noindent such that the universal property of products is satisfied.
\end{definition}

\noindent The chosen data is called the product of $A$ and $B$, however, we often refer to the object $A\times B$ itself as the product.  The above notation denoting the chosen product is fixed for the duration of the paper.

Given a pair of $1$-cells $f\maps A\to B$ and $f'\maps A'\to B'$ in $\mathcal{B}$, we can apply the universal property expressed in the following diagram
\[
   \xy
   (-20,8)*+{A}="1";
   (0,8)*+{A\times A'}="2"; 
   (20,8)*+{A'}="3";
   (-20,-10)*+{B}="4";
   (0,-10)*+{B\times B'}="5"; 
   (20,-10)*+{B'}="6";
        {\ar_{\widetilde{\pi}_{A}^{A'}} "2";"1"};
        {\ar^{\widetilde{\pi}_{A'}^{A}} "2";"3"};
        {\ar^{\widetilde{\pi}_{B}^{B'}} "5";"4"};
        {\ar_{\widetilde{\pi}_{B'}^{B}} "5";"6"};
        {\ar_{f} "1";"4"};
        {\ar^{f\times f'} "2";"5"};
        {\ar^{f'} "3";"6"};
\endxy
\]
\noindent to obtain a unique comparison $1$-cell
\[f\times f'\maps A\times A'\to B\times B'.\]

\begin{definition}
Given a pair of $1$-cells $f\maps A\to B$ and $f'\maps A'\to B'$ in $\mathcal{B}$, we define the {\bf product $1$-cell}
\[f\times f'\maps A\times A'\to B\times B',\]
\noindent to be the unique $1$-cell in $\B$ such that
\[\widetilde{\pi}_{B}^{B'}(f\times f') = f\widetilde{\pi}_{A}^{A'}\]
\noindent and
\[\widetilde{\pi}_{B'}^{B}(f\times f') = f'\widetilde{\pi}_{A'}^{A}.\]
\end{definition}

\begin{definition}\label{spanproduct}
Given a pair of spans, or $1$-morphisms
\[
   \xy
   (-25,15)*+{S}="1";
   (-40,0)*+{B}="2"; 
   (-10,0)*+{A}="3";
   (25,15)*+{S'}="1'";
   (10,0)*+{B'}="2'"; 
   (40,0)*+{A'}="3'";
        {\ar_{q} "1";"2"};
        {\ar^{p} "1";"3"};
        {\ar_{q'} "1'";"2'"};
        {\ar^{p'} "1'";"3'"};
\endxy
\]
\noindent in $\Span(\mathcal{B})$, we define the {\bf product of spans}
\[
   \xy
   (-20,0)*+{B\times B'}="1";
   (0,15)*+{S\times S'}="2";
   (20,0)*+{A\times A'}="3";
        {\ar_{q\times q'} "2";"1"};
        {\ar^{p\times p'} "2";"3"};
\endxy
\]
\noindent consisting of the product $1$-cells $p\times p'$ and $q\times q'$.
\end{definition}

\begin{definition}\label{product2cell}
Given a pair of $2$-cells $\alpha\maps f\To g$ and $\alpha'\maps f'\To g'$ in $\mathcal{B}$, we define the {\bf product $2$-cell}
\[\varpi\times\varpi'\maps f\times f'\To g\times g'\]
\noindent to be the unique $2$-cell in $\B$ satisfying
\[\widetilde{\pi}_{B}^{B'}\cdot (\varpi\times\varpi') = \varpi\cdot\widetilde{\pi}_{A}^{A'}\;\;\textrm{ and }\;\;\widetilde{\pi}_{B'}^{B}\cdot (\varpi\times\varpi') = \varpi'\cdot\widetilde{\pi}_{A'}^{A}.\]
\end{definition}

\begin{definition}\label{mapproduct}
Given a pair of maps of spans, or $2$-morphisms
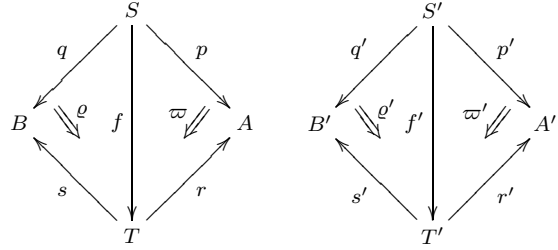
\[
\def\objectstyle{\scriptstyle}
  \def\labelstyle{\scriptstyle}
   \xy
   (-35,0)*+{B}="2";
   (-20,15)*+{S}="1";
   (-20,-15)*+{T}="4";
   (-5,0)*+{A}="3";
   (5,0)*+{B'}="-2";
   (20,15)*+{S'}="-1";
   (20,-15)*+{T'}="-4";
   (35,0)*+{A'}="-3";
        {\ar_{q} "1";"2"};
        {\ar^{p} "1";"3"};
        {\ar^{s} "4";"2"};
        {\ar_{r} "4";"3"};
        {\ar_{f} "1";"4"};
        {\ar_{q'} "-1";"-2"};
        {\ar^{p'} "-1";"-3"};
        {\ar^{s'} "-4";"-2"};
        {\ar_{r'} "-4";"-3"};
        {\ar_{f'} "-1";"-4"};
        {\ar@{=>}^{\scriptstyle \varrho} (-30,2); (-27,-2)};
        {\ar@{=>}_{\scriptstyle \varpi} (-10,2); (-13,-2)};
        {\ar@{=>}^{\scriptstyle  \varrho'} (10,2); (13,-2)};
        {\ar@{=>}_{\scriptstyle \varpi'} (30,2); (27,-2)};
\endxy
\]

\noindent in $\Span(\mathcal{B})$, we define the map of spans
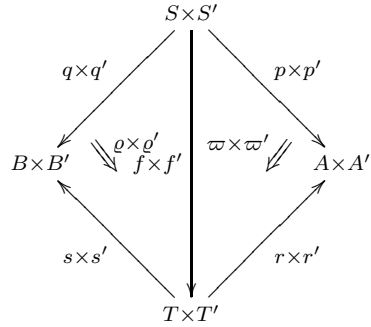
\[
\def\objectstyle{\scriptstyle}
  \def\labelstyle{\scriptstyle}
   \xy
   (-20,0)*+{B\times B'}="2";
   (0,20)*+{S\times S'}="1";
   (0,-20)*+{T\times T'}="4";
   (20,0)*+{A\times A'}="3";
        {\ar_{q\times q'} "1";"2"};
        {\ar^{p\times p'} "1";"3"};
        {\ar^{s\times s'} "4";"2"};
        {\ar_{r\times r'} "4";"3"};
        {\ar_{f\times f'} "1";"4"};
        {\ar@{=>}^{\scriptstyle  \varrho\times  \varrho'} (-13,3); (-10,-1)};
        {\ar@{=>}_{\scriptstyle \varpi\times \varpi'} (13,3); (10,-1)};
\endxy
\]
\noindent consisting of product $1$- and $2$-cells, called the {\bf product of maps of spans}.
\end{definition}

\begin{definition}\label{mapofmapsproduct}
Given a pair of maps of maps of spans, or $3$-morphisms,
\[
   \xy
   (-39,0)*+{\sigma\maps ( \varpi,f, \varrho)}="1"; (-16,0)*+{(\varsigma,g,\varphi)}="2";
   (22,0)*+{\sigma'\maps (\varpi',f', \varrho')}="3"; (49,0)*+{(\varsigma',g',\varphi')}="4";
   (0,0)*+{\textrm{and}}="0";
        {\ar@3{->}_{} "1";"2"};
        {\ar@3{->}_{} "3";"4"};
\endxy
\]
\noindent the {\bf product of maps of maps of spans}
\[
   \xy
   (-27,0)*+{\sigma\times \sigma'\maps (\varpi\times\varpi',f\times f', \varrho\times \varrho')}="1"; (27,0)*+{(\varsigma\times \varsigma',g\times g', \varphi\times \varphi')}="2";
        {\ar@3{->}_{} "1";"2"};
\endxy
\]
\noindent is defined by the product $2$-cell
\[\sigma\times\sigma'\maps f\times f'\To g\times g'.\] 
\end{definition}

\noindent The following equations
\[ (r\times r')\cdot(\sigma\times\sigma')(\varpi\times\varpi') = \varsigma\times\varsigma'\]
\noindent and
\[ (s\times s')\cdot(\sigma\times\sigma')(\varrho\times\varrho') = \varphi\times\varphi'\]
\noindent  obtained from uniqueness in the universal property together with the $3$-morphism equations for $\varsigma$ and $\varphi$ verify that the product of maps of maps of spans is well-defined.

\subsection{The Monoidal Structure}

\theorem
Given a strict $2$-category $\mathcal{B}$ with pullbacks and finite products, $\Span(\B)$ has the structure of a monoidal tricategory consisting of
\begin{itemize}
\item a locally strict trifunctor
\[\otimes\maps \Span(\B)\times \Span(\B)\to \Span(\B)\]
\noindent defined in Proposition~\ref{monoidalproduct},

\item a strict trifunctor
\[ I \maps 1 \to \Span(\B), \]
\noindent defined in Proposition~\ref{monoidalunit},

\item biadjoint biequivalences
\begin{itemize}
\item for associativity
\[ \alpha\maps \otimes(\otimes\times 1)\To \otimes(1\times\otimes),\]
\noindent in a tricategory $[\Span(\B)\times\Span(\B)\times\Span(\B),\Span(\B)]$ defined in Proposition~\ref{monoidalassociator},
\item for left units
 \[ \lambda\maps\otimes(I_{\otimes}\times 1)\To 1,\]
 \noindent a tricategory $[\Span(\B)\times\Span(\B),\Span(\B)]$ defined in Proposition~\ref{monoidalleftunitor},
\item and, for right units 
\[ \rho\maps \otimes (1\times I_{\otimes})\To 1,\]
\noindent in a tricategory $[\Span(\B)\times\Span(\B),\Span(\B)]$ defined in Proposition~\ref{monoidalrightunitor},
\end{itemize}

\item adjoint equivalences
\begin{itemize}
\item for associativity
\[
   \xy
   (-14,0)*+{\Pi\maps (1\times\alpha)\alpha(\alpha\times 1)}="1"; (13,0)*+{\alpha\alpha,}="2";
        {\ar@3{->}_{} "1";"2"};
\endxy
\]
\noindent in a bicategory $[((\otimes \times 1)(\otimes\times 1\times 1),\otimes\otimes(1\times\otimes)]$, defined in Proposition~\ref{monoidalpentagonatormodification},
\item for left units
\[
   \xy
   (-10,0)*+{l\maps (1\times \lambda)\alpha}="1"; (9,0)*+{\lambda,}="2";
        {\ar@3{->}_{} "1";"2"};
\endxy
\]
\noindent in a bicategory $[(I\times 1\times 1)(\otimes 1)\otimes, 1 \otimes]$, defined in Proposition~\ref{monoidalleftmediatormodification},
\item for middle units
\[
   \xy
   (-11,0)*+{m\maps (1\times \rho)\alpha}="1"; (11,0)*+{\lambda\times 1,}="2";
        {\ar@3{->}_{} "1";"2"};
\endxy
\]
\noindent in a bicategory $[1 \otimes, 1 \otimes]$, defined in Proposition~\ref{monoidalmiddlemediatormodification},
\item and, for right units
\[
   \xy
   (-8,0)*+{r\maps \rho\alpha}="1"; (8,0)*+{\rho\times 1,}="2";
        {\ar@3{->}_{} "1";"2"};
\endxy
\]
\noindent in a bicategory $[1 \otimes, (1\times 1\times I)(1\times \otimes)\otimes]$, defined in Proposition~\ref{monoidalrightmediatormodification},
\end{itemize}

\item invertible perturbations
\begin{itemize}
\item for associativity
\[
   \xy
   (-22,0)*+{K_{5}\maps \alpha\Pi\Pi(\Pi\times 1)}="1"; (21,0)*+{\Pi\alpha_{\alpha,1,1}\Pi(1\times\Pi)\alpha_{1,\alpha,1},}="2";
        {\ar@3{->}_{} "1";"2"};
\endxy
\]
\noindent defined in Proposition~\ref{associativityperturbation},
\item for $(4,1)$-units
\[
   \xy
   (-16,0)*+{U_{4,1}\maps \alpha_{\rho,1,1}r\Pi}="1"; (15,0)*+{(r\times 1)r\rho_{\alpha}^{-1},}="2";
        {\ar@3{->}_{} "1";"2"};
\endxy
\]
\noindent defined in Proposition~\ref{U41perturbation},
\item for $(4,2)$-units
\[
   \xy
   (-20,0)*+{U_{4,2}\maps \alpha_{\lambda,1,1}m\Pi}="1"; (19,0)*+{(m\times 1)\alpha_{1,\rho,1}(1\times r),}="2";
        {\ar@3{->}_{} "1";"2"};
\endxy
\]
\noindent defined in Proposition~\ref{U42perturbation},
\item for $(4,3)$-units
\[
   \xy
   (-20,0)*+{U_{4,3}\maps m\alpha_{1,1,\rho}\Pi}="1"; (18,0)*+{(l\times 1)\alpha_{1,\lambda,1}(1\times m),}="2";
        {\ar@3{->}_{} "1";"2"};
\endxy
\]
\noindent defined in Proposition~\ref{U43perturbation},
\item and, for $(4,4)$-units 
\[
   \xy
   (-15,0)*+{U_{4,4}\maps l\alpha_{1,1,\rho}\Pi}="1"; (15,0)*+{\lambda^{-1}l(1\times l),}="2";
        {\ar@3{->}_{} "1";"2"};
\endxy
\]
\noindent defined in Proposition~\ref{U44perturbation}.
\end{itemize}
\end{itemize}
\endtheorem

\proof
The structural components are given in the referenced definitions and propositions.

The monoidal structure we define involves non-trivial product cells and perturbations.  This fact necessitates checking the tetracategory axioms.  We verify these coherence equations by explicit calculation in Propositions~\ref{K6Axiom},~\ref{U52Axiom},~\ref{U53Axiom}, and~\ref{U54Axiom}.  The result follows. 
\endproof

In the following sections we construct component cells of the monoidal structure on the tricategory of spans.

\subsection{Monoidal Product}

The monoidal product, which is obtained via the universal property of the product, is a locally strict trifunctor with identity modifications.  

\proposition \label{monoidalproduct}
There is a locally strict trifunctor
\[ \otimes\maps \Span(\B)\times \Span(\B)\to\Span(\B),\]
\noindent consisting of
\begin{itemize}
\item a function 
\[ (A, B)\stackrel{\otimes}\longmapsto A\times B\] 
\noindent on pairs of objects in $\Span(\B)$ defined by the choice of strict product in Definition~\ref{strictproduct},
\item for each two pairs $(A, B), (A', B')$ objects in $\Span(\B)$, a strict functor
\[\otimes := \otimes_{(A,B),(A',B')}\maps \Span(A,B)\times \Span(A',B')\to \Span(A\times A',B\times B')\]
\noindent between strict $\hom$-$2$-categories, consisting of
\begin{itemize} 
\item the product of pairs of spans given in Definition~\ref{spanproduct},
\item the product of pairs of maps of spans given in Definition~\ref{mapproduct},
\item and the product of pairs of maps of maps of spans given in Definition~\ref{mapofmapsproduct},
\end{itemize}
\item for each two triples of objects $(A,B,C),(A',B',C')$ in $\Span(\B)$, a strict adjoint equivalence
\[(\chi,\; \chi^{-1},\;1,\;1)\maps *\cdot(\otimes\times\otimes)\To \otimes(*\times *),\]
\noindent in a strict $2$-category $\Span(\B)(A\times A', C\times C')$, consisting of, for each two pairs of composable spans
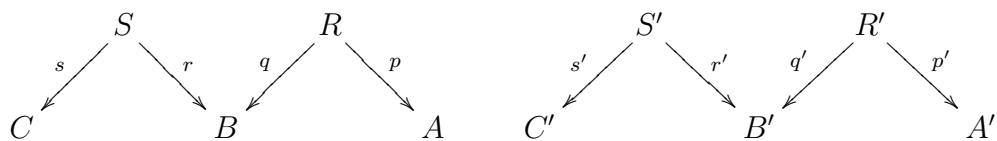
\[
\xymatrix{
& S\ar[dl]_{s}\ar[dr]^{r}&& R\ar[dl]_{q}\ar[dr]^{p}&&& S'\ar[dl]_{s'}\ar[dr]^{r'}&& R'\ar[dl]_{q'}\ar[dr]^{p'}&\\
C&&B&& A & C'&&B'&& A'
}
\]
\begin{itemize}
\item a map of spans
\[
\def\objectstyle{\scriptstyle}
  \def\labelstyle{\scriptstyle}
   \xy
   (-20,0)*+{C\times C'}="2";
   (0,-20)*+{SR\times S'R'}="4";
   (0,20)*+{(S\times S')(R\times R')}="1";
   (20,0)*+{A\times A'}="3";
        {\ar^{(s\times s')(\pi_{S}^{R}\times\pi_{S'}^{R'})} "4";"2"};
        {\ar_{(p\times p')(\pi_{R}^{S}\times\pi_{R'}^{S'})} "4";"3"};
        {\ar_{(s\times s')\pi_{S\times S'}^{R\times R'}} "1";"2"};
        {\ar^{(p\times p')\pi_{R\times R'}^{S\times S'}} "1";"3"};
        {\ar^{\chi} "1";"4"};
        {\ar@{=}^{} (-10,2); (-7,-1)};
        {\ar@{=}_{} (10,2); (7,-1)};
\endxy
\]
\noindent where $\chi := \chi_{(R,S), (R',S')}$ is a $1$-cell in $\B$ satisfying
\[ \pi_{R}^{S}\widetilde{\pi}_{SR}^{S'R'}\chi = \widetilde{\pi}_{R}^{R'}\pi_{R\times R'}^{S\times S'},\;\;\pi_{S}^{R}\widetilde{\pi}_{SR}^{S'R'}\chi = \widetilde{\pi}_{S}^{S'}\pi^{R\times R'}_{S\times S'},\]
\[ \pi_{R'}^{S'}\widetilde{\pi}^{SR}_{S'R'}\chi = \widetilde{\pi}^{R}_{R'}\pi_{R\times R'}^{S\times S'},\;\;\textrm{ and }\;\; \pi_{S'}^{R'}\widetilde{\pi}^{SR}_{S'R'}\chi = \widetilde{\pi}_{S'}^{S}\pi^{R\times R'}_{S\times S'},\] 

\item and an inverse map of spans
\[
\def\objectstyle{\scriptstyle}
  \def\labelstyle{\scriptstyle}
   \xy
   (-20,0)*+{C\times C'}="2";
   (0,20)*+{SR\times S'R'}="4";
   (0,-20)*+{(S\times S')(R\times R')}="1";
   (20,0)*+{A\times A'}="3";
        {\ar_{(s\times s')(\pi_{S}^{R}\times\pi_{S'}^{R'})} "4";"2"};
        {\ar^{(p\times p')(\pi_{R}^{S}\times\pi_{R'}^{S'})} "4";"3"};
        {\ar^{(s\times s')\pi_{S\times S'}^{R\times R'}} "1";"2"};
        {\ar_{(p\times p')\pi_{R\times R'}^{S\times S'}} "1";"3"};
        {\ar^{\chi^{-1}} "4";"1"};
        {\ar@{=}^{} (-10,2); (-7,-1)};
        {\ar@{=}_{} (10,2); (7,-1)};
\endxy
\]
\noindent where $\chi^{-1} := \chi^{-1}_{(R,S), (R',S')}$ is a $1$-cell in $\B$ satisfying
\[\pi_{R\times R'}^{S\times S'}\chi^{-1} = \pi_R^S\times \pi_{R'}^{S'},\;\; \pi_{S\times S'}^{R\times R'}\chi^{-1} = \pi_S^{R}\times \pi_{S'}^{R'},\;\;\textrm{ and }\;\; \kappa_{B\times B'}^{r\times r',q\times q}\cdot\chi^{-1} = \kappa_B^{r,q}\times\kappa_{B'}^{r',q'},\]
\item identity counit and unit isomorphisms of maps of spans
\[ \epsilon_{\chi}\maps \chi\chi^{-1} \Rightarrow 1\;\textrm{ and }\;\eta_{\chi}\maps 1\Rightarrow \chi^{-1}\chi,\]
\end{itemize}
\item  for each pair of objects $A,B$ in $\Span(\B)$, an identity adjoint equivalence
\[ \iota_{A,B}\maps I_{A\times B}\To \otimes(I_A\times I_B),\]
\noindent in the strict $2$-category $\Span(\B)(A\times B,A\times B)$,
\item for $(A,A'),(B,B'),(C,C'),(D,D')\in\Span(\B)\times\Span(\B)$, an identity modification $\omega$ in $\Span(\B)(A\times A', D\times D')$,  consisting of, for each pair of triples of composable spans  $(R,S,T), (R',S',T')$, an equation of maps of spans
\[ (1,\;({\bf a}\times {\bf a})\chi(\chi* 1)\; 1) = (1,\;\chi(1* \chi){\bf a},\; 1),\]
\item for $(A,A'),(B,B')\in\Span(\B)\times\Span(\B)$, an identity modification $\gamma$ in $\Span(\B)(A\times A', B\times B')$,  consisting of, for each pair of spans $R,S$, an equation of maps of spans
\[ ((\kappa^{-1}\times \kappa'^{-1})\cdot\chi\iota, \;({\bf l}\times {\bf l})\chi\iota,\;1) = (\kappa^{-1}, \; {\bf l},\; 1),\]
\item for $(A,A'),(B,B')\in\Span(\B)\times\Span(\B)$, an identity modification $\delta$ in $\Span(\B)(A\times A', B\times B')$,  consisting of, for each pair of spans $R,S$, an equation of maps of spans
\[ (1, \;({\bf r}\times {\bf r})\chi\iota, \;(\kappa\times\kappa')\cdot\chi\iota) = (1, \;{\bf r},\; \kappa).\]
\end{itemize}

\endproposition

\proof
We need to verify functoriality for maps of spans, i.e., that the $2$-functor is strict, and then functoriality for the $2$-functor, which concerns composition of maps of maps of spans.  We define unique auxiliary $1$-cells to give unique definitions of $\chi$ and $\chi^{-1}$ by the universal property, and then verify naturality for the two families of maps of spans.  It is then straightforward to see that these natural transformations together with identity counit and unit define a strict adjoint equivalence.  It is again straightforward to see that $\iota$ uniquely defines a strict adjoint equivalence.  We then have identity modifications and together with the identity modifications of $\Span(\B)$, the trifunctor axioms follow.

For each two pairs of objects $(A,B), (A',B')\in\Span(\B)\times\Span(\B)$, the monoidal product should preserve composition of maps of spans and identity maps of spans.  It is straightforward from definitions to see that identities are preserved.  For composition, consider the two pairs of maps of spans
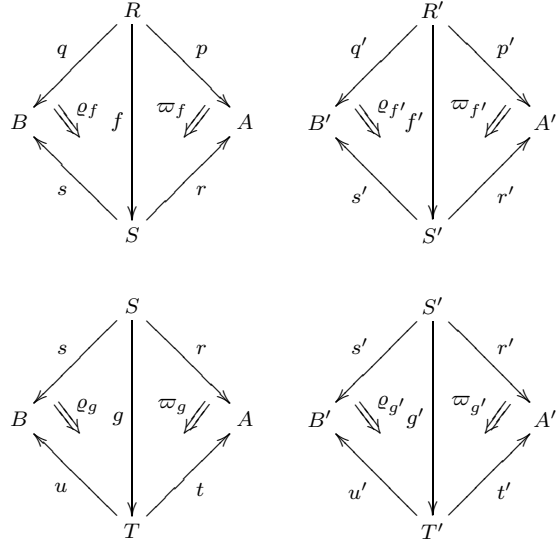
\[
\def\objectstyle{\scriptstyle}
  \def\labelstyle{\scriptstyle}
   \xy
   (-35,0)*+{B}="2";
   (-20,15)*+{R}="1";
   (-20,-15)*+{S}="4";
   (-5,0)*+{A}="3";
   (5,0)*+{B'}="-2";
   (20,15)*+{R'}="-1";
   (20,-15)*+{S'}="-4";
   (35,0)*+{A'}="-3";
        {\ar_{q} "1";"2"};
        {\ar^{p} "1";"3"};
        {\ar^{s} "4";"2"};
        {\ar_{r} "4";"3"};
        {\ar_{f} "1";"4"};
        {\ar_{q'} "-1";"-2"};
        {\ar^{p'} "-1";"-3"};
        {\ar^{s'} "-4";"-2"};
        {\ar_{r'} "-4";"-3"};
        {\ar_{f'} "-1";"-4"};
        {\ar@{=>}^{\scriptstyle \varrho_{f}} (-30,2); (-27,-2)};
        {\ar@{=>}_{\scriptstyle \varpi_{f}} (-10,2); (-13,-2)};
        {\ar@{=>}^{\scriptstyle \varrho_{f'}} (10,2); (13,-2)};
        {\ar@{=>}_{\scriptstyle \varpi_{f'}} (30,2); (27,-2)};
\endxy
\]

\[
\def\objectstyle{\scriptstyle}
  \def\labelstyle{\scriptstyle}
   \xy
   (-35,0)*+{B}="2";
   (-20,15)*+{S}="1";
   (-20,-15)*+{T}="4";
   (-5,0)*+{A}="3";
   (5,0)*+{B'}="-2";
   (20,15)*+{S'}="-1";
   (20,-15)*+{T'}="-4";
   (35,0)*+{A'}="-3";
        {\ar_{s} "1";"2"};
        {\ar^{r} "1";"3"};
        {\ar^{u} "4";"2"};
        {\ar_{t} "4";"3"};
        {\ar_{g} "1";"4"};
        {\ar_{s'} "-1";"-2"};
        {\ar^{r'} "-1";"-3"};
        {\ar^{u'} "-4";"-2"};
        {\ar_{t'} "-4";"-3"};
        {\ar_{g'} "-1";"-4"};
        {\ar@{=>}^{\scriptstyle \varrho_{g}} (-30,2); (-27,-2)};
        {\ar@{=>}_{\scriptstyle \varpi_{g}} (-10,2); (-13,-2)};
        {\ar@{=>}^{\scriptstyle \varrho_{g'}} (10,2); (13,-2)};
        {\ar@{=>}_{\scriptstyle \varpi_{g'}} (30,2); (27,-2)};
\endxy
\]
\noindent Functoriality of maps of spans follows from the equation
\[ (((\varpi_{g}\times\varpi_{g'})\cdot(f\times f'))(\varpi_{f}\times\varpi_{f'}),\; (g\times g')(f\times f'), \;((\varrho_{g}\times\varrho_{g'})\cdot(g\times g'))(\varrho_{f}\times\varrho_{f'})) = \]
\[ ((\varpi_{g}\cdot f)\varpi_{f}\times (\varpi_{g'}\cdot f')\varpi_{f'} ,\; gf\times g'f',\; (\varrho_{g}\cdot f)\varrho_{f} \times (\varrho_{g'}\cdot f')\varrho_{f'}),\]
 \noindent and the obvious preservation of identity maps of spans, so the monoidal product is strict.

To verify functoriality of maps of maps of spans, consider pairs of composable pairs
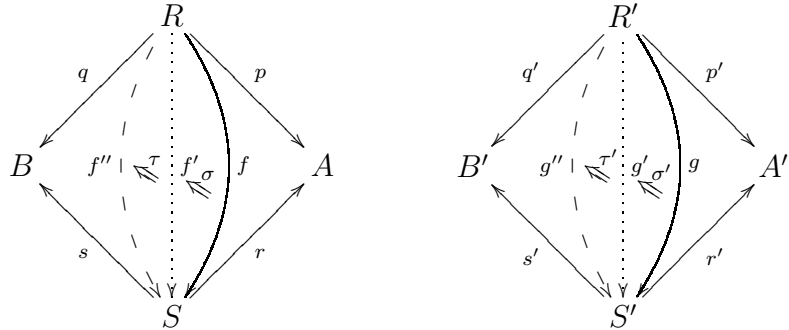
\[
   \xy
   (-30,20)*+{R}="1";
   (-50,0)*+{B}="2";  (-10,0)*+{A}="4";
   (-30,-20)*+{S}="5";
   (-5,0)*+{}="6"; (5,0)*+{}="7";   
   (30,20)*+{R'}="8";
   (10,0)*+{B'}="9"; (50,0)*+{A'}="11";
   (30,-20)*+{S'}="12";
        {\ar_{q} "1";"2"};
        {\ar^{p} "1";"4"};
        {\ar^{s} "5";"2"};
        {\ar_{r} "5";"4"}; 
        {\ar@/^1.8pc/^{f} "1";"5"};
        {\ar@{..>}@/^0pc/^{f'} "1";"5"};
        {\ar@{-->}@/_1.6pc/_{f''} "1";"5"};
        {\ar_{q'} "8";"9"};
        {\ar^{p'} "8";"11"};
        {\ar^{s'} "12";"9"};
        {\ar_{r'} "12";"11"}; 
        {\ar@/^1.8pc/^{g} "8";"12"};
        {\ar@{..>}@/^0pc/^{g'} "8";"12"};
        {\ar@{-->}@/_1.6pc/_{g''} "8";"12"};
        {\ar@{=>}_{\scriptstyle \sigma} (-25,-4); (-28,-2)};
        {\ar@{=>}_{\scriptstyle \tau} (-32,-2); (-35,0)};
       {\ar@{=>}_{\scriptstyle \sigma'} (35,-4); (32,-2)};
        {\ar@{=>}_{\scriptstyle \tau'} (28,-2); (25,0)};
\endxy
\]
\noindent We have the equation of $2$-cells
\[ \tau\sigma\times\tau'\sigma' = (\tau\times\tau')(\sigma\times\sigma')\]
\noindent in $\B$, so composition is preserved.  It is straightforward to see that identity maps of maps of spans are preserved.  It follows that products in $\B$ define a strict functor on $\hom$-$2$-categories.

We define the components of a strict adjoint equivalence $\chi$.  Recall that strict transformations are natural transformations consisting of $1$-morphisms.  For each pair of triples of objects $(A, B, C), (A', B', C')\in\Span(\B)\times\Span(\B)\times\Span(\B)$, we define a natural transformation
\[ \chi_{(A,B,C),(A',B',C')}\maps *(\otimes\times\otimes)\To \otimes(*\times *).\]

We first apply the universal property of pullbacks to obtain an auxiliary pair of maps $\chi_{SR}$ and $\chi_{S'R'}$, which we then use to obtain the $1$-cell $\chi$.  The $1$-cells $\chi_{SR}$ and $\chi_{S'R'}$ are the unique $1$-cells in $\B$ making the diagrams
\[
\xymatrix{
S\times S'\ar[d]_{\widetilde{\pi}_{S}^{S'}} & (S\times S')(R\times R')\ar[r]^<<<{\pi_{R\times R'}^{S\times S'}}\ar[l]_<<<{\pi_{S\times S'}^{R\times R'}}\ar[d]_{\chi_{SR}} & R\times R'\ar[d]^{\widetilde{\pi}_{R}^{R'}}\\
S& SR\ar[r]_{\pi_{R}^{S}}\ar[l]^{\pi_{S}^{R}} & R
}
\]

\[
\xymatrix{
S\times S'\ar[d]_{\widetilde{\pi}_{S'}^{S}} & (S\times S')(R\times R')\ar[r]^<<<{\pi_{R\times R'}^{S\times S'}}\ar[l]_<<<{\pi_{S\times S'}^{R\times R'}}\ar[d]_{\chi_{S'R'}} & R\times R'\ar[d]^{\widetilde{\pi}_{R'}^{R}}\\
S'& S'R'\ar[r]_{\pi_{R'}^{S'}}\ar[l]^{\pi_{S'}^{R'}} & R'
}
\]
\noindent commute and satisfying the equations
\[ \kappa_{B}^{SR}\cdot\chi_{SR} = \widetilde{\pi}_{B}^{B'}\cdot\kappa_{B\times B'}^{(S\times S')(R\times R')}\;\;\;\textrm{ and }\;\;\; \kappa_{B'}^{S'R'}\cdot\chi_{S'R'} = \widetilde{\pi}_{B'}^{B}\cdot\kappa_{B\times B'}^{(S\times S')(R\times R')}.\]

\noindent Applying the universal property of products, we have the unique $1$-cell
\[ \chi_{(R,S),(R',S')}\maps (S\times S')(R\times R') \to SR\times S'R'\]
\noindent such that
\[ \widetilde{\pi}_{SR}^{S'R'}\chi = \chi_{SR}\;\;\textrm{ and }\;\; \widetilde{\pi}^{SR}_{S'R'}\chi = \chi_{S'R'}.\]

\noindent From the above equations we have
\[ \pi_{R}^{S}\widetilde{\pi}_{SR}^{S'R'}\chi = \widetilde{\pi}_{R}^{R'}\pi_{R\times R'}^{S\times S'}, \hspace{10pt} \pi_{S}^{R}\widetilde{\pi}_{SR}^{S'R'}\chi = \widetilde{\pi}_{S}^{S'}\pi^{R\times R'}_{S\times S'},\]
\[ \pi_{R'}^{S'}\widetilde{\pi}^{SR}_{S'R'}\chi = \widetilde{\pi}^{R}_{R'}\pi_{R\times R'}^{S\times S'}, \hspace{10pt} \pi_{S'}^{R'}\widetilde{\pi}^{SR}_{S'R'}\chi = \widetilde{\pi}_{S'}^{S}\pi^{R\times R'}_{S\times S'}.\] 

Similarly, we define $1$-cells $\chi^{-1}$ such that
\[\pi_{R\times R'}^{S\times S'}\chi^{-1} = \pi_R^S\times \pi_{R'}^{S'},\;\; \pi_{S\times S'}^{R\times R'}\chi^{-1} = \pi_S^{R}\times \pi_{S'}^{R'},\;\;\textrm{ and }\;\; \kappa_{B\times B'}^{(S\times S')(R\times R')}\cdot\chi^{-1} = \kappa_B^{SR}\times\kappa_{B'}^{S'R'}.\]

To verify naturality we see that, for each two pairs of horizontally composable maps of spans
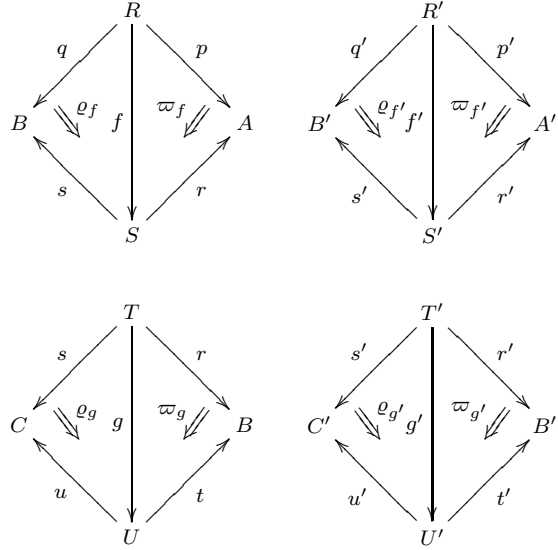
\[
\def\objectstyle{\scriptstyle}
  \def\labelstyle{\scriptstyle}
   \xy
   (-35,0)*+{B}="2";
   (-20,15)*+{R}="1";
   (-20,-15)*+{S}="4";
   (-5,0)*+{A}="3";
   (5,0)*+{B'}="-2";
   (20,15)*+{R'}="-1";
   (20,-15)*+{S'}="-4";
   (35,0)*+{A'}="-3";
        {\ar_{q} "1";"2"};
        {\ar^{p} "1";"3"};
        {\ar^{s} "4";"2"};
        {\ar_{r} "4";"3"};
        {\ar_{f} "1";"4"};
        {\ar_{q'} "-1";"-2"};
        {\ar^{p'} "-1";"-3"};
        {\ar^{s'} "-4";"-2"};
        {\ar_{r'} "-4";"-3"};
        {\ar_{f'} "-1";"-4"};
        {\ar@{=>}^{\scriptstyle \varrho_{f}} (-30,2); (-27,-2)};
        {\ar@{=>}_{\scriptstyle \varpi_{f}} (-10,2); (-13,-2)};
        {\ar@{=>}^{\scriptstyle \varrho_{f'}} (10,2); (13,-2)};
        {\ar@{=>}_{\scriptstyle \varpi_{f'}} (30,2); (27,-2)};
\endxy
\]

\[
\def\objectstyle{\scriptstyle}
  \def\labelstyle{\scriptstyle}
   \xy
   (-35,0)*+{C}="2";
   (-20,15)*+{T}="1";
   (-20,-15)*+{U}="4";
   (-5,0)*+{B}="3";
   (5,0)*+{C'}="-2";
   (20,15)*+{T'}="-1";
   (20,-15)*+{U'}="-4";
   (35,0)*+{B'}="-3";
        {\ar_{s} "1";"2"};
        {\ar^{r} "1";"3"};
        {\ar^{u} "4";"2"};
        {\ar_{t} "4";"3"};
        {\ar_{g} "1";"4"};
        {\ar_{s'} "-1";"-2"};
        {\ar^{r'} "-1";"-3"};
        {\ar^{u'} "-4";"-2"};
        {\ar_{t'} "-4";"-3"};
        {\ar_{g'} "-1";"-4"};
        {\ar@{=>}^{\scriptstyle \varrho_{g}} (-30,2); (-27,-2)};
        {\ar@{=>}_{\scriptstyle \varpi_{g}} (-10,2); (-13,-2)};
        {\ar@{=>}^{\scriptstyle \varrho_{g'}} (10,2); (13,-2)};
        {\ar@{=>}_{\scriptstyle \varpi_{g'}} (30,2); (27,-2)};
\endxy
\]
\noindent there is an identity isomorphism
\[ ((\varpi_{f}\times\varpi_{f'})\cdot\pi_{R\times R'}^{S\times S'},\; \chi'((f\times f')*(g\times g')), \;(\varrho_{g}\times\varrho_{g'})\cdot\pi_{S\times S'}^{R\times R'}) = \]
\[(((\varpi_{f}\cdot\pi_{R}^{S})\times(\varpi_{f'}\cdot\pi_{R'}^{S'}))\chi, \;((f*g)\times(f'*g'))\chi, \;((\varrho_{g}\cdot\pi_{S}^{R})\times(\varrho_{g'}\cdot\pi_{S'}^{R'}))\chi)\]
\noindent between the maps of spans
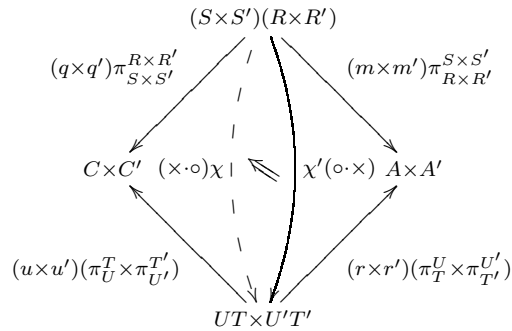
\[
\def\objectstyle{\scriptstyle}
  \def\labelstyle{\scriptstyle}
   \xy
   (-20,0)*+{C\times C'}="2";
   (0,-20)*+{UT\times U'T'}="4";
   (0,20)*+{(S\times S')(R\times R')}="1";
   (20,0)*+{A\times A'}="3";
   (3,-2)*+{}="5";
   (-3,2)*+{}="6";
        {\ar_{(q\times q')\pi_{S\times S'}^{R\times R'}} "1";"2"};
        {\ar^{(m\times m')\pi^{S\times S'}_{R\times R'}} "1";"3"};
        {\ar^{(u\times u')(\pi_U^{T}\times \pi_{U'}^{T'})} "4";"2"};
        {\ar_{(r\times r')(\pi_T^{U}\times \pi_{T'}^{U'})} "4";"3"};
        {\ar@{=>}_{} "5";"6"};
        {\ar@/^1pc/^{\chi'(\circ\cdot\times)} "1";"4"};
        {\ar@/_1pc/@{-->}_{(\times\cdot\circ)\chi} "1";"4"};
\endxy
\]
\noindent It follows that the maps of spans are the components of a natural transformation
\[\chi\maps *(\otimes\times\otimes) \To \otimes (*\times *).\]
\noindent One can check similarly that the collection of maps
\[\chi^{-1}\maps \otimes (*\times *)\To *(\otimes\times\otimes)\]
\noindent is a natural transformation.

The equations of maps of spans
\[ (1,\;\chi\chi^{-1},\;1) = (1,1,1)\;\textrm{ and }\; (1,\;\chi^{-1}\chi,\;1) = (1,1,1),\]
\noindent which are the components of the identity counit and unit verify that $(\chi,\;\chi^{-1},\;1,\;1)$ is a strict adjoint equivalence.  The axioms are immediate.

The existence of the identity adjoint equivalence $\iota$ and the identity modifications is straightforward.  All axioms are immediate.  We have the desired locally strict trifunctor.
\endproof

\subsection{Monoidal Unit}\label{monoidalunitsection}

\proposition \label{monoidalunit}
There is a strict functor between tricategories called the monoidal unit, consisting of the terminal object (or nullary product) in $B$, and the identity morphisms of Definition~\ref{identityspan}, Definition~\ref{identitymap}, and Definition~\ref{identitymapofmaps} for spans, maps of spans, and maps of maps of spans, respectively.  
\endproposition

\proof
Straightforward.
\endproof

\subsection{Monoidal Biadjoint Biequivalences}

The monoidal associativity and unit structures are given as biadjoint biequivalences.  Note that Trimble's definition of tetracategory asks that the tritransformations be equivalences in the appropriate sense at each level~\cite{Tr}.   He calls such a map a triequivalence, which should not be interpreted as a trifunctor that is an equivalence, but rather as a suitable notion of `strong tritransformation'.  Following the definition of Gurski's algebraic tricategory, we replace these structural tritransformations with biadjoint biequivalences, a notion which categorifies that of adjoint equivalence.  See Gurski's thesis~\cite{Gu} for definitions.

The following biadjoint biequivalences are pairs of biadjoint $1$-cells in certain tricategories of trifunctors, tritransformations, trimodifications, and perturbations.  This notion of biadjoint is not to be confused with ambidextrous adjoint pairs, which are sometimes called biadjoints.

Tritransformations, trimodifications, and perturbations are the morphisms of the local tricategories of a tetracategory $\Tricat$ of tricategories.  We need to specify the structure of these tricategories.  Given tricategories $\T$, $\T'$ such that $\T'$ is also a $\Gray$-category, then $\Tricat(T,T')$ is also a $\Gray$-category.~\cite{Gu}  Unfortunately, $\Span(\B)$ is not quite a $\Gray$-category, so there is still work to do in specifying the structure of $\Tricat(\Span(\B)^{n},\Span(\B))$, for $n$ a natural number.  As a corollary of the local bicategory construction in the local tricategories of $\Tricat$, Gurski further shows that if $\T'$ is locally strict, then for trifunctors $F,G\maps \T\to \T'$, the bicategory $\Tricat(\T,\T')(F,G)$ is a strict $2$-category~\cite{Gu}.  Since $\Span(\B)$ is locally strict, so is $\Tricat(\Span(\B)^{n},\Span(\B))$.

To the best of our knowledge the tricategorical structure of $\Tricat(\T,\T')$ does not exist in the literature for an arbitrary tricategory $\T'$.  We expect the details should be straightforward, but we do not have space here to present them here.  This is not too troublesome since we can invoke tricategorical coherence.  The tricategory $\Span(\B)$ is semi-strict and cubical.  If the associator and unit tritransformations were identities, then we could apply Gurski's theorem above and use the $\Gray$-category structure on $\Tricat(\T,\T')$.  Alternatively, we could consider a strictification of the span tricategory $\Span(\B)$ to a triequivalent $\Gray$-category $\Span^{\Gray}(\B)$
\[ {\rm st}\maps\Span(\B)\to\Span^{\Gray}(\B).\]
We can then essentially `whisker' the biadjoint biequivalence structures we define with the strictification maps.  The result being a biadjoint biequivalence in the $\Gray$-category $[\Span(\B)^{n},\Span^{\Gray}(\B)]$, rather than the locally strict tricategory $[\Span(\B)^{n},\Span(\B)]$.  This does not really provide satisfactory resolution to the issue, but instead strongly suggests that there should be numerous solutions to the problem, where in each case the details should be relatively straightforward.  We do not comment further on the tricategory structures in which we define the biadjoint biequivalences, but instead acknowledge this as a missing piece of the construction, which is not likely to trouble the reader to a large extent.

We will need the unit tritransformation
\[ I_{\otimes(\otimes\times 1)}\maps \otimes(\otimes\times 1)\To \otimes(\otimes\times 1)\]
\noindent which consists of
\begin{itemize}
\item for each triple of objects $A, B, C\in\Span(\B)$, a span $I_{ABC}$
\[
   \xy
   (-20,0)*+{(A\times B)\times C}="1";
   (0,15)*+{(A\times B)\times C}="2";
   (20,0)*+{(A\times B)\times C}="3";
        {\ar_{1} "2";"1"};
        {\ar^{1} "2";"3"};
\endxy
\]
\noindent and for each pair of triples of objects $(A,B,C), (A',B',C')$, an adjoint equivalence
\[ (I,I^{\cdot},\epsilon_{I},\eta_{I})\maps (I_{A'B'C'})_*\otimes(\otimes\times 1)\To (I_{ABC})^*\otimes(\otimes\times 1),\]
\noindent consisting of
\begin{itemize}
\item a transformation
\[ I_{(ABC),(A'B'C')}\maps (I_{A'B'C'})_*\otimes(\otimes\times 1)\To (I_{ABC})^*\otimes(\otimes\times 1),\]
\noindent consisting of
\begin{itemize}
\item for each triple of spans $(R,S,T)$, a map of spans $I_{RST}$
\[
\def\objectstyle{\scriptstyle}
  \def\labelstyle{\scriptstyle}
   \xy
   (-20,0)*+{(A'\times B')\times C'}="2";
   (0,-20)*+{((R\times S)\times T)((A\times B)\times C)}="4";
   (0,20)*+{((A'\times B')\times C')((R\times S)\times T)}="1";
   (20,0)*+{(A\times B)\times C}="3";
        {\ar_{\pi_{(A'\times B')\times C'}^{(R\times S)\times T}} "1";"2"};
        {\ar^{((r\times s)\times t)\pi_{(R\times S)\times T}^{(A\times B)\times C}} "1";"3"};
        {\ar^{((r'\times s')\times t')\pi_{(R\times S)\times T}^{(A\times B)\times C}} "4";"2"};
        {\ar_{\pi^{(R\times S)\times T}_{(A\times B)\times C}} "4";"3"};
        {\ar^{I} "1";"4"};
        {\ar@{=}_<<{\scriptstyle } (12,2); (9,-1)};
        {\ar@{=}^<<{\scriptstyle } (-12,2); (-9,-1)};
\endxy
\]
\noindent where $I:=I_{RST}$ is the unique $1$-cell in $\B$ such that
\[ \pi_{(A\times B)\times C}^{(R\times S)\times T}I_{RST} = ((r\times s)\times t)\pi_{(R\times S)\times T}^{(A'\times B')\times C'},\;\;\;\pi_{(R\times S)\times T}^{(A\times B)\times C}I_{RST} = \pi_{(R\times S)\times T}^{(A'\times B')\times C'}\]
\noindent and
\[ \kappa_{(A\times B)\times C}^{(r\times s)\times t,1}\cdot I_{RST} = 1,\]

\item for each pair of triples of spans $(R,S,T), (\bar{R},\bar{S},\bar{T})$, a natural isomorphism
\[ I_{(RST),(\bar{R},\bar{S},\bar{T})}\maps (I_{R,S,T})^*(I_{ABC})^*\otimes(\otimes\times 1) \To (I_{\bar{R},\bar{S},\bar{T}})_{*}(I_{A'B'C'})_*\otimes(\otimes\times 1),\]
\noindent consisting of, for each triple of maps of spans $(f_{R}, f_{S}, f_{T})$, an isomorphism of maps of spans
\[ I_{f_{R}f_{S}f_{T}}\maps (((\varpi_{R}\times \varpi_{S})\times \varpi_{T})\cdot\pi_{(R\times S)\times T}^{(A'\times B')\times C'}, I_{\bar{R},\bar{S},\bar{T}}(((f_{R}\times f_{S})\times f_{T})*1_{(A'\times B')\times C'}),1)\]
\[ \To (1 ,\;(1_{(A\times B)\times C}*((f_{R}\times f_{S})\times f_{T}))I_{R,S,T},\; ((\varrho_{R}\times \varrho_{S})\times \varrho_{T})\cdot\pi_{(R\times S)\times T}^{(A\times B)\times C}I_{R,S,T})\]
\noindent consisting of the unique $2$-cell
\[ I_{f_{R}f_{S}f_{T}}\maps I_{\bar{R},\bar{S},\bar{T}}(((f_{R}\times f_{S})\times f_{T})*1_{(A'\times B')\times C'})\] \[\To (1_{(A\times B)\times C}*((f_{R}\times f_{S})\times f_{T}))I_{R,S,T}\]
\noindent in $\B$, such that
\[ \pi_{(A\times B)\times C}^{(\bar{R}\times\bar{S})\times\bar{T}}\cdot I_{f_{R}f_{S}f_{T}} = \varpi_{(RS)T}^{-1}\cdot\pi_{(R\times S)\times T}^{(A\times B)\times C} \]
\noindent and
\[ \pi^{(A\times B)\times C}_{(\bar{R}\times\bar{S})\times\bar{T}}\cdot I_{f_{R}f_{S}f_{T}} = 1,\]
\end{itemize}

\item a transformation
\[ I^{\cdot}_{(ABC),(A'B'C')}\maps (I_{ABC})^*\otimes(\otimes\times 1) \To  (I_{A'B'C'})_*\otimes(\otimes\times 1),\]
\noindent consisting of
\begin{itemize}
\item for each triple of spans $(R,S,T)$, a map of spans $I^{\cdot}_{RST}$
\[
\def\objectstyle{\scriptstyle}
  \def\labelstyle{\scriptstyle}
   \xy
   (-20,0)*+{(A'\times B')\times C'}="2";
   (0,20)*+{((R\times S)\times T)((A\times B)\times C)}="4";
   (0,-20)*+{((A'\times B')\times C')((R\times S)\times T)}="1";
   (20,0)*+{(A\times B)\times C}="3";
        {\ar^{\pi_{(A'\times B')\times C'}^{(R\times S)\times T}} "1";"2"};
        {\ar_{((r\times s)\times t)\pi_{(R\times S)\times T}^{(A\times B)\times C}} "1";"3"};
        {\ar_{((r'\times s')\times t')\pi_{(R\times S)\times T}^{(A\times B)\times C}} "4";"2"};
        {\ar^{\pi^{(R\times S)\times T}_{(A\times B)\times C}} "4";"3"};
        {\ar^{I^{\cdot}} "4";"1"};
        {\ar@{=}_<<{\scriptstyle } (12,2); (9,-1)};
        {\ar@{=}^<<{\scriptstyle } (-12,2); (-9,-1)};
\endxy
\]
\noindent where $I^{\cdot} := I^{\cdot}_{RST}$ is the unique $1$-cell in $\B$ such that
\[ \pi^{(A'\times B')\times C'}_{(R\times S)\times T}I^{\cdot}_{RST} = \pi_{(R\times S)\times T}^{(A\times B)\times C},\;\;\;\pi^{(R\times S)\times T}_{(A'\times B')\times C'}I^{\cdot}_{RST} =((r'\times s')\times t')\pi_{(R\times S)\times T}^{(A\times B)\times C}\]
\noindent and
\[ \kappa_{(A'\times B')\times C'}^{1,(r'\times s')\times t'}\cdot I^{\cdot}_{RST} = 1,\]

\item for each pair of triples of spans $(R,S,T), (\bar{R},\bar{S},\bar{T})$, a natural isomorphism
\[ I^{\cdot}_{(RST),(\bar{R},\bar{S},\bar{T})}\maps (I^{\cdot}_{R,S,T})^*(I_{ABC})^*\otimes(\otimes\times 1) \To (I_{\bar{R},\bar{S},\bar{T}})_{*}(I^{\cdot}_{A'B'C'})_*\otimes(\otimes\times 1),\]
\noindent consisting of, for each triple of maps of spans $(f_{R}, f_{S}, f_{T})$, an isomorphism of maps of spans
\[ I^{\cdot}_{f_{R}f_{S}f_{T}}\maps (((\varpi_{R}\times \varpi_{S})\times \varpi_{T})\cdot\pi_{(R\times S)\times T}^{(A'\times B')\times C'}, I^{\cdot}_{\bar{R},\bar{S},\bar{T}}(((f_{R}\times f_{S})\times f_{T})*1_{(A'\times B')\times C'}),1)\]
\[ \To (1 ,\;(1_{(A\times B)\times C}*((f_{R}\times f_{S})\times f_{T}))I^{\cdot}_{R,S,T},\; ((\varrho_{R}\times \varrho_{S})\times \varrho_{T})\cdot\pi_{(R\times S)\times T}^{(A\times B)\times C}I^{\cdot}_{R,S,T})\]
\noindent consisting of the unique $2$-cell
\[ I^{\cdot}_{f_{R}f_{S}f_{T}}\maps I^{\cdot}_{\bar{R},\bar{S},\bar{T}}(((f_{R}\times f_{S})\times f_{T})*1_{(A'\times B')\times C'})\]
\[ \To (1_{(A\times B)\times C}*((f_{R}\times f_{S})\times f_{T}))I^{\cdot}_{R,S,T}\]
\noindent in $\B$, such that
\[ \pi^{(A\times B)\times C}_{(\bar{R}\times\bar{S})\times\bar{T}}\cdot I^{\cdot}_{f_{R}f_{S}f_{T}} = 1\]
\noindent and
\[ \pi_{(A\times B)\times C}^{(\bar{R}\times\bar{S})\times\bar{T}}\cdot I^{\cdot}_{f_{R}f_{S}f_{T}} = \varrho_{(f_{R}f_{S})f_{T}}^{-1}\cdot\pi_{(R\times S)\times T}^{(A\times B)\times C},\]

\end{itemize}
\item a modification
\[
   \xy
   (-7,0)*+{\epsilon_{I}\maps II^{\cdot}}="1"; (7,0)*+{1}="2";
        {\ar@3{->}_{} "1";"2"};
\endxy
\]
\noindent consisting of, for each triple of spans $(R,S,T)$, an isomorphism of maps of spans
\[
   \xy
   (-24,0)*+{{\epsilon_{I}}_{RST}\maps (1,\;I_{RST}I_{RST}^{\cdot},\;1)}="1"; (25,0)*+{(1,\;1_{((RS)T)((AB)C)},\; 1),}="2";
        {\ar@3{->}_{} "1";"2"};
\endxy
\]
\noindent consisting of the unique $2$-cell
\[
   \xy
   (-18,0)*+{{\epsilon_{I}}_{RST}\maps I_{RST}I_{RST}^{\cdot}}="1"; (17,0)*+{1_{((RS)T)((AB)C)}}="2";
        {\ar@2{->}_{} "1";"2"};
\endxy
\]
\noindent in $\B$ such that
\[ \pi_{(A\times B)\times C}^{(R\times S)\times T}\cdot {\epsilon_{I}}_{RST} = {\kappa_{(A\times B)\times C}^{(r\times s)\times t,1}}^{-1}\]
\noindent and
\[ \pi^{(A\times B)\times C}_{(R\times S)\times T}\cdot {\epsilon_{I}}_{RST} = 1,\]

\item a modification
\[
   \xy
   (-8,0)*+{\eta_{I}\maps 1}="1"; (7,0)*+{I^{\cdot}I,}="2";
        {\ar@3{->}_{} "1";"2"};
\endxy
\]
\noindent consisting of, for each triple of spans $(R,S,T)$, an isomorphism of maps of spans
\[
   \xy
   (-25,0)*+{{\eta_{I}}_{RST}\maps (1,\;1_{((A'B')C')((RS)T)},\; 1)}="1"; (24,0)*+{(1,\;I_{RST}^{\cdot}I_{RST},\;1),}="2";
        {\ar@3{->}_{} "1";"2"};
\endxy
\]
\noindent consisting of the unique $2$-cell
\[
   \xy
   (-18,0)*+{{\eta_{I}}_{RST}\maps 1_{((A'B')C')((RS)T)}}="1"; (19,0)*+{I_{RST}^{\cdot}I_{RST}}="2";
        {\ar@2{->}_{} "1";"2"};
\endxy
\]
\noindent in $\B$ such that
\[ \pi^{(A'\times B')\times C'}_{(R\times S)\times T}\cdot {\eta_{I}}_{RST} = 1\]
\noindent and
\[ \pi_{(A'\times B')\times C'}^{(R\times S)\times T}\cdot {\eta_{I}}_{RST} = \kappa_{(A'\times B')\times C'}^{1,(r\times s)\times t},\]
\end{itemize}
\item an identity modification
\[ I_{\Pi}\maps (1,({\bf a}\times {\bf a})\chi(\chi*1),1) \To (1,\chi(1*\chi){\bf a},1),\]
\item and an identity modification
\[ I_{M}\maps (1,(1*I)(\iota*1){\bf r}^{-1},1) \To (1,\iota{\bf l}^{-1},1).\]
\end{itemize}

\noindent The remaining unit tritransformations are defined similarly.

\subsubsection*{Monoidal Associativity}

Monoidal associativity is a biadjoint biequivalence consisting of tritransformations, adjoint equivalences of trimodifications and perturbations, and coherence perturbations consisting of isomorphisms of maps of spans.   

The associator for the product of objects $A,B,C\in\B$ is the $1$-cell
\[a_{ABC} \maps (A\times B)\times C\to A\times (B\times C)\]
\noindent in $\B$ defined as the product of $1$-cells
\[ a_{ABC} := \widetilde{\pi}_{A}^{B}\times 1_{B\times C}.\]
\noindent The inverse associator is the $1$-cell
\[ a_{ABC}^{-1}\maps A\times (B\times C)\to (A\times B)\times C\]
\noindent in $\B$ defined as the product of $1$-cells
\[ a^{-1}_{ABC} := 1_{A\times B}\times \widetilde{\pi}_{C}^{B}.\]

\proposition \label{monoidalassociator}
There is a biadjoint biequivalence
\[
   \xy
   (-22,0)*+{(\alpha, \alpha^{\cdot}, \epsilon_{\alpha}, \eta_{\alpha}, \Phi_{\alpha}, \Psi_{\alpha})\maps \otimes(\otimes\times 1)}="1"; (22,0)*+{\otimes(1\times\otimes)}="2";
        {\ar@3{->}_{} "1";"2"};
\endxy
\]
\noindent in a `tricategory' $\Tricat(\Span(\B)^{3},\Span(\B))$, consisting of
\begin{itemize}
\item a tritransformation
\[ \alpha\maps \otimes(\otimes\times 1)\To \otimes(1\times\otimes),\]
\noindent consisting of
\begin{itemize}
\item for each triple $A, B, C$ of objects in $\Span(\B)$, a span $\alpha_{ABC}$
\[
   \xy
   (0,15)*+{(A\times B)\times C}="1";
   (-20,0)*+{A\times (B\times C)}="2"; (20,0)*+{(A\times B)\times C}="3";
        {\ar_{a} "1";"2"};
        {\ar^{1} "1";"3"};
\endxy
\]

\item for each two triples $(A, B, C), (A', B', C')$ of objects in $\Span(\B)$, an adjoint equivalence
\[ ({\alpha}_{\mu}, \alpha_{\mu^{\cdot}}, \alpha_{\epsilon}, \alpha_{\eta}) \maps \Span(\B)(1,\alpha_{A'B'C'})(\otimes(\otimes\times 1))\To \Span(\B)(\alpha_{ABC},1)(\otimes(1\times\otimes)),\]
\noindent in the strict $2$-category
\[\Bicat(\Span(B)^{3}((A,B,C),(A',B',C')),\Span(\B)((A\times B) \times C,A'\times (B'\times C')),\]
\noindent consisting of
\begin{itemize}
\item a strong transformation
\[ {\alpha_{\mu}}_{(A, B, C), (A', B', C')}\maps \Span(\B)(1,\alpha_{A'B'C'})(\otimes(\otimes\times 1))\]
\[\To \Span(\B)(\alpha_{ABC},1)(\otimes(1\times\otimes)),\]
\noindent consisting of
\begin{itemize}
\item for each triple of spans
\[
\xymatrix{
& R\ar[dl]_{r'}\ar[dr]^{r} &  & &S\ar[dl]_{s'}\ar[dr]^{s} &&& T\ar[dl]_{t'}\ar[dr]^{t} &&\\
A' && A & B' && B & C' && C
}
\]
\noindent a map of spans
\[
\def\objectstyle{\scriptstyle}
  \def\labelstyle{\scriptstyle}
   \xy
   (-20,0)*+{A'\times (B'\times C')}="2";
   (0,-20)*+{(R\times (S\times T))((A\times B)\times C)}="1";
   (0,20)*+{((A'\times B')\times C')((R\times S)\times T)}="4";
   (20,0)*+{(A\times B)\times C}="3";
        {\ar^{(r'\times (s'\times t'))\pi^{(A\times B)\times C}_{R\times (S\times T)}} "1";"2"};
        {\ar_{\pi^{R\times (S\times T)}_{(A\times B)\times C}} "1";"3"};
        {\ar_{a\pi^{(R\times S)\times T}_{(A'\times B')\times C'}} "4";"2"};
        {\ar^{((r\times s)\times t)\pi^{(A'\times B')\times C'}_{(R\times S)\times T}} "4";"3"};
        {\ar^{\alpha_{\mu}} "4";"1"};
        {\ar@{=}_<<{\scriptstyle } (12,2); (9,-1)};
        {\ar@{=>}^<<{\scriptstyle a\cdot\kappa^{-1}} (-12,2); (-9,-1)};
\endxy
\]
\noindent where $\kappa := \kappa_{(A'\times B')\times C'}^{1,(r'\times s')\times t'}$ and $\alpha_{\mu}:= {\alpha_{\mu}}_{RST}$ is the unique $1$-cell satisfying
\[ \pi_{(A\times B)\times C}^{R\times (S\times T)}{\alpha_{\mu}}_{RST} = ((r\times s)\times t)\pi_{(R\times S)\times T}^{(A'\times B')\times C'}\;\;\; \pi^{(A\times B)\times C}_{R\times (S\times T)}{\alpha_{\mu}}_{RST} = a\pi_{(R\times S)\times T}^{(A'\times B')\times C'},\]
\noindent and
\[ \kappa^{r\times (s\times t),1}_{A\times (B\times C)}\cdot{\alpha_{\mu}}_{RST} = 1,\]

\item for each pair of triples of spans $(R, S, T), (\bar{R}, \bar{S}, \bar{T})$, a natural isomorphism
\[ {\alpha_{\mu}}_{RST,\bar{R}\bar{S}\bar{T}}\maps ({\alpha_{\mu}}_{\bar{R},\bar{S},\bar{T}})_*\Span(\B)(1,\alpha_{A'B'C'})(\otimes(\otimes\times 1))\To\]
\[ ({\alpha_{\mu}}_{R,S,T})^*\Span(\B)(\alpha_{ABC},1)(\otimes(1\times\otimes)),\]
\noindent consisting of, for each triple of maps of spans
\[
\xymatrix{
& R\ar[dd]^{f_{R}}\ar[dl]_{r'}\ar[dr]^{r} &  & &S\ar[dd]^{f_{S}}\ar[dl]_{s'}\ar[dr]^{s} &&& T\ar[dd]^{f_{T}}\ar[dl]_{t'}\ar[dr]^{t} &&\\
A' && A & B' && B & C' && C\\
& \bar{R}\ar[ul]^{\bar{r}'}\ar[ur]_{\bar{r}} &  & &\bar{S}\ar[ul]^{\bar{s}'}\ar[ur]_{\bar{s}} &&& \bar{T}\ar[ul]^{\bar{t}'}\ar[ur]_{\bar{t}} &&
}
\]
\noindent an isomorphism of maps of spans 
\[{\alpha_{\mu}}_{f_{R},f_{S},f_{T}}\maps (\varpi_{(RS)T}\cdot\pi_{(R\times S)\times T}^{(A'\times B')\times C'},\;{\alpha_{\mu}}_{\bar{R},\bar{S},\bar{T}}(((f_{R}\times f_{S})\times f_{T})*1),\; a\cdot\kappa^{-1}\cdot (((f_{R}\times f_{S})\times f_{T})*1))\]
\[ \To (1,\;(1*(f_{R}\times (f_{S}\times f_{T}))){\alpha_{\mu}}_{R,S,T},\; (\varrho_{R(ST)}\cdot\pi_{R\times (S\times T)}^{(A\times B)\times C}{\alpha_{\mu}}_{R,S,T})(a\cdot\kappa^{-1})),\]
\noindent consisting of the unique $2$-cell
\[  {\alpha_{\mu}}_{f_{R},f_{S},f_{T}}\maps {\alpha_{\mu}}_{\bar{R},\bar{S},\bar{T}}(((f_{R}\times f_{S})\times f_{T})*1) \To (1*(f_{R}\times (f_{S}\times f_{T}))){\alpha_{\mu}}_{R,S,T}\]
\noindent in $\B$ such that
\[ \pi_{(A\times B)\times C}^{\bar{R}\times (\bar{S}\times \bar{T})}\cdot  {\alpha_{\mu}}_{f_{R},f_{S},f_{T}} = ((\varpi_{R}\times\varpi_{S})\times\varpi_{T})^{-1}\cdot \pi_{(R\times S)\times T}^{(A'\times B')\times C'}\]
\noindent and
\[ \pi_{\bar{R}\times (\bar{S}\times\bar{T})}^{(A\times B)\times C}\cdot {\alpha_{\mu}}_{f_{R},f_{S},f_{T}} = 1,\] 
\end{itemize}

\item a strong transformation
\[ {\alpha_{\mu^\cdot}}_{(A, B, C), (A', B', C')}\maps \Span(\B)(\alpha_{ABC},1)(\otimes(1\times\otimes))\]
\[ \To  \Span(\B)(1,\alpha_{A'B'C'})(\otimes(\otimes\times 1)),\]
\noindent consisting of
\begin{itemize}
\item for each triple of spans $R$, $S$, $T$, a map of spans
\[
\def\objectstyle{\scriptstyle}
  \def\labelstyle{\scriptstyle}
   \xy
   (-20,0)*+{A'\times (B'\times C')}="2";
   (0,-20)*+{((A'\times B')\times C')((R\times S)\times T)}="1";
   (0,20)*+{(R\times (S\times T))((A\times B)\times C)}="4";
   (20,0)*+{(A\times B)\times C}="3";
        {\ar^{a\pi_{(A'\times B')\times C'}^{(R\times S)\times T}} "1";"2"};
        {\ar_{((r\times s)\times t)\pi_{(R\times S)\times T}^{(A'\times B')\times C'}} "1";"3"};
        {\ar_{(r'\times(s'\times t'))\pi_{R\times(S\times T)}^{(A\times B)\times C}} "4";"2"};
        {\ar^{\pi_{(A\times B)\times C}^{R\times(S\times T)}} "4";"3"};
        {\ar^{{\alpha_{\mu^\cdot}}} "4";"1"};
        {\ar@{=>}_<<{\scriptstyle a^{-1}\cdot\kappa} (12,2); (9,-1)};
        {\ar@{=}_<<{\scriptstyle } (-12,2); (-9,-1)};
\endxy
\]
\noindent where $\kappa := \kappa_{A\times (B\times C)}^{r\times (s\times t),1}$ and $\alpha_{\mu^{\cdot}}:= {\alpha_{\mu^{\cdot}}}_{RST}$ is the unique $1$-cell satisfying
\[ \pi_{(R\times S)\times T}^{(A'\times B')\times C'}{\alpha_{\mu^{\cdot}}}_{RST} = a^{-1}\pi_{R\times (S\times T)}^{(A\times B)\times C}\;\;\; \pi^{(R\times S)\times T}_{(A'\times B')\times C'}{\alpha_{\mu^{\cdot}}}_{RST} = ((r'\times s')\times t')a^{-1}\pi_{R\times (S\times T)}^{(A\times B)\times C}\]
\noindent and
\[ \kappa^{1,(r'\times s')\times t'}_{(A'\times B')\times C'}\cdot {\alpha_{\mu^{\cdot}}}_{RST} = 1,\]

\item for each pair of triples of spans $(R, S, T), (\bar{R}, \bar{S}, \bar{T})$, a natural isomorphism
\[ {\alpha_{\mu^{\cdot}}}_{RST,\bar{R}\bar{S}\bar{T}} \maps ({\alpha_{\mu^{\cdot}}}_{\bar{R},\bar{S},\bar{T}})_*\Span(\B)(\alpha_{ABC},1)(\otimes(1\times\otimes))\To \]
\[({\alpha_{\mu^{\cdot}}}_{R,S,T})^*\Span(\B)(1,\alpha_{A'B'C'})(\otimes(\otimes\times 1)),\]
\noindent consisting of, for each triple of maps of spans $f_{R}$, $f_{S}$, $f_{T}$, an isomorphism of maps of spans
\[ {\alpha_{\mu^{\cdot}}}_{f_{R},f_{S},f_{T}}\maps (a^{-1}\cdot\kappa\cdot(1*(f_{R}\times (f_{S}\times f_{T}))),\;{\alpha_{\mu^{\cdot}}}_{\bar{R},\bar{S},\bar{T}}(1*(f_{R}\times (f_{S}\times f_{T}))),\; \varrho_{R(ST)}\cdot\pi_{R\times (S\times T)}^{(A\times B)\times C}) \]
\[\To ( a^{-1}\cdot\kappa,\;(((f_{R}\times f_{S})\times f_{T})*1){\alpha_{\mu^{\cdot}}}_{R,S,T},\; (\varrho_{(RS)T}\cdot\pi_{(R\times S)\times T}^{(A'\times B')\times C'})\cdot{\alpha_{\mu^{\cdot}}}_{R,S,T}),\]
\noindent consisting of the unique $2$-cell
\[  {\alpha_{\mu^{\cdot}}}_{f_{R},f_{S},f_{T}}\maps {\alpha_{\mu^{\cdot}}}_{\bar{R},\bar{S},\bar{T}}(1*(f_{R}\times (f_{S}\times f_{T}))) \To (((f_{R}\times f_{S})\times f_{T})*1){\alpha_{\mu^{\cdot}}}_{R,S,T}\]
\noindent in $\B$ such that
\[ \pi^{(A'\times B')\times C'}_{(\bar{R}\times \bar{S})\times \bar{T}}\cdot  {\alpha_{\mu^{\cdot}}}_{f_{R},f_{S},f_{T}} = 1\]
\noindent and
\[ \pi_{(A'\times B')\times C'}^{(\bar{R}\times \bar{S})\times \bar{T}}\cdot {\alpha_{\mu^{\cdot}}}_{f_{R},f_{S},f_{T}} = ((\varrho_{R}\times\varrho_{S})\times\varrho_{T})^{-1}\cdot a^{-1}\pi_{R\times (S\times T)}^{(A\times B)\times C},\] 
\end{itemize}

\item an invertible counit modification
\[\alpha_{\epsilon}\maps {\alpha_{\mu}}{\alpha_{\mu^{\cdot}}}\To 1_{\otimes(1\times\otimes)},\]
\noindent consisting of, for each triple of spans $R$, $S$, $T$, an isomorphism of maps of spans
\[ {\alpha_{\epsilon}}_{RST}\maps  (a^{-1}\cdot\kappa,\;\alpha_{\mu}{\alpha_{\mu^{\cdot}}},\; a\cdot\kappa^{-1}\cdot{\alpha_{\mu^{\cdot}}})\To (1,1_{\otimes(1\times\otimes)},1)\]
\noindent defined by the unique $2$-cell ${\alpha_{\epsilon}}_{RST}$ in $\B$ such that
\[ \pi^{R\times (S\times T)}_{(A\times B)\times C}\cdot{\alpha_{\epsilon}}_{RST}  = a^{-1}\cdot{\kappa^{-1}}_{A\times (B\times C)}^{r\times (s\times t),1}\;\;\;\textrm{ and }\;\;\; \pi_{R\times (S\times T)}^{(A\times B)\times C}\cdot{\alpha_{\epsilon}}_{RST} =  1,\]
\item an invertible unit modification
\[\alpha_{\eta}\maps 1_{\otimes(\otimes\times 1)}\To {\alpha_{\mu^{\cdot}}}{\alpha_{\mu}}\]
\noindent consisting of, for each triple of spans $R$, $S$, $T$, an isomorphism of maps of spans
\[ {\alpha_{\eta}}_{RST}\maps (1,1_{\otimes(\otimes\times 1)},1)\To (a^{-1}\cdot\kappa\cdot\alpha_{\mu},\;{\alpha_{\mu^{\cdot}}}{\alpha_{\mu}},\; a\cdot\kappa^{-1})\]
\noindent defined by the unique $2$-cell ${\alpha_{\eta}}_{RST}$ in $\B$ such that
\[  \pi^{(A'\times B')\times C'}_{(R\times S)\times T}\cdot{\alpha_{\eta}}_{RST} = 1\;\;\;\textrm{ and }\;\;\;  \pi_{(A'\times B')\times C'}^{(R\times S)\times T}\cdot{\alpha_{\eta}}_{RST} = {\kappa^{-1}}_{(A'\times B')\times C'}^{1,(r'\times s')\times t'},\]
\end{itemize}

\item an identity modification $\alpha_{\Pi}$ with component equations of maps of spans
\[ (1,\; (1*\chi)({\alpha_{\mu}}_{RST}*1)(1*{\alpha_{\mu}}_{R'S'T'}),\; a\cdot \kappa^{-1}\cdot \pi) = (1,\;{\alpha_{\mu}}_{(R'R)(S'S)(T'T)}(\chi*1),\; a\cdot \kappa^{-1}\cdot(1*\chi)),\]

\item an identity modification $\alpha_{M}$ with component equations of maps of spans
\[ (1,\; \alpha_{\mu}\iota {\bf r}^{-1},\; a\cdot \kappa^{-1}\cdot\iota {\bf r}^{-1}) = (1,\; \iota {\bf l}^{-1},\; 1),\]
\end{itemize}

\item a tritransformation
\[ \alpha^{\cdot}\maps \otimes(1\times\otimes)\To \otimes(\otimes\times 1),\]
\noindent consisting of
\begin{itemize}
\item for each triple $A, B, C$ of objects in $\Span(\B)$, a span $\alpha^{\cdot}_{ABC}$
\[
\xymatrix{
&(A\times B)\times C\ar[dl]_{1}\ar[dr]^{a}&\\
(A\times B)\times C && A\times (B\times C)
}
\]

\item for each two triples $(A, B, C), (A', B', C')$ of objects in $\Span(\B)$, an adjoint equivalence
\[ ({\alpha}^{\cdot}_{\mu}, \alpha^{\cdot}_{\mu^{\cdot}}, \alpha^{\cdot}_{\epsilon}, \alpha^{\cdot}_{\eta}) \maps \Span(\B)(1,\alpha^{\cdot}_{A'B'C'})(\otimes(1\times\otimes))\To \Span(\B)(\alpha^{\cdot}_{ABC},1)(\otimes(\otimes\times 1)),\]
\noindent in the strict $2$-category
\[ \Bicat(\Span(\B)^{3}((A,B,C),(A',B',C')), \Span(\B)(A\times (B\times C), (A'\times B')\times C')),\]
\noindent consisting of

\begin{itemize}
\item a strong transformation
\[ {\alpha_{\mu}^{\cdot}}_{(A,B,C),(A',B',C')}\maps \Span(\B)(1,\alpha^{\cdot}_{A'B'C'})(\otimes(1\times\otimes))\To \Span(\B)(\alpha^{\cdot}_{ABC},1)(\otimes(\otimes\times 1)),\]
\noindent consisting of
\begin{itemize}
\item for each triple of spans $R$, $S$, $T$, a map of spans
\[
\def\objectstyle{\scriptstyle}
  \def\labelstyle{\scriptstyle}
   \xy
   (-20,0)*+{(A'\times B')\times C'}="2";
   (0,20)*+{((A'\times B')\times C')(R\times (S\times T))}="1";
   (0,-20)*+{((R\times S)\times T)((A\times B)\times C)}="4";
   (20,0)*+{A\times (B\times C)}="3";
        {\ar_{\pi_{(A'\times B')\times C'}^{R\times (S\times T)}} "1";"2"};
        {\ar^{(r\times (s\times t))\pi_{R\times (S\times T)}^{(A'\times B')\times C'}} "1";"3"};
        {\ar^{((r'\times s')\times t')\pi_{(R\times S)\times T}^{(A\times B)\times C}} "4";"2"};
        {\ar_{a\pi_{(A\times B)\times C}^{(R\times S)\times T}} "4";"3"};
        {\ar^{\alpha_{\mu}^{\cdot}} "1";"4"};
        {\ar@{=}^<<{\scriptstyle } (12,2); (9,-1)};
        {\ar@{=>}^<<{\scriptstyle a^{-1}\cdot\kappa} (-12,2); (-9,-1)};
\endxy
\]
\noindent where $\kappa := \kappa_{A'\times (B'\times C')}^{1,r'\times (s'\times t')}$ and $\alpha_{\mu}^{\cdot}:= {\alpha_{\mu}^{\cdot}}_{RST}$ is the unique $1$-cell satisfying
\[ \pi^{(R\times S)\times T}_{(A\times B)\times C}{\alpha_{\mu}^{\cdot}}_{RST} = ((r\times s)\times t)a^{-1}\pi_{R\times (S\times T)}^{(A'\times B')\times C'}\;\;\; \pi_{(R\times S)\times T}^{(A\times B)\times C}{\alpha_{\mu}^{\cdot}}_{RST} = a^{-1}\pi_{R\times (S\times T)}^{(A'\times B')\times C'}\]
\noindent and
\[ \kappa^{(r\times s)\times t,1}_{(A\times B)\times C}\cdot{\alpha_{\mu}^{\cdot}}_{RST} = 1,\]
\item for each pair of triples of spans $(R, S, T), (\bar{R},\bar{S},\bar{T})$, a natural isomorphism
\[ {\alpha_{\mu}^{\cdot}}_{RST,\bar{R}\bar{S}\bar{T}}\maps ({\alpha_{\mu}^{\cdot}}_{\bar{R},\bar{S},\bar{T}})_*\Span(\B)(1,\alpha^{\cdot}_{A',B',C'})(\otimes(1\times\otimes))\To \]
\[({\alpha_{\mu}^{\cdot}}_{R,S,T})^{*}\Span(\B)(\alpha^{\cdot}_{ABC},1)(\otimes(\otimes\times 1))\]
\noindent consisting of, for each triple of maps of spans $f_{R}$, $f_{S}$, $f_{T}$, an isomorphism of maps of spans 
\[{\alpha_{\mu}^{\cdot}}_{f_{R},f_{S},f_{T}}\maps (\varpi_{R(ST)}\cdot\pi_{R\times(S\times T)}^{(A'\times B')\times C'},\;{\alpha^{\cdot}_{\mu}}_{\bar{R},\bar{S},\bar{T}}((f_{R}\times (f_{S}\times f_{T}))*1),\; a^{-1}\cdot\kappa\cdot ((f_{R}\times (f_{S}\times f_{T}))*1))\]
\[\To (1,\;(1*((f_{R}\times f_{S})\times f_{T})){{\alpha}^{\cdot}_{\mu}}_{R,S,T},\; ((\varrho_{(RS)T}\cdot\pi_{(R\times S)\times T}^{(A'\times B')\times C'})\cdot{\alpha^{\cdot}_{\mu}}_{R,S,T})(a^{-1}\cdot\kappa)),\]

\noindent consisting of the unique $2$-cell
\[  {\alpha_{\mu}^{\cdot}}_{f_{R},f_{S},f_{T}}\maps {\alpha^{\cdot}_{\mu}}_{\bar{R},\bar{S},\bar{T}}((f_{R}\times (f_{S}\times f_{T}))*1) \To (1*((f_{R}\times f_{S})\times f_{T})){{\alpha}^{\cdot}_{\mu}}_{R,S,T}\]
\noindent in $\B$ such that
\[ \pi^{(A\times B)\times C}_{(\bar{R}\times \bar{S})\times \bar{T}}\cdot  {\alpha_{\mu}^{\cdot}}_{f_{R},f_{S},f_{T}} = 1\]
\noindent and
\[ \pi_{(A\times B)\times C}^{(\bar{R}\times \bar{S})\times \bar{T}}\cdot {\alpha_{\mu}^{\cdot}}_{f_{R},f_{S},f_{T}} = ((\varrho_{R}\times\varrho_{S})\times\varrho_{T})^{-1}\cdot a^{-1}\pi_{R\times (S\times T)}^{(A'\times B')\times C'},\] 
\end{itemize}

\item a strong transformation
\[ \alpha_{\mu^{\cdot}}^{\cdot}\maps \Span(\B)(\alpha^{\cdot}_{ABC},1)(\otimes(\otimes\times 1))\To \Span(\B)(1,\alpha_{A'B'C'}^{\cdot})(\otimes(1\times\otimes)),\]
\noindent consisting of, 
\begin{itemize}
\item for each triple of spans $R$, $S$, $T$, a map of spans
\[
\def\objectstyle{\scriptstyle}
  \def\labelstyle{\scriptstyle}
   \xy
   (-20,0)*+{(A'\times B')\times C'}="2";
   (0,-20)*+{((A'\times B')\times C')(R\times (S\times T))}="1";
   (0,20)*+{((R\times S)\times T)((A\times B)\times C)}="4";
   (20,0)*+{A\times (B\times C)}="3";
        {\ar^{\pi_{(A'\times B')\times C'}^{R\times (S\times T)}} "1";"2"};
        {\ar_{(r\times (s\times t))\pi_{R\times (S\times T)}^{(A'\times B')\times C'}} "1";"3"};
        {\ar_{((r'\times s')\times t')\pi_{(R\times S)\times T}^{(A\times B)\times C}} "4";"2"};
        {\ar^{a\pi_{(A\times B)\times C}^{(R\times S)\times T}} "4";"3"};
        {\ar_{{\alpha^{\cdot}_{\mu^{\cdot}}}} "4";"1"};
        {\ar@{=>}_<<{\scriptstyle a\cdot\kappa} (12,2); (9,-1)};
        {\ar@{=}^<<{\scriptstyle } (-12,2); (-9,-1)};
\endxy
\]
\noindent where $\kappa := \kappa_{(A\times B)\times C}^{(r\times s)\times t,1}$ and ${\alpha^{\cdot}_{\mu^{\cdot}}}:= {\alpha^{\cdot}_{\mu^{\cdot}}}_{RST}$ is the unique $1$-cell satisfying
\[ \pi^{(A'\times B')\times C'}_{R\times (S\times T)}{\alpha^{\cdot}_{\mu^{\cdot}}}_{RST} =  a\pi_{(R\times S)\times T}^{(A\times B)\times C}\;\;\; \pi_{(A'\times B')\times C'}^{R\times (S\times T)}{\alpha^{\cdot}_{\mu^{\cdot}}}_{RST} = ((r'\times s')\times t')\pi_{(R\times S)\times T}^{(A\times B)\times C},\]
\noindent and
\[ \kappa^{1,r'\times (s'\times t')}_{A'\times (B'\times C')}\cdot {\alpha^{\cdot}_{\mu^{\cdot}}}_{RST} = 1,\]

\item for each pair of triples of spans $(R,S,T), (\bar{R},\bar{S},\bar{T})$, a natural isomorphism
\[ {\alpha^{\cdot}_{\mu^{\cdot}}}_{RST,\bar{R}\bar{S}\bar{T}}\maps ({\alpha_{\mu^{\cdot}}^{\cdot}}_{\bar{R},\bar{S},\bar{T}})_* \Span(\B)(\alpha^{\cdot}_{ABC},1)(\otimes(\otimes\times 1))\To\]
\[  ({\alpha_{\mu^{\cdot}}^{\cdot}}_{R,S,T})^*\Span(\B)(1,\alpha_{A'B'C'}^{\cdot})(\otimes(1\times\otimes)),\]
\noindent consisting of, for each triple of maps of spans $f_{R}$, $f_{S}$, $f_{T}$, an isomorphism of maps of spans
\[ {\alpha_{\mu^{\cdot}}^{\cdot}}_{f_{R},f_{S},f_{T}} \maps ((a\cdot\kappa\cdot(1*((f_{R}\times f_{S})\times f_{T}))),\; {\alpha^{\cdot}_{\mu^{\cdot}}}_{\bar{R},\bar{S},\bar{T}}(1*((f_{R}\times f_{S})\times f_{T})),\; \varrho_{(RS)T}\cdot \pi_{(R\times S)\times T}^{(A\times B)\times C})\]
\[\To ((\varpi_{R(ST)}\cdot \pi_{R\times(S\times T)}^{(A\times B)\times C}{\alpha^{\cdot}_{\mu^{\cdot}}}_{R,S,T})(a\cdot\kappa),\; ((f_{R}\times (f_{S}\times f_{T}))*1){\alpha^{\cdot}_{\mu^{\cdot}}}_{R,S,T},\; 1),\]
\noindent consisting of the unique $2$-cell
\[{\alpha_{\mu^{\cdot}}^{\cdot}}_{f_{R},f_{S},f_{T}}\maps{\alpha^{\cdot}_{\mu^{\cdot}}}_{\bar{R},\bar{S},\bar{T}}(1*((f_{R}\times f_{S})\times f_{T}))\To ((f_{R}\times (f_{S}\times f_{T}))*1){\alpha^{\cdot}_{\mu^{\cdot}}}_{R,S,T}\]
\noindent in $\B$ such that
\[ \pi_{\bar{R}\times (\bar{S}\times\bar{T})}^{A'\times (B'\times C')}\cdot {\alpha_{\mu^{\cdot}}^{\cdot}}_{f_{R},f_{S},f_{T}} = 1 \]
\noindent and
\[ \pi^{\bar{R}\times (\bar{S}\times\bar{T})}_{A'\times (B'\times C')}\cdot {\alpha_{\mu^{\cdot}}^{\cdot}}_{f_{R},f_{S},f_{T}} =  ((\varrho_{R}\times\varrho_{S})\times\varrho_{T})^{-1}\cdot\pi_{(R\times S)\times T}^{(A\times B)\times C},\]
\end{itemize}

\item an invertible counit modification
\[\alpha^{\cdot}_{\epsilon_{\mu}}\maps {\alpha^{\cdot}_{\mu}}{\alpha^{\cdot}_{\mu^{\cdot}}}\To 1_{\otimes(\otimes\times 1)}\]
\noindent consisting of, for each triple of spans $R$, $S$, $T$, an isomorphism of maps of spans
\[ {\alpha^{\cdot}_{\epsilon_{\mu}}}_{RST}\maps  (a^{-1}\cdot\kappa,\;\alpha^{\cdot}_{\mu}{\alpha^{\cdot}_{\mu^{\cdot}}},\; a\cdot\kappa^{-1}\cdot{\alpha^{\cdot}_{\mu^{\cdot}}})\To (1,1_{\otimes(\otimes\times 1)},1)\]
\noindent defined by the unique $2$-cell ${\alpha^{\cdot}_{\epsilon_{\mu}}}_{RST}$ in $\B$ such that
\[ \pi^{(R\times S)\times T}_{(A\times B)\times C}\cdot{\alpha^{\cdot}_{\epsilon_{\mu}}}_{RST} = a^{-1}\cdot{\kappa^{-1}}_{A\times (B\times C)}^{r\times (s\times t),1}\;\;\;\textrm{ and }\;\;\; \pi_{(R\times S)\times T}^{(A\times B)\times C}\cdot{\alpha^{\cdot}_{\epsilon_{\mu}}}_{RST} =  1,\]
\item an invertible unit modification
\[\alpha^{\cdot}_{\eta_{\mu}}\maps 1_{\otimes(1\times\otimes)}\To {\alpha^{\cdot}_{\mu^{\cdot}}}{\alpha^{\cdot}_{\mu}}\]
\noindent consisting of, for each triple of spans $R$, $S$, $T$, an isomorphism of maps of spans
\[ {\alpha^{\cdot}_{\eta_{\mu}}}_{RST}\maps (1,1_{\otimes(1\times\otimes)},1)\To (a^{-1}\cdot\kappa\cdot\alpha^{\cdot}_{\mu},\;{\alpha^{\cdot}_{\mu^{\cdot}}}{\alpha^{\cdot}_{\mu}},\; a\cdot\kappa^{-1})\]
\noindent defined by the unique $2$-cell ${\alpha^{\cdot}_{\eta_{\mu}}}_{RST}$ in $\B$ such that
\[  \pi^{(A'\times B')\times C'}_{(R\times S)\times T}\cdot{\alpha^{\cdot}_{\eta_{\mu}}}_{RST} = 1\;\;\;\textrm{ and }\;\;\;  \pi_{(A'\times B')\times C'}^{(R\times S)\times T}\cdot{\alpha^{\cdot}_{\eta_{\mu}}}_{RST} = \kappa_{(A'\times B')\times C'}^{1,(r\times s')\times t'},\]
\end{itemize}

\item an identity modification $\alpha^{\cdot}_{\Pi}$ with component equations of maps of spans
\[ (1,\; (1*\chi)({\alpha^{\cdot}_{\mu}}_{RST}*1)(1*{\alpha_{\mu}^{\cdot}}_{R'S'T'}),\; a^{-1}\cdot \kappa\cdot \pi) = (1,\;{\alpha^{\cdot}_{\mu}}_{(R'R)(S'S)(T'T)}(\chi*1),\; a^{-1}\cdot \kappa\cdot(1*\chi)),\]

\item an identity modification $\alpha^{\cdot}_{M}$ with component equations of maps of spans
\[ (1,\; \alpha_{\mu}^{\cdot}\iota {\bf r}^{-1},\; a^{-1}\cdot \kappa\cdot\iota {\bf r}^{-1}) = (1,\; \iota {\bf l}^{-1},\; 1),\]
\end{itemize}

\item a strict adjoint equivalence $(\epsilon_{\alpha}, \epsilon^{\cdot}_{\alpha}, \epsilon_{\epsilon_{\alpha}}, \eta_{\epsilon_{\alpha}})$, in the strict $2$-category
\[\Bicat(\Span(B)^{3}((A,B,C),(A',B',C')),\Span(\B)(A\times (B\times C),A'\times (B'\times C'))),\]
\noindent consisting of
\begin{itemize}
\item a trimodification
\[ \epsilon_{\alpha}\maps  \alpha\alpha^{\cdot} \Rightarrow I_{\otimes(1\times\otimes)}\]
\noindent consisting of
\begin{itemize}
\item for each triple of objects $A, B, C\in\Span(\B)$, a map of spans
\[
\def\objectstyle{\scriptstyle}
  \def\labelstyle{\scriptstyle}
   \xy
   (20,0)*+{A\times (B\times C)}="2";
   (0,20)*+{((A\times B)\times C)((A\times B)\times C)}="1";
   (0,-20)*+{A\times (B\times C)}="4";
   (-20,0)*+{A\times (B\times C)}="3";
        {\ar^{a\pi_{(A\times B)\times C}^{(A\times B)\times C}} "1";"2"};
        {\ar_{a\pi_{(A\times B)\times C}^{(A\times B)\times C}} "1";"3"};
        {\ar_{1} "4";"2"};
        {\ar^{1} "4";"3"};
        {\ar_{{\epsilon_{\alpha}}} "1";"4"};
        {\ar@{=}_<<{\scriptstyle } (12,2); (9,-1)};
        {\ar@{=}^<<{\scriptstyle } (-12,2); (-9,-1)};
\endxy
\]
\noindent where ${\epsilon_{\alpha}} := {\epsilon_{\alpha}}_{A,B,C} = a\pi_{(A\times B)\times C}^{(A\times B)\times C}$,

\item for each pair of triples of objects $(A, B, C), (A', B', C')$, an identity modification
\[ m_{\epsilon_{\alpha}}\maps ({{\epsilon_{\alpha}}}_{(A,B,C)})^*\alpha\alpha^{\cdot}\To I_{\otimes(1\times\otimes)}({{\epsilon_{\alpha}}}_{(A',B',C')})_*\]
\noindent consisting of, for each triple of spans $R$, $S$, $T$, an equation of maps of spans
\[ (1,\;({{\epsilon_{\alpha}}}_{ABC}*1_{R(ST)})(1_{(AB)C}*{\alpha_{\mu}}_{RST})({\alpha^{\cdot}_{\mu}}_{RST}*1_{(A'B')C'}),\; a\cdot\kappa^{-1}\cdot\pi({\alpha^{\cdot}_{\mu}}_{RST}*1_{(A'B')C'})) \]
\[ = (1,\;I_{R(ST)}(1_{R(ST)}*{{\epsilon_{\alpha}}}_{A'B'C'}),\; 1),\]

\end{itemize}

\item a trimodification
\[ \epsilon^{\cdot}_{\alpha}\maps I_{\otimes(1\times\otimes)}\Rightarrow \alpha\alpha^{\cdot}\]
\noindent consisting of
\begin{itemize}
\item for each triple of objects $A, B, C\in\Span(\B)$, a map of spans
\[
\def\objectstyle{\scriptstyle}
  \def\labelstyle{\scriptstyle}
   \xy
   (20,0)*+{A\times (B\times C)}="2";
   (0,20)*+{A\times (B\times C)}="1";
   (0,-20)*+{((A\times B)\times C)((A\times B)\times C)}="4";
   (-20,0)*+{A\times (B\times C)}="3";
        {\ar_{a\pi_{(A\times B)\times C}^{(A\times B)\times C}} "4";"2"};
        {\ar^{a\pi_{(A\times B)\times C}^{(A\times B)\times C}} "4";"3"};
        {\ar^{1} "1";"2"};
        {\ar_{1} "1";"3"};
        {\ar_{{\epsilon^{\cdot}_{\alpha}}} "1";"4"};
        {\ar@{=}_<<{\scriptstyle } (12,2); (9,-1)};
        {\ar@{=}^<<{\scriptstyle } (-12,2); (-9,-1)};
\endxy
\]
\noindent where ${\epsilon^{\cdot}_{\alpha}} := {\epsilon^{\cdot}_{\alpha}}_{A,B,C}$ is the unique $1$-cell in $\B$ satisfying
\[ \pi^{(A\times B)\times C}_{(A\times B)\times C}{\epsilon^{\cdot}_{\alpha}}_{A,B,C} =  a^{-1},\]

\item for each pair of triples of objects $(A, B, C), (A', B', C')$, an identity modification
\[ m_{\epsilon^{\cdot}_{\alpha}}\maps ({m_{\epsilon^{\cdot}_{\alpha}}}_{(A,B,C)})^*I_{\otimes(1\times\otimes)}  \To \alpha\alpha^{\cdot}({m_{\epsilon^{\cdot}_{\alpha}}}_{A',B',C'})_*\]
\noindent consisting of, for each triple of spans $R$, $S$, $T$, an equation of maps of spans
\[  (1,\;({{\epsilon^{\cdot}_{\alpha}}}_{ABC}*1_{R(ST)})I_{R(ST)},\; 1) = \]
\[  (1,\;   (1_{(AB)C}*{\alpha_{\mu}}_{RST})({\alpha^{\cdot}_{\mu}}_{RST}*1_{(A'B')C'})(1_{R(ST)}*{{\epsilon^{\cdot}_{\alpha}}}_{A'B'C'}),\; a\cdot\kappa^{-1}\cdot\pi({\alpha^{\cdot}_{\mu}}_{RST}*1_{(A'B')C'})),\]
\end{itemize}

\item an identity perturbation
\[
   \xy
   (-12,0)*+{\epsilon_{\epsilon_{\alpha}}\maps \epsilon_{\alpha}\epsilon_{\alpha}^{\cdot}}="1"; (13,0)*+{1_{I_{\otimes(1\times\otimes)}}}="2";
        {\ar@3{->}_{} "1";"2"};
\endxy
\]
\noindent consisting of, for each triple of objects $A$, $B$, $C$, an equation of maps of spans
\[
   \xy
   (-19,0)*+{{\epsilon_{\epsilon_{\alpha}}}_{ABC}\maps {\epsilon_{\alpha}}_{ABC}{\epsilon_{\alpha}^{\cdot}}_{ABC}}="1"; (19,0)*+{{1_{I_{\otimes(1\times\otimes)}}}_{ABC},}="2";
        {\ar@3{->}_{} "1";"2"};
\endxy
\]

\item an identity perturbation
\[
   \xy
   (-11,0)*+{\eta_{\epsilon_{\alpha}}\maps 1_{\alpha\alpha^{\cdot}}}="1"; (12,0)*+{\epsilon^{\cdot}_{\alpha}\epsilon_{\alpha}}="2";
        {\ar@3{->}_{} "1";"2"};
\endxy
\]
\noindent consisting of, for each triple of objects $A$, $B$, $C$, an equation of maps of spans
\[
   \xy
   (-19,0)*+{{\eta_{\epsilon_{\alpha}}}_{ABC}\maps 1_{\alpha_{ABC}\alpha^{\cdot}_{ABC}}}="1"; (19,0)*+{{\epsilon^{\cdot}_{\alpha}}_{ABC}{\epsilon_{\alpha}}_{ABC},}="2";
        {\ar@3{->}_{} "1";"2"};
\endxy
\]
\end{itemize}

\item a strict adjoint equivalence $(\eta_{\alpha}, \eta^{\cdot}_{\alpha}, \epsilon_{\eta_{\alpha}}, \eta_{\eta_{\alpha}})$, in the  strict $2$-category
\[\Bicat(\Span(B)^{3}((A,B,C),(A',B',C')),\Span(\B)((A\times B)\times C,(A'\times B')\times C')),\]
\noindent consisting of
\begin{itemize}
\item a trimodification
\[
   \xy
   (-11,0)*+{\eta_{\alpha}\maps  \alpha^{\cdot}\alpha}="1"; (11,0)*+{I_{\otimes(\otimes\times 1)}}="2";
        {\ar@3{->}_{} "1";"2"};
\endxy
\]
\noindent consisting of
\begin{itemize}
\item for each triple of objects $A, B, C\in\Span(\B)$, a map of spans
\[
\def\objectstyle{\scriptstyle}
  \def\labelstyle{\scriptstyle}
   \xy
   (20,0)*+{(A\times B))\times C}="2";
   (0,20)*+{((A\times B)\times C)((A\times B)\times C)}="1";
   (0,-20)*+{(A\times B)\times C}="4";
   (-20,0)*+{(A\times B)\times C}="3";
        {\ar^{\pi_{(A\times B)\times C}^{(A\times B)\times C}} "1";"2"};
        {\ar_{\pi_{(A\times B)\times C}^{(A\times B)\times C}} "1";"3"};
        {\ar_{1} "4";"2"};
        {\ar^{1} "4";"3"};
        {\ar_{{\eta_{\alpha}}} "1";"4"};
        {\ar@{=}_<<{\scriptstyle } (12,2); (9,-1)};
        {\ar@{=}^<<{\scriptstyle } (-12,2); (-9,-1)};
\endxy
\]
\noindent where ${\eta_{\alpha}} := {\eta_{\alpha}}_{A,B,C} = \pi_{(A\times B)\times C}^{(A\times B)\times C}$,

\item for each pair of triples of objects $(A, B, C), (A', B', C')$, an identity modification
\[ m_{\eta_{\alpha}}\maps ({{\eta_{\alpha}}}_{(A,B,C)})^*\alpha^{\cdot}\alpha\To I_{\otimes(\otimes\times 1)}({{\eta_{\alpha}}}_{A',B',C'})_*\]
\noindent consisting of, for each triple of spans $R$, $S$, $T$, an equation of maps of spans
\[ (1,\;({{\eta_{\alpha}}}_{ABC}*1_{(RS)T})(1_{(AB)C}*{\alpha^{\cdot}_{\mu}}_{RST})({\alpha_{\mu}}_{RST}*1_{(A'B')C'}),\; a^{-1}\cdot\kappa\cdot\pi({\alpha_{\mu}}_{RST}*1_{(A'B')C'})) \]
\[ = (1,\;I_{(RS)T}(1_{(RS)T}*{{\eta_{\alpha}}}_{A'B'C'}),\; 1),\]

\end{itemize}

\item a trimodification
\[
   \xy
   (-11,0)*+{\eta^{\cdot}_{\alpha}\maps I_{\otimes(\otimes\times 1)}}="1"; (11,0)*+{\alpha^{\cdot}\alpha}="2";
        {\ar@3{->}_{} "1";"2"};
\endxy
\]
\noindent consisting of
\begin{itemize}
\item for each triple of objects $A, B, C\in\Span(\B)$, a map of spans
\[
\def\objectstyle{\scriptstyle}
  \def\labelstyle{\scriptstyle}
   \xy
   (20,0)*+{(A\times B)\times C}="2";
   (0,20)*+{(A\times B)\times C}="1";
   (0,-20)*+{((A\times B)\times C)((A\times B)\times C)}="4";
   (-20,0)*+{(A\times B)\times C}="3";
        {\ar_{\pi_{(A\times B)\times C}^{(A\times B)\times C}} "4";"2"};
        {\ar^{\pi_{(A\times B)\times C}^{(A\times B)\times C}} "4";"3"};
        {\ar^{1} "1";"2"};
        {\ar_{1} "1";"3"};
        {\ar_{{\eta^{\cdot}_{\alpha}}} "1";"4"};
        {\ar@{=}_<<{\scriptstyle } (12,2); (9,-1)};
        {\ar@{=}^<<{\scriptstyle } (-12,2); (-9,-1)};
\endxy
\]
\noindent where ${\eta^{\cdot}_{\alpha}} := {\eta^{\cdot}_{\alpha}}_{A,B,C}$ is the unique $1$-cell in $\B$ satisfying
\[ \pi^{(A\times B)\times C}_{(A\times B)\times C}{\eta^{\cdot}_{\alpha}}_{A,B,C} = 1,\]

\item for each pair of triples of objects $(A, B, C), (A', B', C')$, an identity modification
\[ m_{\eta^{\cdot}_{\alpha}}\maps ({m_{\eta^{\cdot}_{\alpha}}}_{(A,B,C)})^*I_{\otimes(\otimes\times 1)} \To \alpha^{\cdot}\alpha({m_{\eta^{\cdot}_{\alpha}}}_{A',B',C'})_*\]
\noindent consisting of, for each triple of spans $R$, $S$, $T$, an equation of maps of spans
\[ (1,\;({{\eta^{\cdot}_{\alpha}}}_{ABC}*1_{(RS)T})I_{(RS)T},\; 1)= \]
\[ (1,\; (1_{(AB)C}*{\alpha^{\cdot}_{\mu}}_{RST})({\alpha_{\mu}}_{RST}*1_{(A'B')C'})(1_{(RS)T}*{{\eta^{\cdot}_{\alpha}}}_{A'B'C'}),\; a^{-1}\cdot\kappa\cdot\pi({\alpha_{\mu}}_{RST}*1_{(A'B')C'})),\]
\end{itemize}

\item an identity perturbation
\[
   \xy
   (-11,0)*+{\epsilon_{\eta_{\alpha}}\maps \eta_{\alpha}\eta_{\alpha}^{\cdot}}="1"; (11,0)*+{1_{I_{\otimes(\otimes\times 1)}}}="2";
        {\ar@3{->}_{} "1";"2"};
\endxy
\]
\noindent consisting of, for each triple of objects $A$, $B$, $C$, an equation of maps of spans
\[
   \xy
   (-18,0)*+{{\epsilon_{\eta_{\alpha}}}_{ABC}\maps {\eta_{\alpha}}_{ABC}{\eta_{\alpha}^{\cdot}}_{ABC}}="1"; (18,0)*+{{1_{I_{\otimes(1\times\otimes)}}}_{ABC},}="2";
        {\ar@3{->}_{} "1";"2"};
\endxy
\]

\item an identity perturbation
\[
   \xy
   (-11,0)*+{\eta_{\eta_{\alpha}}\maps 1_{I_{\otimes(\otimes\times 1)}}}="1"; (11,0)*+{\eta_{\alpha}\eta^{\cdot}_{\alpha}}="2";
        {\ar@3{->}_{} "1";"2"};
\endxy
\]
\noindent consisting of, for each triple of objects $A$, $B$, $C$, an equation of maps of spans
\[
   \xy
   (-18,0)*+{{\eta_{\eta_{\alpha}}}_{ABC}\maps 1_{\alpha_{ABC}\alpha^{\cdot}_{ABC}}}="1"; (18,0)*+{{\eta^{\cdot}_{\alpha}}_{ABC}{\eta_{\alpha}}_{ABC},}="2";
        {\ar@3{->}_{} "1";"2"};
\endxy
\]

\end{itemize}

\item an identity perturbation
\[
   \xy
   (-33,0)*+{\Phi_{\alpha}\maps (1,\;{\bf l}(1*\epsilon_{\alpha}){\bf a}(\eta_{\alpha}^{\cdot}*1){\bf r}^{-1},\; \kappa^{-1}\cdot (1*\epsilon_{\alpha}){\bf a}(\eta_{\alpha}^{\cdot}*1){\bf r}^{-1})}="1"; (33,0)*+{(1,1,1)}="2";
        {\ar@3{->}_{} "1";"2"};
\endxy
\]
\noindent consisting of, for each triple of objects $A, B, C\in \Span(\B)$, an equation of maps of spans,

\item and, an identity perturbation
\[
   \xy
   (-33,0)*+{\Psi_{\alpha}\maps (\kappa\cdot (\epsilon_{\alpha}*1){\bf a}^{-1}(1*\eta_{\alpha}^{\cdot}){\bf l}^{-1},\; {\bf r}(\epsilon_{\alpha}*1){\bf a}^{-1}(1*\eta_{\alpha}^{\cdot}){\bf l}^{-1},\; 1)}="1"; (33,0)*+{(1,1,1)}="2";
        {\ar@3{->}_{} "1";"2"};
\endxy
\]
\noindent consisting of, for each triple of objects $A, B, C\in \Span(\B)$, an equation of maps of spans.
\end{itemize}
\endproposition

\proof
We check that $\alpha$ and $\alpha^{\cdot}$ are tritransformations, that $\epsilon_{\alpha}$ and $\eta_{\alpha}$ are adjoint equivalences, and that $\Phi_{\alpha}$ and $\Psi_{\alpha}$ are invertible perturbations.  Finally, we verify the axioms of a biadjoint biequivalence.

We first check that $\alpha$ is a tritransformation.  The proof for $\alpha^{\cdot}$ follows similarly.  We need to verify that $(\alpha_{\mu},\;\alpha_{\mu^{\cdot}},\;{\alpha}_{\epsilon},\;{\alpha}_{\eta})$ is an adjoint equivalence of strong transformations and modifications.

Naturality for ${\alpha_{\mu}}_{(R,S,T),(\bar{R},\bar{S},\bar{T})}$ is a straightforward calculation showing that, for each triple of maps of maps of spans $(\sigma_{R},\sigma_{S},\sigma_{T})$, the equation
\[ {\alpha_{\mu}}_{g_{R},g_{S},g_{T}}{1_{\alpha_{\mu}}}(\sigma_{(RS)T}*1) = (1*\sigma_{R(ST)}){1_{\alpha_{\mu}}}{{\alpha_{\mu}}_{f_{R},f_{S},f_{T}}}\]
\noindent holds.

Since the monoidal product is locally strict, the transformation axioms simplify.  One is a simple calculation and the other is immediate.  It follows that $\alpha_{\mu}$ is a strong transformation.  Similarly,  $\alpha_{\mu^{\cdot}}$ is a strong transformation.

The $2$-cell equation
\[ ((r\times(s\times t))\cdot (\pi\cdot{\alpha_{\epsilon}}_{RST}))(\kappa_{A\times (B\times C)}^{r\times (s\times t),1}\cdot \alpha_{\mu}\alpha_{\mu^{\cdot}}) = (\kappa_{A\times (B\times C)}^{r\times (s\times t),1}\cdot 1)(1\cdot (\pi\cdot{\alpha_{\epsilon}}_{RST}))\]
\noindent allows us to apply the universal property in defining the component isomorphisms of spans of the counit.  A similar equation gives the unit isomorphism.  The modification axiom for $\alpha_{\epsilon}$ is the following equation of $2$-cells
\[ 1_{f_{R}\times (f_{S}\times f_{T})}({\alpha_{\epsilon}}_{\bar{R},\bar{S},\bar{T}}1_{1*(f_{R}\times (f_{S}\times f_{T}))}) =\]
\[ (1_{1*(f_{R}\times (f_{S}\times f_{T}))}{\alpha_{\epsilon}}_{R,S,T})({\alpha_{\mu}}_{f_{R},f_{S},f_{T}}1_{{\alpha_{\mu^{\cdot}}}_{R,S,T}})(1_{{\alpha_{\mu}}_{R,S,T}}{\alpha_{\mu^{\cdot}}}_{f_{R},f_{S},f_{T}})\]
\noindent which is easily verified by definitions.  A similar modification axiom can be checked for $\alpha_{\eta}$.

The transformations and modifications form an adjoint equivalence in a strict $2$-category and the axioms reduce to the equations
\[ (1_{\alpha_{\mu^{\cdot}}}\alpha_{\epsilon})(\alpha_{\eta}1_{\alpha_{\mu^{\cdot}}}) = 1_{\alpha_{\mu^{\cdot}}}\]
\noindent and
\[ (\alpha_{\eta}1_{\alpha_{\mu}})(1_{\alpha_{\mu}}\alpha_{\epsilon}) = 1_{\alpha_{\mu}},\]
\noindent which are verified by simple calculations.

The modification axioms are immediate for the collections of identity cells $\alpha_{\Pi}$ and $\alpha_{M}$.

The tritransformation axioms are immediate since all modifications cells from $\Span(\B)$, the monoidal product, and $\alpha$ are identities.  It follows that $\alpha$ and similarly $\alpha^{\cdot}$ are tritransformations.

Next we check that $\epsilon_{\alpha}$ and $\epsilon^{\cdot}_{\alpha}$ are trimodifications and the $1$-cells of an adjoint equivalence.  Similar results will hold for $\eta_{\alpha}$ and $\eta^{\cdot}_{\alpha}$. 

The modification axiom is immediate since the $2$-cells ${m_{\epsilon_{\alpha}}}_{RST}$ are identities.  The trimodification axioms are also immediately satisfied since these $2$-cells and the modification components $\alpha_{\Pi}$, $\alpha_{M}$, and the analogous modifications for $\alpha^{\cdot}$ and $I_{\otimes(1\times\otimes)}$ are all identities.  It follows that $\epsilon_{\alpha}$, and similarly $\epsilon_{\alpha}^{\cdot}$ are trimodifications.  

The counit and unit perturbations $\epsilon_{\epsilon_{\alpha}}$ and $\eta_{\epsilon_{\alpha}}$ are each identities and thus trivially satisfy the perturbation axioms.  Further the adjoint equivalence axioms are immediately satisfied.  It follows that 
$(\epsilon_{\alpha}, \epsilon^{\cdot}_{\alpha}, \epsilon_{\epsilon_{\alpha}}, \eta_{\epsilon_{\alpha}})$, and similarly $(\eta_{\alpha}, \eta^{\cdot}_{\alpha}, \epsilon_{\eta_{\alpha}}, \eta_{\eta_{\alpha}})$, are adjoint equivalences.

Finally, the axioms of a biadjoint biequivalence will be satisfied since the perturbations $\Phi_{\alpha}$ and $\Psi_{\alpha}$ are identities and the modifications of $\Span(\B)$ are all identities.  To check these equations explicitly we need to specify the tricategory structure on $\Tricat(\Span^{3}(\B),\Span(\B))$, but we do not include these details here.
\endproof

\subsubsection*{Monoidal Left Unitor}

The monoidal left unitor is a biadjoint biequivalence consisting of tritransformations, adjoint equivalences of trimodifications and perturbations, and coherence perturbations consisting of isomorphisms of maps of spans.   

The left unitor for the product of an object $A\in\B$ and the unit $1\in\B$ is the $1$-cell
\[ \lambda_{A}\maps 1\times A\to A\]
\noindent in $\B$ defined as the projection
\[\widetilde{\pi}^{1}_{A}\maps 1\times A\to A\]
\noindent with $1$-cell
\[ \widetilde{\pi}_{1\times A}^{1}\maps A\to 1\times A\]
\noindent defined as the unique inverse $1$-cell in $\B$.

\proposition \label{monoidalleftunitor}
There is a biadjoint biequivalence
\[ (\lambda,\;\lambda^{\cdot},\; \epsilon_{\lambda},\;\eta_{\lambda}\;\Phi_{\lambda},\Psi_{\lambda})\maps \otimes(I\times 1)\To \otimes 1,\]
\noindent in a `tricategory' $\Tricat(\Span(\B),\Span(\B))$, consisting of
\begin{itemize}
\item a tritransformation
\[ \lambda\maps \otimes(I\times 1)\To \otimes 1,\]
\noindent consisting of
\begin{itemize}
\item for each object $A\in\Span(\B)$, a span $\lambda_{A}$
\[
   \xy
   (-20,0)*+{A}="1";
   (0,15)*+{1\times A}="2";
   (20,0)*+{1\times A}="3";
        {\ar_{\widetilde{\pi}_{A}^{1}} "2";"1"};
        {\ar^{1} "2";"3"};
\endxy
\]

\item for each pair of objects $A, B$ in $\Span(\B)$, an adjoint equivalence
\[ (\lambda_{\mu},\lambda_{\mu^{\cdot}},\lambda_{\epsilon}, \lambda_{\eta})\maps \Span(\B)(1, \lambda_{B})(\otimes(I\times 1))\To \Span(\B)(\lambda_{A}, 1)1,\]
\noindent in the strict $2$-category
\[\Bicat(\Span(\B)(A,B),\Span(\B)(1\times A,B))\]
\noindent consisting of

\begin{itemize}
\item a strong transformation
\[ {\lambda_{\mu}}_{A,B} \maps \Span(\B)(1, \lambda_{B})(\otimes(I\times 1))\To \Span(\B)(\lambda_{A}, 1)1,\]
\noindent consisting of
\begin{itemize}
\item for each span
\[
   \xy
   (-15,0)*+{B}="1";
   (0,15)*+{R}="2";
   (15,0)*+{A}="3";
        {\ar_{q} "2";"1"};
        {\ar^{p} "2";"3"};
\endxy
\]
\noindent a map of spans
\[
\def\objectstyle{\scriptstyle}
  \def\labelstyle{\scriptstyle}
   \xy
   (-20,0)*+{B}="2";
   (0,20)*+{(1\times B)(1 \times R)}="1";
   (0,-20)*+{R(1\times A)}="4";
   (20,0)*+{1\times A}="3";
        {\ar_{\widetilde{\pi}_{B}^{1}\pi_{1\times B}^{1\times R}} "1";"2"};
        {\ar^{(1\times p)\pi_{1\times R}^{1\times B}} "1";"3"};
        {\ar^{q\pi_{R}^{1\times A}} "4";"2"};
        {\ar_{\pi_{1\times A}^{R}} "4";"3"};
        {\ar^{\lambda_{\mu}} "1";"4"};
        {\ar@{=}_<<{\scriptstyle } (12,2); (9,-1)};
        {\ar@{=>}^<<{\scriptstyle \widetilde{\pi}_{B}^{1}\cdot\kappa^{-1}} (-12,2); (-9,-1)};
\endxy
\]
\noindent where $\widetilde{\pi}_{B}^{1}\cdot\kappa := \widetilde{\pi}_{B}^{1}\cdot\kappa_{1\times B}^{1,1\times q}$ and $\lambda_{\mu} := {\lambda_{\mu}}_{R}$ is the unique $1$-cell satisfying
\[ \pi_{1\times A}^{R}{\lambda_{\mu}}_{R} =  (1\times p)\pi_{1\times R}^{1\times B} \;\;\;  \pi^{1\times A}_{R}{\lambda_{\mu}}_{R} =   \widetilde{\pi}_{R}^{1}\pi_{1\times R}^{1\times B}\]
\noindent and
\[ \kappa^{p,\widetilde{\pi}_{A}^{1}}_{A}\cdot{\lambda_{\mu}}_{R} = 1,\]

\item for each pair of spans $R$, $\bar{R}$, a natural isomorphism
\[ {\lambda_{\mu}}_{R,\bar{R}}\maps ({\lambda_{\mu}}_{\bar{R}})_*\Span(\B)(1,\lambda_{B})(\otimes(I\times 1)) \To ({\lambda_{\mu}}_{R})^*\Span(\B)(\lambda_{A},1)1,\]
\noindent consisting of, for each map of spans $f_{R}$, an isomorphism of maps of spans
\[ {\lambda_{\mu}}_{f_{R}}\maps ((1\times\varpi_{R})\cdot\pi,\; {\lambda_{\mu}}_{\bar{R}}((1\times f_{R})*1),\; \widetilde{\pi}\cdot\kappa^{-1}\cdot ((1\times f_{R})*1))\]
\[ \To  (1,\; (1* f_{R}){\lambda_{\mu}}_{R},\; (\varrho_{R}\cdot\pi{\lambda_{\mu}}_{R})(\widetilde{\pi}\cdot\kappa^{-1})),\]
\noindent consisting of the unique $2$-cell
\[ {\lambda_{\mu}}_{f_{R}}\maps {\lambda_{\mu}}_{\bar{R}}((1\times f_{R})*1) \To (1*f_{R}){\lambda_{\mu}}_{R}\]
\noindent in $\B$ such that
\[ \pi_{1\times A}^{\bar{R}}\cdot {\lambda_{\mu}}_{f_{R}} = (1\times\varpi_{R})^{-1}\cdot\pi_{1\times R}^{1\times B}\;\;\;\textrm{ and }\;\;\; \pi^{1\times A}_{\bar{R}}\cdot {\lambda_{\mu}}_{f_{R}} = 1, \]
\end{itemize}

\item a strong transformation 
\[ {\lambda_{\mu^{\cdot}}}_{A,B} \maps \Span(\B)(\lambda_{A}, 1)1 \To \Span(\B)(1, \lambda_{B})(\otimes(I\times 1)),\]
\noindent consisting of
\begin{itemize}
\item for each span $R$, a map of spans
\[
\def\objectstyle{\scriptstyle}
  \def\labelstyle{\scriptstyle}
   \xy
   (-20,0)*+{B}="2";
   (0,-20)*+{(1\times B)(1\times R)}="1";
   (0,20)*+{R(1\times A)}="4";
   (20,0)*+{1\times A}="3";
        {\ar^{\widetilde{\pi}_{B}^{1}\pi_{1\times B}^{1\times R}} "1";"2"};
        {\ar_{(1\times p)\pi_{1\times R}^{1\times B}} "1";"3"};
        {\ar_{q\pi_{R}^{1\times A}} "4";"2"};
        {\ar^{\pi_{1\times A}^{R}} "4";"3"};
        {\ar^{\lambda_{\mu^{\cdot}}} "4";"1"};
        {\ar@{=>}_<<{\scriptstyle \widetilde{\pi}\cdot\kappa} (12,2); (9,-1)};
        {\ar@{=}_<<{\scriptstyle } (-12,2); (-9,-1)};
\endxy
\]
\noindent where $\widetilde{\pi}\cdot\kappa :=\widetilde{\pi}_{1\times A}^{A}\cdot{\kappa}_{A}^{p,\widetilde{\pi}_{A}^{1}}$ and $\lambda_{\mu^{\cdot}} := {\lambda_{\mu^{\cdot}}}_{R}$ is the unique $1$-cell satisfying
\[ \pi_{1\times R}^{1\times B}{\lambda_{\mu^{\cdot}}}_{R} = \widetilde{\pi}_{1\times R}^{R}\pi_{R}^{1\times A}\;\;\;  \pi^{1\times R}_{1\times B}{\lambda_{\mu^{\cdot}}}_{R} = (1\times q)\widetilde{\pi}_{1\times R}^{R}\pi_{R}^{1\times A}\]
\noindent and
\[ \kappa^{1,1\times q}_{1\times B}\cdot{\lambda_{\mu^{\cdot}}}_{R} = 1,\]

\item for each pair of spans $R$, $\bar{R}$, a natural isomorphism
\[ {\lambda_{\mu^{\cdot}}}_{R,\bar{R}}\maps ({\lambda_{\mu^{\cdot}}}_{\bar{R}})_* \Span(\B)(\lambda_{A},1)1 \To ({\lambda_{\mu^{\cdot}}}_{R})^*\Span(\B)(1,\lambda_{B})(\otimes(I\times 1)),\]
\noindent consisting of, for each map of spans $f_{R}$, an isomorphism of maps of spans
\[ {\lambda_{\mu^{\cdot}}}_{f_{R}}\maps (\widetilde{\pi}\cdot\kappa\cdot(1* f_{R}),\; {\lambda_{\mu^{\cdot}}}_{\bar{R}}(1*f_{R}),\; \varrho_{R}\cdot\pi) \To \]
\[ (((1\times\varpi_{R})\cdot\pi{\lambda_{\mu^{\cdot}}}_{R})(\widetilde{\pi}\cdot\kappa),\; ((1\times f_{R})*1){\lambda_{\mu^{\cdot}}}_{R},\; 1),\]
\noindent consisting of the unique $2$-cell
\[ {\lambda_{\mu^{\cdot}}}_{f_{R}}\maps {\lambda_{\mu^{\cdot}}}_{\bar{R}}(1*f_{R})\To ((1\times f_{R})*1){\lambda_{\mu^{\cdot}}}_{R}\]
\noindent in $\B$ such that
\[ \pi_{1\times \bar{R}}^{1\times B}\cdot {\lambda_{\mu^{\cdot}}}_{f_{R}} = 1\;\;\;\textrm{ and }\;\;\; \pi_{1\times B}^{1\times\bar{R}}\cdot {\lambda_{\mu^{\cdot}}}_{f_{R}} = \widetilde{\pi}_{1\times B}^{B}\cdot\varrho_{R}^{-1}\cdot\pi_{R}^{1\times A}, \]
\end{itemize}

\item an invertible counit modification
\[ {\lambda_{\epsilon}}\maps \lambda_{\mu}\lambda_{\mu^{\cdot}} \To 1\]
\noindent consisting of, for each span $R$, an isomorphism of maps of spans
\[ {\lambda_{\epsilon}}_{R}\maps (\widetilde{\pi}\cdot\kappa ,\; \lambda_{\mu}\lambda_{\mu^{\cdot}},\; \widetilde{\pi}\cdot{\kappa}^{-1}\cdot\lambda_{\mu^{\cdot}}) \To (1,1,1)\]
\noindent defined by the unique $2$-cell ${\lambda_{\epsilon}}_{R}$ in $\B$ such that
\[ \pi^{R}_{1\times A}\cdot{\lambda_{\epsilon}}_{R} = \widetilde{\pi}_{1\times A}^{A}\cdot{\kappa_{A}^{p,\widetilde{\pi}_{A}^{1}}}^{-1} \;\;\;\textrm{ and }\;\;\; \pi_{R}^{1\times A}\cdot{\lambda_{\epsilon}}_{R} = 1,\]

\item an invertible unit modification
\[ {\lambda_{\eta}}\maps 1\To \lambda_{\mu^{\cdot}}\lambda_{\mu}\]
\noindent consisting of, for each span $R$, an isomorphism of maps of spans
\[ {\lambda_{\eta}}_{R}\maps (1,1,1)\To (\kappa\cdot{\lambda_{\mu}},\;\lambda_{\mu^{\cdot}}\lambda_{\mu},\; \widetilde{\pi}\cdot\kappa^{-1})\]
\noindent defined by the unique $2$-cell ${\lambda_{\eta}}_{R}$ in $\B$ such that
\[ \pi_{1\times R}^{1\times B}\cdot {\lambda_{\eta}}_{R} = 1\;\;\;\textrm{ and }\;\;\;  \pi^{1\times R}_{1\times B}\cdot {\lambda_{\eta}}_{R} = {\kappa_{1\times B}^{1,1\times q}}^{-1},\]
\end{itemize}
\item an identity modification $\lambda_{\Pi}$ with component equations of maps of spans
\[ (1,\;\chi({\lambda_{\mu}}_{R}*1)(1*{\lambda_{\mu}}_{S}),\; \widetilde{\pi}\cdot {\kappa}^{-1}\cdot \pi) = (1,{\lambda_{\mu}}_{SR}(\chi*1), \widetilde{\pi}\cdot{\kappa}^{-1}\cdot (\chi*1)),\]
\item an identity modification $\lambda_{M}$ with component equations of maps of spans
\[ (1,\; {\lambda_{\mu}}_{A}\iota {\bf r}^{-1},\; \widetilde{\pi}\cdot {\kappa}^{-1}\cdot\iota {\bf r}^{-1}) =  (1,\; \iota {\bf l}^{-1},\; 1),\]
\end{itemize}

\item a tritransformation
\[ \lambda^{\cdot}\maps 1\To\otimes(I\times 1),\]
\noindent consisting of
\begin{itemize}
\item for each object $A\in\Span(\B)$, a span $\lambda^{\cdot}_{A}$
\[
   \xy
   (-20,0)*+{1\times A}="1";
   (0,15)*+{1\times A}="2";
   (20,0)*+{A}="3";
        {\ar_{1} "2";"1"};
        {\ar^{\widetilde{\pi}_{A}^{1}} "2";"3"};
\endxy
\]
\item for each pair of objects $A, B$ in $\Span(\B)$, an adjoint equivalence
\[ (\lambda^{\cdot}_{\mu},\lambda^{\cdot}_{\mu^{\cdot}},\lambda^{\cdot}_{\epsilon}, \lambda^{\cdot}_{\eta}) \maps \Span(\B)(1, \lambda^{\cdot}_{B})1\To \Span(\B)(\lambda^{\cdot}_{A}, 1)(\otimes(I\times 1)),\]
\noindent in the strict $2$-category
\[\Bicat(\Span(\B)(A,B),\Span(\B)(A,1\times B)),\] 
\noindent consisting of
\begin{itemize}
\item a strong transformation
\[ {\lambda^{\cdot}_{\mu}}_{A,B} \maps \Span(\B)(1, \lambda^{\cdot}_{B})1\To \Span(\B)(\lambda^{\cdot}_{A}, 1)(\otimes(I\times 1)),\]
\noindent consisting of
\begin{itemize}
\item for each span $R$, a map of spans
\[
\def\objectstyle{\scriptstyle}
  \def\labelstyle{\scriptstyle}
   \xy
   (-20,0)*+{1\times B}="2";
   (0,20)*+{(1\times B)R}="1";
   (0,-20)*+{(1\times R)(1\times A)}="4";
   (20,0)*+{A}="3";
        {\ar_{\pi_{1\times B}^{R}} "1";"2"};
        {\ar^{p\pi_{R}^{1\times B}} "1";"3"};
        {\ar^{(1\times q)\pi_{1\times R}^{1\times A}} "4";"2"};
        {\ar_{\widetilde{\pi}_{A}^{1}\pi_{1\times A}^{1\times R}} "4";"3"};
        {\ar^{\lambda^{\cdot}_{\mu}} "1";"4"};
        {\ar@{=}_<<{\scriptstyle } (12,2); (9,-1)};
        {\ar@{=>}^<<{\scriptstyle \widetilde{\pi}_{1\times B}^{B}\cdot\kappa^{-1}} (-14,2); (-11.5,-1)};
\endxy
\]
\noindent where $\widetilde{\pi}_{1\times B}^{B}\cdot\kappa := \widetilde{\pi}_{1\times B}^{B}\cdot{\kappa_{B}^{1,q}}$ and $\lambda^{\cdot}_{\mu} := {\lambda^{\cdot}_{\mu}}_{R}$ is the unique $1$-cell satisfying
\[ \pi_{1\times A}^{1\times R}{\lambda^{\cdot}_{\mu}}_{R} =  (1\times p)\widetilde{\pi}_{1\times R}^{R}\pi_{R}^{1\times B} \;\;\;  \pi^{1\times A}_{1\times R}{\lambda^{\cdot}_{\mu}}_{R} =   \widetilde{\pi}_{1\times R}^{R}\pi_{R}^{1\times B}\]
\noindent and
\[ \kappa^{1\times p,1}_{1\times A}\cdot{\lambda^{\cdot}_{\mu}}_{R} = 1,\]

\item for each pair of spans $R$, $\bar{R}$, a natural isomorphism
\[ {\lambda^{\cdot}_{\mu}}_{R,\bar{R}}\maps ({\lambda^{\cdot}_{\mu}}_{\bar{R}})_*\Span(\B)(1,\lambda^{\cdot}_{B})1 \To ({\lambda^{\cdot}_{\mu}}_{R})^*\Span(\B)(\lambda^{\cdot}_{A},1)(\otimes(I\times 1)),\]
\noindent consisting of, for each map of spans $f_{R}$, an isomorphism of maps of spans
\[ {\lambda^{\cdot}_{\mu}}_{f_{R}}\maps (\varpi_{R}\cdot\pi,\; {\lambda^{\cdot}_{\mu}}_{\bar{R}}(f_{R}*1),\; \widetilde{\pi}\cdot\kappa^{-1}\cdot (f_{R}*1)) \To\]
\[  (1,\; (1*(1\times f_{R})){\lambda^{\cdot}_{\mu}}_{R},\; ((1\times\varrho_{R})\cdot\pi{\lambda^{\cdot}_{\mu}}_{R})(\widetilde{\pi}\cdot\kappa^{-1})),\]
\noindent consisting of the unique $2$-cell
\[ {\lambda^{\cdot}_{\mu}}_{f_{R}}\maps{\lambda^{\cdot}_{\mu}}_{\bar{R}}(f_{R}*1) \To (1*(1\times f_{R})){\lambda^{\cdot}_{\mu}}_{R}\]
\noindent in $\B$ such that
\[ \pi_{1\times A}^{1\times \bar{R}}\cdot {\lambda^{\cdot}_{\mu}}_{f_{R}} = \widetilde{\pi}_{1\times A}^{A}\cdot\varpi_{R}\cdot\pi_{R}^{1\times B}\;\;\;\textrm{ and }\;\;\; \pi^{1\times A}_{1\times \bar{R}}\cdot {\lambda^{\cdot}_{\mu}}_{f_{R}} = 1, \]
\end{itemize}

\item a strong transformation 
\[ {\lambda^{\cdot}_{\mu^{\cdot}}}_{A,B} \maps \Span(\B)(\lambda^{\cdot}_{A}, 1)(\otimes(I\times 1)) \To \Span(\B)(1, \lambda^{\cdot}_{B})1,\]
\noindent consisting of
\begin{itemize}
\item for each span $R$, a map of spans
\[
\def\objectstyle{\scriptstyle}
  \def\labelstyle{\scriptstyle}
   \xy
   (-20,0)*+{1\times B}="2";
   (0,-20)*+{(1\times B)R}="1";
   (0,20)*+{(1\times R)(1\times A)}="4";
   (20,0)*+{A}="3";
        {\ar^{\pi_{1\times B}^{R}} "1";"2"};
        {\ar_{p\pi_{R}^{1\times B}} "1";"3"};
        {\ar_{(1\times q)\pi_{1\times R}^{1\times A}} "4";"2"};
        {\ar^{\widetilde{\pi}_{A}^{1}\pi_{1\times A}^{1\times R}} "4";"3"};
        {\ar^{\lambda^{\cdot}_{\mu^{\cdot}}} "4";"1"};
        {\ar@{=>}_<<{\scriptstyle \widetilde{\pi}_{A}^{1}\cdot\kappa} (12,2); (9,-1)};
        {\ar@{=}_<<{\scriptstyle } (-12,2); (-9,-1)};
\endxy
\]
\noindent where $\widetilde{\pi}_{A}^{1}\cdot\kappa := \widetilde{\pi}_{A}^{1}\cdot\kappa_{1\times A}^{p,1}$ and $\lambda^{\cdot}_{\mu^{\cdot}} := {\lambda^{\cdot}_{\mu^{\cdot}}}_{R}$ is the unique $1$-cell satisfying
\[ \pi_{R}^{1\times B}{\lambda^{\cdot}_{\mu^{\cdot}}}_{R} = \widetilde{\pi}_{R}^{1}\pi_{1\times R}^{1\times A}\;\;\;  \pi^{R}_{1\times B}{\lambda^{\cdot}_{\mu^{\cdot}}}_{R} = (1\times q)\pi_{1\times R}^{1\times A}\]
\noindent and
\[ \kappa^{\widetilde{\pi}_{B}^{1},q}_{B}\cdot{\lambda^{\cdot}_{\mu^{\cdot}}}_{R} = 1,\]

\item for each pair of spans $R$, $\bar{R}$, a natural isomorphism
\[ {\lambda^{\cdot}_{\mu^{\cdot}}}_{R,\bar{R}}\maps ({\lambda^{\cdot}_{\mu^{\cdot}}}_{\bar{R}})_* \Span(\lambda^{\cdot}_{A},1)(\otimes(I\times 1)) \To ({\lambda^{\cdot}_{\mu^{\cdot}}}_{R})^*\Span(1,\lambda_{B})1,\]
\noindent consisting of, for each map of spans $f_{R}$, an isomorphism of maps of spans
\[ {\lambda^{\cdot}_{\mu^{\cdot}}}_{f_{R}}\maps (\widetilde{\pi}\cdot\kappa\cdot(1* (1\times f_{R})),\; {\lambda^{\cdot}_{\mu^{\cdot}}}_{\bar{R}}(1*(1\times f_{R})),\; (1\times\varrho_{R})\cdot\pi{\lambda^{\cdot}_{\mu^{\cdot}}}_{\bar{R}})\] \[\To  ((\varpi_{R}\cdot\pi{\lambda^{\cdot}_{\mu^{\cdot}}}_{R})(\widetilde{\pi}\cdot\kappa),\; (f_{R}*1){\lambda^{\cdot}_{\mu^{\cdot}}}_{R},\; 1),\]
\noindent consisting of the unique $2$-cell
\[ {\lambda^{\cdot}_{\mu^{\cdot}}}_{f_{R}}\maps {\lambda^{\cdot}_{\mu^{\cdot}}}_{\bar{R}}(1*(1\times f_{R}))\To (f_{R}*1){\lambda^{\cdot}_{\mu^{\cdot}}}_{R}\]
\noindent in $\B$ such that
\[ \pi_{\bar{R}}^{1\times B}\cdot {\lambda^{\cdot}_{\mu^{\cdot}}}_{f_{R}} = 1\;\;\;\textrm{ and }\;\;\; \pi_{1\times B}^{\bar{R}}\cdot {\lambda^{\cdot}_{\mu^{\cdot}}}_{f_{R}} = (1\times \varrho_{R})^{-1}\cdot\pi_{1\times R}^{1\times A}, \]
\end{itemize}

\item an invertible counit modification
\[ {\lambda^{\cdot}_{\epsilon}}\maps \lambda^{\cdot}_{\mu}\lambda^{\cdot}_{\mu^{\cdot}} \To 1_{1}\]
\noindent consisting of, for each span $R$, an isomorphism of maps of spans
\[ {\lambda^{\cdot}_{\epsilon}}_{R}\maps (\widetilde{\pi}\cdot\kappa ,\; \lambda^{\cdot}_{\mu}\lambda^{\cdot}_{\mu^{\cdot}},\; \widetilde{\pi}\cdot{\kappa}^{-1}\cdot\lambda^{\cdot}_{\mu^{\cdot}}) \To (1,1_{1},1)\]
\noindent defined by the unique $2$-cell ${\lambda^{\cdot}_{\epsilon}}_{R}$ in $\B$ such that
\[ \pi^{1\times R}_{1\times A}\cdot{{\lambda}^{\cdot}_{\epsilon}}_{R} = {\kappa_{1\times A}^{1\times p,1}} \;\;\;\textrm{ and }\;\;\; \pi_{1\times R}^{1\times A}\cdot{{\lambda}^{\cdot}_{\epsilon}}_{R} = 1,\]

\item an invertible unit modification
\[ {\lambda^{\cdot}_{\eta}}\maps 1\To \lambda^{\cdot}_{\mu^{\cdot}}\lambda^{\cdot}_{\mu}\]
\noindent consisting of, for each span $R$, an isomorphism of maps of spans
\[ {\lambda^{\cdot}_{\eta}}_{R}\maps (1,1,1)\To (\widetilde{\pi}\cdot\kappa\cdot{\lambda^{\cdot}_{\mu}},\;\lambda^{\cdot}_{\mu^{\cdot}}\lambda^{\cdot}_{\mu},\; \widetilde{\pi}\cdot\kappa^{-1})\]
\noindent defined by the unique $2$-cell ${\lambda^{\cdot}_{\eta}}_{R}$ in $\B$ such that
\[ \pi_{R}^{1\times B}\cdot {\lambda^{\cdot}_{\eta}}_{R} = 1\;\;\;\textrm{ and }\;\;\;  \pi^{R}_{1\times B}\cdot {\lambda^{\cdot}_{\eta}}_{R} = \widetilde{\pi}_{1\times B}^{B}\cdot{\kappa_{B}^{\widetilde{\pi}_{B}^{1},q}}^{-1},\]
\end{itemize}
\item an identity modification $\lambda^{\cdot}_{\Pi}$ with component equations of maps of spans
\[ (1,\;\chi({\lambda^{\cdot}_{\mu}}_{R}*1)(1*{\lambda^{\cdot}_{\mu}}_{S}),\; \widetilde{\pi}\cdot {\kappa}^{-1}\cdot \pi) = (1,{\lambda^{\cdot}_{\mu}}_{SR}(\chi*1), \widetilde{\pi}\cdot{\kappa}^{-1}\cdot (\chi*1)),\]
\item an identity modification $\lambda^{\cdot}_{M}$ with component equations of maps of spans
\[ (1,\;{\lambda^{\cdot}_{\mu}}_{A}\iota {\bf r}^{-1},\; \widetilde{\pi}\cdot {\kappa}^{-1}\cdot\iota {\bf r}^{-1}) =  (1,\; \iota {\bf l}^{-1},\; 1),\]
\end{itemize}

\item a strict adjoint equivalence $(\epsilon_{\lambda},\;\epsilon_{\lambda}^{\cdot},\epsilon_{\epsilon_{\lambda}},\;\eta_{\epsilon_{\lambda}})$
\noindent consisting of
\begin{itemize}
\item a trimodification
\[ \epsilon_{\lambda}\maps \lambda\lambda^{\cdot} \Rightarrow I_{1},\]
\noindent consisting of
\begin{itemize}
\item for each object $A$ in $\Span(\B)$, a map of spans
\[
\def\objectstyle{\scriptstyle}
  \def\labelstyle{\scriptstyle}
   \xy
   (-20,0)*+{A}="2";
   (0,20)*+{(1\times A)(1\times A)}="1";
   (0,-20)*+{A}="4";
   (20,0)*+{A}="3";
        {\ar_{\widetilde{\pi}_{A}^{1}\pi_{1\times A}^{1\times A}} "1";"2"};
        {\ar^{\widetilde{\pi}_{A}^{1}\pi_{1\times A}^{1\times A}} "1";"3"};
        {\ar^{1} "4";"2"};
        {\ar_{1} "4";"3"};
        {\ar^{\epsilon_{\lambda}} "1";"4"};
        {\ar@{=}_<<{\scriptstyle } (12,2); (9,-1)};
        {\ar@{=}^<<{\scriptstyle } (-12,2); (-9,-1)};
\endxy
\]
\noindent where $\epsilon_{\lambda} := {\epsilon_{\lambda}}_{A} = \widetilde{\pi}_{A}^{1}\pi_{1\times A}^{1\times A}$,

\item and, for each pair $A, B \in\Span(\B)$, an identity modification
\[ m_{\epsilon_{\lambda}}\maps \left({\epsilon_{\lambda}}_{A}\right)^*\lambda\lambda^{\cdot} \Rightarrow I_{1}\left({\epsilon_{\lambda}}_{B}\right)_*,\]
\noindent consisting of, for each span, an equation of maps of spans,
\[ (1,\;(\epsilon_{\lambda} *1)(1*\lambda_{\mu})({\lambda_{\mu}^{\cdot}}*1),\; \widetilde{\pi}\cdot\kappa^{-1}\cdot \pi({\lambda_{\mu}^{\cdot}}*1)) = (1,\; I_{1}(1*\epsilon_{\lambda}),\; 1)\]
\end{itemize}

\item a trimodification
\[ \epsilon_{\lambda}^{\cdot}\maps I_{1} \Rightarrow \lambda_{\mu}\lambda_{\mu}^{\cdot},\]
\noindent consisting of
\begin{itemize}
\item for each object $A$ in $\Span(\B)$, a map of spans
\[
\def\objectstyle{\scriptstyle}
  \def\labelstyle{\scriptstyle}
   \xy
   (-20,0)*+{A}="2";
   (0,20)*+{A}="1";
   (0,-20)*+{(1\times A)(1\times A)}="4";
   (20,0)*+{A}="3";
        {\ar^{\widetilde{\pi}_{A}^{1}\pi_{1\times A}^{1\times A}} "4";"2"};
        {\ar_{\widetilde{\pi}_{A}^{1}\pi_{1\times A}^{1\times A}} "4";"3"};
        {\ar_{1} "1";"2"};
        {\ar^{1} "1";"3"};
        {\ar^{\epsilon_{\lambda}^{\cdot}} "1";"4"};
        {\ar@{=}_<<{\scriptstyle } (12,2); (9,-1)};
        {\ar@{=}^<<{\scriptstyle } (-12,2); (-9,-1)};
\endxy
\]
\noindent where $\epsilon_{\lambda}^{\cdot} := {\epsilon_{\lambda}^{\cdot}}_{A}$ is the unique $1$-cell in $\B$ such that
\[ \pi_{1\times A}^{1\times A}\cdot\mu_{\lambda}^{\cdot} = \widetilde{\pi}_{1\times A}^{A}\;\;\;\textrm{ and }\;\;\;\kappa_{1\times A}^{1,1}\cdot\mu_{\lambda}^{\cdot} = 1,\]
\item and, for each pair $A, B \in\Span(\B)$, an identity modification
\[ m_{\epsilon^{\cdot}_{\lambda}}\maps \left({\epsilon^{\cdot}_{\lambda}}_{A}\right)^*I_{1} \Rightarrow \left({\epsilon^{\cdot}_{\lambda}}_{B}\right)_*\lambda\lambda^{\cdot},\]
\noindent consisting of, for each span, an equation of maps of spans,
\[  (1,\;({\epsilon_{\lambda}^{\cdot}}_{A} *1)I_{1},\; 1) = (1,\; (1*\lambda_{\epsilon})({\lambda_{\epsilon}^{\cdot}}*1)(1*\epsilon_{\lambda}^{\cdot}),\; \widetilde{\pi}\cdot{\kappa}^{-1}\cdot \pi({\lambda_{\epsilon}^{\cdot}}*1)(1*{\epsilon_{\lambda}^{\cdot}}_{B})),\]

\end{itemize}

\item an identity counit perturbation
\[
   \xy
   (-9,0)*+{\epsilon_{\epsilon_{\lambda}}\maps \epsilon_{\lambda}\epsilon_{\lambda}^{\cdot}}="1"; (9,0)*+{1_{I_{1}}}="2";
        {\ar@3{->}_{} "1";"2"};
\endxy
\]
\noindent consisting of, for each object $A\in\Span(\B)$, an equation of maps of spans
\[
   \xy
   (-12,0)*+{{\epsilon_{\epsilon_{\lambda}}}_{A}\maps {\epsilon_{\lambda}}_{A}{\epsilon_{\lambda}^{\cdot}}_{A}}="1"; (10,0)*+{{1_{I_{1}}}_{A}}="2";
        {\ar@3{->}_{} "1";"2"};
\endxy
\]

\item and, an identity unit perturbation
\[
   \xy
   (-11,0)*+{\eta_{\epsilon_{\lambda}}\maps 1_{\lambda\lambda^{\cdot}}}="1"; (10,0)*+{\mu_{\lambda}^{\cdot}\mu_{\lambda}}="2";
        {\ar@3{->}_{} "1";"2"};
\endxy
\]
\noindent consisting of, for each object $A\in\Span(\B)$, an equation of maps of spans
\[
   \xy
   (-12,0)*+{{\eta_{\epsilon_{\lambda}}}_{A}\maps {1_{\lambda\lambda^{\cdot}}}_{A}}="1"; (12,0)*+{{\epsilon_{\lambda}^{\cdot}}_{A}{\epsilon_{\lambda}}_{A}}="2";
        {\ar@3{->}_{} "1";"2"};
\endxy
\]
\end{itemize}

\item a strict adjoint equivalence $(\eta_{\lambda},\;\eta_{\lambda}^{\cdot},\epsilon_{\eta_{\lambda}},\;\eta_{\eta_{\lambda}})$, consisting of
\begin{itemize}
\item a trimodification
\[
   \xy
   (-11,0)*+{\eta_{\lambda}\maps \lambda^{\cdot}\lambda}="1"; (10,0)*+{1_{\otimes(I\times 1)}}="2";
        {\ar@3{->}_{} "1";"2"};
\endxy
\]
\noindent consisting of
\begin{itemize}
\item for each object $A$ in $\Span(\B)$, a map of spans
\[
\def\objectstyle{\scriptstyle}
  \def\labelstyle{\scriptstyle}
   \xy
   (-20,0)*+{1\times A}="2";
   (0,20)*+{(1\times A)(1\times A)}="1";
   (0,-20)*+{1\times A}="4";
   (20,0)*+{1\times A}="3";
        {\ar_{\pi_{1\times A}^{1\times A}} "1";"2"};
        {\ar^{\pi_{1\times A}^{1\times A}} "1";"3"};
        {\ar^{1} "4";"2"};
        {\ar_{1} "4";"3"};
        {\ar^{\eta_{\lambda}} "1";"4"};
        {\ar@{=}_<<{\scriptstyle } (12,2); (9,-1)};
        {\ar@{=}^<<{\scriptstyle } (-12,2); (-9,-1)};
\endxy
\]
\noindent where $\eta_{\lambda} := {\eta_{\lambda}}_{A} = \pi_{1\times A}^{1\times A}$, 
\item and, for each pair of objects $A, B\in\Span(\B)$, an identity modification
\[ m_{\eta_{\lambda}}\maps ({\eta_{\lambda}}_{A})^*\lambda^{\cdot}\lambda \To ({\eta_{\lambda}}_{B})_*I_{\otimes(I\times 1)},\]
\noindent consisting of, for each span, an equation of maps of spans
\[  (1,\;(\eta_{\lambda} *1)(1*{\lambda_{\mu}^{\cdot}})(\lambda_{\mu}*1),\; \widetilde{\pi}\cdot\kappa^{-1}\cdot \pi({\lambda_{\mu}}*1)) = (1,\; I(1*\eta_{\lambda}),\; 1),\]

\end{itemize}

\item a trimodification
\[
   \xy
   (-11,0)*+{\eta_{\lambda}^{\cdot}\maps 1_{\otimes(I\times 1)}}="1"; (10,0)*+{\lambda^{\cdot}\lambda,}="2";
        {\ar@3{->}_{} "1";"2"};
\endxy
\]
\noindent consisting of
\begin{itemize}
\item for each object $A$ in $\Span(\B)$, a map of spans
\[
\def\objectstyle{\scriptstyle}
  \def\labelstyle{\scriptstyle}
   \xy
   (-20,0)*+{1\times A}="2";
   (0,20)*+{1\times A}="1";
   (0,-20)*+{(1\times A)(1\times A)}="4";
   (20,0)*+{1\times A}="3";
        {\ar^{\pi_{1\times A}^{1\times A}} "4";"2"};
        {\ar_{\pi_{1\times A}^{1\times A}} "4";"3"};
        {\ar_{1} "1";"2"};
        {\ar^{1} "1";"3"};
        {\ar^{\eta_{\lambda}^{\cdot}} "1";"4"};
        {\ar@{=}_<<{\scriptstyle } (12,2); (9,-1)};
        {\ar@{=}^<<{\scriptstyle } (-12,2); (-9,-1)};
\endxy
\]
\noindent where $\eta_{\lambda}^{\cdot} := {\eta_\lambda}^{\cdot}_{A}$ is the unique $1$-cell in $\B$ such that
\[ \pi_{1\times A}^{1\times A}\cdot\eta_{\lambda}^{\cdot} = 1\;\;\;\textrm{ and }\;\;\;\kappa_{1\times A}^{1,1}\cdot\eta_{\lambda}^{\cdot} = 1,\]

\item and, for each pair of objects $A, B\in\Span(\B)$, an identity modification
\[ m_{\eta^{\cdot}_{\lambda}}\maps ({\eta^{\cdot}_{\lambda}}_{A})_*I_{\otimes(I\times 1)}\To ({\eta^{\cdot}_{\lambda}}_{B})^* \lambda^{\cdot}\lambda,\]
\noindent consisting of, for each span, an equation of maps of spans
\[ (1,\;I(\eta_{\lambda}^{\cdot} *1),\; 1) = (1,\; (1*\eta_{\lambda}^{\cdot})(1*{\lambda_{\mu}}^{\cdot})({\lambda_{\mu}}*1),\; \widetilde{\pi}\cdot\kappa^{-1}\cdot(\lambda_{\mu}*1))\]
\noindent consisting of, for each span, an equation of maps of spans,

\end{itemize}

\item an identity counit perturbation
\[
   \xy
   (-9,0)*+{\epsilon_{\eta_{\lambda}}\maps \eta_{\lambda}\eta_{\lambda}^{\cdot}}="1"; (9,0)*+{1_{\lambda^{\cdot}\lambda}}="2";
        {\ar@3{->}_{} "1";"2"};
\endxy
\]
\noindent consisting of, for each object $A\in\Span(\B)$, an equation of maps of spans,

\item an identity unit perturbation
\[
   \xy
   (-12,0)*+{\eta_{\eta_{\lambda}}\maps 1_{\otimes(I\times 1)}}="1"; (10,0)*+{\eta_{\lambda}^{\cdot}\eta_{\lambda}}="2";
        {\ar@3{->}_{} "1";"2"};
\endxy
\]
\noindent consisting of, for each object $A\in\Span(\B)$, an equation of maps of spans,
\end{itemize}

\item an identity perturbation
\[
   \xy
   (-30,0)*+{\Phi_{\lambda}\maps (1,\;{\bf l}(1*\eta_{\lambda})(\eta_{\lambda}^{\cdot}*1){\bf r}^{-1},\; \kappa^{-1}\cdot (1*\eta_{\lambda})(\eta_{\lambda}^{\cdot}*1){\bf r}^{-1})}="1"; (31,0)*+{(1,1,1)}="2";
        {\ar@3{->}_{} "1";"2"};
\endxy
\]
\noindent consisting of, for each object $A\in\Span(\B)$, an equation of maps of spans,

\item and, for each triple of objects $A, B, C\in \Span(\B)$, an identity isomorphism of maps of spans
\[
   \xy
   (-30,0)*+{\Psi_{\lambda}\maps (\kappa\cdot (\eta_{\lambda}*1)(1*\eta_{\lambda}^{\cdot}){\bf l}^{-1},\; {\bf r}(\eta_{\lambda}*1)(1*\eta_{\lambda}^{\cdot}){\bf l}^{-1},\; 1)}="1"; (31,0)*+{(1,1,1).}="2";
        {\ar@3{->}_{} "1";"2"};
\endxy
\]
\end{itemize}

\endproposition

\proof
We check that $\lambda$ and $\lambda^{\cdot}$ are tritransformations, that $\epsilon_{\lambda}$ and $\eta_{\lambda}$ are adjoint equivalences, and that $\Phi_{\lambda}$ and $\Psi_{\lambda}$ are invertible perturbations.  Finally, we verify the axioms of a biadjoint biequivalence.

We first check that $\lambda$ is a tritransformation.  The proof for $\lambda^{\cdot}$ follows similarly.  We need to verify that $(\lambda_{\mu},\;\lambda_{\mu^{\cdot}},\;{\lambda}_{\epsilon},\;{\lambda}_{\eta})$ is an adjoint equivalence of strong transformations and modifications.

Naturality for ${\lambda_{\mu}}_{R,\bar{R}}$ is a straightforward calculation showing that, for each map of maps of spans $\sigma_{R}$, the equation
\[ {\lambda_{\mu}}_{g_{R}}((1\times \sigma_{R})*1_{{\lambda_{\mu}}_{\bar{R}}}) = (1_{{\lambda_{\mu}}_{R}}*(1*\sigma_{R})){{\lambda_{\mu}}_{f_{R}}}\]
\noindent holds.

Since the monoidal product is locally strict, the transformation axioms simplify.  One is a simple calculation and the other is immediate.  It follows that $\lambda_{\mu}$ is a strong transformation.  Similarly,  $\lambda_{\mu^{\cdot}}$ is a strong transformation.

The $2$-cell equation
\[ (p\cdot (\pi\cdot{\lambda_{\epsilon}}_{R}))(\kappa_{A}^{p,\widetilde{\pi}^{1}_{A}}\cdot \lambda_{\mu}\lambda_{\mu^{\cdot}}) = (\kappa_{A}^{p,\widetilde{\pi}_{A}^{1}}\cdot 1)(\widetilde{\pi}_{A}^{1}\cdot (\pi\cdot{\lambda_{\epsilon}}_{R}))\]
\noindent allows us to apply the universal property in defining the component isomorphisms of spans of the counit.  A similar equation gives the unit isomorphism.  The modification axiom for $\lambda_{\epsilon}$ is the following equation of $2$-cells
\[ 1_{f_{R}}({\lambda_{\epsilon}}_{\bar{R}}1_{1*f_{R}}) = (1_{1*f_{R}}{\lambda_{\epsilon}}_{R})({\lambda_{\mu}}_{f_{R}}1_{{\lambda_{\mu^{\cdot}}}_{R}})(1_{{\lambda_{\mu}}_{R}}{\lambda_{\mu^{\cdot}}}_{f_{R}})\]
\noindent which is easily verified by definitions.  A similar modification axiom can be checked for $\lambda_{\eta}$.

The transformations and modifications form an adjoint equivalence in a strict $2$-category and the axioms reduce to the equations
\[ (1_{\lambda_{\mu^{\cdot}}}\lambda_{\epsilon})(\lambda_{\eta}1_{\lambda_{\mu^{\cdot}}}) = 1_{\lambda_{\mu^{\cdot}}}\]
\noindent and
\[ (\lambda_{\eta}1_{\lambda_{\mu}})(1_{\lambda_{\mu}}\lambda_{\epsilon}) = 1_{\lambda_{\mu}},\]
\noindent which are verified by simple calculations.

The modification axioms are immediate for the collections of identity cells $\lambda_{\Pi}$ and $\lambda_{M}$.

The tritransformation axioms are immediate since all modifications cells from $\Span(\B)$, the monoidal product, and $\lambda$ are identities.  It follows that $\lambda$ and similarly $\lambda^{\cdot}$ are tritransformations.

Next we check that $\epsilon_{\lambda}$ and $\epsilon^{\cdot}_{\lambda}$ are trimodifications and the $1$-cells of an adjoint equivalence.  Similar results will hold for $\eta_{\lambda}$ and $\eta^{\cdot}_{\lambda}$. 

The modification axiom is immediate since the $2$-cells ${m_{\epsilon_{\lambda}}}_{R}$ are identities.  The trimodification axioms are also immediately satisfied since these $2$-cells and the modification components $\lambda_{\Pi}$, $\lambda_{M}$, and the analogous modifications for $\lambda^{\cdot}$ and $I_{1}$ are all identities.  It follows that $\epsilon_{\lambda}$, and similarly $\epsilon_{\lambda}^{\cdot}$ are trimodifications.  

The counit and unit perturbations $\epsilon_{\epsilon_{\lambda}}$ and $\eta_{\epsilon_{\lambda}}$ are each identities and thus trivially satisfy the perturbation axioms.  Further the adjoint equivalence axioms are immediately satisfied.  It follows that 
$(\epsilon_{\lambda}, \epsilon^{\cdot}_{\lambda}, \epsilon_{\epsilon_{\lambda}}, \eta_{\epsilon_{\lambda}})$, and similarly $(\eta_{\lambda}, \eta^{\cdot}_{\lambda}, \epsilon_{\eta_{\lambda}}, \eta_{\eta_{\lambda}})$, are adjoint equivalences.

Finally, the axioms of a biadjoint biequivalence will be satisfied since the perturbations $\Phi_{\lambda}$ and $\Psi_{\lambda}$ are identities and the modifications of $\Span(\B)$ are all identities.  To check these equations explicitly we need to specify the tricategory structure on $\Tricat(\Span(\B),\Span(\B))$, but we do not include these details here.
\endproof

\subsubsection*{Monoidal Right Unitor}

The monoidal right unitor is a biadjoint biequivalence
\[ \rho\maps \otimes (1\times I)\Rightarrow 1\]
\noindent consisting of transformations, adjoint equivalences of modifications and perturbations, and coherence perturbations consisting of isomorphisms of maps of spans.   

\proposition \label{monoidalrightunitor}
There is a biadjoint biequivalence $(\rho,\;\rho^{\cdot},\; \epsilon_{\rho},\;\eta_{\rho},\;\Phi_{\rho},\Psi_{\rho})$ consisting of
\begin{itemize}
\item a tritransformation
\[ \rho\maps \otimes (1\times I)\Rightarrow 1,\]
\noindent consisting of
\begin{itemize}
\item for each object $A\in\Span(\B)$, the span $\rho_{A}$
\[
   \xy
   (-20,0)*+{A}="1";
   (0,15)*+{A\times 1}="2";
   (20,0)*+{A\times 1}="3";
        {\ar_{\widetilde{\pi}_{A}^{1}} "2";"1"};
        {\ar^{1} "2";"3"};
\endxy
\]
\item for each pair of objects $A, B$ in $\Span(\B)$, an adjoint equivalence
\[ (\rho_{\mu},\rho_{\mu^{\cdot}},\rho_{\epsilon},\rho_{\eta})\maps \hom_{\Span}(1,\rho_{B})(\otimes(1\times I))\To \hom_{\Span}(\rho_{A},1)1,\]
\noindent in a strict $2$-category of transformations and modifications, consisting of

\begin{itemize}
\item a strong transformation
\[ \rho_{\mu_{A,B}}\maps \hom_{\Span}(1,\rho_{B})(\otimes(1\times I))\To \hom_{\Span}(\rho_{A},1)1,\]
\noindent consisting of
\begin{itemize}
\item for each span $R$, a map of spans
\[
\def\objectstyle{\scriptstyle}
  \def\labelstyle{\scriptstyle}
   \xy
   (-20,0)*+{B}="2";
   (0,20)*+{(B\times 1)(R \times 1)}="1";
   (0,-20)*+{R(A\times 1)}="4";
   (20,0)*+{A\times 1}="3";
        {\ar_{\widetilde{\pi}_{B}^{1}\pi_{B\times 1}^{R\times 1}} "1";"2"};
        {\ar^{(1\times p)\pi_{R\times 1}^{B\times 1}} "1";"3"};
        {\ar^{q\pi_{R}^{A\times 1}} "4";"2"};
        {\ar_{\pi_{A\times 1}^{R}} "4";"3"};
        {\ar^{\rho_{\mu}} "1";"4"};
        {\ar@{=}_<<{\scriptstyle } (12,2); (9,-1)};
        {\ar@{=>}^<<{\scriptstyle \widetilde{\pi}_{B}^{1}\cdot\kappa^{-1}} (-12,2); (-9,-1)};
\endxy
\]
\noindent where $\kappa := \kappa_{B\times 1}^{(B\times 1)(R\times 1)}$ and $\rho_{\mu}:= {\rho_{\mu}}_{R}$ is the unique $1$-cell satisfying
\[ \pi_{A\times 1}^{R}{\rho_{\mu}}_{R} =  (p\times 1)\pi_{R\times 1}^{B\times 1} \;\;\;  \pi^{A\times 1}_{R}{\rho_{\mu}}_{R} =   \widetilde{\pi}_{R}^{1}\pi_{R\times 1}^{B\times 1},\]
\noindent and
\[ \kappa^{R(A\times 1)}_{A}\cdot{\rho_{\mu}}_{R} = 1,\]

\item for each pair of spans $R$, $\bar{R}$, a natural isomorphism
\[ {\rho_{\mu}}_{R,\bar{R}}\maps \left({\rho_{\mu}}_{\bar{R}}\right)_*\hom_{\Span}(1,\rho_{B})(\otimes(1\times I)) \To \left({\rho_{\mu}}_{R}\right)^*\hom_{\Span}(\rho_{A},1)1,\]
\noindent consisting of, for each map of spans $f_{R}$, an isomorphism of maps of spans
\[ {\rho_{\mu}}_{f_{R}}\maps ((\alpha_{R}\times 1)\cdot\pi,\; {\rho_{\mu}}_{\bar{R}}((f_{R}\times 1)*1),\;\widetilde{\pi}\cdot\kappa^{-1}\cdot((f_{R}\times 1)*1))\]
\[ \To (1,\; (1*f_{R}){\rho_{\mu}}_{R},\; (\beta_{R}\cdot\pi{\rho_{\mu}}_{R})(\widetilde{\pi}\cdot\kappa^{-1})),\]
\noindent consisting of the unique $2$-cell
\[ {\rho_{\mu}}_{f_{R}}\maps  {\rho_{\mu}}_{\bar{R}}((f_{R}\times 1)*1) \To (1*f_{R}){\rho_{\mu}}_{R}\]
\noindent in $\B$ such that
\[ \pi_{A\times 1}^{\bar{R}}\cdot{\rho_{\mu}}_{f_{R}} = (\alpha_{R}^{-1}\times 1)\cdot \pi_{R\times 1}^{B\times 1}\;\;\;\textrm{ and }\;\;\; \pi^{A\times 1}_{\bar{R}}\cdot{\rho_{\mu}}_{f_{R}} = 1,\]

\end{itemize}

\item a strong transformation
\[ {\rho_{\mu^{\cdot}}}_{A,B}\maps \hom_{\Span}(\rho_{A},1)1 \To \hom_{\Span}(1,\rho_{B})(\otimes(1\times I)),\]
\noindent consisting of

\begin{itemize}

\item for each span $R$, a map of spans
\[
\def\objectstyle{\scriptstyle}
  \def\labelstyle{\scriptstyle}
   \xy
   (-20,0)*+{B}="2";
   (0,-20)*+{(B\times 1)(R\times 1)}="1";
   (0,20)*+{R(A\times 1)}="4";
   (20,0)*+{A\times 1}="3";
        {\ar^{\widetilde{\pi}_{B}^{1}\pi_{B\times 1}^{R\times 1}} "1";"2"};
        {\ar_{(p\times 1)\pi_{R\times 1}^{B\times 1}} "1";"3"};
        {\ar_{q\pi_{R}^{A\times 1}} "4";"2"};
        {\ar^{\pi_{A\times 1}^{R}} "4";"3"};
        {\ar^{\rho_{\mu^{\cdot}}} "4";"1"};
        {\ar@{=>}_<<{\scriptstyle \widetilde{\pi}\cdot\kappa} (12,2); (9,-1)};
        {\ar@{=}_<<{\scriptstyle } (-12,2); (-9,-1)};
\endxy
\]
\noindent where $\widetilde{\pi}\cdot\kappa := \widetilde{\pi}_{A\times 1}^{A}\cdot{\kappa_{A}^{R(A\times 1)}}$ and ${\rho_{\mu^{\cdot}}}:= {\rho_{\mu^{\cdot}}}_{R}$ is the unique $1$-cell satisfying
\[ \pi_{R\times 1}^{B\times 1}{\rho_{\mu^{\cdot}}}_{R} = \widetilde{\pi}_{R\times 1}^{R}\pi_{R}^{A\times 1}\;\;\;  \pi^{R\times 1}_{B\times 1}{\rho_{\mu^{\cdot}}}_{R} = (q\times 1)\widetilde{\pi}_{R\times 1}^{R}\pi_{R}^{A\times 1}\]
\noindent and
\[ \kappa^{(B\times 1)(R\times 1)}_{B\times 1}\cdot{\rho_{\mu^{\cdot}}}_{R} = 1,\]

\item for each pair of spans $R$, $\bar{R}$, a natural isomorphism
\[ {\rho_{\mu^{\cdot}}}_{R,\bar{R}}\maps \left({\rho_{\mu^{\cdot}}}_{\bar{R}}\right)_*\hom_{\Span}(\rho_{A},1)1 \To  \left({\rho_{\mu^{\cdot}}}_{R}\right)^*\hom_{\Span}(1,\rho_{B})(\otimes(1\times I)),\]
\noindent consisting of, for each map of spans $f_{R}$, an isomorphism of maps of spans
\[ {\rho_{\mu^{\cdot}}}_{f_{R}}\maps (\widetilde{\pi}\cdot\kappa\cdot(1*f_{R}),\; {\rho_{\mu^{\cdot}}}_{\bar{R}}(1*f_{R}),\; \beta_{R}\cdot\pi)\]
\[ \To (((\alpha_{R}\times 1)\cdot\pi{\rho_{\mu^{\cdot}}}_{R})(\widetilde{\pi}\cdot\kappa),\; ((f_{R}\times 1)*1){\rho_{\mu^{\cdot}}}_{R},\; 1)\]
\noindent consisting of the unique $2$-cell
\[ {\rho_{\mu^{\cdot}}}_{f_{R}}\maps {\rho_{\mu^{\cdot}}}_{\bar{R}}(1*f_{R}) \To ((f_{R}\times 1)*1){\rho_{\mu^{\cdot}}}_{R}\]
\noindent in $\B$ such that
\[ \pi_{\bar{R}\times 1}^{B\times 1}\cdot {\rho_{\mu^{\cdot}}}_{f_{R}} = 1\;\;\;\textrm{ and }\;\;\; \pi^{\bar{R}\times 1}_{B\times 1}\cdot {\rho_{\mu^{\cdot}}}_{f_{R}} = \widetilde{\pi}_{B\times 1}^{B}\beta_{R}^{-1}\cdot\pi_{R}^{A\times 1}, \]
\end{itemize}

\item an invertible counit modification
\[ \rho_{\epsilon}\maps \rho_{\mu}\rho_{\mu^{\cdot}} \To 1\]
\noindent consisting of, for each span $R$, an isomorphism of maps of spans
\[ {\rho_{\epsilon}}_{R}\maps (\widetilde{\pi}\cdot\kappa,\; \rho_{\mu}\rho_{\mu^{\cdot}},\; \widetilde{\pi}_{B}^{1}\cdot\kappa^{-1}\cdot\rho_{\mu^{\cdot}}) \To (1,\; 1,\; 1)\]
\noindent defined by the unique $2$-cell ${\rho_{\epsilon}}_{R}$ in $\B$ such that
\[ \pi_{A\times 1}^{R}\cdot {\rho_{\epsilon}}_{R} = \widetilde{\pi}_{A\times 1}^{A}\cdot {\kappa_{A}^{p,\widetilde{\pi}_{A}^{1}}}^{-1}\;\;\;\textrm{ and }\;\;\;  \pi^{A\times 1}_{R}\cdot {\rho_{\epsilon}}_{R} = 1,\]

\item an invertible unit modification
\[ \rho_{\eta}\maps 1\To \rho_{\mu^{\cdot}}\rho_{\mu}\]
\noindent consisting of, for each span $R$, an isomorphism of maps of spans
\[ {\rho_{\eta}}_{R}\maps (1,1,1) \To (\widetilde{\pi}\cdot\kappa\cdot\rho_{\mu},\;\rho_{\mu^{\cdot}}\rho_{\mu},\; \widetilde{\pi}\cdot\kappa^{-1})\]
\noindent defined by the unique $2$-cell ${\rho_{\eta}}_{R}$ in $\B$ such that
\[ \pi_{R\times 1}^{B\times 1}\cdot {\rho_{\eta}}_{R} = 1 \;\;\;\textrm{ and }\;\;\; \pi_{B\times 1}^{R\times 1}\cdot {\rho_{\eta}}_{R} =  {\kappa_{B\times 1}^{1,q\times 1}}^{-1},\]

\end{itemize}

\item an identity modification with component equations of maps of spans
\[ (1,\;\chi({\rho_{\mu}}_{R}*1)(1*{\rho_{\mu}}_{S}),\; \widetilde{\pi}\cdot {\kappa}^{-1}\cdot \pi) = (1,\;{\rho_{\mu}}_{SR}(\chi*1), \widetilde{\pi}\cdot {\kappa}^{-1}\cdot (\chi*1)),\]
\item an identity modification with component equations of maps of spans
\[ (1,\; {\rho_{\mu}}_{A}\iota {\bf r}^{-1},\; \widetilde{\pi}\cdot {\kappa}^{-1}\cdot\iota {\bf r}^{-1}) =  (1,\; \iota {\bf l}^{-1},\; 1),\]
\end{itemize}

\item a tritransformation
\[ \rho^{\cdot}\maps 1\Rightarrow \otimes (1\times I),\]
\noindent consisting of
\begin{itemize}
\item for each object $A\in\Span(\B)$, the span $\rho^{\cdot}_{A}$
\[
   \xy
   (-20,0)*+{A\times 1}="1";
   (0,15)*+{A\times 1}="2";
   (20,0)*+{A}="3";
        {\ar_{1} "2";"1"};
        {\ar^{\widetilde{\pi}_{A}^{1}} "2";"3"};
\endxy
\]
%\[
%\xymatrix{
%&A\times 1\ar[dl]_{1}\ar[dr]^{\widetilde{\pi}_{A}^{1}}&\\
%A\times 1 && A
%}
%\]

\item for each pair of objects $A, B$ in $\Span(\B)$, an adjoint equivalence
\[ (\rho^{\cdot}_{\mu},\;\rho_{\mu^{\cdot}}^{\cdot},\;\rho^{\cdot}_{\epsilon},\;\rho^{\cdot}_{\eta})\maps \hom_{\Span}(1,\rho^{\cdot}_{B})1\To \hom_{\Span}(\rho_{A}^{\cdot}, 1)(\otimes(1\times I)),\]
\noindent in a strict $2$-category of transformations and modifications, consisting of
\begin{itemize}
\item a strong transformation
\[ {\rho^{\cdot}_{\mu}}_{A,B}\maps \hom_{\Span}(1, \rho^{\cdot}_{B})1\To \hom_{\Span}(\rho_{A}^{\cdot}, 1)(\otimes(1\times I)),\]
\noindent consisting of

\begin{itemize}
\item for each span $R$, a map of spans
\[
\def\objectstyle{\scriptstyle}
  \def\labelstyle{\scriptstyle}
   \xy
   (-20,0)*+{B\times 1}="2";
   (0,20)*+{(B\times 1)R}="1";
   (0,-20)*+{(R\times 1)(A\times 1)}="4";
   (20,0)*+{A}="3";
        {\ar_{\pi_{B\times 1}^{R}} "1";"2"};
        {\ar^{p\pi_{R}^{B\times 1}} "1";"3"};
        {\ar^{(q\times 1)\pi_{R\times 1}^{A\times 1}} "4";"2"};
        {\ar_{\widetilde{\pi}_{A}^{1}\pi_{A\times 1}^{R\times 1}} "4";"3"};
        {\ar^{{\rho_{\mu}^{\cdot}}} "1";"4"};
        {\ar@{=}_<<{\scriptstyle } (12,2); (9,-1)};
        {\ar@{=>}^<<{\scriptstyle \widetilde{\pi}\cdot\kappa^{-1}} (-12,2); (-9,-1)};
\endxy
\]
\noindent where $\widetilde{\pi}\cdot\kappa^{-1} := \widetilde{\pi}_{B\times 1}^{B}\cdot{\kappa_{B}^{\widetilde{\pi}_{B\times 1}^{B},q}}^{-1}$ and ${\rho_{\mu}^{\cdot}}:= {\rho_{\mu}^{\cdot}}_{R}$ is the unique $1$-cell satisfying
\[ \pi_{A\times 1}^{R\times 1}{\rho^{\cdot}_{\mu}}_{R} = (p\times 1)\widetilde{\pi}_{R\times 1}^{R}\pi_{R}^{B\times 1}\;\;\; \pi^{A\times 1}_{R\times 1}{\rho^{\cdot}_{\mu}}_{R} = \widetilde{\pi}_{R\times 1}^{R}\pi_{R}^{B\times 1}\]
\noindent and
\[ \kappa^{p\times 1,1}_{A\times 1}\cdot{\rho^{\cdot}_{\mu}}_{R} = 1,\]

\item for each pair of spans $R$, $\bar{R}$, a natural isomorphism
\[ {\rho^{\cdot}_{\mu}}_{R,\bar{R}}\maps \left({\rho^{\cdot}_{\mu}}_{\bar{R}}\right)_*\hom_{\Span}(1,\rho^{\cdot}_{B})1 \To \left({\rho^{\cdot}_{\mu}}_{R}\right)^*\hom_{\Span}(\rho^{\cdot}_{A},1)(\otimes(1\times I)),\]
\noindent consisting of, for each map of spans $f_{R}$, an isomorphism of maps of spans
\[ {\rho^{\cdot}_{\mu}}_{f_{R}}\maps (\alpha_{R}\cdot\pi,\; {\rho^{\cdot}_{\mu}}_{\bar{R}}(f_{R}* 1),\; \widetilde{\pi}\cdot\kappa^{-1}\cdot(f_{R}*1))\]
\[ \To  (1,\; (1*(f_{R}\times 1)){\rho^{\cdot}_{\mu}}_{R},\; ((\beta_{R}\times 1)\cdot\pi{\rho^{\cdot}_{\mu}}_{R})(\widetilde{\pi}\cdot\kappa^{-1})),\]
\noindent consisting of the unique $2$-cell
\[ {\rho^{\cdot}_{\mu}}_{f_{R}}\maps {\rho^{\cdot}_{\mu}}_{\bar{R}}(f_{R}* 1) \To (1*(f_{R}\times 1)){\rho^{\cdot}_{\mu}}_{R}\]
\noindent in $\B$ such that
\[ \pi_{\bar{R}\times 1}^{A\times 1}\cdot{\rho^{\cdot}_{\mu}}_{f_{R}} = 1 \;\;\;\textrm{ and }\;\;\;  \pi^{\bar{R}\times 1}_{A\times 1}\cdot{\rho^{\cdot}_{\mu}}_{f_{R}} = (\alpha_{R}^{-1}\times 1)\cdot\widetilde{\pi}_{R\times 1}^{R}\pi_{R}^{B\times 1}, \]
\end{itemize}

\item a strong transformation
\[ {\rho^{\cdot}_{\mu^{\cdot}}}_{A,B}\maps \hom_{\Span}(\rho^{\cdot}_{A},1)(\otimes(1\times I))\To \hom_{\Span}(1,\rho_{B}^{\cdot})1,\]
\noindent consisting of

\begin{itemize}
\item for each span $R$, a map of spans
\[
\def\objectstyle{\scriptstyle}
  \def\labelstyle{\scriptstyle}
   \xy
   (-20,0)*+{B\times 1}="2";
   (0,-20)*+{(B\times 1)R}="1";
   (0,20)*+{(R\times 1)(A\times 1)}="4";
   (20,0)*+{A}="3";
        {\ar^{\pi^{R}_{B\times 1}} "1";"2"};
        {\ar_{p\pi_{R}^{B\times 1}} "1";"3"};
        {\ar_{(q\times 1)\pi_{R\times 1}^{A\times 1}} "4";"2"};
        {\ar^{\widetilde{\pi}_{A}^{1}\pi^{R\times 1}_{A\times 1}} "4";"3"};
        {\ar_{{\rho_{\mu^{\cdot}}^{\cdot}}} "4";"1"};
        {\ar@{=>}_<<{\scriptstyle \widetilde{\pi}_{A}^{1}\cdot\kappa} (12,2); (9,-1)};
        {\ar@{=}_<<{\scriptstyle } (-12,2); (-9,-1)};
\endxy
\]
\noindent where $\kappa := \kappa_{A\times 1}^{p\times 1,1}$ and ${\rho^{\cdot}_{\mu^{\cdot}}} := {\rho^{\cdot}_{\mu^{\cdot}}}_{R}$ is the unique $1$-cell satisfying
\[ \pi_{R}^{B\times 1}{\rho^{\cdot}_{\mu^{\cdot}}}_{R} = \widetilde{\pi}_{R}^{1}\pi_{R\times 1}^{A\times 1}\;\;\; \pi^{R}_{B\times 1}{\rho^{\cdot}_{\mu^{\cdot}}}_{R} = (q\times 1)\pi_{R\times 1}^{A\times 1}\]
\noindent and
\[ \kappa^{\widetilde{\pi}_{R}^{1},q}_{B}\cdot{\rho^{\cdot}_{\mu^{\cdot}}}_{R} = 1,\]

\item for each pair of spans $R$, $\bar{R}$, a natural isomorphism
\[ {\rho^{\cdot}_{\mu^{\cdot}}}_{R,\bar{R}}\maps \left({\rho^{\cdot}_{\mu^{\cdot}}}_{\bar{R}}\right)_*\hom_{\Span}(1,\rho^{\cdot}_{B})(\otimes(1\times I)) \To \left({\rho^{\cdot}_{\mu^{\cdot}}}_{R}\right)^*\hom_{\Span}(\rho^{\cdot}_{A},1)1,\]
\noindent consisting of, for each map of spans $f_{R}$, an isomorphism of maps of spans
\[ {\rho^{\cdot}_{\mu^{\cdot}}}_{f_{R}}\maps (\widetilde{\pi}\cdot\kappa\cdot(1*(f_{R}\times 1)),\; {\rho^{\cdot}_{\mu^{\cdot}}}_{\bar{R}}(1*(f_{R}\times 1)),\; (\beta_{R}\times 1)\cdot\pi) \To \]
\[ (((\alpha_{R}\times 1)\cdot\pi{\rho^{\cdot}_{\mu^{\cdot}}}_{R})(\widetilde{\pi}\cdot\kappa),\; (f_{R}*1){\rho^{\cdot}_{\mu^{\cdot}}}_{R},\; 1),\]
\noindent consisting of the unique $2$-cell
\[ {\rho^{\cdot}_{\mu^{\cdot}}}_{f_{R}}\maps {\rho^{\cdot}_{\mu^{\cdot}}}_{\bar{R}}(1*(f_{R}\times 1)) \To (f_{R}* 1){\rho^{\cdot}_{\mu^{\cdot}}}_{R}\]
\noindent in $\B$ such that
\[ \pi_{\bar{R}}^{B\times 1}\cdot{\rho^{\cdot}_{\mu^{\cdot}}}_{f_{R}} = 1 \;\;\;\textrm{ and }\;\;\;  \pi^{\bar{R}}_{B\times 1}\cdot{\rho^{\cdot}_{\mu^{\cdot}}}_{f_{R}} = (\beta_{R}^{-1}\times 1)\cdot\pi_{R\times 1}^{A\times 1}, \]
\end{itemize}

\item an invertible counit modification
\[ {\rho^{\cdot}_{\epsilon}}\maps \rho^{\cdot}_{\mu}\rho^{\cdot}_{\mu^{\cdot}} \To 1\]
\noindent consisting of, for each span $R$, an isomorphism of maps of spans
\[ {\rho^{\cdot}_{\epsilon}}_{R}\maps (\widetilde{\pi}\cdot\kappa,\; \rho^{\cdot}_{\mu}\rho^{\cdot}_{\mu^{\cdot}},\; \widetilde{\pi}\cdot\kappa^{-1}\cdot\rho^{\cdot}_{\mu^{\cdot}}) \To (1,\; 1,\; 1)\]
\noindent defined by the unique $2$-cell ${\rho^{\cdot}_{\epsilon}}_{R}$ in $\B$ such that
\[ \pi_{A\times 1}^{R\times 1}\cdot {\rho^{\cdot}_{\epsilon}}_{R} = {\kappa_{A\times 1}^{p\times 1,1}}^{-1}\;\;\;\textrm{ and }\;\;\;  \pi^{A\times 1}_{R\times 1}\cdot {\rho^{\cdot}_{\epsilon}}_{R} = 1,\]

\item an invertible unit modification
\[ \rho^{\cdot}_{\eta}\maps 1\To \rho^{\cdot}_{\mu^{\cdot}}\rho^{\cdot}_{\mu}\]
\noindent consisting of, for each span $R$, an isomorphism of maps of spans
\[ {\rho^{\cdot}_{\eta}}_{R}\maps (1,1,1) \To (\widetilde{\pi}\cdot\kappa\cdot\rho^{\cdot}_{\mu},\;\rho^{\cdot}_{\mu^{\cdot}}\rho^{\cdot}_{\mu},\; \widetilde{\pi}\cdot\kappa^{-1})\]
\noindent defined by the unique $2$-cell ${\rho_{\eta}}_{R}$ in $\B$ such that
\[ \pi_{R}^{B\times 1}\cdot {\rho^{\cdot}_{\eta}}_{R} = 1 \;\;\;\textrm{ and }\;\;\; \pi_{B\times 1}^{R}\cdot {\rho^{\cdot}_{\eta}}_{R} =  \widetilde{\pi}_{B\times 1}^{B}\cdot{\kappa_{B}^{\widetilde{\pi}_{B}^{1},q}}^{-1},\]

\end{itemize}

\item an identity modification with component equations of maps of spans
\[ (1,\;\chi({\rho^{\cdot}_{\mu}}_{R}*1)(1*{\rho^{\cdot}_{\mu}}_{S}),\; \widetilde{\pi}\cdot {\kappa}^{-1}\cdot \pi) = (1,\;{\rho^{\cdot}_{\mu}}_{SR}(\chi*1), \widetilde{\pi}\cdot {\kappa}^{-1}\cdot (\chi*1)),\]
\item an identity modification with component equations of maps of spans
\[ (1,\; {\rho^{\cdot}_{\mu}}_{A}\iota {\bf r}^{-1},\; \widetilde{\pi}\cdot {\kappa}^{-1}\cdot\iota {\bf r}^{-1}) =  (1,\; \iota {\bf l}^{-1},\; 1),\]
\end{itemize}

\item a strict adjoint equivalence $(\mu_{\rho},\;\mu_{\rho}^{\cdot},\epsilon_{\mu_{\rho}},\;\eta_{\mu_{\rho}})$
\noindent consisting of
\begin{itemize}
\item a trimodification
\[ \mu_{\rho}\maps \rho\rho^{\cdot} \Rightarrow I_{1},\]
\noindent consisting of
\begin{itemize}
\item for each object $A$ in $\Span(\B)$, a map of spans
\[
\def\objectstyle{\scriptstyle}
  \def\labelstyle{\scriptstyle}
   \xy
   (-20,0)*+{A}="2";
   (0,20)*+{(A\times 1)(A\times 1)}="1";
   (0,-20)*+{A}="4";
   (20,0)*+{A}="3";
        {\ar_{\widetilde{\pi}_{A}^{1}\pi_{A\times 1}^{A\times 1}} "1";"2"};
        {\ar^{\widetilde{\pi}_{A}^{1}\pi_{A\times 1}^{A\times 1}} "1";"3"};
        {\ar^{1} "4";"2"};
        {\ar_{1} "4";"3"};
        {\ar^{\mu_{\rho}} "1";"4"};
        {\ar@{=}_<<{\scriptstyle } (12,2); (9,-1)};
        {\ar@{=}^<<{\scriptstyle } (-12,2); (-9,-1)};
\endxy
\]
\noindent where $\mu_{\rho} := {\mu_{\rho}}_{A} = \widetilde{\pi}_{A}^{1}\pi_{A\times 1}^{A\times 1}$,

\item and, for each pair $A, B \in\Span(\B)$, an identity modification
\[ m_{\mu_{\rho}}\maps (1,\;(\mu_{\rho} *1)(1*\rho_{\mu})({\rho_{\mu}^{\cdot}}*1),\; \widetilde{\pi}\cdot\kappa^{-1}\cdot \pi(\rho^{\cdot}_{\mu}*1)) \To (1,\; I_{1}(1*\mu_{\rho}),\; 1),\]
\noindent consisting of, for each span, an equation of maps of spans,
\end{itemize}

\item a trimodification
\[ \mu^{\cdot}_{\rho}\maps I_{1} \Rightarrow \rho\rho^{\cdot},\]
\noindent consisting of
\begin{itemize}
\item for each object $A$ in $\Span(\B)$, a map of spans
\[
\def\objectstyle{\scriptstyle}
  \def\labelstyle{\scriptstyle}
   \xy
   (-20,0)*+{A}="2";
   (0,20)*+{A}="1";
   (0,-20)*+{(A\times 1)(A\times 1)}="4";
   (20,0)*+{A}="3";
        {\ar^{\widetilde{\pi}_{A}^{1}\pi_{A\times 1}^{A\times 1}} "4";"2"};
        {\ar_{\widetilde{\pi}_{A}^{1}\pi_{A\times 1}^{A\times 1}} "4";"3"};
        {\ar_{1} "1";"2"};
        {\ar^{1} "1";"3"};
        {\ar^{\mu_{\rho}^{\cdot}} "1";"4"};
        {\ar@{=}_<<{\scriptstyle } (12,2); (9,-1)};
        {\ar@{=}^<<{\scriptstyle } (-12,2); (-9,-1)};
\endxy
\]
\noindent where $\mu_{\rho}^{\cdot}$ is the unique $1$-cell in $\B$ such that
\[ \pi_{A\times 1}^{A\times 1}\cdot\mu_{\rho}^{\cdot} = \widetilde{\pi}_{A\times 1}^{A}\;\;\;\textrm{ and }\;\;\;\kappa_{A\times 1}^{(A\times 1)(A\times 1)}\cdot\mu_{\rho}^{\cdot} = 1,\]
\item and, for each pair $A, B \in\Span(\B)$, an identity modification
\[ m_{\mu_{\rho}^{\cdot}}\maps (1,\;(\mu_{\rho}^{\cdot} *1)I,\; 1) \To (1,\; (1*\rho_{\mu})({\rho_{\mu}^{\cdot}}*1)(1*\mu_{\rho}^{\cdot}),\; \widetilde{\pi}\cdot{\kappa}^{-1}\cdot \pi({\rho_{\mu}^{\cdot}}*1)(1*\mu_{\rho}^{\cdot})),\]
\noindent consisting of, for each span, an equation of maps of spans,
\end{itemize}

\item an identity counit perturbation
\[
   \xy
   (-9,0)*+{\epsilon_{\mu_{\rho}}\maps \mu_{\rho}\mu_{\rho}^{\cdot}}="1"; (10,0)*+{1_{I_{1}}}="2";
        {\ar@3{->}_{} "1";"2"};
\endxy
\]
\noindent consisting of, for each object $A\in\Span(\B)$, an equation of maps of spans,
\item and, an identity unit perturbation
\[
   \xy
   (-10,0)*+{\eta_{\mu_{\rho}}\maps 1_{\rho\rho^{\cdot}}}="1"; (10,0)*+{\mu_{\rho}^{\cdot}\mu_{\rho}}="2";
        {\ar@3{->}_{} "1";"2"};
\endxy
\]
\noindent consisting of, for each object $A\in\Span(\B)$, an equation of maps of spans,
\end{itemize}

\item a strict adjoint equivalence $(\eta_{\rho},\;\eta_{\rho}^{\cdot},\epsilon_{\eta_{\rho}},\;\eta_{\eta_{\rho}})$, consisting of
\begin{itemize}
\item a trimodification
\[
   \xy
   (-10,0)*+{\eta_{\rho}\maps \rho^{\cdot}\rho}="1"; (9,0)*+{1_{1\times I}}="2";
        {\ar@3{->}_{} "1";"2"};
\endxy
\]
\begin{itemize}
\item for each object $A\in\Span(\B)$, a map of spans
\[
\def\objectstyle{\scriptstyle}
  \def\labelstyle{\scriptstyle}
   \xy
   (-20,0)*+{A\times 1}="2";
   (0,20)*+{(A\times 1)(A\times 1)}="1";
   (0,-20)*+{A\times 1}="4";
   (20,0)*+{A\times 1}="3";
        {\ar_{\pi_{A\times 1}^{A\times 1}} "1";"2"};
        {\ar^{\pi_{A\times 1}^{A\times 1}} "1";"3"};
        {\ar^{1} "4";"2"};
        {\ar_{1} "4";"3"};
        {\ar^{\eta_{\rho}} "1";"4"};
        {\ar@{=}_<<{\scriptstyle } (12,2); (9,-1)};
        {\ar@{=}^<<{\scriptstyle } (-12,2); (-9,-1)};
\endxy
\]
\noindent where $\eta_{\rho} := {\eta_{\rho}}_{A} = \pi_{A\times 1}^{A\times 1}$,
\item and, for each pair of objects $A, B\in\Span(\B)$, an identity modification
\[ m_{\eta_{\rho}}\maps (1,\;(\eta_{\rho} *1)(1*{\rho_{\mu}^{\cdot}})(\rho_{\mu}*1),\; \widetilde{\pi}\cdot\kappa^{-1}\cdot \pi(\rho_{\mu}*1)) = (1,\; I(1*\eta_{\rho}),\; 1),\]
\noindent consisting of, for each span, an equation of maps of spans,

\end{itemize}

\item a trimodification
\[
   \xy
   (-11,0)*+{\eta_{\rho}^{\cdot}\maps 1_{\otimes(1\times I)}}="1"; (10,0)*+{\rho^{\cdot}\rho,}="2";
        {\ar@3{->}_{} "1";"2"};
\endxy
\]
\noindent consisting of
\begin{itemize}
\item for each object $A\in\Span(\B)$, a map of spans
\[
\def\objectstyle{\scriptstyle}
  \def\labelstyle{\scriptstyle}
   \xy
   (-20,0)*+{A\times 1}="2";
   (0,20)*+{A\times 1}="1";
   (0,-20)*+{(A\times 1)(A\times 1)}="4";
   (20,0)*+{A\times 1}="3";
        {\ar^{\pi_{A\times 1}^{A\times 1}} "4";"2"};
        {\ar_{\pi_{A\times 1}^{A\times 1}} "4";"3"};
        {\ar_{1} "1";"2"};
        {\ar^{1} "1";"3"};
        {\ar^{\eta_{\rho}^{\cdot}} "1";"4"};
        {\ar@{=}_<<{\scriptstyle } (12,2); (9,-1)};
        {\ar@{=}^<<{\scriptstyle } (-12,2); (-9,-1)};
\endxy
\]
\noindent where $\eta_{\rho}^{\cdot} := {\eta_\rho}^{\cdot}_{A}$ is the unique $1$-cell in $\B$ such that
\[ \pi_{A\times 1}^{A\times 1}\cdot\eta_{\rho}^{\cdot} = 1\;\;\;\textrm{ and }\;\;\;\kappa_{A\times 1}^{(A\times 1)(A\times 1)}\cdot\eta_{\rho}^{\cdot} = 1,\]

\item and, for each pair of objects $A, B\in\Span(\B)$, an identity modification
\[ m_{\eta^{\cdot}_{\rho}}\maps (1,\;(\eta_{\rho}^{\cdot} *1)I,\; 1) = (1,\; (1*{\rho_{\mu}^{\cdot}})({\rho_{\mu}}*1)(1*\eta_{\rho}^{\cdot}),\; \widetilde{\pi}\cdot\kappa^{-1}\cdot\pi(\rho_{\mu}*1)(1*\eta_{\rho}^{\cdot}))\]
\noindent consisting of, for each span, an equation of maps of spans,

\end{itemize}

\item an identity counit perturbation
\[
   \xy
   (-10,0)*+{\epsilon_{\eta_{\rho}}\maps \eta_{\rho}\eta_{\rho}^{\cdot}}="1"; (10,0)*+{1_{\rho^{\cdot}\rho}}="2";
        {\ar@3{->}_{} "1";"2"};
\endxy
\]
\noindent consisting of, for each object $A\in\Span(\B)$, an equation of maps of spans,
\item and, an identity unit perturbation
\[
   \xy
   (-10,0)*+{\eta_{\eta_{\rho}}\maps 1_{\otimes(1\times I)}}="1"; (10,0)*+{\eta_{\rho}^{\cdot}\eta_{\rho}}="2";
        {\ar@3{->}_{} "1";"2"};
\endxy
\]
\noindent consisting of, for each object $A\in\Span(\B)$, an equation of maps of spans,
\end{itemize}

\item an identity perturbation
\[
   \xy
   (-31,0)*+{\Phi_{\rho}\maps (1,\;{\bf l}(1*\mu_{\rho})(\eta_{\rho}^{\cdot}*1){\bf r}^{-1},\; \kappa^{-1}\cdot (1*\mu_{\rho})(\eta_{\rho}^{\cdot}*1){\bf r}^{-1})}="1"; (32,0)*+{(1,1,1),}="2";
        {\ar@3{->}_{} "1";"2"};
\endxy
\]
\noindent consisting of, for each object $A\in\Span(\B)$, an equation of maps of spans,

\item an identity perturbation
\[
   \xy
   (-31,0)*+{\Psi_{\rho}\maps (\kappa\cdot (\mu_{\rho}*1)(1*\eta_{\rho}^{\cdot}){\bf l}^{-1},\; {\bf r}(\mu_{\rho}*1)(1*\eta_{\rho}^{\cdot}){\bf l}^{-1},\; 1)}="1"; (30,0)*+{(1,1,1)}="2";
        {\ar@3{->}_{} "1";"2"};
\endxy
\]
\noindent consisting of, for each object $A\in\Span(\B)$, an equation of maps of spans.
\end{itemize}
\endproposition

\proof
The proof is similar to that of the previous proposition for the biadjoint biequivalence $\lambda$.  We omit the details.
\endproof

\subsection{Monoidal Adjoint Equivalences}

\subsubsection*{Monoidal Pentagonator Modification}

We define a strict adjoint equivalence
\[
   \xy
   (-14,0)*+{\Pi\maps (1\times\alpha)\alpha(\alpha\times 1)}="1"; (13,0)*+{\alpha\alpha}="2";
        {\ar@3{->}_{} "1";"2"};
\endxy
\]
\noindent in the strict $2$-category of tritransformations, trimodifications, and perturbations.

\proposition \label{monoidalpentagonatormodification}
There is a strict adjoint equivalence $(\Pi,\;\Pi^*,\;1,\;1)$ consisting of
\begin{itemize}
\item a trimodification
\[
   \xy
   (-14,0)*+{\Pi\maps (1\times\alpha)\alpha(\alpha\times 1)}="1"; (13,0)*+{\alpha\alpha}="2";
        {\ar@3{->}_{} "1";"2"};
\endxy
\]
\noindent consisting of
\begin{itemize}
\item for each $4$-tuple of objects $A, B, C, D$ in $\Span(\B)$, a map of spans
\[
\def\objectstyle{\scriptstyle}
  \def\labelstyle{\scriptstyle}
   \xy
   (-20,0)*+{A\times (B\times (C\times D))}="2";
   (0,20)*+{(A\times ((B\times C)\times D))(((A\times (B\times C))\times D)(((A\times B)\times C)\times D))}="1";
   (0,-20)*+{((A\times B)\times (C\times D))(((A\times B)\times C)\times D)}="4";
   (20,0)*+{((A\times B)\times C)\times D}="3";
        {\ar^{a_{AB(CD)}\pi_{(A\times B)\times (C\times D)}^{((A\times B)\times C)\times D}} "4";"2"};
        {\ar_{\pi^{(A\times B)\times (C\times D)}_{((A\times B)\times C)\times D}} "4";"3"};
        {\ar_{(1_A\times a_{BCD})\pi_{A\times ((B\times C)\times D)}^{((A\times (B\times C))\times D)(((A\times B)\times C)\times D)}} "1";"2"};
        {\ar^{\pi_{((A\times B)\times C)\times D}^{(A\times (B\times C))\times D}\pi^{A\times ((B\times C)\times D)}_{((A\times (B\times C))\times D)(((A\times B)\times C)\times D)}} "1";"3"};
        {\ar^{\Pi} "1";"4"};
        {\ar@{=}_<<{\scriptstyle } (9,2); (7,-1)};
        {\ar@{=}^<<{\scriptstyle } (-9,2); (-7,-1)};
\endxy
\]
\noindent where $\Pi := \Pi_{ABCD}$ is the unique $1$-cell in $\B$ such that
\[\pi_{((A\times B)\times C)\times D}^{(A\times B)\times (C\times D)}\cdot\Pi_{ABCD} =  \pi_{((A\times B)\times C)\times D}^{(A\times (B\times C))\times D}\pi_{((A\times (B\times C))\times D)(((A\times B)\times C)\times D)}^{A\times ((B\times C)\times D)},\]
\[\pi^{((A\times B)\times C)\times D}_{(A\times B)\times (C\times D)}\cdot\Pi_{ABCD} =  a_{AB(CD)}^{-1}(1_{A}\times a_{BCD})\pi_{A\times ((B\times C)\times D)}^{((A\times (B\times C))\times D)(((A\times B)\times C)\times D)},\]
\noindent and
\[\kappa_{(A\times B)\times (C\times D)}^{1,a}\cdot\Pi_{ABCD} =  1,\]

\item and, for each two $4$-tuples $(A, B, C, D), (A', B', C', D')$ of objects in $\Span(\B)$, an identity modification
\[ m_{\Pi}\maps (1,\;(\Pi*1)(1*(1\times\alpha))((1*\alpha)*1)(((\alpha\times 1)*1)*1),\; (1\times(a\cdot\kappa^{-1}))\cdot\pi(1*\alpha*1)((\alpha\times 1)*1^2))\]
\[ \To (1,\; (1*\alpha)(\alpha*1)(1*\Pi),\; a\cdot\kappa^{-1}\cdot\pi(\alpha*1)(1*\Pi)),\]
\noindent consisting of, for each four spans, an equation of maps of spans,
\end{itemize}

\item a trimodification
\[
   \xy
   (-14,0)*+{\Pi^{*}\maps\alpha\alpha}="1"; (13,0)*+{(1\times\alpha)\alpha(\alpha\times 1),}="2";
        {\ar@3{->}_{} "1";"2"};
\endxy
\]
\noindent consisting of
\begin{itemize}
\item for each $4$-tuple of objects $A, B, C, D$ in $\Span(\B)$, a map of spans
\[
\def\objectstyle{\scriptstyle}
  \def\labelstyle{\scriptstyle}
   \xy
   (-20,0)*+{A\times (B\times (C\times D))}="2";
   (0,20)*+{((A\times B)\times (C\times D))(((A\times B)\times C)\times D)}="1";
   (0,-20)*+{(A\times ((B\times C)\times D))(((A\times (B\times C))\times D)(((A\times B)\times C)\times D))}="4";
   (20,0)*+{((A\times B)\times C)\times D}="3";
        {\ar_{a_{AB(CD)}\pi_{(A\times B)\times (C\times D)}^{((A\times B)\times C)\times D}} "1";"2"};
        {\ar^{\pi^{(A\times B)\times (C\times D)}_{((A\times B)\times C)\times D}} "1";"3"};
        {\ar^{(1_A\times a_{BCD})\pi_{A\times ((B\times C)\times D)}^{((A\times (B\times C))\times D)(((A\times B)\times C)\times D)}} "4";"2"};
        {\ar_{\pi_{((A\times B)\times C)\times D}^{(A\times (B\times C))\times D}\pi^{A\times ((B\times C)\times D)}_{((A\times (B\times C))\times D)(((A\times B)\times C)\times D)}} "4";"3"};
        {\ar^{\Pi^*} "1";"4"};
        {\ar@{=}_<<{\scriptstyle } (9,2); (7,-1)};
        {\ar@{=}^<<{\scriptstyle } (-9,2); (-7,-1)};
\endxy
\]
\noindent where and $\Pi^{*}:=\Pi^{*}_{ABCD}$ is a $1$-cell in $\B$ such that
\[\pi^{((A\times (B\times C))\times D)}_{(((A\times B)\times C)\times D)}\pi_{((A\times (B\times C))\times D)(((A\times B)\times C)\times D)}^{A\times ((B\times C)\times D)}\cdot\Pi^{*}_{ABCD} = \pi_{(((A\times B)\times C)\times D}^{(A\times B)\times (C\times D)},\]
\[\pi_{((A\times (B\times C))\times D)}^{(((A\times B)\times C)\times D)}\pi_{((A\times (B\times C))\times D)(((A\times B)\times C)\times D)}^{A\times ((B\times C)\times D)}\cdot\Pi^{*}_{ABCD} = (a\times 1)a^{-1}\pi^{(((A\times B)\times C)\times D}_{(A\times B)\times (C\times D)},\]
\[\pi^{((A\times (B\times C))\times D)(((A\times B)\times C)\times D)}_{A\times ((B\times C)\times D)}\cdot\Pi^{*}_{ABCD} = (1\times a^{-1})a\pi_{(A\times B)\times (C\times D)}^{((A\times B)\times C)\times D},\]
\noindent and
\[\kappa_{A\times ((B\times C)\times D)}^{1,a\pi_{(A\times (B\times C))\times D}^{((A\times B)\times C)\times D}}\cdot\Pi^{*}_{ABCD} = ((a\times 1)a^{-1})\cdot\kappa_{(A\times B)\times (C\times D)}^{1,a},\]

\item and, for each two $4$-tuples $(A, B, C, D), (A', B', C', D')$ of objects in $\Span(\B)$, an identity modification
\[m_{\Pi^*}\maps (1,\; (\Pi^**1)(1*\alpha)(\alpha*1),\; a\cdot\kappa^{-1}\cdot(\alpha*1))\To\]
\[ (1,\; (1*(1\times\alpha))((1*\alpha*1)((\alpha\times 1)*1^2)(1*\Pi^*),\; (1\times a\cdot\kappa^{-1})\cdot\pi(1*\alpha*1)((\alpha\times 1)*1^2)(1*\Pi^*)),\]
\noindent consisting of, for each three spans, an equation of maps of spans,
\end{itemize}

\item and, identity counit and unit perturbations.
\end{itemize}
\endproposition

\proof
Since $m_{\Pi}$ and the modification components of $\alpha$ are all identities, $\Pi$ is a trimodification.

To define $\Pi^*$ uniquely we need an auxiliary map, which is the unique $1$-cell
\[ \Pi^0\maps ((A\times B)\times (C\times D))(((A\times B)\times C)\times D)\to ((A\times (B\times C))\times D)(((A\times B)\times C)\times D)\]
\noindent in $\B$ such that
\[ \pi_{((A\times B)\times C)\times D}^{(A\times (B\times C))\times D}\cdot\Pi^0 = \pi_{((A\times B)\times C)\times D}^{(A\times B)\times (C\times D)},\;\;\pi^{((A\times B)\times C)\times D}_{(A\times (B\times C))\times D}\cdot\Pi^0 = (a\times 1)a^{-1}\pi^{((A\times B)\times C)\times D}_{(A\times B)\times (C\times D)},\]
\noindent and
\[ \kappa_{(A\times (B\times C))\times D}^{1,a\times 1}\cdot \Pi^0 = (a\times 1)a^{-1}\cdot\kappa_{(A\times B)\times (C\times D)}^{((A\times B)\times (C\times D))(((A\times B)\times C)\times D)}.\]
\noindent It is then straightforward to see that $\Pi^*$ is a modification, and that the pair of modifications define a strict adjoint equivalence.
\endproof

\subsubsection*{Monoidal Left Mediator Unit Modification}

We define an adjoint equivalence
\[
   \xy
   (-10,0)*+{l\maps \lambda\alpha}="1"; (9,0)*+{(\lambda\times 1)}="2";
        {\ar@3{->}_{} "1";"2"};
\endxy
\]
\noindent in the strict $2$-category of tritransformations, trimodifications, and perturbations.

\proposition \label{monoidalleftmediatormodification}
There is an adjoint equivalence $(l,\;l^{*},\;1,\;\eta_{l})$, consisting of
\begin{itemize}
\item a trimodification
\[
   \xy
   (-9,0)*+{l\maps \lambda\alpha}="1"; (10,0)*+{(\lambda\times 1),}="2";
        {\ar@3{->}_{} "1";"2"};
\endxy
\]
\noindent consisting of
\begin{itemize}
\item for each pair of objects $A, B$ in $\Span(\B)$, a map of spans
\[
\def\objectstyle{\scriptstyle}
  \def\labelstyle{\scriptstyle}
   \xy
   (-20,0)*+{A\times B}="2";
   (0,20)*+{(1\times (A\times B))((1\times A)\times B)}="1";
   (0,-20)*+{(1\times A)\times B}="4";
   (20,0)*+{(1\times A)\times B}="3";
        {\ar^{\widetilde{\pi}_{A}^{1}\times 1_{B}} "4";"2"};
        {\ar_{1} "4";"3"};
        {\ar_{\widetilde{\pi}_{A\times B}^{1}\pi_{1\times (A\times B)}^{(1\times A)\times B}} "1";"2"};
        {\ar^{\pi^{1\times (A\times B)}_{(1\times A)\times B}} "1";"3"};
        {\ar^{l_{AB}} "1";"4"};
        {\ar@{=}_<<{\scriptstyle } (14,2); (11,-1)};
        {\ar@{=>}^<<{\scriptstyle \kappa^{-1}} (-14,2); (-11,-1)};
\endxy
\]
\noindent where $\kappa :=  \kappa_{1\times (A\times B)}^{1,a}$ and $l_{AB} := \pi_{(1\times A)\times B}^{1\times (A\times B)}$,

\item and, for each two pairs $(A, B), (A', B')$ of objects in $\Span(\B)$, an identity modification
\[ m_{l}\maps (1, (l*1)(1*\lambda)(\alpha*1),\; (\widetilde{\pi}\cdot\kappa^{-1})\cdot \pi(\alpha*1))\]
\[ \To (1,\;(\lambda\times 1)(1*l),\; (1\times \widetilde{\pi}\cdot\kappa^{-1}\cdot (1*l))(\widetilde{\pi}\cdot\kappa\cdot\pi)),\]
\noindent consisting of, for each two spans, an equation of maps of spans,

\end{itemize}

\item a trimodification
\[
   \xy
   (-9,0)*+{l^*\maps \lambda\times 1}="1"; (10,0)*+{\lambda\alpha,}="2";
        {\ar@3{->}_{} "1";"2"};
\endxy
\]
\noindent consisting of

\begin{itemize}
\item for each pair of objects $A, B$ in $\Span(\B)$, a map of spans
\[
\def\objectstyle{\scriptstyle}
  \def\labelstyle{\scriptstyle}
   \xy
   (-20,0)*+{(A\times B)}="2";
   (0,-20)*+{(1\times (A\times B))((1\times A)\times B)}="1";
   (0,20)*+{(1\times A)\times B}="4";
   (20,0)*+{(1\times A)\times B}="3";
        {\ar_{\widetilde{\pi}_{A}^{1}\times 1_{B}} "4";"2"};
        {\ar^{1} "4";"3"};
        {\ar^{\widetilde{\pi}_{A\times B}^{1}\pi_{1\times (A\times B)}^{(1\times A)\times B}} "1";"2"};
        {\ar_{\pi^{1\times (A\times B)}_{(1\times A)\times B}} "1";"3"};
        {\ar^{l_{AB}^{*}} "4";"1"};
        {\ar@{=}_<<{\scriptstyle } (14,2); (11,-1)};
        {\ar@{=}^<<{\scriptstyle } (-14,2); (-11,-1)};
\endxy
\]
\noindent where $l_{AB}^{*}$ is the unique $1$-cell in $\B$ such that
\[ \pi_{(1\times A)\times B}^{1\times (A\times B)}\cdot l_{AB}^{*} = 1,\;\;\;\pi^{(1\times A)\times B}_{1\times (A\times B)}\cdot l_{AB}^{*} = a,\;\;\textrm{ and }\;\;\kappa_{1\times (A\times B)}^{1,a}\cdot l_{AB}^{*} = 1,\]

\item and, for each two pairs $(A, B), (A', B')$ of objects in $\Span(\B)$, an identity modification
\[m_{l^{*}}\maps (1,\; (l^{*}*1)(\lambda\times 1),\; \widetilde{\pi}_{A'\times B'}^{1}\cdot {\kappa^{-1}}_{(A'\times B')\times 1}^{1,(p'\times q')\times 1}) \To \]
\[(1,\; (1*\lambda)(\alpha*1)(1*l^{*}),\; (\widetilde{\pi}_{A'\times B'}^{1}\cdot{\kappa^{-1}}_{(A'\times B')\times 1}^{1,(p'\times q')\times 1})\cdot\pi (\alpha*1)(1*l^*)),\] 
\noindent consisting of, for each two spans, an equation of maps of spans,

\end{itemize}

\item an identity counit perturbation,

\item and, a unit perturbation
\[ \eta_{l}\maps (1,1,1)\To (1,\; l^{*}l,\; 1),\]
\noindent consisting of, for each pair of objects $A, B\in\Span(\B)$, an isomorphism of maps of spans
\[ \eta_{l_{AB}}\maps (1,1,1)\To (1,\; l_{AB}^{*}l_{AB},\; 1),\]
\noindent where $\eta_{l_{AB}}$ is the unique $2$-cell in $\B$ such that
\[ \pi_{(1\times A)\times B}^{1\times (A\times B)}\cdot\eta_{l_{AB}} = 1\;\;\;\textrm{ and } \pi^{(1\times A)\times B}_{1\times (A\times B)}\cdot\eta_{l_{AB}} = {\kappa^{-1}}_{1\times (A\times B)}^{1,a}.\]
\end{itemize}
\endproposition

\proof
The modification axioms hold since $m_{l}$, $m_{l^*}$ and the modification components of $\alpha$ and $\lambda$ are all identities.  The adjoint equivalence axioms follow.
\endproof

\subsubsection*{Monoidal Middle Mediator Unit Modification}

We define an adjoint equivalence
\[
   \xy
   (-11,0)*+{m\maps (1\times \lambda)\alpha}="1"; (12,0)*+{\rho\times 1}="2";
        {\ar@3{->}_{} "1";"2"};
\endxy
\]
\noindent in the strict $2$-category of tritransformations, trimodifications, and perturbations.

\proposition \label{monoidalmiddlemediatormodification}
There is an adjoint equivalence $(m,\;m^{*},\;1,\;\eta_{m})$, consisting of
\begin{itemize}
\item a trimodification
\[
   \xy
   (-11,0)*+{m\maps  (1\times \lambda)\alpha}="1"; (12,0)*+{\rho\times 1,}="2";
        {\ar@3{->}_{} "1";"2"};
\endxy
\]
\noindent consisting of
\begin{itemize}
\item for each pair of objects $A, B$ in $\Span(\B)$, a map of spans
\[
\def\objectstyle{\scriptstyle}
  \def\labelstyle{\scriptstyle}
   \xy
   (-20,0)*+{(A\times B)}="2";
   (0,20)*+{(A\times (1\times B))((A\times 1)\times B)}="1";
   (0,-20)*+{(A\times 1)\times B}="4";
   (20,0)*+{(A\times 1)\times B}="3";
        {\ar^{\rho_{A}\times 1_{B}} "4";"2"};
        {\ar_{1} "4";"3"};
        {\ar_{(1_{A}\times\lambda_{B})\pi_{A\times (1\times B)}^{(A\times 1)\times B}} "1";"2"};
        {\ar^{\pi^{A\times (1\times B)}_{(A\times 1)\times B}} "1";"3"};
        {\ar^{m_{AB}} "1";"4"};
        {\ar@{=}_<<{\scriptstyle } (14.5,2); (11,-1)};
        {\ar@{=>}^<<{\scriptstyle (1\times \widetilde{\pi})\cdot\kappa^{-1}} (-14.5,2); (-11,-1)};
\endxy
\]
\noindent where $(1\times \widetilde{\pi})\cdot\kappa^{-1} := (1_{A}\times \widetilde{\pi}_{B}^{1})\cdot{\kappa^{-1}}_{A\times (1\times B)}^{1,a}$ and $m_{AB} := \pi_{(A\times 1)\times B}^{A\times (1\times B)}$,

\item and, for each two pairs $(A, B), (A', B')$ of objects in $\Span(\B)$, an identity modification
\[m_{m}\maps (1,\; (m*1)(1*(1\times \lambda))(\alpha*1),\; ((1\times\lambda)\cdot\kappa^{-1})\cdot\pi (\alpha*1)) \To \]
\[(1,\; (\rho\times 1)(1*m),\; ((\widetilde{\pi}\cdot\kappa^{-1}\times 1)\cdot \pi(1*m))((1\times \widetilde{\pi})\cdot\kappa^{-1} \cdot\pi)),\]
\noindent consisting of, for each two spans, an equation of maps of spans,

\end{itemize}

\item a trimodification
\[
   \xy
   (-14,0)*+{m^{*}\maps (\rho\times 1)}="1"; (12,0)*+{(1\times \lambda)\alpha,}="2";
        {\ar@3{->}_{} "1";"2"};
\endxy
\]
\noindent consisting of

\begin{itemize}
\item for each pair of objects $A, B$ in $\Span(\B)$, a map of spans
\[
\def\objectstyle{\scriptstyle}
  \def\labelstyle{\scriptstyle}
   \xy
   (-20,0)*+{(A\times B)}="2";
   (0,20)*+{(A\times 1)\times B}="1";
   (0,-20)*+{(A\times (1\times B))((A\times 1)\times B)}="4";
   (20,0)*+{(A\times 1)\times B}="3";
        {\ar^{(1_{A}\times\widetilde{\pi}_{B}^{1})\pi_{A\times (1\times B)}^{(A\times 1)\times B}} "4";"2"};
        {\ar_{\pi^{A\times (1\times B)}_{(A\times 1)\times B}} "4";"3"};
        {\ar_{\widetilde{\pi}_{A}^{1}\times 1_{B}} "1";"2"};
        {\ar^{1} "1";"3"};
        {\ar^{m^{*}_{AB}} "1";"4"};
        {\ar@{=}_<<{\scriptstyle } (14,2); (11,-1)};
        {\ar@{=}^<<{\scriptstyle } (-14,2); (-11,-1)};
\endxy
\]
\noindent where $m^{*}_{AB}$ is the unique $1$-cell such that
 \[ \pi_{(A\times 1)\times B}^{A\times (1\times B)}\cdot m^{*} = 1\;\;\;\textrm{ and } \pi_{A\times (1\times B)}^{(A\times 1)\times B}\cdot m^{*} = a,\]
 \noindent and
 \[ \kappa_{A\times (1\times B)}^{1,a}\cdot m^{*} = 1,\]

\item and, for each two pairs $(A, B), (A', B')$ of objects in $\Span(\B)$, an identity modification
\[ m_{m^{*}}\maps (1,\; (m^{*}*1)(\rho\times 1),\; (\widetilde{\pi}_{B}^{1}\cdot\kappa^{-1}\times 1)\cdot\pi )\]
\[ \To (1,\; (1*(1\times \lambda))(\alpha*1)(1*m^{*}),\; (1\times \widetilde{\pi}\cdot\kappa^{-1})\cdot \pi(\alpha *1)(1*m^{*})),\]
\noindent consisting of, for each two spans, an equation of maps of spans,

\end{itemize}

\item an identity counit perturbation,

\item and, a unit perturbation
\[ \eta_{m}\maps (1,1,1)\To (1,\; m^{*}m,\; 1),\]
\noindent consisting of, for each pair of objects $A, B\in\Span(\B)$, an isomorphism of maps of spans
\[ \eta_{m_{AB}}\maps (1,1,1)\To (1,\; m_{AB}^{*}m_{AB},\; 1),\]
\noindent where $\eta_{m_{AB}}$ is the unique $2$-cell in $\B$ such that
\[ \pi_{(A\times 1)\times B}^{A\times (1\times B)}\cdot \eta_{m_{AB}} = 1\;\;\;\textrm{ and} \;\;\;\pi^{(A\times 1)\times B}_{A\times (1\times B)}\cdot \eta_{m_{AB}} = {\kappa^{-1}}_{A\times (1\times B)}^{1,a}.\]
\end{itemize}
\endproposition

\proof
The modification axioms hold since $m_{m}$, $m_{m^*}$ and the modification components of $\alpha$, $\lambda$, and $\rho$ are all identities.  The adjoint equivalence axioms follow.
\endproof

\subsubsection*{Monoidal Right Mediator Unit Modification}

We define an adjoint equivalence
\[
   \xy
   (-9,0)*+{r\maps (1\times\rho)\alpha}="1"; (9,0)*+{\rho}="2";
        {\ar@3{->}_{} "1";"2"};
\endxy
\]
\noindent in the $2$-category of transformations, modifications, and perturbations.

\proposition \label{monoidalrightmediatormodification}
There is an adjoint equivalence $(r,\;r^{*},\;1,\;\eta_{r})$, consisting of
\begin{itemize}
\item a trimodification
\[
   \xy
   (-10,0)*+{r\maps (1\times\rho)\alpha}="1"; (9,0)*+{\rho,}="2";
        {\ar@3{->}_{} "1";"2"};
\endxy
\]
\noindent consisting of
\begin{itemize}
\item for each pair of objects $A, B$ in $\Span(\B)$, a map of spans
\[
\def\objectstyle{\scriptstyle}
  \def\labelstyle{\scriptstyle}
   \xy
   (-20,0)*+{A\times B}="2";
   (0,20)*+{(A\times (B\times 1))((A\times B)\times 1)}="1";
   (0,-20)*+{(A\times B)\times 1}="4";
   (20,0)*+{(A\times B)\times 1}="3";
        {\ar^{\rho_{A}\times 1_{B}} "4";"2"};
        {\ar_{1} "4";"3"};
        {\ar_{(1\times \rho_{B})\pi_{A\times (B\times 1)}^{(A\times B)\times 1}} "1";"2"};
        {\ar^{\pi^{A\times (B\times 1)}_{(A\times B)\times 1}} "1";"3"};
        {\ar^{r_{AB}} "1";"4"};
        {\ar@{=}_<<{\scriptstyle } (14,2); (11,-1)};
        {\ar@{=>}^<<{\scriptstyle (1\times \rho)\cdot\kappa^{-1}} (-14.5,2); (-11,-1)};
\endxy
\]
\noindent where $(1\times\rho)\cdot\kappa^{-1} :=  (\widetilde{\pi}_{A}^{1}\times 1_{B})\cdot{\kappa^{-1}}_{A\times (B\times 1)}^{1,a}$ and $r_{AB} := \pi_{(A\times B)\times 1}^{A\times (B\times 1)}$,

\item and, for each two pairs $(A, B), (A', B')$ of objects in $\Span(\B)$, an identity modification
\[m_{r}\maps(1,\; (r*1)(1*(1\times\rho))(\alpha*1),\; \widetilde{\pi}\cdot\kappa^{-1}\cdot \pi (\alpha*1)) \To (1,\;\rho(1*r),\; \rho\cdot\kappa^{-1}\cdot\pi),\]
\noindent consisting of, for each two spans, an equation of maps of spans,

\end{itemize}

\item a trimodification
\[
   \xy
   (-10,0)*+{r^{*}\maps \rho}="1"; (9,0)*+{(1\times \rho)\alpha,}="2";
        {\ar@3{->}_{} "1";"2"};
\endxy
\]
\noindent consisting of

\begin{itemize}
\item for each pair of objects $A, B$ in $\Span(\B)$, a map of spans
\[
\def\objectstyle{\scriptstyle}
  \def\labelstyle{\scriptstyle}
   \xy
   (-20,0)*+{A\times B}="2";
   (0,20)*+{(A\times B)\times 1}="1";
   (0,-20)*+{(A\times (B\times 1))((A\times B)\times 1)}="4";
   (20,0)*+{(A\times B)\times 1}="3";
        {\ar^{(1_{A}\times \rho_{B})\pi_{A\times (B\times 1)}^{(A\times B)\times 1}} "4";"2"};
        {\ar_{\pi_{(A\times B)\times 1}^{A\times (B\times 1)}} "4";"3"};
        {\ar_{\rho_{A\times B}} "1";"2"};
        {\ar^{1} "1";"3"};
        {\ar^{r^{*}_{AB}} "1";"4"};
        {\ar@{=}_<<{\scriptstyle } (12,2); (9,-1)};
        {\ar@{=}^<<{\scriptstyle } (-12,2); (-9,-1)};
\endxy
\]
\noindent where $r^{*}_{AB}$ is the unique $1$-cell in $\B$ such that
\[ \pi_{(A\times B)\times 1}^{A\times (B\times 1)}r^{*}_{AB} = 1\,\;\;\;\pi_{A\times (B\times 1)}^{(A\times B)\times 1}r^{*}_{AB} = a,\;\;\textrm{ and }\;\;\; \kappa_{A\times (B\times 1)}^{1,a}\cdot r^{*}_{AB} = 1,\]

\item and, for each two pairs $(A, B), (A', B')$ of objects in $\Span(\B)$, an identity modification
\[m_{r^{*}}\maps (1,\; (r^{*}*1)(\rho*1),\; 1) \To (1,\; (1*(1\times\rho))(\alpha*1)(1*r^{*}),\;\widetilde{\pi}\cdot\kappa^{-1}\cdot \pi(\alpha*1)(1*r^{*})),\]
\noindent consisting of, for each two spans, an equation of maps of spans,

\end{itemize}

\item an identity counit perturbation,

\item and, a unit perturbation
\[ \eta_{r}\maps (1,1,1)\To (1,\; r^{*}r,\; 1),\]
\noindent consisting of, for each pair of objects $A, B\in\Span(\B)$, an isomorphism of maps of spans
\[ \eta_{r_{AB}}\maps (1,1,1)\To (1,\; r^{*}_{AB}r_{AB},\; 1)\]
\noindent where $\eta_{r_{AB}}$ is the unique $2$-cell in $\B$ such that
\[ \pi_{(A\times B)\times 1}^{A\times (B\times 1)}\cdot  \eta_{r_{AB}} = 1\;\;\;\textrm{ and }\;\;\; \pi_{A\times (B\times 1)}^{(A\times B)\times 1}\cdot  \eta_{r_{AB}} = {\kappa^{-1}}_{A\times (B\times 1)}^{1,a}.\]
\end{itemize}
\endproposition

\proof
The modification axioms hold since $m_{r}$, $m_{r^*}$ and the modification components of $\alpha$ and $\rho$ are all identities.  The adjoint equivalence axioms follow.
\endproof

\subsection{Monoidal Perturbations}

\subsubsection*{$K^5$ Perturbation}

\proposition \label{associativityperturbation}
There is a perturbation consisting of, for each $5$-tuple $A, B, C, D, E$ of objects in $\Span(\B)$, an isomorphism of maps of spans
\tiny
\[  K^{5}_{ABCDE}\maps (1,\; (\Pi*1)(1^2*{\alpha_{\mu^{\cdot}}}_{1,1,\alpha})(1*\Pi*1)((\Pi\times 1)*1^3),\; 1) \Rightarrow \] 
\[((a^{-1}\cdot\kappa)\cdot(1*\Pi*1)(1^3*(1\times\Pi))(1^2*{\alpha_{\mu^{\cdot}}}_{1,\alpha,1}*1^2),\; (1*\Pi)({\alpha_{\mu^{\cdot}}}_{\alpha,1,1}*1^2)(1*\Pi*1)(1^3*(1\times\Pi))(1^2*{\alpha_{\mu^{\cdot}}}_{1,\alpha,1}*1^2),\; 1),\]
\normalsize 
\noindent defined by the unique $2$-cell
\tiny
\[ K^{5}_{ABCDE}\maps (\Pi*1)(1^2*{\alpha_{\mu^{\cdot}}}_{1,1,\alpha})(1*\Pi*1)((\Pi\times 1)*1^3) \Rightarrow (1*\Pi)({\alpha_{\mu^{\cdot}}}_{\alpha,1,1}*1^2)(1*\Pi*1)(1^3*(1\times\Pi))(1^2*{\alpha_{\mu^{\cdot}}}_{1,\alpha,1}*1^2),\]
\normalsize 
\noindent in $\B$ such that
\[\pi_{(((AB)C)D)E}^{((AB)C)(DE)}\pi_{(((AB)C)(DE))((((AB)C)D)E)}^{(AB)(C(DE))}\cdot K^5 = \kappa_{(((AB)C)D)E}^{1,(a\times 1)\times 1}, \]
\tiny
\[\pi^{(((AB)C)D)E}_{((AB)C)(DE)}\pi_{(((AB)C)(DE))((((AB)C)D)E)}^{(AB)(C(DE))}\cdot K^5 = 
{\kappa^{-1}}_{((A(BC))D)E}^{1,a\times 1}\circ {\kappa^{-1}}_{(A((BC)D))E}^{1,(1\times a)\times 1} \circ {\kappa^{-1}}_{(A(B(CD)))E}^{1,a} \circ {\kappa^{-1}}_{A((B(CD))E))}^{1,1\times a},\]
\normalsize 
\noindent and
\[ \pi^{(((AB)C)(DE))((((AB)C)D)E)}_{(AB)(C(DE))}\cdot K^5 = 1.\] 
\endproposition

\proof
It is straightforward to check that the equation of $2$-cells holds so that we can apply the universal property, which yields the component $2$-cell isomorphism.  Checking the perturbation axiom is then a straightforward calculation.
\endproof

\subsubsection*{$U_{4,1}$ Perturbation}

\proposition \label{U41perturbation}
There is a perturbation consisting of, for each triple $A, B, C$ of objects in $\Span(\B)$, an isomorphism of maps of spans
\[
   \xy
   (-30,3)*+{U^{4,1}_{ABC}\maps (a^{-1}\cdot\kappa\cdot (1*l)(\Pi *1),\; {\alpha_{\mu^{\cdot}}}_{\lambda,1,1}(1*l)(\Pi *1),\; \kappa^{-1}\cdot\pi(\Pi *1))}="1"; (-54,-3)*+{}="2"; (-6,-3)*+{(1,\; ((l\times 1)*1)(1*l*1)(1^2*{\lambda_{\alpha}}),\; \widetilde{\pi}\cdot\kappa^{-1}\cdot\pi)}="3";
        {\ar@3{->}_{} "2";"3"};
\endxy
\]
\noindent defined by the unique $2$-cell
\[ U^{4,1}_{ABC}\maps  {\alpha_{\mu^{\cdot}}}_{\lambda,1,1}(1*l)(\Pi *1)\To ((l\times 1)*1)(1*l*1)(1^2*{\lambda_{\alpha}})\]
\noindent in $\B$ such that
\[\pi_{((1A)B)C}^{(AB)C}\cdot U_{4,1} = {\kappa^{-1}}_{((1A)B)C}^{1,a\times 1}\circ {\kappa^{-1}}_{(1(AB))C}^{1,a}\]
\noindent and
\[ \pi^{((1A)B)C}_{(AB)C}\cdot U_{4,1} = {\kappa}_{1((AB)C)}^{1,a\times 1}.\] 
\endproposition

\proof
The equation of $2$-cells in $\B$
\[ (1\cdot  (\pi^{((1A)B)C}_{(AB)C}\cdot U_{4,1}))(\kappa_{(AB)C}^{1,(\lambda_{A}\times 1)\times 1}\cdot  (\alpha^{-1}_{\lambda,1,1}(1*l)(\Pi *1))) =\]
\[ (\kappa_{(AB)C}^{1,(\lambda_{A}\times 1)\times 1}\cdot ((((l\times 1)*1)(1*l*1)(1^2*\lambda_\alpha))((\lambda_{A}\times 1_{B})\times 1_{C}))\cdot (\pi_{((1A)B)C}^{(AB)C}\cdot U_{4,1}))\]
\noindent allows us to apply the universal property to obtain the the $2$-cell $U_{4,1}$.  The perturbation axiom is satisfied since the components of the relevant modifications are identities.
\endproof

\subsubsection*{$U_{4,2}$ Perturbation}

\proposition \label{U42perturbation}
There is a perturbation consisting of, for each triple $A, B, C$ of objects in $\Span(\B)$, an isomorphism of maps of spans
\[
   \xy
   (-30,3)*+{U^{4,2}_{ABC}\maps (a^{-1}\cdot\kappa\cdot (1*m)(\Pi*1),\; {\alpha_{\mu^{\cdot}}}_{\rho,1,1}(1*m)(\Pi*1),\; (1\times \widetilde{\pi})\cdot\kappa^{-1}\cdot\pi(\Pi*1))}="1"; 
   (-85,-3)*+{}="2"; 
   (-77,-3)*+{}="4"; 
   (-26,-3)*+{(1,\; ((m\times 1)*1)(1*\alpha^{-1}_{1,\lambda,1})(1^2*(1\times l)),\; (1\times \kappa^{-1})\cdot\pi)}="3";
        {\ar@3{->}_{} "2";"4"};
\endxy
\]

\noindent defined by the unique $2$-cell
\[ U^{4,2}_{ABC}\maps  \alpha^{-1}_{\rho,1,1}(1*m)(\Pi*1) \To ((m\times 1)*1)(1*\alpha^{-1}_{1,\lambda,1})(1^2*(1\times l))\]
\noindent in $\B$ such that
\[\pi_{((A1)B)C}^{(AB)C}\cdot U_{4,2} = {\kappa^{-1}}_{((A1)B)C}^{1,a\times 1}\circ {\kappa^{-1}}_{(A(1B))C}^{1,a}\]
\noindent and
\[ \pi^{((AB)1)C}_{(AB)C}\cdot U_{4,2} = 1.\] 
\endproposition

\proof
The equation of $2$-cells in $\B$
\[ (1\cdot  (\pi^{((1A)B)C}_{(AB)C}\cdot U_{4,2}))(\kappa_{(AB)C}^{1,(\lambda_{A}\times 1)\times 1}\cdot  (\alpha_{\rho,1,1}(1*m)(1*\Pi))) =\]
\[ (\kappa_{(AB)C}^{1,(\lambda_{A}\times 1)\times 1}\cdot (((m\times 1)*1)(1*\alpha^{-1}_{1,\lambda,1})(1^2*(1\times l)))((\lambda_{A}\times 1_{B})\times 1_{C}))\cdot (\pi_{((1A)B)C}^{(AB)C}\cdot U_{4,2}))\]
\noindent allows us to apply the universal property to obtain the the $2$-cell $U_{4,2}$.  The perturbation axiom is satisfied since the components of the relevant modifications are identities.
\endproof

\subsubsection*{$U_{4,3}$ Perturbation}

\proposition \label{U43perturbation}
There is a perturbation consisting of, for each triple $A, B, C$ of objects in $\Span(\B)$, an isomorphism of maps of spans
\[
   \xy
   (-23,0)*+{U^{4,3}_{ABC}\maps (a^{-1}\cdot\kappa\cdot (\Pi *1),\; (m*1)(1*{\alpha_{\mu^{\cdot}}}_{(1,1,\rho)})(\Pi*1),\; 1)}="1"; 
   (-97,-6)*+{}="3"; 
   (-89,-6)*+{}="4"; 
   (-25,-6)*+{(1,\; ((l\times 1)*1)(1* {\alpha_{\mu^{\cdot}}}_{1,\lambda,1})(1*(1*(1\times m))),\; (1\times ((1\times \widetilde{\pi})\cdot\kappa^{-1}))\cdot\pi)}="2";
        {\ar@3{->}_{} "3";"4"};
\endxy
\]
\noindent defined by the unique $2$-cell
\[ U^{4,3}_{ABC}\maps  (m*1)(1*{\alpha_{\mu^{\cdot}}}_{(1,1,\rho)})(\Pi*1)\To ((l\times 1)*1)(1* {\alpha_{\mu^{\cdot}}}_{1,\lambda,1})(1*(1*(1\times m)))\]
\noindent in $\B$ such that
\[\pi_{((AB)1)C}^{(AB)C}\cdot U_{4,3} = 1\]
\noindent and
\[ \pi^{((AB)1)C}_{(AB)C}\cdot U_{4,3} = {\kappa^{-1}}_{A((B1)C)}^{1,1\times\alpha}.\] 
\endproposition

\proof
The equation of $2$-cells in $\B$
\[ (1\cdot  (\pi^{((AB)1)C}_{(AB)C}\cdot U_{4,3}))(\kappa_{(AB)C}^{1,\rho_{AB}\times 1}\cdot  (m*1)(1*\alpha^{-1}_{1,1,\lambda})(\Pi*1)) =\]
\[ (\kappa_{(AB)C}^{1,\rho_{AB}\times 1}\cdot (((r\times 1)*1)(1* \alpha^{-1}_{1,\lambda,1})(1^2*(1\times m))))((\rho_{AB}\times 1_{C})\cdot (\pi^{((AB)1)C}_{(AB)C}\cdot U_{4,3})).\]
\noindent allows us to apply the universal property to obtain the the $2$-cell $U_{4,3}$.  The perturbation axiom is satisfied since the components of the relevant modifications are identities.
\endproof

\subsubsection*{$U_{4,4}$ Perturbation}

\proposition \label{U44perturbation}
There is a perturbation consisting of, for each triple $A, B, C$ of objects in $\Span(\B)$, an isomorphism of maps of spans
\[
\def\objectstyle{\scriptstyle}
  \def\labelstyle{\scriptstyle}
   \xy
   (-23,0)*+{U^{4,4}_{ABC}\maps (1,\; (r*1)(1*{\alpha_{\mu^{\cdot}}}_{1,1,\rho})(\Pi *1),\; 1)}="1"; 
   (12,0)*+{}="3";
   (-15,-8)*+{(\widetilde{\pi}\cdot\kappa^{-1}\cdot(1*\rho)(1^2*(1\times r)))(\widetilde{\pi}\cdot\kappa^{-1}\cdot\pi(1^2*(1\times r)))(((1\times ((1\times\rho)\cdot\kappa^{-1})\cdot\pi))(1,\;\rho_{\alpha}(1*\rho)(1^2*(1\times r)),}="2";
        {\ar@3{->}_{} "1";"3"};
\endxy
\]

\noindent defined by the unique $2$-cell
\[ U^{4,4}_{ABC}\maps  \lambda(1*l)(1^2*(1\times l))((1\times 1)\times \lambda_{(AB)C})\To \lambda^{-1}(1*l)(1^2*(1\times l))\]
\noindent in $\B$ such that
\[\pi_{((AB)C)1}^{(AB)C}\cdot U_{4,4} = 1\]
\noindent and
\[ \pi^{((AB)C)1}_{(AB)C}\cdot U_{4,4} = {\kappa^{-1}}_{A((B1)C)}^{1,1\times a}.\] 
\endproposition

\proof
The equation of $2$-cells in $\B$
\[ (1\cdot  (\pi^{((AB)C)1}_{(AB)C}\cdot U_{4,4}))(\kappa_{(AB)C}^{1,\rho_{(AB)C}}\cdot (r*1)(1*\alpha^{-1}_{1,1,\rho})(\Pi *1)) =\]
\[ ((\kappa_{(AB)C}^{1,\rho_{(AB)C}}\cdot \rho_{\alpha}(1*\rho)(1^2*(1\times r))((1\times 1)\times \lambda_{(AB)C})\cdot (\pi_{((AB)C)1}^{(AB)C}\cdot U_{4,4}))\]
\noindent allows us to apply the universal property to obtain the the $2$-cell $U_{4,4}$.  The perturbation axiom is satisfied since the components of the relevant modifications are identities.
\endproof

\subsection{Monoidal Tricategory Axioms}

There are four axioms that a monoidal tricategory must satisfy.  These were written down by Trimble in diagram form over the span of forty-six pages~\cite{Tr}.  These diagrams are reproduced in Definition~\ref{tetracategorydefinition} along with clarifications and explanations of the geometric $2$- and $3$-cells comprising the diagrams.  Checking the axioms is not terribly difficult in our construction; however, certain ``twistings" sometimes make opaque the constituent cells of the geometric $3$-cell domains and codomains.

\subsubsection*{The $K_{6}$ Associativity Axiom}

\proposition \label{K6Axiom}
The $K_{6}$ axiom for the monoidal structure on $\Span(\B)$ holds.
\endproposition

\proof

We have

\[ \pi_{((((AB)C)D)E)F}\cdot (1*K_{5})(1*\alpha *1)(1*K_{5}*1)(1*\Pi*1)(1*\alpha*1)(1*K_{5}) = \]
\[ 1^2 \circ {\kappa_{((A((BC)D))E)F}^{1,((1\times a)\times 1)\times 1}} \circ 1^2 \circ \left({\kappa_{(((A(BC))D)E)F}^{1,(a\times 1)\times 1}} \circ {\kappa_{((((AB)C)D)E)F}^{1,((a\times 1)\times 1)\times 1}}\right) \]
\noindent and

\[ \pi_{((((AB)C)D)E)F}\cdot (1*\Pi*1)(1*(1\times K_{5}*1)(1*\alpha*1)\]\[(1*K_{5}*1)(1*\Pi*1)(1*\Pi*1)(1*K_{5}*1)((K_{5}\times 1)*1) = \]
\[ 1^3 \circ \left({\kappa_{((A((BC)D))E)F}^{1,((1\times a)\times 1)\times 1}} \circ {\kappa_{(((A(BC))D)E)F}^{1,(a\times 1)\times 1}}\right) \circ  1^3\circ {\kappa_{((((AB)C)D)E)F}^{1,((a\times 1)\times 1)\times 1}} \]
\noindent up to whiskering by associator $1$-cells $a$ and projection $1$-cells $\widetilde{\pi}_{A}^{1}$.

We have

\[ \pi_{(((AB)C)D)(EF)}\cdot (1*K_{5})(1*\alpha *1)(1*K_{5}*1)(1*\Pi*1)(1*\alpha*1)(1*K_{5}) = \]
\[ \left.\circ {\kappa^{-1}}_{(A(B((CD)E)))F}^{1,(1\times (1\times a))\times 1} \circ {\kappa^{-1}}_{(A(B(C(DE))))F)}^{1,a} \circ {\kappa^{-1}}_{A((B(C(DE)))F)}^{1,1\times (a\times 1)} \circ {\kappa^{-1}}_{A(B((C(DE))F))}^{1,1\times (1\times a)}\right) \circ 1^2 \]
\[ 1 \circ {\kappa_{((A((BC)D))E)F}^{1,((1\times a)\times 1)\times 1}} \circ 1 \left(\circ {\kappa^{-1}}_{((A((BC)D))E)F}^{1,((1\times a)\times 1)\times 1}
\circ {\kappa^{-1}}_{((A(B(CD))E)F}^{1,a\times 1}
\circ {\kappa^{-1}}_{(A((B(CD))E)F}^{1,a\times 1}\right.\]

\noindent and

\[  \pi_{(((AB)C)D)(EF)}\cdot (1*\Pi*1)(1*(1\times K_{5}*1)(1*\alpha*1)(1*K_{5}*1)(1*\Pi*1)(1*\Pi*1)(1*K_{5}*1)((K_{5}\times 1)*1)  =\]
\[\left. \circ {\kappa^{-1}}_{(A(B(C(DE))))F)}^{1,a} 
\circ {\kappa^{-1}}_{A((B(C(DE)))F)}^{1,1\times (a\times 1)}
 \circ {\kappa^{-1}}_{A(B((C(DE))F))}^{1,1\times (1\times a)}\right)
\circ 1^4 \]
\[ 1^3 \circ \left({\kappa^{-1}}_{((A(B(CD))E)F}^{1,a\times 1}
\circ {\kappa^{-1}}_{(A((B(CD))E)F}^{1,a\times 1}
\circ {\kappa^{-1}}_{(A(B((CD)E)))F}^{1,(1\times (1\times a))\times 1}\right.\]

\noindent up to whiskering by associator $1$-cells $a$ and projection $1$-cells $\widetilde{\pi}_{A}^{1}$.

We have

\[ \pi_{((AB)C)(D(EF))}\cdot (1*K_{5})(1*\alpha *1)(1*K_{5}*1)(1*\Pi*1)(1*\alpha*1)(1*K_{5}) = \]
\[ \left. \circ {\kappa^{-1}}_{(A(B(C(DE))))F)}^{1,a} \circ {\kappa^{-1}}_{A((B(C(DE)))F)}^{1,1\times (a\times 1)} \circ {\kappa^{-1}}_{A(B((C(DE))F))}^{1,1\times (1\times a)}\right)   \circ 1^5 \]
\[ \left({\kappa^{-1}}_{((A(B(CD))E)F}^{1,a\times 1}
\circ {\kappa^{-1}}_{(A((B(CD))E)F}^{1,a\times 1} \circ {\kappa^{-1}}_{(A(B((CD)E)))F}^{1,(1\times (1\times a))\times 1}\right.\]

\noindent and

\[ \pi_{((AB)C)(D(EF))}\cdot (1*\Pi*1)(1*(1\times K_{5}*1)(1*\alpha*1)(1*K_{5}*1)(1*\Pi*1)(1*\Pi*1)(1*K_{5}*1)((K_{5}\times 1)*1) = \]
\[ \left.\circ {\kappa^{-1}}_{(A(B(C(DE))))F)}^{1,a} \circ {\kappa^{-1}}_{A((B(C(DE)))F)}^{1,1\times (a\times 1)}\circ {\kappa^{-1}}_{A(B((C(DE))F))}^{1,1\times (1\times a)}\right) \circ 1^4 \]
\[ 1^3 \circ \left({\kappa^{-1}}_{((A(B(CD))E)F}^{1,a\times 1}  \circ {\kappa^{-1}}_{(A(B((CD))E))F}^{1,(1\times a)\times 1} \circ {\kappa^{-1}}_{(A(B((CD)E)))F}^{1,(1\times (1\times a))\times 1}\right.\]

\noindent up to whiskering by associator $1$-cells $a$ and projection $1$-cells $\widetilde{\pi}_{A}^{1}$.

We have

\[ \pi_{(AB)(C(D(EF)))}\cdot (1*K_{5})(1*\alpha *1)(1*K_{5}*1)(1*\Pi*1)(1*\alpha*1)(1*K_{5}) = 1^6\]
\noindent and
\[  \pi_{(AB)(C(D(EF)))}\cdot (1*\Pi*1)(1*(1\times K_{5}*1)(1*\alpha*1)(1*K_{5}*1)(1*\Pi*1)(1*\Pi*1)(1*K_{5}*1)((K_{5}\times 1)*1) = 1^8\]

\noindent up to whiskering by associator $1$-cells $a$ and projection $1$-cells $\widetilde{\pi}_{A}^{1}$.

The axiom then holds by the universal property.
\endproof

\subsubsection*{The $U_{5,2}$ Unit Axiom}

\proposition \label{U52Axiom}
The $U_{5,2}$ axiom for the monoidal structure on $\Span(\B)$ holds.
\endproposition

\proof

It is straightforward to verify the axiom.   We have

\[  \pi_{(((A1)B)C)D}\cdot (1*(1\times U_{4,1}))(1*\Pi *1)(1*\alpha)(1*U_{4,2}*1)(K_{5}*1) = \]
\[ 1^2\circ {\kappa^{-1}}_{(((A1)B)C)D}^{1,(a\times 1)\times 1} \circ \left({\kappa^{-1}}_{((A(1B))C)D}^{1,a\times 1} \circ {\kappa^{-1}}_{(A((1B)C))D}^{1,a}\circ {\kappa^{-1}}_{A(((1B)C)D)}^{1,(1\times a)\times 1}\circ {\kappa^{-1}}_{A((1(BC))D)}^{1,1\times a}\right)\circ 1\]
\noindent and

\[ \pi_{(((A1)B)C)D}\cdot (1*\alpha *1)(1*(U_{4,1}\times 1)*1)(1* U_{4,2} *1)(1* m *1)(1*\Pi) = \]
\[  1\circ \left({\kappa^{-1}}_{(((A1)B)C)D}^{1,(a\times 1)\times 1} \circ {\kappa^{-1}}_{((A(1B))C)D}^{1,a\times 1}\right) \circ \left({\kappa^{-1}}_{(A((1B)C))D}^{1,a} \circ {\kappa^{-1}}_{A(((1B)C)D)}^{1,(1\times a)\times 1}\right)\circ {\kappa^{-1}}_{A((1(BC))D)}^{1,1\times a}\circ 1\]
\noindent up to whiskering by structural $1$-cells.

We have

\[  \pi_{((AB)C)D}\cdot (1*(1\times U_{4,1}))(1*\Pi *1)(1*\alpha)(1*U_{4,2}*1)(K_{5}*1) = \]
\[ \circ \left({\kappa^{-1}}_{((A(1B))C)D}^{1,a\times 1}
\circ {\kappa^{-1}}_{(A((1B)C))D}^{1,a} \circ {\kappa^{-1}}_{A(((1B)C)D)}^{1,(1\times a)\times 1}
\circ {\kappa^{-1}}_{A((1(BC))D)}^{1,1\times a}\right) \circ 1\]
\[  \left({\kappa_{A(((1B)C)D)}^{1,(1\times a)\times 1}}^{-1}\circ {\kappa_{A((1(BC))D)}^{1,1\times a}}^{-1}\right)\circ 1 \circ \left({\kappa_{A((1(BC))D)}^{1,1\times a}} \circ {\kappa_{A(((1B)C)D)}^{1,(1\times a)\times 1}} \circ {\kappa_{(A((1B)C))D}^{1,a}} \circ {\kappa_{((A(1B))C)D}^{1,a\times 1}} \right)  \]
\[ \left({\kappa_{A(((1B)C)D)}^{1,(1\times a)\times 1}}^{-1}\circ {\kappa_{A((1(BC))D)}^{1,1\times a}}^{-1}\right)\circ 1^2 = \]

\noindent and

\[ \pi_{((AB)C)D}\cdot (1*\alpha *1)(1*(U_{4,1}\times 1)*1)(1* U_{4,2} *1)(1* m *1)(1*\Pi) = \]
\[ 1^2 \circ {\kappa^{-1}}_{A(((1B)C)D)}^{1,(1\times a)\times 1} \circ {\kappa^{-1}}_{A((1(BC))D)}^{1,1\times a}\circ 1\]
\noindent up to whiskering by structural $1$-cells.

We have

\[ \pi_{(AB)(CD)}\cdot (1*(1\times U_{4,1}))(1*\Pi *1)(1*\alpha)(1*U_{4,2}*1)(K_{5}*1) = 1^5\]
\noindent and
\[ \pi_{(AB)(CD)}\cdot (1*\alpha *1)(1*(U_{4,1}\times 1)*1)(1* U_{4,2} *1)(1* m *1)(1*\Pi) = 1^5\]
\noindent up to whiskering by structural $1$-cells.
The axiom then holds by the universal property.
\endproof

\proposition \label{U53Axiom}
The $U_{5,3}$ axiom for the monoidal structure on $\Span(\B)$ holds.
\endproposition

\proof

We have

\[ \pi_{(((AB)1)C)D}\cdot (1*\alpha *1)(1*(U_{4,3}\times 1)*1)(1*(1\times\Pi)*1)(1*\alpha *1)(1*(1\times U_{4,2})) = \]
\[ 1^4\circ \left({\kappa^{-1}}_{(((AB)1)C)D}^{1,(a\times 1)\times 1} \circ {\kappa^{-1}}_{((A(B1))C)D}^{1,a\times 1} \circ {\kappa^{-1}}_{(A((B1)C))D}^{1,a}\circ {\kappa^{-1}}_{A(((B1)C)D)}^{1,(1\times a)\times 1}\circ {\kappa^{-1}}_{A((B(1C))D)}^{1,1\times a}\right)\]
\noindent and

\[ \pi_{(((AB)1)C)D}\cdot ((1\times U_{4,2})*1)(1*\Pi)(1* \alpha)(1* U_{4,3} *1)(K_{5}*1) = \]
\[ 1^2\circ  {\kappa^{-1}}_{(((AB)1)C)D}^{1,(a\times 1)\times 1} \circ 1\circ \left({\kappa^{-1}}_{((A(B1))C)D}^{1,a\times 1} \circ {\kappa^{-1}}_{(A((B1)C))D}^{1,a}\circ {\kappa^{-1}}_{A(((B1)C)D)}^{1,(1\times a)\times 1}\circ {\kappa^{-1}}_{A((B(1C))D)}^{1,1\times a}\right)\]
\noindent up to whiskering by associator $1$-cells $a$ and projection $1$-cells $\widetilde{\pi}_{A}^{1}$.

We have

\[ \pi_{((AB)C)D}\cdot (1*\alpha *1)(1*(U_{4,3}\times 1)*1)(1*(1\times\Pi)*1)(1*\alpha *1)(1*(1\times U_{4,2}))) = \]
\[ 1\circ \left({\kappa^{-1}}_{(((AB)1)C)D}^{1,(a\times 1)\times 1} \circ {\kappa^{-1}}_{((A(B1))C)D}^{1,a\times 1} \circ {\kappa^{-1}}_{(A((B1)C))D}^{1,a}\circ {\kappa^{-1}}_{A(((B1)C)D)}^{1,(1\times a)\times 1}\right)  \circ {\kappa^{-1}}_{A((B(1C))D)}^{1,1\times a}\circ 1^2\]
\noindent and

\[ \pi_{((AB)C)D}\cdot ((1\times U_{4,2})*1)(1*\Pi)(1* \alpha)(1* U_{4,3} *1)(K_{5}*1)  = \]
\[  1^2\circ {\kappa^{-1}}_{(((AB)1)C)D}^{1,(a\times 1)\times 1} \circ 1\circ \left({\kappa^{-1}}_{((A(B1))C)D}^{1,a\times 1} \circ {\kappa^{-1}}_{(A((B1)C))D}^{1,a}\circ {\kappa^{-1}}_{A(((B1)C)D)}^{1,(1\times a)\times 1}\circ {\kappa^{-1}}_{A((B(1C))D)}^{1,1\times a}\right)\]
\noindent up to whiskering by associator $1$-cells $a$ and projection $1$-cells $\widetilde{\pi}_{A}^{1}$.

We have

\[ \pi_{(AB)(CD)}\cdot (1*\alpha *1)(1*(U_{4,3}\times 1)*1)(1*(1\times\Pi)*1)(1*\alpha *1)(1*(1\times U_{4,2})) = \]
\[ 1^3 \circ {\kappa^{-1}}_{A(B((1C)D))}^{1,1\times (1\times a)} \circ 1\]
\noindent and

\[ \pi_{(AB)(CD)}\cdot ((1\times U_{4,2})*1)(1*\Pi)(1* \alpha)(1* U_{4,3} *1)(K_{5}*1) = \]
\[ 1^3 \circ {\kappa^{-1}}_{A(B((1C)D))}^{1,1\times (1\times a)} \circ 1\]
\noindent up to whiskering by associator $1$-cells $a$ and projection $1$-cells $\widetilde{\pi}_{A}^{1}$.

The axiom then holds by the universal property.
\endproof

\proposition \label{U54Axiom}
The $U_{5,4}$ axiom for the monoidal structure on $\Span(\B)$ holds.
\endproposition

\proof

We have

\[ \pi_{(((AB)C)1)D}\cdot (1*\alpha *1)(1*(U_{4,4}\times 1)*1)(1*\Pi*1)(1*\alpha *1)(1*U_{4,3}) = 1^5 \]
\noindent and

\[ \pi_{(((AB)C)1)D}\cdot ((1\times U_{4,3})*1)(1*U_{4,3}*1)(1* m *1)(1* \Pi *1)(K_{5}*1) = \]
\[ 1^2 \circ {\kappa^{-1}}_{(((AB)C)1)D}^{1,(a\times 1)\times 1} \circ 1 \circ {\kappa_{(((AB)C)1)D}^{1,(a\times 1)\times 1}} \]
\noindent up to whiskering by associator $1$-cells $a$ and projection $1$-cells $\widetilde{\pi}_{A}^{1}$.

We have

\[ \pi_{((AB)C)D}\cdot (1*\alpha *1)(1*(U_{4,4}\times 1)*1)(1*\Pi*1)(1*\alpha *1)(1*U_{4,3}) = \]

\[ 1\circ \left({\kappa^{-1}}_{(((AB)C)1)D}^{1,(a\times 1)\times 1} \circ {\kappa^{-1}}_{((A(BC))1)D}^{1,a\times 1} \circ {\kappa^{-1}}_{(A((BC)1))D}^{1,a}\circ {\kappa^{-1}}_{A(((BC)1)D)}^{1,(1\times a)\times 1} \circ {\kappa^{-1}}_{A((B(C1))D)}^{1,1\times a}\right) \circ 1^3\]
\noindent and

\[ \pi_{((AB)C)D}\cdot ((1\times U_{4,3})*1)(1*U_{4,3}*1)(1* m *1)(1* \Pi *1)(K_{5}*1)  = \]
\[  1^2\circ \left({\kappa^{-1}}_{(((AB)C)1)D}^{1,(a\times 1)\times 1} \circ {\kappa^{-1}}_{((A(BC))1)D}^{1,a\times 1} \circ {\kappa^{-1}}_{(A((BC)1))D}^{1,a}\circ {\kappa^{-1}}_{A(((BC)1)D)}^{1,(1\times a)\times 1}\circ {\kappa^{-1}}_{A((B(C1))D)}^{1,1\times a}\right) \circ 1^2\]
\noindent up to whiskering by associator $1$-cells $a$ and projection $1$-cells $\widetilde{\pi}_{A}^{1}$.

We have

\[ \pi_{(AB)(CD)}\cdot (1*\alpha *1)(1*(U_{4,4}\times 1)*1)(1*\Pi*1)(1*\alpha *1)(1*U_{4,3}) = \]
\[ 1^3 \circ {\kappa^{-1}}_{A(B((C1)D))}^{1,1\times (1\times a)} \circ 1\]
\noindent and

\[ \pi_{(AB)(CD)}\cdot ((1\times U_{4,3})*1)(1*U_{4,3}*1)(1* m *1)(1* \Pi *1)(K_{5}*1) = \]
\[ {\kappa^{-1}}_{A(B((C1)D))}^{1,1\times (1\times a)} \circ 1^4\]
\noindent up to whiskering by associator $1$-cells $a$ and projection $1$-cells $\widetilde{\pi}_{A}^{1}$.

The axiom then holds by the universal property.
\endproof

\appendix
\section{Trimble's Remarks}\label{appendix}

It seems appropriate to preserve, along with the tetracategory definition, Trimble's original remarks accompanying his work.  The following section taken from John Baez's website~\cite{Tr} are Trimble's remarks on his motivation, techniques, and hurdles in writing down the definition.

\begin{quote}
In $1995$, at Ross Street's request, I gave a very explicit description of weak $4$-categories, or tetracategories as I called them then, in terms of nuts-and-bolts pasting diagrams, taking advantage of methods I was trying to develop then into a working definition of weak $n$-category. Over the years various people have expressed interest in seeing what these diagrams look like -- for a while they achieved a certain notoriety among the few people who have actually laid eyes on them (Ross Street and John Power may still have copies of my diagrams, and on occasion have pulled them out for visitors to look at, mostly for entertainment I think).

Despite their notorious complexity, there seems to be some interest in having these diagrams publicly available; John Baez has graciously offered to put them on his website. The theory of (weak) n-categories has come a long way since the time I first drew these diagrams up, and it may be wondered what there is to gain by re-examining the apparently primitive approach I was then taking. I don't have any definitive answer to that, but I do know that there are interesting combinatorics lurking within the methods I was using, which I suspect may hold clues to certain coherence-theoretic aspects of weak n-categories; therefore in these notes I wanted to give some idea of the methods used to construct these diagrams (as well as to provide explanatory notes in case the diagrams are too baffling), in case anyone wants to push them further.

Disclaimer: there is an ``operadic" definition of weak $n$-category, described by me in a talk at Cambridge in $1999$, and subsequently written up by Tom Leinster, Eugenia Cheng and Aaron Lauda. The notion of tetracategory given here is not the operadic one in dimension $4$. The truth be told, I have never considered the operadic version as anything other than a somewhat strictified case of the ``true" notion of n-category; it was invented primarily to give a precise meaning to Grothendieck's ``fundamental n-groupoid of a topological space", and is conjectured to be an algebraic characterization of the n-type of a space. (The operadic definition does have the advantage of being short, almost a ``two-line definition", and it does embody a great deal of structure one expects for weak n-categories, which may make it appealing for didactic purposes. There are also some conjectures about fundamental $n$-groupoids which I consider important and which, to my knowledge, have not been completely addressed: one is that they do characterize $n$-types; another is a van Kampen theorem.)

There is a strong connection with operads, however: as already noticed by a number of people in the early $90$'s, there is a strong resemblance between the ``higher associativity conditions" in the notions of bicategories and tricategories, and the convex polyhedra $K_n$ called ``associahedra", introduced $30$ years earlier by Stasheff in his work on the homotopy types of loop spaces, and which together form an operad, whose formal definition was given later by May.

My early efforts to define weak n-categories (in $1994-1995$) were really just an attempt to take this visual resemblance seriously and build a theory around it. Thus, the higher associativities could all be described by the polytopes $K_n$, outfitted with orientations or directions on each cellular face $t$, which divide the boundary of t into negative and positive parts. The negative boundary cells are oriented in such a way that they paste together sensibly (according to one or another formal framework, e.g., Street's parity complexes, Mike Johnson's pasting schemes, etc.) to form the ``domain" of $t$; similarly the positive boundary cells paste together to form the codomain of $t$.

In addition to the associahedra used for the higher associativities, one needs some allied polyhedra for the higher unit conditions; in short, I needed a suitable collection of ``monoidahedra" to capture the combinatorial structure of the $n$-category data and axioms. This turned out to be easy; a minor tweak on the machine used to define associahedra does the trick. There are also ``functoriahedra" for describing the data and equations occurring for (weak) $n$-functors. All these polyhedra are in fact describable by means of bar constructions, as I've tried to describe informally to a few people, and especially to Baez. (In part this has long been known: the ubiquity of abstract bar constructions probably first hit me on the head while reading May's Geometry of Iterated Loop Spaces [where operads were first introduced]. The specific application to functoriahedra is, I would guess, not generally known.)

The next step is more technical: putting parity complex structures on these polyhedra. Here I received vital input from Dominic Verity: using some ideas from the ``surface diagrams" he and Street had begun working on, together with a general position argument, we had a theorem which guarantees that such parity complex structures exist. (I think it doesn't matter which parity structure is chosen, i.e. that one should get the same notion of n-category regardless, or that different choices of parity structure correspond to different presentations of the same theory.)

In practice, it would be nice to have a systematic way of choosing parity structure to make the pictures as pretty as possible, but for that I have just a handful of tricks and ideas, without any general theory. I believe Jack Duskin told me he had a nice way of getting parity structures on general associahedra, but I don't have details on that.

The final step in my $1995$ attempt to define $n$-categories was to interpret these pasting schemes as cells in local $\hom(n-1)$-categories. Here is where the project got bogged down. With the wealth of $n$-categorical theory which has appeared in the meantime, it should be possible, I think, to fix this part of the program (perhaps by adapting Penon's definition?). At the time, it seemed as if I needed an ($n-1$)-categorical coherence theorem in order to interpret these pastings, which I didn't have of course. We did have coherence of tricategories, due to Gordon-Power-Street (GPS), and this enabled me to give a rigorous definition of tetracategory.

I still feel that despite the technical difficulties which blocked my progress (some of which may be spurious), the program may be worth resurrecting, because in fact in dimensions $2$ and $3$ this approach does yield exactly the classical notions of bicategory and tricategory, without further tweakings or encryptions, and also because the underlying ideas and methods are really quite simple -- a fact which is obscured by the dazzling complexity of the diagrams which were produced at Street's request.  The fact is that these diagrams were produced merely by patiently turning the crank of a simple machine -- anyone with time to kill can produce these things automatically. The general machine, and what one might be able to say with its help, was the real point behind my failed attempts back in $1995$.
\end{quote}

\section{Limits in $2$-Categories}\label{sec-Limits}

We give a short exposition on $2$-limits for the readers benefit and reference.  Limits can be described by {\em cones}, {\em representable functors}, or {\em adjoint functors}.  We will discuss limits simultaneously as a mixture of the first two.  The adjoint functor definition is the most succinct, however, we find cones to provide the most intuitive approach and representable functors to be most amenable to generalizations from $1$-dimensional limits to our setting of weighted and $2$-dimensional limits.  In this section we will describe conical limits, weighted limits, strict $2$-limits, and pseudo limits, looking to `pullbacks' and `products', for examples.

Intuitively, we think of a limit as the `best' solution to a particular problem.  The problem is as follows.  Suppose we are given an indexing category $D$, which we will call the {\em shape} of the limit, and a functor $F\maps D\to\mathcal{C}$, which we will call a {\em diagram of shape $D$}.  Our first task is to find an object, $``\Lim\; F"\in\mathcal{C}$, along with a collection of morphisms $\{g_d\maps ``\Lim\; F"\to F(d)\}$, called a {\em cone}, such that for any morphism $\delta\maps d\to d'$ in $D$ the triangle
\[
   \xy
   (0,20)*+{``\Lim\; F"}="1";
   (-20,0)*+{F(d)}="3";
   (20,0)*+{F(d')}="5";
        {\ar_{g_{d}} "1";"3"};
        {\ar^{g_{d'}} "1";"5"};
        {\ar_{F(\delta)} "3";"5"};
\endxy
\]
\noindent commutes.  Notice the quotes around our potential solution.  We mean to signify that this is a possible solution, and when we find the best possible solution, in the sense which we will now describe, then we will remove the quotes and say that we have found {\em the} limit.

A solution is considered to be the best solution, and thus, the limit, if it satisfies a certain {\em universal property}.  We say a solution ($\Lim\; F,\{f_d\}$) is the limit if, for every other possible solution ($``\Lim\; F",\{g_d\}$), there exists a unique comparison morphism $h\maps ``\Lim\; F"\to \Lim\; F$ in $\mathcal{C}$ such that for each $d\in D$, the triangle
\[
   \xy
   (0,-20)*+{F(d)}="1";
   (-20,0)*+{``\Lim\; F"}="3";
   (20,0)*+{ \Lim\; F}="5";
        {\ar_{g_{d}} "3";"1"};
        {\ar^{f_{d}} "5";"1"};
        {\ar^{h} "3";"5"};
\endxy
\]
\noindent commutes. 
Note that our use of the definite article `the' is an allowable abuse of language, so long as we remember that the best solution, or limit, is unique only up to the proper notion of equivalence.  In a category, this notion of equivalence is {\em isomorphism}.

Given this informal definition of limits in terms of families of morphisms called cones, we recall the more abstract definition of limits in terms of {\em representable presheaves}.  The notion of representability is closely related to the famous Yoneda embedding.  Representability is a property of a presheaf, which jumping ahead just for a moment, will correspond to the existence of a limit.  So, we will say a limit exists when a certain presheaf is representable.   However, in practice we often want to work with a specific chosen limit, i.e., a {\em representation} of a presheaf.  Recall, the following definition of a limit of a presheaf $F$ as a pair consisting of a representing object $\Lim\; F$ and a natural isomorphism $\phi$.  Here $\Delta 1\maps D\to \Set$ is the {\em constant point functor} $d\mapsto \{*\}$, which takes every object in $D$ to the singleton set in $\Set$.

\begin{definition}
Let $\mathcal{C},D$ be categories and $F\maps $D$\to\mathcal{C}$ a functor, i.e., a diagram of shape $D$.  A {\bf limit} of $F$ is a pair $(\Lim\; F\in\mathcal{C},\phi)$, where
\[ \phi \maps[D,\Set](\Delta 1(\--),\mathcal{C}(\--,F(\--)))\To \mathcal{C}(\--,\Lim\; F)\]
\noindent is a natural isomorphism in $[\mathcal{C}^{\op},\Set]$.
\end{definition}

We are now in a position to generalize from conical limits (weighted by the constant point functor) in categories to more general weighted and $2$-categorical limits.

\subsection{Strict and Pseudo Limits in $2$-Categories}

Limits in $1$-categories are defined using {\em natural transformations} in the functor category $[D,\Set]$ as cones.  Natural transformations are the $2$-morphisms in the $2$-category of categories, functors, and natural transformations.  In this section we need to generalize our cones to {\em strict transformations} and {\em strong transformations}, which are the $2$-morphisms in the $3$-categories of $2$-categories, homomorphisms, (strict or strong) transformations, and modifications, respectively.  This distinction between strict and strong transformations is the primary difference between strict and pseudo limits.  However, pseudo limits can be imitated by strict limits using weighted limits.  We will demonstrate this in our definition of iso-comma object.

Let $D$ be a (strict) $2$-category, then $[D,\Cat]$ is the $2$-category of (strict) $2$-functors, strict transformations, and modifications.  We denote by $2\Cat(D,\Cat)$ the strict $2$-category (which is strict since the codomain is a strict $2$-category) of homomorphisms, strong transformations, and modifications.  Then ${\rm Ps}(D,\Cat)$ will denote the full sub-$2$-category of (strict) $2$-functors in $2\Cat(D,\Cat)$.

All $2$-categorical limits will be considered in the generalized context of weighted limits.  This means we replace the constant point functor $\Delta 1\maps D\to\Set$ used in defining conical limits with an arbitrary (strict) $2$-functor $W\maps D\to\Cat$, called the {\em weight}, and write $\Lim\;(W,F)$ in place of $\Lim\; F$.  However, many important examples are still conical limits, i.e., limits weighted by the strict $2$-functor $\Delta 1\maps D\to\Cat$ that assigns each object of $D$ the terminal category and all $1$- and $2$-morphisms to the respective identity morphisms.  Note that weighted limits do show up in the context of $1$-categorical limits when extending to limits in enriched category theory.  Of course, strict $2$-categories are $\Cat$-enriched categories, hence the $\Cat$-valued weights in the following definitions.

Strict $2$-limits and pseudo limits have a common property as $2$-categorical limits, which differentiates them from the more general $2$-dimensional limits.  That is, the universal properties are stated as {\em isomorphisms} of categories, rather than {\em equivalences}.  Disregarding some of the nuances suggested above, which arise when considering weighted limits, the main difference between strict and pseudo limits is that for strict limits each cone commutes with the diagram ``on the nose", whereas for a pseudo limit, commutativity holds only up to specified isomorphisms.  This distinction arises by replacing $[D,\Cat]$ in the definition of strict limits with ${\rm Ps}(D,\Cat)$ in the definition of pseudo limits.

We now recall the definitions.

\begin{definition}
Let $\mathcal{C},D$ be (strict) $2$-categories, $F\maps $D$\to\mathcal{C}$ a diagram, and $W\maps D\to\Cat$ a weight, both (strict) $2$-functors.  A {\bf strict} $W$-{\bf weighted $2$-limit} of $F$ is a pair $(\Lim\;(W,F)\in\mathcal{C},\phi)$, where
\[ \phi \maps[D,\Cat](W(\--),\mathcal{C}(\--,F(\--)))\To \mathcal{C}(\--,\Lim\;(W,F))\]
\noindent is an invertible strict transformation in $[\mathcal{C},\Cat]$.
\end{definition}

\begin{definition}
Let $\mathcal{C},D$ be (strict) $2$-categories, $F\maps $D$\to\mathcal{C}$ a diagram, and $W\maps D\to\Cat$ a weight, both (strict) $2$-functors.  A $W$-{\bf weighted pseudo limit} of $F$ is a pair $({\rm pslim}(W,F)\in\mathcal{C},\phi)$, where
\[ \phi \maps {\rm Ps}(D,\Cat)(W(\--),\mathcal{C}(\--,F(\--)))\To \mathcal{C}(\--,{\rm pslim}(W,F))\]
\noindent is an invertible strict transformation in $[\mathcal{C},\Cat]$.
\end{definition}

\noindent The components of $\phi$ are isomorphisms of categories in each of the preceding definitions.

\subsection{Some Finite Limits in $2$-Categories}

\subsubsection{Terminal Objects}

Let $\mathcal{C}$ be a (strict) $2$-category and $D$ be the initial $2$-category, i.e., the $2$-category that has an empty set of objects.  For the unique diagram $F\maps D\to\mathcal{C}$, the limit is an object ${\bf 1}\in\mathcal{C}$ with an invertible strict transformation
\[ \phi\maps [D,\Cat](\Delta 1(\--),\mathcal{C}(\--,F(\--)))\to \mathcal{C}(\--,{\bf 1}).\]

\noindent For every object $A\in\mathcal{C}$, there is an isomorphism of categories
\[ \phi_A\maps [D,\Cat](\Delta 1(\--),\mathcal{C}(A,F(\--)))\to \mathcal{C}(A,{\bf 1}),\]
\noindent where the domain has exactly one transformation (since $D$ has no objects), which gets mapped to a unique morphism $A\to {\bf 1}$ in $\mathcal{C}$.

\noindent It is not difficult to check that strict terminal objects as we have defined them are also pseudo limits.

\subsubsection{Products}

Let $\mathcal{C}$ be a (strict) $2$-category and $D$ the (strict) $2$-category with exactly two objects and only identity $1$- and $2$-morphisms.  For a diagram $F\maps D\to\mathcal{C}$, the limit is an object $\Lim\; F\in\mathcal{C}$ with an invertible strict transformation

\[ \phi\maps [D,\Cat](\Delta 1(\--),\mathcal{C}(\--,F(\--)))\to \mathcal{C}(\--,\Lim\; F).\]

\noindent Note that, in this case, a diagram $F(D)$ is a pair of objects $X,Y\in\mathcal{C}$ with just identity $1$- and $2$-morphisms, so the limiting object $\Lim\; F$ would, in practice, be written as $X\times Y$.\\

\noindent For every object $A\in\mathcal{C}$, there is an isomorphism of categories
\[ \phi_A\maps [D,\Cat](\Delta 1(\--),\mathcal{C}(A,F(\--)))\to \mathcal{C}(A,\Lim\; F).\]

\noindent The pair of $1$-morphisms in the limiting cone
\[ \phi_{\Lim\; F}^{-1}(1_{\Lim\; F})\maps \Delta 1(\--)\To \mathcal{C}(\Lim\; F,F(\--)))\]
\noindent are the ``projections"
\[Y\leftarrow \Lim\; F \rightarrow X.\]

The universal property is straightforward to verify.  

$\bullet$ Since the components of $\phi$ are isomorphisms, for each object $A\in\mathcal{C}$, a cone
\[ \sigma\maps \Delta 1(\--)\to \mathcal{C}(A,F(\--)),\]
which consists of a pair of $1$-morphisms
\[Y\leftarrow A\rightarrow X.\]
is mapped to a unique comparison map 
\[ \phi_A(\sigma)\maps A\to \Lim\; F,\]
and we see that
\[
   \xy
   (0,20)*+{A}="1";
   (-20,0)*+{Y}="3";
   (20,0)*+{X}="5";
   (0,0)*+{\Lim\; F}="7";
        {\ar^{} "1";"3"};
        {\ar^{} "1";"5"};
        {\ar|{\phi_A(\sigma)} "1";"7"};
        {\ar^{} "7";"3"};
        {\ar_{} "7";"5"};
\endxy
\]
commutes by naturality of $\phi$.  This is the one-dimensional aspect of the universal property.

$\bullet$ Now, consider an object $A\in\mathcal{C}$ and a pair of $1$-morphisms $h,k\maps A\to \Lim\; F$ such that for each $d\in D$, there is a $2$-morphism 
\[ \phi_{\Lim\; F}^{-1}(1_{\Lim\; F})_dh\To \phi_{\Lim\; F}^{-1}(1_{\Lim\; F})_dk.\]
Since there are no non-trivial morphisms in $D$, it is immediate that this collection of $2$-morphisms is a modification in $[D,\Cat](\Delta 1(\--),\mathcal{C}(A,F(\--)))$, i.e., a map of cones.  Since $\phi_A$ is an isomorphism of categories, and thus, fully faithful, this modification maps to a unique $2$-morphism $h\To k$ in $\mathcal{C}$.  This is the two-dimensional aspect of the universal property.

\noindent As with terminal objects, strict products as we have defined them are also pseudo limits.

\subsubsection{Pullbacks}

In $2$-categories, pullbacks are more interesting as limits than terminal objects and products.  The conical strict $2$-limit with the usual diagram for pullbacks is not a pseudo limit.   This is closely related to the fact that this limit requires diagrams of $1$-morphisms to commute on the nose, which is generally not a property imposed in a $2$-category.  We usually ask for a diagram to commute up to a $2$-morphism, so we work instead with the pseudo pullback.

We can define the conical strict pullback, which is a representation of the $2$-functor
\[ [D,\Set](\Delta 1(\--),\mathcal{C}(\--,F(\--)))\maps \mathcal{C}^{\op}\to \Cat\]
\noindent where $D$ is the $2$-category
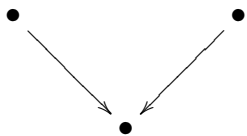
\[
   \xy
   (-15,0)*+{\bullet}="3";
   (15,0)*+{\bullet}="5";
   (0,-15)*+{\bullet}="7";
        {\ar^{} "3";"7"};
        {\ar^{} "5";"7"};
\endxy
\]
\noindent with no non-identity $2$-morphisms.  Unravelling this definition gives the pullback as a $1$-categorical limit with a $2$-categorical universal property.

Since we will not use the strict pullback, we will not give anymore details.  Instead, we will define the pseudo pullback as a conical pseudo limit.   This is a good opportunity to give give some intuition for non-conical limits as well.   Any conical pseudo limit can be defined as a non-conical strict $2$-limit.  We will first define the pseudo pullback as a conical, or $\Delta 1$-weighted, pseudo limit.

Let $\mathcal{C}$ be a (strict) $2$-category and $D$ as above.  For a (strict) $2$-functor, or diagram, $F\maps D\to\mathcal{C}$, the pseudo limit is an object $\Lim\; F\in\mathcal{C}$ with an invertible strict transformation
\[ \phi\maps {\rm Ps}(D,\Cat)(\Delta 1(\--),\mathcal{C}(\--,F(\--)))\to \mathcal{C}(\--,\Lim\; F)\]
\noindent in $[\mathcal{C},\Cat]$.

Recall, a diagram $F(D)$ is a cospan
\[
   \xy
   (-15,0)*+{Y}="3";
   (15,0)*+{X}="5";
   (0,-15)*+{Z}="7";
        {\ar^{} "3";"7"};
        {\ar^{} "5";"7"};
\endxy
\]
in $\mathcal{C}$, so again the limiting object $\Lim\; F$ would, in practice, be written as $X\times_Z Y$.

\noindent For every object $A\in\mathcal{C}$, there is an isomorphism of categories
\[ \phi_A\maps {\rm Ps}(D,\Cat)(\Delta 1(\--),\mathcal{C}(A,F(\--)))\to \mathcal{C}(A,\Lim\; F).\]

\noindent There are three morphisms in the limiting cone
\[ \phi_{\Lim\; F}^{-1}(1_{\Lim\; F})\maps \Delta 1(\--)\To \mathcal{C}(\Lim\; F,F(\--)))\]
\noindent coming from the three objects of $D$.  These are the ``projections"
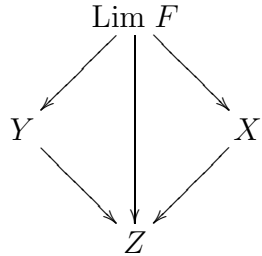
\[
   \xy
   (0,15)*+{\Lim\; F}="1";
   (-15,0)*+{Y}="3";
   (15,0)*+{X}="5";
   (0,-15)*+{Z}="7";
        {\ar^{} "1";"3"};
        {\ar^{} "1";"5"};
        {\ar_{} "1";"7"};
        {\ar^{} "3";"7"};
        {\ar_{} "5";"7"};
\endxy
\]
\noindent Since the cone is now a strong transformation rather than the strict transformations defining cones for strict $2$-limits, each $1$-morphism in $D$ gives us an invertible $2$-morphism in $\mathcal{C}$.  It follows that the pseudo pullback has $2$-cells as pictured
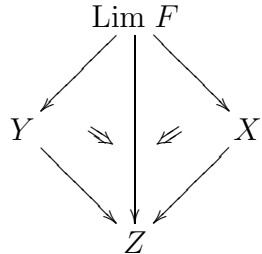
\[
   \xy
   (0,15)*+{\Lim\; F}="1";
   (-15,0)*+{Y}="2"; (15,0)*+{X}="3";
   (0,-15)*+{Z}="4";
        {\ar_{} "1";"2"};
        {\ar^{} "1";"3"};
        {\ar_{} "1";"4"};
        {\ar^{} "2";"4"};
        {\ar^{} "3";"4"};
        {\ar@{=>}_<<{\scriptstyle } (-6,0); (-3,-2)};
        {\ar@{=>}_<<{\scriptstyle } (6,0); (3,-2)};
\endxy
\]

The universal property is straightforward to verify.

$\bullet$  For each object $A\in\mathcal{C}$, an object
\[ \sigma\maps \Delta 1(\--)\to \mathcal{C}(A,F(\--)),\]
in ${\rm Ps}(D,\Cat)(\Delta 1(\--),\mathcal{C}(A,F(\--)))$ is a cone
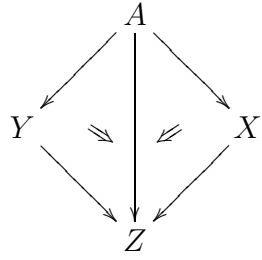
\[
   \xy
   (0,15)*+{A}="1";
   (-15,0)*+{Y}="2"; (15,0)*+{X}="3";
   (0,-15)*+{Z}="4";
        {\ar_{} "1";"2"};
        {\ar^{} "1";"3"};
        {\ar_{} "1";"4"};
        {\ar^{} "2";"4"};
        {\ar^{} "3";"4"};
        {\ar@{=>}_<<{\scriptstyle } (-6,0); (-3,-2)};
        {\ar@{=>}_<<{\scriptstyle } (6,0); (3,-2)};
\endxy
\]
We obtain a unique comparison map 
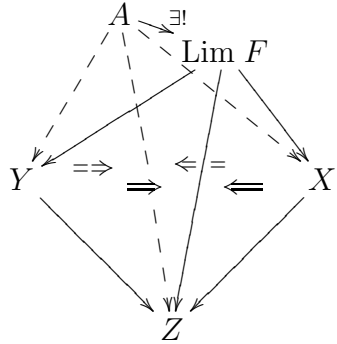
\[ \phi_A(\sigma)\maps A\to \Lim\; F,\]
\[
   \xy
   (-3,1)*+{}="5";
   (-2.5,-3)*+{}="6";
   (-7,22)*+{A}="1";
   (7,17)*+{\Lim\; F}="-1";
   (-20,0)*+{Y}="2"; (20,0)*+{X}="3";
   (0,-20)*+{Z}="4";
        {\ar^{\exists !} "1";"-1"};
        {\ar@{-->}_{} "1";"2"};
        {\ar@{-->}^{} "1";"3"};
        {\ar@{--}_{} "1";"5"};
        {\ar@{-->}_{} "6";"4"};
        {\ar_{} "-1";"2"};
        {\ar^{} "-1";"3"};
        {\ar_{} "-1";"4"};
        {\ar^{} "2";"4"};
        {\ar^{} "3";"4"};
        {\ar@{=>}_<<{\scriptstyle } (-6,-1); (-1.5,-1)};
        {\ar@{=>}_<<{\scriptstyle } (12,-1); (7,-1)};
        {\ar@{==>}_<<{\scriptstyle } (-13,1.5); (-8,1.5)};
        {\ar@{==>}_<<{\scriptstyle } (7,2); (0.5,2)};
\endxy
\]
and the diagram commutes by naturality of $\phi$.  That is, since $\phi$ is a strict invertible transformation, the cone at vertex $A$ must be equal to the limiting cone precomposed with the $1$-morphism $\phi_A(\sigma)$.  This amounts to a pair of $1$-cell equations and a pair of $2$-cell equations, which appear on the left and right sides of the above diagram.

$\bullet$ Now, consider an object $A\in\mathcal{C}$ and a pair of $1$-morphisms $h,k\maps A\to \Lim\; F$ such that for each $d\in D$, there is a $2$-morphism 
\[ M_d\maps \phi_{\Lim\; F}^{-1}(1_{\Lim\; F})_dh\To \phi_{\Lim\; F}^{-1}(1_{\Lim\; F})_dk\]
satisfying for each $\delta\maps d\to d'$ in $D$, an equation
\[ \phi_{\Lim\; F}^{-1}(1_{\Lim\; F})(\delta)\circ F(\delta)\cdot M_{d} = M_{d'}\circ\phi_{\Lim\; F}^{-1}(1_{\Lim\; F})(\delta)\]
between $2$-morphisms in $\mathcal{C}$.  These are the equations expressing that the collection of $2$-morphisms $\{M_d\}_{d\in D}$ is a modification, i.e., a $1$-morphism in ${\rm Ps}(D,\Cat)(\Delta 1(\--),\mathcal{C}(A,F(\--)))$.  Since $\phi_A$ is an isomorphism of categories, and thus, fully faithful, this modification maps to a unique $2$-morphism $\gamma\maps h\To k$ in $\mathcal{C}$ such that for each $d\in D$
\[ \phi_{\Lim\; F}^{-1}(1_{\Lim\; F})_{d}\cdot \gamma = M_d.\]
\noindent This is the two-dimensional aspect of the universal property.

In the definition that follows, we do not require the modification axioms as an extra condition.  We circumvent this requirement by defining the pullback as a weighted strict $2$-limit rather than a pseudo $2$-limit.  The analogous condition to the equation relating $M_d$ and $M_{d'}$ above, will be the naturality equation, which we will see is non-trivial in weighted limits.

Since the $2$-cells of the cones are invertible, one often inverts one of the $2$-cells, composes the resulting pair, and then discards the projection to $Z$ when defining the pseudo pullback by cones.  The resulting limit is called the {\em iso-comma object}.
\[
   \xy
   (0,15)*+{\Lim\; F}="1";
   (-15,0)*+{Y}="2"; (15,0)*+{X}="3";
   (0,-15)*+{Z}="4";
        {\ar_{} "1";"2"};
        {\ar^{} "1";"3"};
        {\ar^{} "2";"4"};
        {\ar^{} "3";"4"};
        {\ar@{=>}_<<{\scriptstyle } (3,0); (-3,0)};
\endxy
\]
We will give an explicit working definition of the iso-comma construction in Section~\ref{explicit}.
As suggested at the beginning of this section, pseudo pullbacks can also be defined as non-conical strict $2$-limits.  This definition is actually closer to that of the iso-comma object, since there are only $2$ ``distinct'' projection maps and one invertible $2$-morphism.

We take the same definition for the $2$-category $D$
\[
   \xy
   (-15,15)*+{d_1}="1";
   (0,0)*+{d_3}="2";
   (15,15)*+{d_2}="3";
        {\ar_{\delta_1} "1";"2"};
        {\ar^{\delta_2} "3";"2"};
\endxy
\]
but now define a non-conical weighting $2$-functor $W\maps D\to\Cat$.  The idea is that while the cones for strict $2$-limits do not contain $2$-morphisms, this weakened form of commutativity can be introduced by a non-conical weighting functor.  Let {\bf 1} be the terminal $2$-category.  We make the assignments
\[
   \xy
   (-10,0)*+{d_1 \mapsto W(d_1) :=}="1";
   (14,0)*+{\bullet}="2";
   (7,0)*+{id}="3";
        {\ar@/_1pc/@{-}^{} "2";"3"};
        {\ar@/_1pc/^{} "3";"2"};
\endxy
\]
\[
   \xy
      (-10,0)*+{d_2 \mapsto W(d_2) :=}="1";
   (14,0)*+{\bullet}="2";
   (7,0)*+{id}="3";
        {\ar@/_1pc/@{-}^{} "2";"3"};
        {\ar@/_1pc/^{} "3";"2"};
\endxy
\]
\[
   \xy
   (-10,0)*+{d_3 \mapsto W(d_3) :=}="1";
   (14,0)*+{\bullet}="2";
   (7,0)*+{id}="3";
   (48,0)*+{id}="4";
   (41,0)*+{\bullet}="5";
        {\ar@/_1pc/@{-}^{} "2";"3"};
        {\ar@/_1pc/^{} "3";"2"};
        {\ar@/_1pc/@{-}^{} "5";"4"};
        {\ar@/_1pc/^{} "4";"5"}; 
        {\ar@/_1pc/^{} "2";"5"};
        {\ar@/_1pc/^{} "5";"2"};
\endxy
\]
where the two non-identity arrows are an invertible pair.  Then the image of $W$ in $\Cat$ is a pair of functors mapping the terminal categories $W(d_1)$ and $W(d_2)$ each to one of the pair of objects in $W(d_3)$
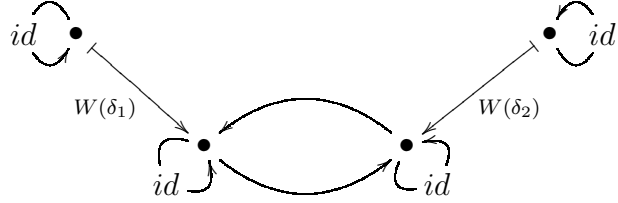
\[
   \xy
   (-3,15)*+{\bullet}="-2";
   (-10,15)*+{id}="-3";
   (67,15)*+{id}="-4";
   (60,15)*+{\bullet}="-5";   
   (14,0)*+{\bullet}="2";
   (9,-5)*+{id}="3";
   (45,-5)*+{id}="4";
   (41,0)*+{\bullet}="5";
        {\ar@/_1pc/@{-}^{} "2";"3"};
        {\ar@/_1pc/^{} "3";"2"};
        {\ar@/_1pc/@{-}^{} "5";"4"};
        {\ar@/_1pc/^{} "4";"5"}; 
        {\ar@/_1.5pc/^{} "2";"5"};
        {\ar@/_1.5pc/^{} "5";"2"};
        {\ar@/_1pc/@{-}^{} "-2";"-3"};
        {\ar@/_1pc/^{} "-3";"-2"};
        {\ar@/_1pc/@{-}^{} "-5";"-4"};
        {\ar@/_1pc/^{} "-4";"-5"}; 
        {\ar@{|->}_{W(\delta_1)} "-2";"2"};
        {\ar@{|->}^{W(\delta_2)} "-5";"5"};        
\endxy
\]

For a strict $2$-functor $F\maps D\to\C$, the $W$-weighted strict $2$-pullback of shape $D$ is an object $\Lim\;(W,F)$ together with a strict invertible transformation
\[ \phi \maps[D,\Cat](W(\--),\mathcal{C}(\--,F(\--)))\To \mathcal{C}(\--,\Lim\;(W,F))\]
\noindent in $[\mathcal{C},\Cat]$.  The strict transformation
\[ \phi_{\Lim\;(W,F)}^{-1}(1_{\Lim\;(W,F)})\maps W(\--) \to \mathcal{C}(\Lim\;(W,F),F(\--)).\] 
\noindent is the limiting cone
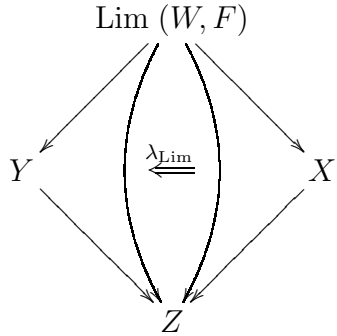
\[
   \xy
   (0,20)*+{\Lim\;(W,F)}="1";
   (-20,0)*+{Y}="2"; (20,0)*+{X}="3";
   (0,-20)*+{Z}="4";
        {\ar_{} "1";"2"};
        {\ar^{} "1";"3"};
        {\ar@/^1.5pc/_{} "1";"4"};
        {\ar@/_1.5pc/_{} "1";"4"};
        {\ar^{} "2";"4"};
        {\ar^{} "3";"4"};
        {\ar@{=>}_{\scriptstyle \lambda_{\Lim\;}} (3,0); (-3,0)};
\endxy
\]
\noindent The morphisms to $X$ and $Y$ come from the single object of $W(d_1)$ and $W(d_2)$, respectively.  Due to the non-conical weighting $W$, there are two morphisms in the cone with codomain $Z$, each coming from one of the objects in $W(d_3)$.  The invertible pair of morphisms in $W(d_3)$  yield the invertible $2$-morphism in the diagram and thereby ``weaken" the usual strict pullback.

The universal property is as follows.

$\bullet$ For each object $A\in\mathcal{C}$, an object
\[ \sigma\maps W(\--)\to \mathcal{C}(A,F(\--)),\]
in $[D,\Cat](W(\--),\mathcal{C}(A,F(\--)))$ is a cone
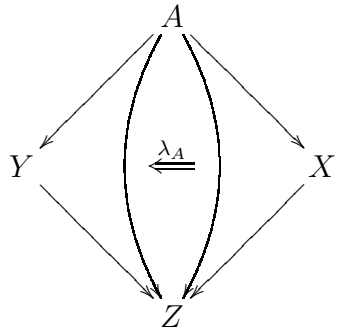
\[
   \xy
   (0,20)*+{A}="1";
   (-20,0)*+{Y}="2"; (20,0)*+{X}="3";
   (0,-20)*+{Z}="4";
        {\ar_{} "1";"2"};
        {\ar^{} "1";"3"};
        {\ar@/^1.5pc/_{} "1";"4"};
        {\ar@/_1.5pc/_{} "1";"4"};
        {\ar^{} "2";"4"};
        {\ar^{} "3";"4"};
        {\ar@{=>}_{\scriptstyle \lambda_A} (3,0); (-3,0)};
\endxy
\]
Since $\phi_A$ is an isomorphism, we obtain a unique comparison map 
\[ \phi_A(\sigma)\maps A\to \Lim\; F,\]
such that 
\[ \phi_{\Lim\;(W,F)}^{-1}(1_{\Lim\;(W,F)})_d\phi_A(\sigma) = \sigma_d\]
and
\[ \lambda_{\Lim\;}\cdot \phi_A(\sigma) = \lambda_A.\]
This follows from $\phi$ being a strict invertible transformation.  This is the one-dimensional aspect of the universal property.

$\bullet$ Now, consider an object $A\in\mathcal{C}$ and a pair of $1$-morphisms $h,k\maps A\to \Lim\;(W,F)$ such that for each object $d\in D$, there is a natural transformation 
\[ M_d\maps \phi_{\Lim\;(W,F)}^{-1}(1_{\Lim\;(W,F)})_dh\To \phi_{\Lim\;(W,F)}^{-1}(1_{\Lim\;(W,F)})_dk.\]
The naturality equation is analogous to the modification equation required in our previous discussion of the pseudo pullback.

The maps $\{M_d\;|\;d\in D\}$ are the components of a modification.   Since $\phi_A$ is an isomorphism of categories, and thus, fully faithful, this modification maps to a unique $2$-morphism $\gamma\maps h\To k$ in $\mathcal{C}$ such that for each $d\in D$,
\[ \phi_{\Lim\;(W,F)}^{-1}(1_{\Lim\;(W,F)})_{d}\cdot \gamma = M_d.\]
This is the two-dimensional aspect of the universal property.

\end{document}